\documentclass[12pt,oneside,british,a4wide]{amsart}
\usepackage{amsmath, amssymb,verbatim,appendix}
\usepackage[pagebackref=true]{hyperref}
\usepackage[mathscr]{eucal}
\usepackage{amscd}
\usepackage{amsthm}
\usepackage{stmaryrd}
\usepackage{enumitem}
\usepackage{comment}
\usepackage[latin1]{inputenc} 
\usepackage{tikz}
\usetikzlibrary{shapes,arrows}
\usepackage{cite}
\usepackage{url}
\usepackage{float}
\usepackage[normalem]{ulem}
\usepackage{xcolor}

\usepackage{setspace,a4wide,color,xcolor,graphicx}
 
\def\showauthornotes{1}

\ifnum\showauthornotes=1
\newcommand{\Authornote}[2]{{\sf\small\color{red}{[#1: #2]}}}
\else
\newcommand{\Authornote}[2]{}
\fi

\newtheorem{theorem}{Theorem}[section]
\newtheorem{lemma}[theorem]{Lemma}
\newtheorem{corollary}[theorem]{Corollary}
\newtheorem{proposition}[theorem]{Proposition}
\newtheorem{notation}[theorem]{Notation}
\newtheorem{question}[theorem]{Question}
\newtheorem{conjecture}[theorem]{Conjecture}

\newtheorem{observation}[theorem]{Observation}
\newtheorem{claim}[theorem]{Claim}
\newtheorem{fact}[theorem]{Fact}

\numberwithin{equation}{section}
\setlist[enumerate]{font={\rmfamily}}
\setlist[enumerate,1]{label={(\roman*)}}

\theoremstyle{definition}
\newtheorem{example}[theorem]{Example}
\newtheorem{definition}[theorem]{Definition}
\newtheorem{remark}[theorem]{Remark}

\def\E{\mathbb{E}}
\def\Z{\mathbb{Z}}
\def\R{\mathbb{R}}

\def\C{\mathbb{C}}

\def\F{\mathbb{F}}

\def\OP{\mathrm{OP}}
\def\NOP{\mathrm{NOP}}
\def\FOP{\mathrm{FOP}}
\def\NFOP{\mathrm{NFOP}}
\def\IP{\mathrm{IP}}
\def\NIP{\mathrm{NIP}}
\def\HOP{\mathrm{HOP}}
\def\NHOP{\mathrm{NHOP}}

\def\VC{\mathrm{VC}}
\def\GS{\mathrm{GS}}
\def\QGS{\mathrm{QGS}}
\def\IFOP{\mathrm{IFOP}}

\def\e{\epsilon}

\def\ra{\rightarrow}

\newcommand{\abar}{\bar{a}}
\newcommand{\bbar}{\bar{b}}
\newcommand{\cbar}{\bar{c}}
\newcommand{\dbar}{\bar{d}}
\newcommand{\ebar}{\bar{e}}
\newcommand{\fbar}{\bar{f}}
\newcommand{\gbar}{\bar{g}}
\newcommand{\ibar}{\bar{i}}
\newcommand{\jbar}{\bar{j}}
\newcommand{\kbar}{\bar{k}}

\newcommand{\rbar}{\bar{r}}
\newcommand{\sbar}{\bar{s}}
\newcommand{\tbar}{\bar{t}}
\newcommand{\ubar}{\bar{u}}
\newcommand{\vbar}{\bar{v}}
\newcommand{\wbar}{\bar{w}}
\newcommand{\xbar}{\bar{x}}
\newcommand{\ybar}{\bar{y}}

\newcommand{\calL}{\mathcal{L}}

\newcommand{\calH}{\mathcal{H}}
\newcommand{\calC}{\mathcal{C}}
\newcommand{\calM}{\mathcal{M}}
\newcommand{\calN}{\mathcal{N}}

\newcommand{\calB}{\mathcal{B}}

\newcommand{\calP}{\mathcal{P}}

\newcommand{\calI}{\mathcal{I}}

\newcommand{\calQ}{\mathcal{Q}}
\newcommand{\calR}{\mathcal{R}}
\newcommand{\calS}{\mathcal{S}}

\newcommand{\calX}{\mathcal{X}}
\newcommand{\calY}{\mathcal{Y}}

\newcommand{\codim}{{\mathrm{codim}}}

\newcommand{\Red}{{\mathrm{Red}}}
\newcommand{\Span}{{\mathrm{Span}}}
\newcommand{\At}{{\mathrm{At}}}
\newcommand{\rk}{{\mathrm{rk}}}

\def\disc{\operatorname{disc}}
\def\Red{\operatorname{Red}}
\def\dev{\operatorname{dev}}

\newenvironment{proofof}[1]{\indent{\itshape Proof of #1}.\;}{\qed}

\begin{document}
\title[Higher-order generalizations of stability and arithmetic regularity]{Higher-order generalizations of stability and arithmetic regularity}

\author{C. Terry}

\author{J. Wolf}

\address{Department of Mathematics, Statistics, and Computer Science, University of Illinois Chicago, Chicago IL 60607, USA}

\email{caterry@uic.edu}

\address{Department of Pure Mathematics and Mathematical Statistics, Centre for Mathematical Sciences, Wilberforce Road, Cambridge CB3 0WB, UK}

\email{julia.wolf@dpmms.cam.ac.uk}

\date{}

\begin{abstract}We define a natural notion of higher order stability and show that subsets of $\F_p^n$ that are tame in this sense can be approximately described by a union of low-complexity quadratic varieties, up to linear error. This generalizes the arithmetic regularity lemma for stable subsets of $\F_p^n$, proved in earlier work of the authors, to the realm of higher-order Fourier analysis.  

This result is strictly stronger than the structure theorem for sets of bounded $\VC_2$-dimension, first proved by the authors in earlier versions of this paper \cite{Terry.2021av2} and now available as a separate manuscript \cite{Terry.2021d}. Taken together, these results provide group theoretic analogues of results obtained for 3-uniform hypergraphs in \cite{Terry.2021b}.
\end{abstract}

\maketitle

\tableofcontents

\section{Introduction}

This paper uses arithmetic regularity lemmas as a framework for investigating higher order generalizations of model theoretic stability in the setting of groups.  We begin with a non-technical account of the history and motivation for our results, before providing a more detailed introduction.

\subsection{Stability and graph regularity}\label{subsec:stabilitygraphs}

Stability theory, first developed in depth by Shelah \cite{Shelah.1990o5n}, is a foundational area within model theory. It continues to be an active area of research, while also serving as a touchstone for work in more general contexts, such as NIP, simple, and $\textrm{NTP}_2$ theories (see, for example, \cite{Simon.2015,Kim.1997,Chernikov.2014,Kaplan.2020}).  On the other hand, Szemer\'{e}di's regularity lemma \cite{Szemeredi.1975} is a key tool in combinatorics and theoretical computer science, which at first glance appears completely unrelated to model theoretic stability.  Roughly speaking, the regularity lemma states that the vertex set of any finite graph can be partitioned into a bounded number of parts, so that almost all pairs of parts are ``regular" in the sense that they behave approximately like bipartite random graphs.

Szemer\'{e}di's regularity lemma is often viewed as an ingredient in proofs of other theorems, as is suggested by the word ``lemma" in its name. However, results of the last twenty years have shown it to be an interesting theorem in its own right, with deep connections to model theory. These results typically take the form of ``tame" regularity lemmas, meaning stronger regularity lemmas that hold under additional combinatorial assumptions.  In the first such work, Alon-Fischer-Newman  \cite{Alon.2007} and Lov\'{a}sz-Szegedy \cite{Lovasz.2007}  proved strong  regularity lemmas for graphs of bounded $\VC$-dimension (corresponding to NIP in model theory).  A few years later, Malliaris-Shelah \cite{Malliaris.2014} proved even stronger versions in the setting of stable graphs.   The main difference between the  NIP and stable regularity lemmas lies in the behavior of the so-called \emph{irregular pairs}.  In particular,  Malliaris and Shelah showed that stable graphs admit regular partitions with \emph{no irregular pairs}, a fact known to be false in general for NIP graphs.     Model theoretic stability theory is central to all known proofs that irregular pairs can be eliminated in the stable setting. The original paper of Malliaris and Shelah uses Shelah 2-rank and finite indiscernible sequences, among other ideas. Subsequent proofs of Malliaris-Pillay \cite{Malliaris.2016} and Pillay \cite{Pillay.2020} rely on the existence of unique non-forking extensions and generic compact domination, respectively.

While stability theory was developed well before connections to tame regularity lemmas were understood, it is important for the perspective of this paper to note that these developments could have taken place in the reverse order. In particular, examining tame regularity lemmas would inevitably have led to the discovery of several important aspects of stability theory.

\subsection{Stable group theory and linear arithmetic regularity}\label{subsec:linearinformal}

Beginning with the Baldwin-Saxl theorem \cite{Baldwin.1976} in the 1980s, model theoretic stability has  been shown to interact with group structure in deep and meaningful ways.  The study of such interactions, now known as \emph{stable group theory}, has guided several important research programs, for example the long-running quest to classify finite simple groups of finite Morley rank \cite{Cherlin.1979}, ground-breaking work on stability of the free group \cite{Sela.2013}, beautiful results about groups definable in $o$-minimal expansions of the real field \cite{Peterzil.2007}, and connections to topological dynamics found in NIP groups \cite{Pillay.2013}.

Szemer\'{e}di's regularity lemma has several generalizations beyond the setting of graphs, including analogues tailored to studying subsets of groups. The first such result was introduced by Green \cite{Green.2005}, who termed it an \emph{arithmetic regularity lemma}. In the context of the elementary abelian $p$-group $\F_p^n$ (a group of considerable interest in additive combinatorics), Green's arithmetic regularity lemma roughly states that given an arbitrary subset $A$, the group can be decomposed into the cosets of a bounded index subgroup, such that the set $A$ behaves approximately like a random set on almost every coset. We refer to such cosets as  ``regular with respect to the set $A$," and to such a coset partition as  a ``linear arithmetic regularity lemma for $A$." The word \emph{linear} distinguishes this style of arithmetic regularity lemma from the higher order versions we will discuss later.  Arithmetic regularity lemmas have since found numerous applications in arithmetic combinatorics and theoretical computer science (see, for example, Green's lecture notes \cite{Green.2007} and the surveys \cite{Gowers.2010,Hatami.2019}). 

It was Hrushovski \cite{Hrushovski.2012} who first realized the potential connections between model theory and arithmetic combinatorics, using tools from stable group theory to prove a non-abelian analogue of Freiman's theorem. However, it was not until several years later that the first ``tame" arithmetic regularity lemma was proved. Specifically, the authors showed in \cite{Terry.2019} that stable subsets of $\F_p^n$ admit efficient linear arithmetic regularity lemmas with no irregular cosets (see Section \ref{subsec:introlinear} for a precise statement).  This theorem was extended to more general groups in \cite{Terry.2020, Conant.2017p4, Conant.2020}, with the results of Conant-Pillay-Terry \cite{Conant.2017p4} demonstrating remarkably close connections between stable linear arithmetic regularity lemmas and classical tools from stable group theory, specifically generic types and model theoretic connected components.

After the first stable arithmetic regularity lemma appeared in \cite{Terry.2019}, analogous theorems were proved in the more general setting of bounded VC-dimension (i.e. the NIP setting in model theory) \cite{Sisask.2018,Alon.2018is,Conant.2018zd,CT}.   In $\F_p^n$, these results exhibit similar quantitative and qualitative improvements to the stable versions with one important caveat.  While in the stable setting, one can guarantee \emph{no irregular cosets}, this is (provably) not the case in the NIP setting.  All NIP arithmetic regularity lemmas to date employ approximate versions of the so-called stabilizer, an important tool from both stable and NIP group theory.  In \cite{Conant.2018zd}, direct connections to more sophisticated tools are revealed,  including Borel definability of types and generic compact domination.

In sum, the pursuit of tame linear arithmetic regularity lemmas has uncovered deep connections to stable and NIP group theory, with the two settings being distinguished in $\F_p^n$ by the behavior of the irregular cosets.   As in the case of graphs, a reverse approach, beginning with tame linear arithmetic regularity lemmas, would very likely have led to the development of many important aspects of stable and NIP group theory.

\subsection{The search for higher arity stability} Both stability and NIP can be considered to be ``binary" notions, in the sense that their definitions deal with definable  \emph{binary} relations (i.e. bi-partitioned formulas). About ten years ago, Shelah \cite{Shelah.2017} defined a ``higher arity" analogue of NIP, the $k$-ary version of which is denoted $\NIP_{k-1}$. This notion has been shown to be a robust higher arity analogue of NIP, including in various algebraic settings  \cite{Shelah.2017,Chernikov.2021, Hempel.2014}. In light of this work, the following question gained currency in the model theory community: what is the corresponding ``higher arity" analogue of stability? A suitable definition, along with evidence of its utility, remained elusive, despite regular developments around $\NIP_k$. Part and parcel of this question is the quest for higher arity analogues of stable group theory.

Several years before this question came to prominence, new ``higher order" versions of the regularity lemma were developed in combinatorics to overcome the limitations of the classical version.  Specifically, in the 2000's, more complicated hypergraph regularity lemmas were developed by several authors, which were able to produce finitary proofs  of the multidimensional Szemer\'{e}di theorem \cite{Gowers.2007,Frankl.2002, Nagle.2006, Rodl.2004}. Contemporaneous and highly quantitative work of Gowers on Szemer\'edi's theorem in the integers \cite{Gowers.1998, Gowers.2001} gave rise, over the course of the subsequent decade, to an analogous hierarchy of higher order \emph{arithmetic} regularity lemmas. The first level of this hierarchy is occupied by the linear arithmetic regularity lemma, and the next level is known as a \emph{quadratic arithmetic regularity lemma}.  In $\F_p^n$, the quadratic arithmetic regularity lemma  informally says that for any subset $A$, there is a partition of $\F_p^n$ into fibres of boundedly many linear and quadratic forms,  so that $A$ is locally ``uniform" on almost all fibres, in the sense of a local version of the Gowers $U^3$ norm.  We will refer to such a partition as a ``quadratic arithmetic regularity lemma for $A$" (see \cite{Terry.2021d} for a gentle introduction).

The goal of this paper and its companion \cite{Terry.2021b} is to study higher arity stability by taking the  ``reverse" approach. Specifically, our goal is to identify meaningful definitions of higher arity stability by proving tame versions of higher arity and higher order regularity lemmas.  The companion paper \cite{Terry.2021b} takes this approach from the point of view of hypergraph regularity, while this paper focuses on the setting of groups, specifically quadratic arithmetic regularity lemmas in elementary abelian $p$-groups. 

In view of our earlier discussion, a true higher arity analogue of stability should have implications for the ``irregular" parts of a higher order regularity lemma. The main contribution of this paper and its companion is the identification of a new ternary analogue of stability, called $\NFOP_2$, which does just that  in both the hypergraph and arithmetic settings.  Our companion paper \cite{Terry.2021b} shows that $\NFOP_2$ characterizes a particular behavior of irregular triads in regular decompositions of 3-uniform hypergraphs.  The main theorem of the present paper is a corresponding higher order arithmetic regularity lemma for $\NFOP_2$ subsets of elementary abelian $p$-groups (see Theorem \ref{thm:fop}). We show that such sets admit quadratic arithmetic regularity lemmas in which the irregular components exhibit additional algebraic structure (specifically, they are controlled by lower order parts of the decomposition).  A key ingredient in the proof is an approximate structure theorem for $\NIP_2$ subsets of elementary abelian $p$-groups (see Theorem \ref{thm:vc2finite}).  While the first proof of Theorem \ref{thm:vc2finite} originally appeared in early versions of this paper \cite{Terry.2021av2}, it has now been separated out from the present manuscript and given a new, simplified proof in \cite{Terry.2021d}.   Together, these results connect $\NFOP_2$ and $\NIP_2$ to \emph{quadratic} structure in groups, in analogy to the earlier results discussed above, which connect stability and NIP to \emph{linear} structure in groups.

 Our results represent the first ``higher order" tame regularity lemmas in the setting of groups. Further, while higher arity versions of NIP  have previously been studied in infinite groups \cite{Chernikov.2021, Hempel.2014}, our results constitute the first theorems about any notion of higher arity \emph{stability} in groups. To mark this distinction, we will give an example to show that the conclusion regarding the irregular components in our $\NFOP_2$ quadratic arithmetic regularity lemma is truly ``special," in the sense that it does not hold for all sets that are $\NIP_2$.  

The remainder of this introduction gives a more technical account of our results and the relevant history.

\subsection{Formal introduction to tame linear arithmetic regularity}\label{subsec:introlinear}

In this subsection we flesh out the discussion of the tame linear arithmetic regularity lemma initiated in Section \ref{subsec:linearinformal}.  To streamline the exposition, we state most results in $\F_p^n$, and discuss extensions to other classes of groups at the end.  

We begin by defining the order property and stable subsets of groups.

\begin{definition}[Order property and stability for subsets of groups]\label{def:stable}
Given $k\geq 1$, a group $G$ and a set $A\subseteq G$, we say that $A$ \emph{has the $k$-order property} (or \emph{has $k$-$\OP$}) if there exist $a_1,\ldots, a_k,b_1,\ldots, b_k\in G$ such that $a_i\cdot b_j\in A$ if and only if $i\leq j$.  

If $A$ does not have the $k$-order property in $G$, then we say that $A$ is \emph{$k$-stable}. 
\end{definition}

A basic example of a stable set is any coset of a subgroup. Indeed, a set $A\subseteq G$ is $2$-stable if and only if $A$ is a coset of a subgroup of $G$ (see \cite{Terry.2019}).  Moreover, any union of $k$ cosets of a subgroup is still $(k+1)$-stable (see \cite[Lemma 1.5]{Sanders.2018}). The first tame arithmetic regularity lemma, proved by the authors in the setting of elementary abelian $p$-groups \cite{Terry.2019}, tells us that these examples are canonical, in the sense that all stable sets are well approximated by unions of cosets.

\begin{theorem}[Arithmetic regularity lemma for stable sets \cite{Terry.2019}]\label{thm:stablegps}
For all primes $p$, integers $k\geq 1$ and $\e\in (0,1)$, there exists $M=\exp(\e^{-O_{p,k}(1)})$ such that for all sufficiently large $n$, the following holds.  For any $k$-stable $A\subseteq G= \F_p^n$, there exists $H\leqslant G$ of index at most $M$ such that  
\begin{enumerate}[label=\normalfont(\roman*)] 
\item for all cosets $C$ of $H$, $|A\cap C|/|C|\in [0,\e)\cup (1-\e,1]$;
\item there is a union $Y$ of cosets of $H$ such that $|A\Delta Y|\leq \e|G|$. 
\end{enumerate}
\end{theorem}

After Theorem \ref{thm:stablegps} appeared, several analogous results were proved under the more general assumption of bounded VC-dimension, which we now define. 

\begin{definition}[$\VC$-dimension of a subset of a group]\label{def:vcdim}
Given $k\geq 1$, a finite group $G$ and a subset $A\subseteq G$, we say that $A$ \emph{has $\VC$-dimension at least $k$} (or \emph{has $k$-$\IP$})\footnote{IP stands for ``independence property" in model theory.}  if there exist elements $\{a_i: i\in [k]\}\cup \{b_S:S\subseteq [k]\}\subseteq G$ such that $a_i\cdot b_S\in A$ if and only if $i\in S$.  We define the \emph{$\VC$-dimension of $A$}, denoted $\VC(A)$, to be the largest $k$ such that $A$ has $\VC$-dimension at least $k$.

If $\VC(A)<k$, we also say  that $A$ is $k$-NIP.
\end{definition}

It is not difficult to show that $k$-stable subsets are also $k$-NIP, but the converse does not hold in general (see \cite[Section 5]{Alon.2018is} or Section \ref{subsec:gs}). Nonetheless, it is possible to obtain a weaker version  of Theorem \ref{thm:stablegps} for sets of bounded VC-dimension, as first proved  by Alon, Fox, and Zhao in \cite{Alon.2018is}.

\begin{theorem}[Arithmetic regularity lemma for sets of bounded $\VC$-dimension \cite{Alon.2018is}]\label{thm:vc}
For all primes $p$, integers $k\geq 1$ and $\e\in (0,1)$, there exists $M=\e^{-O_{k}(1)}$ such that for any $A\subseteq G= \F_p^n$ with $\VC(A)<k$, there exists a subgroup $H\leqslant G$ of index at most $M$ and a set $\Sigma$ of cosets of $H$ such that 
\begin{enumerate}[label=\normalfont(\roman*)] 
\item   $|\Sigma|\leq \e |G/H|$ and for all cosets $C\notin \Sigma$, $|A\cap C|/|C|\in [0,\e)\cup (1-\e,1]$; 
\item there is a union $Y$ of cosets of $H$ such that $|A\Delta Y|\leq \e|G|$.   
\end{enumerate}
\end{theorem}

In the context of graph regularity lemmas, the hallmark difference between stability and NIP is the absence of irregular pairs under the assumption of stability.  Analogously, in the context of linear arithmetic regularity lemmas in $\F_p^n$, the difference between stability and NIP is the absence of ``irregular" cosets in the stable setting.   The following remark is required to see how this is a consequence of Theorems \ref{thm:stablegps} and \ref{thm:vc}.

\begin{remark}\label{rem:vc01}
When $A\subseteq \F_p^n$ has bounded VC-dimension, saying that a coset $C$ is regular with respect to $A$ in the sense of Green's linear arithmetic regularity lemma is roughly equivalent to saying that the density of $A$ on $C$ is close to $0$ or $1$ (up to polynomial changes in parameters).  Indeed, Lemma 5 in \cite{Terry.2019} shows that given an arbitrary set $A\subseteq \F_p^n$, if $A$ has density near $0$ or $1$ on $C$, then $C$ is regular with respect to $A$.  On the other hand, an easy adaptation of Lemma 6 from the same paper shows that given any $A\subseteq \F_p^n$ of bounded VC-dimension, and any coset $C$ which is regular  with respect to $A$, the density of $A$ must be near $0$ or $1$ on $C$. 
\end{remark}

 In light of Remark \ref{rem:vc01}, we see that  Theorem \ref{thm:stablegps} guarantees that $A$ is regular on \emph{all} cosets, while Theorem \ref{thm:vc} allows for a small collection of irregular cosets, labelled by $\Sigma$.  Further, one can show that this difference is essential, in the sense that  there exist sets of bounded VC-dimension for which the conclusions of Theorem \ref{thm:stablegps} fail.

\begin{definition}[Green-Sanders example]\label{def:GS}
Let $p>2$ be a prime and let $n\geq 1$ be an integer. Given $1\leq i\leq n$,  let $e_i$ denote the $i$th standard basis vector in $\F_p^n$ and let $H_i:=\{x\in \mathbb{F}_p^n: x_1=\ldots=x_i=0\}$. The \emph{Green-Sanders $p$-example of dimension $n$} is the subset of $\F_p^n$ defined by
\[A_{\GS}(p,n)=\bigcup_{i=1}^n (H_i+e_i).\]
\end{definition}

While the set $A_{\GS}(p,n)$ is simple to describe, it resists certain types of algebraic decompositions, making it useful as a counterexample.  It  first appeared   in a paper of Green and Sanders \cite{Green.2015qy4}, who showed that in any linear arithmetic regularity lemma for $A_{\GS}(p,n)$,   there is always an irregular coset.   These examples are relevant to the present discussion because they also have very small VC-dimension, as we will prove.

\begin{proposition}[Green-Sanders example has bounded $\VC$-dimension]\label{prop:gs3vc}
For all primes $p>2$ and integers $n\geq 1$, $A_{\GS}(p,n)$ has $\VC$-dimension at most $3$.\footnote{When $p=3$, we show this is tight. After the original version of this paper appeared, Gladkova discovered a better bound and shorter proof in the case $p>3$  \cite{Gladkova.2024}. For this reason, we have omitted proofs in the case  $p>3$ from the current version.}
\end{proposition}
 
The Green-Sanders result from  \cite{Green.2015qy4} implies that any subspace of $\F_p^n$ will always yield a coset which is irregular with respect to $A_{\GS}(p,n)$. Combining this with Remark \ref{rem:vc01}, we can deduce that there will always be a coset on which the density of $A_{\GS}(p,n)$ is bounded away from $0$ or $1$. We will provide a direct proof of this as well.

\begin{proposition}[Green-Sanders example requires an irregular coset]\label{prop:gsvc}
Let $p>2$ be a prime and let $n\geq 1$ be an integer. For any $0<\e<1/p$ and any subgroup $H\leqslant \F_p^n$, there is a coset $C$ of $H$ with $\e|C|<|A_{\GS}(p,n)\cap C|< (1-\e)|C|$.   
\end{proposition}
 
Together, Propositions \ref{prop:gs3vc} and \ref{prop:gsvc} show that Theorem \ref{thm:stablegps} cannot hold in general for subsets of $\F_p^n$ of bounded VC-dimension.    As we shall see, the Green-Sanders example turns out to play a key role in the higher order setting also.

We end this subsection by discussing generalizations of Theorems \ref{thm:stablegps} and \ref{thm:vc} to other classes of groups, and connections of  these results to model theory.  In the stable setting, a stronger  version of Theorem \ref{thm:stablegps}  in arbitrary finite groups was proved in \cite{Conant.2017p4}. The proof there uses an ultraproduct contruction, which, while giving ineffective bounds, is able to show that the subgroup $H$ in Theorem \ref{thm:stablegps} arises directly from the local connected component, $G_{\varphi}^0$, a central object in stable group theory.  Shortly after, a quantitative extension of Theorem \ref{thm:stablegps} to the abelian setting was proved in \cite{Terry.2020} using Shelah 2-rank (also a key ingredient in the proof of Theorem \ref{thm:stablegps} in \cite{Terry.2019}).  In the best iteration of these results to date, Conant \cite{Conant.2020} gave a quantitatively efficient proof of the stronger results from \cite{Conant.2017p4} for arbitrary finite groups, using approximate versions of the $G_{\varphi}^0$. We state the version of this for $\F_p^n$ below.

\begin{theorem}[Strong arithmetic regularity for stable sets \cite{Conant.2017p4}]\label{thm:stablestructure}
For all $k\geq 1$ and increasing functions $\psi:\mathbb{N}\rightarrow \mathbb{R}^+$,  there exists  $M=M(k,\psi)$  such that for any  $k$-stable $A\subseteq G=\F_p^n$, there is $H\leqslant G$ of codimension $m\leq M$ such that 
\begin{enumerate}[label=\normalfont(\roman*)] 
\item for all cosets $C$ of $H$,  $|A\cap C|/|C|\in [0,p^{-\psi(m)})\cup (1-p^{-\psi(m)},1]$;
\item there is a union $Y$ of cosets of $H$ so that $|A\Delta Y|\leq p^{-\psi(m)}|G|$.
\end{enumerate} 
\end{theorem} 

Turning to sets of bounded VC-dimension, we first mention that  Theorem \ref{thm:vc} was proved by  Alon, Fox, and Zhao in the broader context of abelian groups of bounded exponent \cite{Alon.2018is}.  Analogous results\footnote{It is worth noting here that outside the bounded exponent setting, sets of small VC-dimension are approximated by \emph{Bohr sets} rather than cosets of a subgroup.} to Theorem \ref{thm:vc} were proved around the same time by Sisask in the abelian setting \cite{Sisask.2018}, and in the setting of arbitrary finite groups by Conant, Pillay, and the first author \cite{Conant.2018zd}.  Both \cite{Alon.2018is} and \cite{Sisask.2018} obtain efficient quantitative bounds.  In \cite{Conant.2018zd}, Fourier analytic tools from additive combinatorics are replaced by structure theorems for compact groups, allowing for extension to the non-abelian setting. This shift requires an ultraproduct construction, which produces ineffective bounds while also giving rise to qualitative improvements, including the statement below for NIP subsets of $\F_p^n$.

\begin{theorem}[Strong arithmetic regularity for sets of bounded $\VC$-dimension \cite{Conant.2018zd}]\label{thm:vcstrong}
For all integers $k\geq 1$, reals $\e\in (0,1)$,  and increasing functions $\psi:\mathbb{N}\rightarrow \mathbb{R}^+$,  there exists $M=M(k,\e,\psi)$ such that for any $A\subseteq G=\F_p^n$ with $\VC(A)<k$, there exist  a subgroup $H\leqslant  G$ of codimension $m\leq M$, and a set of cosets $\Sigma$ of $H$ such that 
\begin{enumerate}[label=\normalfont(\roman*)] 
\item $|\Sigma|\leq \e |G/H|$; 
\item for all cosets $C\notin \Sigma$, $|A\cap C|/|C|\in [0,p^{-\psi(m)})\cup (1-p^{-\psi(m)},1]$; 
\item there is a union $Y$ of cosets of $H$ so that $|(A\Delta Y)\setminus Z|\leq p^{-\psi(m)}|G|$, where $Z=\bigcup_{D\in \Sigma}D$.   
\end{enumerate}
\end{theorem}

 A fully quantitative version of the results in \cite{Conant.2018zd} remains open in the non-abelian setting, and is closely related to the problem of finding a fully quantitative proof of the Breuillard-Green-Tao theorem \cite{BGT}.  Even in $\F_p^n$, there exists to date no quantitative proof of  Theorem \ref{thm:vcstrong}.

\subsection{Main results}
In this subsection we will state our main results.  We will first introduce Theorem \ref{thm:vc2finite}, which is a quadratic analogue of Theorem \ref{thm:vc} for a higher arity version of VC-dimension.\footnote{   This result was first proved in earlier versions of this paper \cite{Terry.2021av2}, but has now been moved to a separate manuscript with a simplified proof \cite{Terry.2021d}.} We will then turn to explaining our main result, Theorem \ref{thm:fop}, which is a quadratic analogue of  Theorem \ref{thm:stablegps} for a higher arity version of stability that we shall define. 

We begin by recasting our understanding of Theorems \ref{thm:stablegps} and \ref{thm:vc} in terms of linear factors. This terminology has its origin in ergodic theory but is now standard in additive combinatorics, and in particular in the study of higher order Fourier analysis on elementary abelian $p$-groups. A \emph{linear factor} is simply a set $\calL=\{r_1,\ldots, r_{\ell}\}$ of linearly independent vectors in $\F_p^n$. Any linear factor generates a partition of $\F_p^n$ into \emph{atoms}, meaning sets of the form 
\begin{align}\label{linear}
L(a_1,\ldots, a_{\ell})=\{x \in \F_p^n: x^T r_1=a_1,x^T r_2=a_2,\dots, x^T r_\ell=a_\ell \},
\end{align}
for some $a_1,\ldots, a_{\ell}\in \F_p$.   We let $\At(\calL)$ denote the set of atoms of $\calL$, and refer to $\ell$ as its \emph{complexity}.  The atom of $\calL$ containing $0\in \F_p^n$ is simply the codimension-$\ell$ subgroup $H= \{x\in \F_p^n: x^Tr_1=\ldots =x^Tr_{\ell}=0\}$, i.e. the atom $L(0)$ with label $0\in \F_p^{\ell}$. The partition $\At(\calL)$ then consists of the cosets of $L(0)$.  Conversely, any subgroup $H\leqslant \F_p^n$ of codimension $\ell$ naturally gives rise to a linear factor $\calL$ of complexity $\ell$, whose atoms are the cosets of $H$.  Theorems \ref{thm:stablegps} and \ref{thm:vc} are therefore statements about approximating stable and NIP sets by unions of atoms of \emph{linear} factors, making them tame \emph{linear} arithmetic regularity lemmas. On the other hand, the main result of this paper and its companion \cite{Terry.2021d} deal with higher order analogues known as \emph{quadratic} arithmetic regularity lemmas.  
 
 We will refrain from giving here a full description of a general quadratic arithmetic regularity lemma, instead referring the reader to \cite[Section 2]{Terry.2021d} for a comprehensive introduction to quadratic Fourier analysis. For the purposes of this introduction, it suffices to give a brief description of \emph{quadratic factors}, which constitute the ``structured" component in a quadratic arithmetic regularity lemma.

Roughly speaking, a quadratic factor partitions the group $\F_p^n$ into the fibers of a quadratic form.  More concretely, a quadratic factor will consist of a pair $\calB=(\calL,\calQ)$ where $\calL=\{r_1,\ldots, r_{\ell}\}$ is a linear factor and $\calQ=\{M_1,\ldots, M_q\}$ is a set of symmetric $n\times n$ matrices over $\F_p$.  An \emph{atom of $\calB$} is a set of the form
\begin{align}\label{quad}
B(a;b)=L(a)\cap \{x \in \F_p^n: x^T M_1 x =b_1, x^T M_2 x =b_2,\dots, x^T M_q x =b_q \},
\end{align}
for some $a=(a_1,\dots,a_\ell)\in \F_p^{\ell}$ and $b=(b_1,\ldots,b_q)\in \F_p^q$ (here we are using the notation for linear atoms introduced in (\ref{linear})).

We refer to the tuple $(a,b)\in \F_p^{\ell}\times \F_p^q$ as the \emph{label} of the atom $B(a;b)$ and denote  the set of atoms of $\calB$  by $\At(\calB)$.    The \emph{complexity} of $\calB$ is the pair $(\ell,q)$. For the purpose of most applications, there will be a \emph{rank} requirement on the factor (see Section \ref{sec:fourierreview}) but this may be safely disregarded at the level of this introduction.

Our higher order analogue of Theorem \ref{thm:vc} is stated in terms of the following higher arity analgoue of VC-dimension,  first defined by Shelah \cite{Shelah.2017}.

\begin{definition}[$\VC_2$-dimension of a subset of a group]\label{def:vc2dim}
Given $k\geq 1$, a finite group $G$ and a subset $A\subseteq G$, we say that $A$ \emph{has $\VC_2$-dimension at least $k$}  (or \emph{has $k$-$\IP_2$}) if there exist 
$$
\{a_S: S\subseteq [k]\times [k]\}\cup \{b_{i}:i\in [k]\}\cup \{c_i: i\in [k]\}\subseteq G
$$
 such that $b_i\cdot c_j\cdot a_S\in A$ if and only if $(i,j)\in S$.  The \emph{$\VC_2$-dimension of $A$}, denoted $\VC_2( A)$, is defined to be the largest $k$ such that $A$ has $\VC_2$-dimension at least $k$.
 
If $\VC_2(A)<k$, we also say that $A$ is $k$-$\NIP_2$.
\end{definition}

An atom of a quadratic factor turns out to be a typical example of a set with bounded $\VC_2$-dimension. Indeed, although such an atom may have unbounded VC-dimension, it will always have $\VC_2$-dimension at most $1$ (see \cite[Lemma 3.9]{Terry.2021d} and Section \ref{subsec:atomsprops} respectively).  Combining this with a result by Chernikov-Palac\'{i}n-Takeuchi \cite{Chernikov.2019}, one can deduce that a union of $k$ quadratic atoms has $\VC_2$-dimension bounded by some function of $k$.  Therefore,  unions of small numbers of quadratic atoms are basic examples of sets of bounded $\VC_2$-dimension. Theorem \ref{thm:vc2finite} below says that all sets of bounded $\VC_2$-dimension resemble these examples, in that they are well approximated by atoms of a bounded complexity quadratic factor  (see \cite[Theorem 1.9 and Corollary 1.10]{Terry.2021d}). 

In what follows, a \emph{growth function} refers to any increasing function from $\mathbb{R}^+\rightarrow \mathbb{R}^+$, and given two growth functions $\omega, \omega'$, we write $\omega\geq \omega'$ if $\omega(x)\geq \omega'(x)$ for all $x\in \mathbb{R}^+$. 

\begin{theorem}[Arithmetic regularity lemma for sets of bounded $\VC_2$-dimension \cite{Terry.2021d}]\label{thm:vc2finite}
For all primes\footnote{There are some technicalities in quadratic Fourier analysis in characteristic 2, which we chose to avoid in \cite{Terry.2021d} for the sake of clarity. However, Theorem \ref{thm:vc2finite} should hold in characteristic 2, and can likely be obtained as in \cite{Terry.2021d} with a little more attention to detail.} $p>2$, integers $k\geq 1$ and reals $\mu\in (0,1)$, there exists a growth function $\omega_0=\omega_0(p,k,\mu)$ such that for all growth functions $\omega\geq\omega_0$, there exists a constant $C_0=C_0(p,k,\mu,\omega)$ such that the following holds  for all sufficiently large $n$.  

Suppose that $A\subseteq G=\F_p^n$ satisfies $\VC_2(A)<k$.  Then there are integers $\ell,q\geq 0$, a quadratic factor $\calB=(\calL,\calQ)$ on $G$ of complexity $(\ell,q)$, and a set $\Sigma \subseteq \At(\calB)$ such that 
\begin{enumerate}[label=\normalfont(\roman*)] 
\item $\ell,q\leq C_0$;
\item $\calB$ has rank at least $\omega(\ell+q)$;
\item $|\Sigma|\leq \e |\At(\calB)|$; 
\item for all $B\in \At(\calB)\setminus \Sigma$, $|A\cap B|/|B|\in [0,\mu)\cup (1-\mu,1]$;
\item there is a union $Y$ of atoms of $\calB$ such that $|A\Delta Y|\leq \mu|G|$. 
\end{enumerate}
\end{theorem}

 In other words, Theorem \ref{thm:vc2finite} tells us that given a subset $A\subseteq \F_p^n$ of bounded $\VC_2$-dimension, there exists a high-rank quadratic factor of bounded complexity, so that $A$ has density near $0$ or $1$ on all atoms of $\calB$ \emph{outside a small collection of exceptional atoms}, labeled by $\Sigma$.  This implies that the set $A$ is approximately equal to a union  of atoms of said quadratic factor, as stated in part (v) above.  Theorem \ref{thm:vc2finite} should be viewed as a strong version of the quadratic arithmetic regularity lemma, as the following remark elucidates.
 
 \begin{remark}\label{rem:vc2}
Given a set $A\subseteq \F_p^n$ with bounded $\VC_2$-dimension and an atom $B$ of a high rank quadratic factor, saying that $B$ is ``regular" with respect to $A$ in the sense of quadratic arithmetic regularity is equivalent to saying that $A$ has density near $0$ or $1$ on $B$ (up to polynomial changes in the parameters involved).       Indeed, given a high rank quadratic factor $\calB$,  if an arbitrary set $A\subseteq \F_p^n$  has density near $0$ or $1$ on an atom $B$ of $\calB$, then \cite[Proposition 3.20]{Terry.2021d} shows that $B$ is  regular  with respect to $A$ in this sense.  Conversely, if $A$ has bounded $\VC_2$-dimension, then  \cite[Proposition 4.3]{Terry.2021d} shows that $A$ must have density near $0$ or $1$ on any atom of $\calB$ which is regular with respect to $A$  in this sense.  
 \end{remark}

In light of Remark \ref{rem:vc2}, the set $\Sigma$ of exceptional atoms appearing in Theorem \ref{thm:vc2finite} is also the set of irregular  atoms, in the sense of the quadratic arithmetic regularity lemma. We will implicitly use Remark \ref{rem:vc2} throughout the rest of the introduction to conflate ``regularity" with ``density near $0$ or $1$" when the set in question has bounded $\VC_2$-dimension (as will be the case for all sets under discussion).

Theorem \ref{thm:vc2finite} represents a quadratic version of Theorem \ref{thm:vc} in which VC-dimension has been replaced by the higher arity analogue, $\VC_2$-dimension. The main goal of the present paper is to prove a quadratic version of Theorem \ref{thm:stablegps}, in which \emph{stability} is replaced by a higher arity analogue.  To introduce these results, we begin by defining a new higher arity generalization of the order property, which we call the \emph{functional order property}.\footnote{This definition is related to, but distinct from, a previous unpublished definition of Takeuchi called ``$\OP_2$." See Section \ref{subsec:fopbasic} for further discussion.}

\begin{definition}[Functional order property ($\FOP_2$)]\label{def:fop}
Given $\ell\geq 1$, a group $G$ and a subset $A$, we say that $A$ \emph{has $\ell$-$\FOP_2$} if there exist elements $y_1, \dots,y_\ell, z_1, \dots,z_\ell$ in $G$, and for any function $f:[\ell]^2\rightarrow [\ell]$, there are $x^f_1,\dots,x^f_\ell$ in $G$ such that $x^f_i\cdot y_j\cdot z_k \in A$ holds if and only if $k\leq f(i,j)$. 

If $A$ does not have $\ell$-$\FOP_2$, we say that $A$ is $\ell$-$\NFOP_2$ .
\end{definition}

There are certain contexts in which $\FOP_2$ is undeniably the correct answer to the question: IP is to $\IP_2$ as OP is to what?  An example of such a result appears in Section \ref{ss:reduced} of this paper, and plays a crucial role in our proofs (see also Section 6.1 of \cite{Terry.2021b}).   Before stating our main theorem, we shall give some examples of $\NFOP_2$ sets and discuss what we might expect our main result to say.  A more in-depth look at Definition \ref{def:fop} and its motivation appears in Section \ref{sec:context}.

To build context,  we first consider how $\NFOP_2$ interacts with the other combinatorial definitions mentioned above (we postpone all proofs to Section \ref{subsec:fopbasic}). First, it is not difficult to show that any $\NFOP_2$ set has bounded $\VC_2$-dimension. 

\begin{lemma}[$\NFOP_2\Rightarrow \NIP_2$]\label{lem:nip2fop}
Let $G$ be a group and suppose that $A\subseteq G$. For any $\ell\geq 1$, if $A$ does not have $\ell$-$\FOP_2$, then $\VC_2(A)< \ell$.
\end{lemma}

Thus, any $\ell$-$\NFOP_2$ subset of a group is also $\ell$-$\NIP_2$.  On the other hand, the collection of $\ell$-$\NFOP_2$ subsets contains all $\ell$-$\NIP$ sets.

\begin{lemma}[$\NIP\Rightarrow \NFOP_2$]\label{lem:vcfop}
Let $G$ be a group and suppose that $A\subseteq G$. For any\footnote{The restriction $\ell \geq 2$ is required due to a discrepancy in the degenerate cases: a $1$-$\NFOP_2$ set must be empty, while a $1$-NIP set may be empty or the whole group.} $\ell\geq 2$, if $A$ is $\ell$-$\NIP$, then $A$ is $\ell$-$\NFOP_2$. 
\end{lemma}

We have already seen that quadratic atoms are basic examples of $\NIP_2$ sets. It turns out that they are also basic examples of $\NFOP_2$ sets. In particular, we shall show that any atom of a quadratic factor has no $2$-$\FOP_2$ (see Section \ref{subsec:atomsprops}). We will further show that unions of quadratic atoms admit only short functional order properties (see Theorem \ref{thm:unionsatoms}).\footnote{We prove general closure properties for $\NFOP_2$ formulas in our companion paper \cite{Terry.2021b}, but include here a more explicit proof for unions of quadratic atoms in $\F_p^n$. See Section \ref{subsec:atomsprops} for further discussion.}  Consequently, unions of small numbers of quadratic atoms provide natural examples of $\NFOP_2$ sets. On the other hand, since $\NFOP_2$ sets are also $\NIP_2$ (by Lemma \ref{lem:nip2fop}),  Theorem \ref{thm:vc2finite} tells us any $\NFOP_2$ set $A\subseteq\F_p^n$ can be well approximated by atoms of a high rank quadratic factor of bounded complexity, in the sense that $A$ has density near $0$ or $1$ on all atoms of said factor, \emph{outside a small set $\Sigma$ of irregular atoms}.  

 The main result of this paper is a stronger version of Theorem \ref{thm:vc2finite} for $\NFOP_2$ sets, in which special properties hold of the set $\Sigma$ of irregular atoms. In analogy with Theorem \ref{thm:stablegps}, the reader might guess that we intend to prove that $\NFOP_2$ sets satisfy a version of Theorem \ref{thm:vc2finite} in which there are \emph{no} irregular atoms.  However,  this type of statement is too strong. In fact, the set $A_{\GS}(p,n)$ of Definition \ref{def:GS} provides a counterexample.  
 
 First, we recall from Proposition \ref{prop:gs3vc} that $A_{\GS}(p,n)$ has VC-dimension at most $3$, and therefore, by Lemma \ref{lem:vcfop}, $A_{\GS}(p,n)$ has no $4$-$\FOP_2$. On the other hand, we will prove Proposition \ref{prop:gslinear} below, which implies that any quadratic arithmetic regularity lemma for $A_{\GS}(p,n)$ will include at least one  irregular atom (in fact several).   We note that Remark \ref{rem:vc2} and Proposition \ref{prop:gsvc} allows us to conflate irregular atoms with those of intermediate density when it comes to the set $A_{\GS}(p,n)$. For this and several statements below, we remind the reader of the notation appearing in (\ref{linear}) and (\ref{quad}).

\begin{proposition}[Green-Sanders example requires an irregular atom]\label{prop:gslinear}
For all primes $p>2$ and all $\e\in (0,1)$, there exists a growth function $\omega=\omega(p,\e)$ such that for all $\ell,q\geq 0$ and all sufficiently large $n$ the following holds. 

For  any quadratic factor $\calB=(\calL,\calQ)$ on $\mathbb{F}_p^n$ of complexity $(\ell,q)$ and rank at least $\omega(\ell+q)$, there is an atom $B\in \At(\calB)$ such that 
$$
\frac{|B\cap A_{\GS}(p,n)|}{|B|}\in \Big(\frac{1-\e}{p},\frac{1+\e}{p}\Big).  
$$
In fact, the above holds for all atoms of the form $B=B(0;b)$ with $b\in \F_p^q$.  
\end{proposition}

In light of Remark \ref{rem:vc2} and Proposition \ref{prop:gslinear},  it is not possible to show that all $\NFOP_2$ sets admit quadratic arithmetic regularity lemmas with  \emph{no} irregular atoms.  On the other hand, this makes $A_{\GS}(p,n)$ useful as a case study in what the correct theorem for $\NFOP_2$ sets should say.  Indeed, $A_{\GS}(p,n)$ admits a rather special type of quadratic decomposition where the irregular atoms are highly constrained.   

\begin{proposition}[Special quadratic decompositions for Green-Sanders examples]\label{prop:gslinear1}
For all primes $p>2$ and all $\e\in (0,1)$, there exists a growth function $\omega=\omega(p,\e)$ such that for all integers $\ell\geq 1$ and $q\geq 0$, and all sufficiently large $n$, the following holds. 

There exists a quadratic factor $\calB=(\calL, \calQ)$ on $\mathbb{F}_p^n$ of complexity $(\ell,q)$ and rank at least $\omega(\ell+q)$ with the property that for every $B\in \At(\calB)$, either $B\subseteq L(0)$, in which case 
\[\frac{|B\cap A_{\GS}(p,n)|}{|B|}\in \left(\frac{1-\e}{p},1-\frac{1-\e}{p}\right),\] 
or $B\cap L(0)=\emptyset$, in which case $|B\cap A_{\GS}(p,n)|/|B|\in \{0,1\}$.  
\end{proposition}

Proposition \ref{prop:gslinear1} tells us that while $A_{\GS}(p,n)$ requires some  irregular atoms in any quadratic arithmetic regularity lemma, these can be made to take a  very special form. Specifically, the irregular atoms of $\calB$ are all contained in a single atom $L(0)$ of the linear factor $\calL$. On all $\calB$-atoms avoiding $L(0)$, the density of $A_{\GS}(p,n)$ is either $0$ or $1$.  It turns out that this is approximately the right type of statement for general $\NFOP_2$ sets.  In particular, we show that any $\NFOP_2$ set satisfies a stronger version of Theorem \ref{thm:vc2finite} in which all the irregular atoms are contained in a small number of purely linear atoms.

\begin{theorem}[Arithmetic regularity lemma for $\NFOP_2$ sets]\label{thm:fop}
For all primes $p>2$, integers $k\geq 1$, and increasing functions $\psi:\mathbb{N}\rightarrow \mathbb{R}^+$, there exists $D_0=D_0(p, k,\psi)$ such that for every growth function $\omega$, there exists $D_1=D_1(p, k,\psi,\omega)$ so that the following holds for all sufficiently large $n$.\footnote{We note that items (iii) and (v) are stronger than in earlier versions of this paper,  but do follow from the original proof in \cite{Terry.2021av2}.}

Suppose that $A\subseteq G=\F_p^n$ does not have $k$-$\FOP_2$. Then there exist integers $\ell\geq 1$ and $q\geq 0$,   a quadratic factor $\calB=(\calL,\calQ)$ on $G$ of complexity $(\ell,q)$, and a set $\Sigma\subseteq \At(\calL)$ such that 
\begin{enumerate}[label=\normalfont(\roman*)] 
\item $\ell\leq D_1$ and $q\leq D_0$;
\item $\calB$ has rank at least $\omega(\ell+q)$; 
\item $|\Sigma|\leq p^{-\psi(q)}|\At(\calL)|$;
\item for all $L\in \At(\calL)\setminus \Sigma$ and all $Q\in \At(\calQ)$, 
\[\frac{|A\cap L\cap Q|}{|L\cap Q|}\in [0,p^{-\psi(q)})\cup (1-p^{-\psi(q)},1];\] 
\item for some union $Y$ of atoms of $\calB$, $|A\Delta Y |\leq p^{-\psi(q)}|G|$.
\end{enumerate}
\end{theorem}

Let us emphasize the differences between Theorem \ref{thm:fop} and Theorem \ref{thm:vc2finite}.  First and foremost,   in Theorem \ref{thm:fop}, there is the additional condition that all irregular atoms are contained in a small collection $\Sigma$ of \emph{purely linear atoms}. Moreover, all approximations are in terms of a \emph{much smaller error}, namely $p^{-\psi(q)}$, which can be made arbitrarily small compared with the complexity of the quadratic part of the factor.

We will show that the conclusions of Theorem \ref{thm:fop} do not hold for all sets with bounded $\VC_2$-dimension.  Consequently, Theorems \ref{thm:vc2finite} and \ref{thm:fop} can be used to distinguish $\NIP_2$, a higher arity version of VC-dimension, from $\NFOP_2$, a higher arity version of stability.  We will achieve this using a quadratic analogue of the Green-Sanders example introduced in Definition \ref{def:GS}.  To define this example, we will use the fact that for odd $n$, there exists a set $\{M_1,\ldots, M_n\}$  of symmetric $n\times n$ matrices over $\F_p$, all non-trivial linear combinations of which have rank $n$. The existence of such a set of matrices follows from \cite[Lemma 3]{Dumas.2011}, for example.

\begin{definition}[Quadratic Green-Sanders example]\label{def:QGS}
 Let $p>2$ be a prime and let $n\geq 1$ be an odd integer. Let $(M_i)_{i=1}^n$ be as described in the preceding paragraph.  Given $1\leq i\leq n$, let
\[Q_i(e_i)=\{x\in \F_p^n:x^TM_1x=x^TM_2x=\dots = x^TM_{i-1}x=0, x^TM_ix=1\}.\]
The \emph{quadratic Green-Sanders $p$-example of dimension $n$}, denoted $A_{\QGS}(p,n)$, is the subset of $\mathbb{F}_p^n$ defined by
\[A_{\QGS}(p,n)=\bigcup_{i=1}^n Q_i(e_i).\]
\end{definition}

We show that the conclusions of Theorem \ref{thm:fop} fail for this example in a rather strong way.

\begin{theorem}\label{thm:ex1}
There exists an increasing function $\psi:\mathbb{N}\rightarrow \mathbb{R}^+$ so that for all primes $p>2$, there is a growth function $\omega=\omega(p)$ such that for all integers $\ell\geq 1$ and $q\geq 0$,  all sufficiently large odd integers $n$, and any quadratic factor $\calB=(\calL,\calQ)$ on $\F_p^n$ of complexity $(\ell,q)$ and rank at least $\omega(\ell+q)$, the following hold.
\begin{enumerate}[label=\normalfont(\roman*)] 
\item For all $L\in \At(\calL)$, there is $Q\in \At(\calQ)$ such that 
\[p^{-\psi(q)}< \frac{|A_{\QGS}(p,n)\cap L\cap Q|}{|L\cap Q|}< 1-p^{-\psi(q)}.\]
\item For any union $Y$ of elements of $\At(\calB)$, $|A\Delta Y|>p^{-4\psi(q)}|G|$. 
\end{enumerate}
\end{theorem}

Theorem \ref{thm:ex1} part (i) shows that $A_{\QGS}(p,n)$ fails to satisfy conclusions (iii)-(iv) of Theorem \ref{thm:fop} with respect to the growth function $\psi$ and \emph{any} $\Sigma\subseteq  \At(\calL)$.  Part (ii) of Theorem \ref{thm:ex1} shows $A_{\QGS}(p,n)$ fails to satisfy conclusion (v) of Theorem \ref{thm:fop} with respect to the growth function $4\psi$. 

The quadratic Green-Sanders examples were first defined by the authors in earlier versions of this paper \cite{Terry.2021av2}, where they were shown to satisfy Theorem \ref{thm:ex1} above  and   conjectured to have bounded $\VC_2$-dimension.  This conjecture was subsequently answered in the affirmative by Gladkova \cite{Gladkova.2024}.

\begin{theorem}[Quadratic Green-Sanders example has bounded $\VC_2$-dimension \cite{Gladkova.2024}]\label{thm:Gladkova}
For all primes $p>2$ and odd integers $n\geq 1$, the set $A_{\QGS}(p,n)$ of Definition \ref{def:QGS} has $\VC_2$-dimension at most $501$.
\end{theorem}

Combining Theorems \ref{thm:ex1} and \ref{thm:Gladkova}, the sets $A_{\QGS}(p,n)$ demonstrate that the conclusions of Theorem \ref{thm:fop} cannot be extended to the $\NIP_2$ setting.  This aligns with the results of our companion paper on tame hypergraph regularity \cite{Terry.2021b}, where we show that $\NFOP_2$ is \emph{characterized} by an analogous hypergraph regularity lemma that is strictly stronger than the type of regularity lemma characterizing $\NIP_2$.  We remark that Theorem \ref{thm:ex1} combined with Theorem \ref{thm:fop} implies that for every $k \geq 1$, a sufficiently  large quadratic Green-Sanders example must have $k$-$\FOP_2$.  We will also prove this explicitly using the machinery developed for the proof of Theorem \ref{thm:fop} (see Corollary \ref{cor:qgscor}).   

We end this subsection by conjecturing a stronger version of Theorem \ref{thm:vc2finite}, suggested by the analogous results in the linear case. In particular,  we conjecture that one can  improve  the density bounds on the regular atoms, as in Theorem \ref{thm:fop}.

\begin{conjecture}\label{con:strucVC2}
 For all primes $p>2$, integers $k\geq 1$, reals $\e\in (0,1)$ and increasing function $\psi:\mathbb{N}\rightarrow \mathbb{R}^+$, there exists $D_0=D_0(p, k,\e,\psi)$ such that for every growth function $\omega$, there exists $D_1=D_1(p, k,\e,\psi,\omega)$  such that for all sufficiently large $n$ the following holds. 

Suppose that $A\subseteq G=\F_p^n$ has $\VC_2(A)<k$.  Then there exist integers $\ell\geq 1$ and $q\geq 0$, a quadratic factor $\calB=(\calL,\calQ)$ on $G$ of complexity $(\ell,q)$, and $\Sigma\subseteq \At(\calB)$ such that
\begin{enumerate}[label=\normalfont(\roman*)]
\item $\ell\leq D_1$ and  $q\leq D_0$;
\item  $\calB$ has rank at least $\omega(\ell+q)$; 
 \item $|\Sigma|\leq \e |\At(\calB)|$;
\item for all $B\in \At(\calB)\setminus \Sigma$, $\frac{|A\cap B|}{|B|}\in [0,p^{-\psi(q)})\cup (1-p^{-\psi(q)},1]$;
\item for some union $Y$ of atoms of $\calB$,  $|(A\Delta Y)\setminus Z|\leq p^{-\psi(q)}|G|$, where $Z=\bigcup_{X\in \Sigma}X$.
\end{enumerate}
\end{conjecture}

We note that the conclusion of Conjecture \ref{con:strucVC2} is weaker than that of Theorem \ref{thm:fop}, since the set of irregular atoms is not assumed to be constrained to a small number of linear atoms.  Further, it claims a weaker bound on the number of irregular atoms, namely $\e |\At(\calB)|$ rather than the bound implied by part (iii) of Theorem \ref{thm:fop} (which has the form $p^{-\psi(q)}|\At(\calB)|$). This weaker bound is necessary. Indeed, if the bound on $\Sigma$ in Conjecture \ref{con:strucVC2} were replaced with one of the form $p^{-\psi(q)}|\At(\calB)|$, one could deduce all conclusions of Theorem \ref{thm:fop} for sets of bounded $\VC_2$-dimension.  However, we know this is impossible due to the properties of $A_{\QGS}(p,n)$ described above.  On the other hand, it is straightforward to see that $A_{\QGS}(p,n)$ is not a counterexample to Conjecture \ref{con:strucVC2}.

\subsection{Connections to work on hypergraphs}\label{subsec:discussion}
The results in this paper are intimately connected with parallel work of the authors \cite{Terry.2021b} in the setting of 3-uniform hypergraphs. Indeed, the shape of the group theoretic results obtained here motivated several of the statements in \cite{Terry.2021b}, and the aforementioned Green-Sanders example of Definition \ref{def:GS} plays a pivotal role there. Conversely, the results in \cite{Terry.2021b} help contextualize the results in the present paper. To explain these connections, we must first give an informal explanation of hypergraph regularity (see Appendix \ref{ss:hypergraphs} for details). 

Recall that a $3$-uniform hypergraph is a pair $H=(V,E)$, where $V$ is a vertex set and $E\subseteq{V\choose 3}$.  In this setting, a regular decomposition $\calP$ consists of a partition $V=V_1\cup \ldots \cup V_t$ of the vertex set, along with, for each $ij\in {[t]\choose 2}$, a partition of the pairs of vertices between $V_i$ and $V_j$, denoted by $K_2[V_i, V_j]=\bigcup_{\alpha\leq \ell}P_{ij}^{\alpha}$. Roughly speaking, such a decomposition yields a partition of ${V\choose 3}$ into the so-called \emph{triads} of $\calP$, each of which is a $3$-partite graph of the form $G_{ijk}^{\alpha\beta\gamma}=(V_i\cup V_j\cup V_k, P_{ij}^{\alpha}\cup P_{jk}^{\beta}\cup P_{ik}^{\gamma})$.  

A triad $G_{ijk}^{\alpha\beta\gamma}$ of $\calP$ is said to be \emph{regular with respect to $H$} if the bipartite graphs $P_{ij}^{\alpha}$, $P_{jk}^{\beta}$, and $P_{ik}^{\gamma}$ are each sufficiently quasirandom, and the edges of the hypergraph $H$ are ``uniformly distributed" across the triples of vertices forming  triangles in $G_{ijk}^{\alpha\beta\gamma}$ (see Definition \ref{def:regtrip}). A decomposition $\calP$ is then considered a \emph{regular decomposition of $H$} if almost all triples of vertices belong to a triad of $\calP$ which is regular with respect to $H$ (see Definition \ref{def:regptn}).  The main problem motivating \cite{Terry.2021b} is to characterize when hereditary properties of $3$-uniform hypergraphs admit regular decompositions in which the irregular  triads satisfy certain special properties.

The special property relevant to the results of this paper is as follows.  A hereditary property of $3$-uniform hypergraphs $\calH$ is said to \emph{admit linear $\disc_{2,3}$-error} \cite[Definition 2.11(3)]{Terry.2021b} if every sufficiently large element of $\calH$ has a regular decomposition $\calP$ such that for some small set $\Omega\subseteq {[t]\choose 3}$, every irregular triad of $\calP$ has vertex set $V_i\cup V_j\cup V_k$ for some $ijk\in \Omega$.   The name ``linear error" was chosen due to an explicit connection to Theorem \ref{thm:fop} (see Appendix \ref{ss:hypergraphs}).

We show in \cite[Theorem 1.7]{Terry.2021b} that a hereditary property of $3$-uniform hypergraphs admits linear $\disc_{2,3}$-error if and only if it does not have the functional order property.  We will show in Appendix \ref{ss:hypergraphs} that Theorem \ref{thm:fop}  matches up in a precise way with this theorem about 3-uniform hypergraphs.  Furthermore, we will explain how Theorem \ref{thm:vc2finite} aligns with an analogous result about $3$-uniform hypergraphs of bounded $\VC_2$-dimension.

\subsection{Future directions}\label{subsec:future}
It is straightforward to imagine higher-order generalizations of the work presented in this paper. Indeed, it is easy to write down a definition of a higher-order analogue of the functional order property, denoted by $\FOP_k$, where $k\geq 3$ (see \cite[Definition 2.61]{Terry.2021b}). We conjecture that subsets of $\F_p^n$ which are $\NFOP_k$ are essentially unions of atoms of a bounded complexity factor of degree $k$, up to an error that is contained in a small number of atoms of degree $k-1$. 

In a similar vein, it is natural to attempt to generalize the results in this paper beyond elementary abelian $p$-groups. From the point of view of arithmetic combinatorics, an obvious candidate are the prime cyclic groups, where much of the higher-order Fourier analysis required in this paper is already well developed. Having said this, the technical effort required for a quantitative generalization to this class of groups is likely to be substantial, and the difficulties of moving beyond the case of prime cyclic groups are probably at present insurmountable. 

It is here that model theory is likely to make a major contribution. In previous work on the linear case, in particular \cite{Conant.2017p4,Conant.2018zd}, it was established that strong linear structure is intimately related to the so-called \emph{connected component} (in the model theoretic sense) of the group. We expect an analogous model theoretic object to lie at the heart of strong quadratic decomposition theorems, which will almost certainly be related to so-called ``nilstructures" as laid out, for example, in \cite{Green.2012xge, Candela.2019}.

\medskip

\textbf{Structure of the paper.}
In Section \ref{sec:fourierreview}, we present the definitions of linear and quadratic factors along with some other preliminaries. The reader is referred to \cite{Terry.2021d} for a more comprehensive introduction to quadratic Fourier analysis and the associated regularity lemmas.  In Section \ref{sec:context}, we establish foundational facts about the functional order property, and study the key examples discussed in the introduction. Section \ref{sec:sufficient} uses a counting lemma from \cite{Terry.2021d} to prove the first ingredient for our main theorem, namely a sufficient condition for finding $k$-$\FOP_2$ based on the existence of certain order properties in an auxiliary $3$-coloring of the configuration (or label) space of a quadratic factor.  An analogous condition is given for $k$-$\IP_2$ as well.  Section \ref{sec:FOP} contains the proof of our main theorem, Theorem \ref{thm:fop}.  This proof relies on the application of a generalization of the stable linear arithmetic regularity lemma of \cite{Terry.2019} (i.e. Theorem \ref{thm:stablegps}), applied to an auxiliary $3$-coloring of the configuration space.  The proof of said generalization of \cite{Terry.2019} is the topic of Section \ref{sec:stablereg}. Section \ref{sec:FOP} makes heavy use of a standard translation between the order property and Shelah 2-rank (discussed in Section \ref{sec:FOP} in terms of trees), which also played a crucial role in earlier work of the authors \cite{Terry.2019, Terry.2020}.  Finally, Appendix \ref{ss:hypergraphs} describes how the results of this paper relate to hypergraph regularity and the main results of \cite{Terry.2021b}, Appendix \ref{app:easy} contains proofs of elementary lemmas needed for Section \ref{sec:FOP}, and Appendix \ref{app:generalprops} contains the rather technical proofs of several results stated in Section \ref{sec:context}.

\medskip

\textbf{Note on earlier versions of this work.} An earlier version of this manuscript \cite{Terry.2021av2} contained the first proofs of Theorem \ref{thm:vc2finite}, and its key ingredient, Proposition \ref{prop:countinggp}. The proofs of both results have been transferred and completely reworked into the separate manuscript \cite{Terry.2021d} to allow the present version of this paper to focus exclusively on higher-order stability. Some further refinement and reorganization has occurred between the earlier version \cite{Terry.2021av2} and this one, with the aim of sharpening its focus on the main result, Theorem \ref{thm:fop}.

\medskip

\textbf{Acknowledgements.} The first author would like to thank Maryanthe Malliaris for many useful mathematical conversations.  She also would like to thank Gabriel Conant for crucial observations about an early version of the proof of the main theorem, as well as for the diagrams appearing in the appendix. The second author is grateful to Tim Gowers for several helpful conversations during the preparation of this manuscript. She is also indebted to Robert Calderbank for bringing \cite{Delsarte.1975} to her attention during a visit in 2018. Both authors would like to thank Val Gladkova for pointing out an error in the proof of Proposition \ref{prop:qgsbad} in a previous version of the manuscript.

The first author has been partially supported by NSF grants DMS-2115518 and DMS-2239737, and a Sloane Research Fellowship. During the final stages of preparing this manuscript, the second author was supported by an Open Fellowship from the UK Engineering and Physical Sciences Research Council (EP/Z53352X/1). The authors' collaboration has been supported by several travel grants over the years, including by the London Mathematical Society, the Simons Foundation, and the Association for Women in Mathematics.

\section{Preliminaries}\label{sec:fourierreview}

This section contains some notational conventions used throughout the rest of the paper, followed by preliminaries on linear and quadratic factors.

\subsection{Notation}

For the convenience of the reader, we collect here some notation that appears later in the paper. 

As is common in discrete mathematics, given an integer $k\geq 1$, we write $[k]=\{1,2,\dots,k\}$. Given a set $X$ and an integer $k\geq 1$, we write ${X\choose k}=\{Y\subseteq X: |Y|=k\}$.  By convention, the natural numbers start at $0$ in this paper.

We will occasionally have need to refer to graphs and hypergraphs.  A \emph{$k$-uniform hypergraph} (or \emph{$k$-graph} for short) is a pair $G=(V,E)$ where $V$ is the vertex set of $G$, and $E\subseteq {V\choose k}$ is its edge set.  Given a $k$-graph $G$, $V(G)$ denotes its vertex set, and $E(G)$ denotes its edge set. When $k=2$, we call this simply a \emph{graph}.  A \emph{bipartite graph} is a graph $G=(V,E)$ for which there exists a partition $V=U\cup W$ so that for all $e\in E$, $|e\cap U|= |e\cap W|=1$.  We will write $G=(U\cup W, E)$ to denote $G$ is bipartite with a distinguished compatible partition $V=U\cup W$.  Similarly, a \emph{tripartite $3$-graph} is a $3$-graph $G=(V,E)$ for which there exists a partition $V=U\cup W\cup Z$ so that for all $e\in E$, $|e\cap U|= |e\cap W|=|e\cap Z|=1$.  We will write $G=(U\cup W\cup Z, E)$ to denote $G$ is tripartite with a distinguished compatible partition $V=U\cup W\cup Z$.  When discussing edges of graphs or hypergraphs, we will use the notation $xy$ to denote a set of two distinct vertices $\{x,y\}$, and similarly, $xyz$ to denote a set of three distinct vertices $\{x,y,z\}$.

Given $a\in \F_p^k$ and $i\in [k]$, we let $a_i$ denote the $i$th coordinate of $a$.  We will frequently write $0$ to mean a vector consisting of all $0$'s, e.g. $0\in \F_p^m$.  In certain instances, where it adds clarity,  we will move to a convention of putting bars over vectors, but we will make clear when this convention is used.

We will use the following asymptotic notation. For a function $f$ and positive-valued function $g$ on a domain $D$, we write $f = O(g)$ if there exists a constant $C$ such that $|f(x)| \leq Cg(x)$ for all $x\in D$.  Similarly, we write $f=(1+O(g))h$ to mean that  for some constant $C$, $|f(x)-h(x)|\leq Cg(x)h(x)$ for all $x\in D$. When $C$ depends on another parameter, this is (usually) indicated by a subscript on $O$. We also use the notation $x=1\pm \e$ to mean that $|1-x|\leq \e$.

Given a group $G$ and a set $A\subseteq G$, we denote by $\neg A$ the complement of $A$ in $G$.  When convenient, we will sometimes write $A^1$ to denote the original set $A$, and write $A^0$ to denote its complement in $G$.

Finally, we shall assume throughout Sections 4 and 5 that the prime $p$ of the base field is greater than 2. The case of characteristic 2 should be amenable to similar treatment with a little more work, but we have chosen to focus on the case $p>2$ here for the sake of clarity.

\subsection{Linear and quadratic factors}

We begin by defining linear factors and their atoms.

\begin{definition}[Linear factor]\label{def:linfac}
Given an integer $\ell\geq 1$, a \emph{linear factor $\calL$ on $\F_p^n$ of complexity $\ell$} is a set $\calL=\{v_1,\ldots, v_{\ell}\}\subseteq \F_p^n$ of linearly independent vectors. The \emph{atoms of $\calL$} are the sets of the  form
$$
L(a)=\{x\in \F_p^n: x^T v_1=a_1,\ldots, x^T v_{\ell}=a_{\ell}\},
$$
for some $a=(a_1,\ldots, a_{\ell})\in \F_p^{\ell}$.
The tuple $a$ is called the \emph{label} of the atom $L(a)$. 

We denote the set of atoms of $\calL$ by $\At(\calL)=\{L(a): a\in \F_p^{\ell}\}$.
\end{definition}

\begin{definition}[Purely quadratic factor]\label{def:quadfac}
Given an integer $q\geq 1$, a \emph{purely quadratic factor on $\F_p^n$ of complexity $q$} is a set $\calQ=\{M_1,\ldots, M_q\}$ of symmetric $n\times n$ matrices with entries in $\F_p$.  The \emph{atoms of $\calB$} are the sets of the form 
$$
Q(b)=\{x\in \F_p^n: x^TM_1x=b_1,\ldots, x^TM_qx=b_q\},
$$
for some $b=(b_1,\ldots, b_q)\in \F_p^q$.  The tuple $b$ is the \emph{label} of the atom $Q(b)$. 

We denote the set of atoms of $\calQ$ by $\At(\calQ)=\{Q(b): b\in \F_p^n\}$. 
\end{definition}

This paper focuses on the following combinations of linear and purely quadratic factors.

\begin{definition}[Quadratic factor]\label{def:quadfac}
Given integers $\ell,q\geq 1$, a \emph{quadratic factor $\calB$ on $\F_p^n$ of complexity $(\ell,q)$} is a pair $(\calL,\calQ)$ where $\calL=\{v_1,\ldots, v_{\ell}\}$ is a linear factor on $\F_p^n$, and $\calQ=\{M_1,\ldots, M_q\}$ is a set of symmetric $n\times n$ matrices with entries in $\F_p$.  The \emph{atoms of $\calB$} are the sets of form 
$$
B(a;b)=\{x\in \F_p^n: x^T v_1=a_1,\ldots, x^T v_{\ell}=a_{\ell}, x^TM_1x=b_1,\ldots, x^TM_qx=b_q\},
$$
for some $a=(a_1,\ldots, a_{\ell})\in \F_p^{\ell}$ and $b=(b_1,\ldots, b_q)\in \F_p^q$. The tuple $(a,b)$ is the \emph{label} of the atom $B(a;b)$. 

 We denote the set of atoms of $\calB$ by  $\At(\calB)=\{B(a;b): (a,b)\in \F_p^{\ell}\times \F_p^q\}$.
\end{definition}
Observe that, in the notation of Definition \ref{def:quadfac}, we have 
$$
B(a;b)=L(a)\cap Q(b).
$$

We next extend  the above definitions and notation to factors with trivial complexities.  For this, we let $\F_p^0$  denote the unique $0$-dimensional $\F_p$-vector space consisting of the single element $0\in \F_p$.  

\begin{definition}[Trivial linear and purely quadratic factors]
The \emph{linear factor on $\F_p^n$ of complexity $0$} is defined to be $\calL=\emptyset$. It has a unique atom, $\F_p^n$, denoted by $L(0)$ where $0\in \F_p^0$.  

The \emph{purely quadratic factor on $\F_p^n$ of complexity $0$} is defined to be $\calQ=\emptyset$. It has a unique atom, $\F_p^n$, denoted by $Q(0)$ where $0\in \F_p^0$.   
\end{definition}
These conventions allow us to define quadratic factors on $\F_p^n$ with complexity    $(\ell,q)$, where one or both of $\ell$,$q$ are $0$ as follows. 

\begin{definition}[Quadratic factors with trivial complexities]
 Given an integer $q\geq 0$, a  \emph{quadratic factor $\calB$ on $\F_p^n$ of complexity $(0,q)$} is a pair $\calB=(\emptyset,\calQ)$ where $\calQ$ is a purely quadratic factor on $\F_p^n$ of complexity $q$. In this case, we   set $\At(\calB)=\At(\calQ)$.  
 
 Given an integer $\ell\geq 0$, a  \emph{quadratic factor $\calB$ on $\F_p^n$ of complexity $(\ell,0)$} is a pair $\calB=(\calL,\emptyset)$ where $\calL$ is a  linear factor on $\F_p^n$ of complexity $\ell$. In this case, we say that $\calB$ is \emph{purely linear} and set $\At(\calB)=\At(\calL)$.  
 \end{definition}

Given a quadratic factor $\calB$ on $\F_p^n$ of complexity $(\ell,q)$ where one or both of $\ell,q$ are $0$,  it now makes sense to define $B(a;b)=L(a)\cap Q(b)$ for each $(a,b)\in \F_p^{\ell}\times \F_p^q$. With these conventions we will then have $\At(\calB)=\{B(a;b): (a,b)\in \F_p^{\ell}\times \F_p^q\}$.  

From here on, a linear or quadratic factor will generally mean one with possibly trivial complexities. However, it is worth noting that large portions of the paper deal exclusively with quadratic factors on $\F_p^n$ of complexity $(\ell,q)$ where $\ell,q\geq 1$.  Given such a factor $\calB$, we we can naturally conflate the space of labels $\F_p^{\ell}\times \F_p^q$ with elements of the group $\F_p^{\ell+q}$, by associating $(d,c)\in \F_p^{\ell}\times \F_p^q$ with $dc\in \F_p^{\ell+q}$ (here by $dc$, we mean $\sum_{i=1}^{\ell}d_ie_i+\sum_{j=1}^qc_je_{\ell+j}$ where $e_1,\ldots, e_{\ell+q}$, are the standard basis vectors in $\F_p^{\ell+q}$).   Using this association,  we let $B(dc)$ denote the atom $B(d;c)$, allowing us to write
$$
\At(\calB)=\{B(a): a\in \F_p^{\ell+q}\}.
$$
In this context, we will write $dc\in \F_p^{\ell+q}$ to implicitly mean $d\in \F_p^{\ell}$ and $c\in \F_p^q$.  
This minor change in notation will be useful in later sections, where we treat the space of labels as an auxiliary group.  We note these conventions will only be applied when both $\ell$ and $q$ are nonzero.

Given two linear factors $\calL,\calL'$ on $\F_p^n$, we say that $\calL'$ is a \emph{refinement of $\calL$}, denoted $\calL'\preceq \calL$, if the partition of $\F_p^n$ given by  $\At(\calL')$ refines the partition given by $\At(\calL)$.  Note that if $\calL\subseteq \calL'$, then $\calL'\preceq \calL$  always holds.  In this special case, we will write $\calL'\preceq_{syn}\calL$, and call $\calL'$ a \emph{syntactic} refinement of $\calL$.

Given two quadratic factors $\calB=(\calL,\calQ)$ and $\calB'=(\calL',\calQ')$ on $\F_p^n$, we say that $\calB'$ is a \emph{refinement of $\calB$}, denoted $\calB'\preceq \calB$, if the partition $\At(\calB')$ of $\F_p^n$ refines the partition given by $\At(\calB)$.  When $\calL\subseteq \calL'$ and $\calQ\subseteq \calQ'$, we say that $\calB'$ is a \emph{syntactic refinement of $\calB$}, and denote this by writing $\calB'\preceq_{syn}\calB$.

In working with quadratic factors, it will be important to consider their rank, as defined below.

\begin{definition}[Rank of a factor]\label{def:rankfac}
Suppose $q\geq 1$, and $\calB=(\calL,\calQ)$ is a quadratic factor on $\F_p^n$ with  $\calQ=\{M_1,\ldots, M_q\}$.  The \emph{rank of $\calB$} is 
$$
\min\{\rk(\lambda_1M_1+\ldots+\lambda_qM_q): \lambda_1,\ldots, \lambda_q\in \F_p\text{ are not all }0\}.
$$
By convention, the rank of a quadratic factor $\calB=(\emptyset, \calL)$ on $\F_p^n$ is $n$.
\end{definition}

 Crucially, when a factor has high rank, all of its atoms have approximately the same size, as the following lemma shows (for a proof, see \cite{Green.2007}).

\begin{lemma}[Size of an atom]\label{lem:sizeofatoms}
Suppose $\ell,q\geq 0$ are integers and $\calB=(\calL,\calQ)$ is a quadratic factor on $\mathbb{F}_p^n$ of complexity $(\ell,q)$ and rank $\rho$.  Then for any $B\in \At(\calB)$,  
\[|B|=(1+O(p^{\ell+q-\rho/2}))p^{n-\ell-q}.\]
\end{lemma}

In practice it will be useful to apply the following immediate corollary. 

\begin{corollary}\label{cor:sizeofatoms}
For any function $\tau:\mathbb{\R}^+\rightarrow \mathbb{R}^+$, there is a growth function $\rho=\rho(\tau)$ such that for all integers $\ell,q\geq 0$, if $\calB$ is a quadratic factor on $\mathbb{F}_p^n$ of complexity $(\ell,q)$ and rank at least $\rho(\ell+q)$, then for any $B\in \At(\calB)$, $|B|=(1\pm p^{-\tau(\ell+q)})p^{n-\ell-q}$.

Given $\e\in (0,1)$, we will write $\rho(\e)$ instead of $\rho(\tau)$ when $\tau$ is the constant function $\log_p(\e^{-1})$ (in which case the estimate above becomes $|B|=(1\pm\e )p^{n-\ell-q}$). 
\end{corollary}

 A standard fact about quadratic factors is that, given a factor $(\calL,\calQ)$, one can always find high rank refinement $(\calL',\calQ')$ whose complexity has not increased too much (see \cite[Lemma 3.11]{Green.2007}).  We will use the following version of this fact, which is marginally more sophisticated than its standard formulation.

\begin{lemma}[High rank refinement]\label{lem:rank}
For any growth function $\rho$, there exists a growth function $\tau=\tau(\rho)$ with the following property. 

Let $\ell\geq 0$ and $q\geq 1$ be integers.  Suppose $(\calL,\calQ)$ is a quadratic factor of complexity $(\ell,q)$, and $\calQ=\{Q_1,\ldots, Q_q\}$ is a fixed enumeration. Then there exists a quadratic factor $(\calL',\calQ')\preceq (\calL,\calQ)$ of complexity $(\ell',q')$ such that $\calL'\supseteq \calL$, $\calQ'\subseteq \calQ$, and 
\begin{enumerate}[label=\normalfont(\roman*)] 
\item $q'\leq q$ and $\ell' \leq \tau(\ell+q)$, and the rank of $(\calL',\calQ')$ is at least $\rho(\ell'+q')$;
\item if $(\calL,\calQ)$ has rank less than $\rho(\ell+q)$, then $\calQ'\subsetneq \calQ$;\label{start}
\item for each $i\in [q]$, if $(\calL,\calQ_i)$ has rank at least $\rho(\ell'+q')$, then $\calQ_i\subseteq \calQ'$, where $\calQ_i=\{M_1,\ldots, M_i\}$.  \label{keep}
\end{enumerate}
\end{lemma}
\begin{proof}
To prove this, one proceeds exactly as in the proof of \cite[Lemma 3.11]{Green.2007} with the added caveat that whenever one encounters a linear combination of matrices of low rank, one eliminates the matrix with the largest index in the fixed enumeration. 
\end{proof}

We remark that one typically only needs  conclusion (i) of Lemma \ref{lem:rank}, but we will need all three for one application, namely the proof of Proposition \ref{prop:qgsbad}.

It will be convenient throughout the paper to have the following more compact terminology for when a set is ``well approximated" by sets in a distinguished partition.  

\begin{definition}[(almost) $\e$-atomic partition]
Let $G$ be a group, let $A\subseteq G$, and let $\e\in (0,1)$. We say that a partition $\calP$ of $G$ is
\begin{itemize}
\item \emph{$\e$-atomic with respect to $A$} if  for all $X\in \calP$, $|A\cap X|/|X|\in [0,\e)\cup (1-\e,1]$;
\item \emph{almost $\e$-atomic with respect to $A$} if there is a set $\Omega\subseteq \calP$ such that $|\bigcup_{X\in \Omega}X|\leq \e|G|$ and for all $X\in \calP\setminus \Omega$, $|A\cap X|/|X|\in [0,\e)\cup (1-\e,1]$.
\end{itemize}
\end{definition}

We extend this definition to partitions generated by linear factors as follows.  
\begin{definition}[(almost) $\e$-atomic linear factor]
Let $A\subseteq \F_p^n$ and let $\e\in (0,1)$. We say that a linear factor $\calL$ on $\F_p^n$ is \emph{(almost) $\e$-atomic with respect to $A$} if the partition $\At(\calL)$ of $\F_p^n$ is (almost) $\e$-atomic with respect to $A$.
\end{definition}

In this terminology,  Theorem \ref{thm:vc} says that given a subset $A\subseteq \F_p^n$ with $\VC(A)<k$, we can find a linear factor of bounded complexity which is almost $\e$-atomic with respect to $A$. Similarly, Theorem \ref{thm:stablegps} says that given a $k$-stable subset $A\subseteq \F_p^n$, we can find a linear factor of bounded complexity which is $\e$-atomic with respect to $A$.  The main theorem of this paper can be analogously rephrased in terms of almost atomic \emph{quadratic factors}, which we now define.

\begin{definition}[almost $\e$-atomic quadratic factor (with linear error)]
Let $\e\in (0,1)$ and let $A\subseteq \F_p^n$. We say that a quadratic factor $\calB=(\calL,\calQ)$ on $\F_p^n$ is 
\begin{itemize}
\item \emph{almost $\e$-atomic with respect to $A$} if the partition $\At(\calB)$ is almost $\e$-atomic with respect to $A$;
\item \emph{almost $\e$-atomic with respect to $A$ with linear error} if the partition $\At(\calL)$ is almost $\e$-atomic with respect to $A$, and moreover, there is a set $\Sigma\subseteq \At(\calB)$ so that for all atoms $L(a)\cap Q(b)\in \At(\calB)$ with  $L(a)\notin \Sigma$,  
\[\frac{|A\cap L(a)\cap Q(b)|}{|L(a)\cap Q(b)|}\in [0,\e)\cup (1-\e,1].\]
\end{itemize}
\end{definition}

In this terminology, we see that Theorem \ref{thm:vc2finite} (now the main result of \cite{Terry.2021d}) says that given a set $A\subseteq \F_p^n$ with $\VC(A)<k$, there exists a high rank quadratic factor of bounded complexity which is almost $\e$-atomic with respect to $A$.  The main theorem of this paper, Theorem \ref{thm:fop}, says that if the set $A$ is moreover $k$-$\NFOP_2$, then we can ensure there exists a high rank quadratic factor of bounded complexity which is almost $\e$-atomic with respect to $A$ \emph{with linear error}.

\section{Basic properties of $\FOP_2$ and important examples}\label{sec:context}

Before giving a detailed outline of this section, we remark that its results are not strictly necessary for the proof of our main result, and the reader primarily interested in the latter may skip straight to Section \ref{sec:sufficient}. 

The goal of this section is threefold: first, we establish relationships between $\FOP_2$ and other notions of tameness, including stability, $\NIP$, and $\NIP_2$. Specifically, we will show that in the context of a distinguished subset $A\subseteq G$, for all $\ell\geq 2$,
\[\ell\mbox{-}\IP_2\Rightarrow \ell\mbox{-}\FOP_2\Rightarrow \ell\mbox{-}\IP\Rightarrow \ell\mbox{-}\OP,\]
so
\[\NOP\Rightarrow \NIP\Rightarrow \NFOP_2\Rightarrow \NIP_2.\]
We shall see that all the implications are strict.

Our second main goal is to introduce several important combinatorial facts about quadratic atoms. It is shown in \cite[Lemmas 3.9 and 3.10]{Terry.2021d} that a single atom in any high rank quadratic factor has unbounded $\VC$-dimension and $\VC_2$-dimension at most 1.\footnote{These facts were first proved in earlier versions of this paper \cite[Example 5.1 and Corollary 5.8]{Terry.2021av2}.} We show in Section \ref{subsec:atomsprops} that a single atom of any quadratic factor (not necessarily of high rank) has no 2-$\FOP_2$. We will use this to prove Theorem \ref{thm:unionsatoms}, which says that for all $k\geq 1$, there is an $\ell\geq 1$ such that a union of $k$ quadratic atoms (even from different factors) has no $\ell$-$\FOP_2$, and consequently $\VC_2$-dimension at most $\ell$.

Finally, we explore the two examples defined in the introduction, namely the linear and quadratic Green-Sanders examples (Definitions \ref{def:GS} and \ref{def:QGS} respectively). We show, among other things, that these provide examples of NIP (respectively $\NIP_2$) sets which do not satisfy the conclusions of Theorem \ref{thm:stablegps} (respectively Theorem \ref{thm:fop}). 

The linear Green-Sanders example, $A_\GS(p,n)$, first appeared in this form in a paper of Green and Sanders  \cite{Green.2015qy4}, where it was shown that $A_\GS(p,n)$ has a large non-trivial Fourier coefficient with respect to the trivial coset of any positive-dimensional subspace of $\F_p^n$.  We will  show that every linear Green-Sanders is an example of a subset of a group with a large order property and very small VC-dimension (another such example in groups of the form $\Z/n\Z$ is given in \cite{Alon.2018is}).  We give an easy proof of Proposition \ref{prop:gsvc}, which says that any linear factor on $\F_p^n$ has an atom on which the density of $A_\GS(p,n)$ is bounded away from 0 and 1. Thus, we can think of $A_\GS(p,n)$ as a natural example of an NIP set for which the conclusion of Theorem  \ref{thm:stablegps} fails.  Turning to quadratic decompositions of $A_{\GS}(p,n)$, we prove Proposition \ref{prop:gslinear}, which shows that, in fact, any high rank quadratic factor of bounded complexity has an atom on which the density of $A_\GS(p,n)$ is bounded away from 0 and 1. On the other hand, we prove Proposition \ref{prop:gslinear1}, which states that  there exists a high rank quadratic factor of bounded complexity with the property that $A_{\GS}(p,n)$ has density equal to 0 or 1 on all atoms avoiding a single linear atom (namely the one labeled by 0).

As part of the proof of Theorem \ref{thm:unionsatoms}, we will define another ternary analogue of the order property, which we call the \emph{hyperplane order property} (see Definition \ref{def:hop}) and which is weaker than the functional order property. We will also show, in Appendix \ref{app:no4hop}, that $A_\GS(p,n)$ omits all instances of the hyperplane order property longer than $4$. Combined with Proposition \ref{prop:gslinear}, this refutes a naturally arising conjecture about which combinatorial restrictions on subsets of $\F_p^n$ ensure quadratic decompositions with no irregular atoms, and also plays an important role in \cite[Section 7]{Terry.2021b}.

The second example is a quadratic analogue of the linear Green-Sanders example, denoted by $A_\QGS(p, n)$ and defined here for the first time. It was proved by Gladkova (Theorem \ref{thm:Gladkova}) that $A_\QGS(p, n)$ has bounded $\VC_2$-dimension.\footnote{This was originally conjectured in earlier versions of this paper \cite{Terry.2021av2}, and subsequently resolved by Gladkova.} In analogy to Proposition \ref{prop:gsvc} concerning the linear Green-Sanders example, we prove Theorem \ref{thm:ex1}, which shows that there does not exist a quadratic factor with respect to which $A_\QGS(p, n)$ satisfies the conclusions of Theorem \ref{thm:fop}. It therefore follows from Theorem  \ref{thm:fop} that $A_\QGS(p, n)$ has an unbounded instance of the functional order property, a fact which we also prove directly later in the paper (see Corollary \ref{cor:qgscor}). For comparison, we showed in \cite{Terry.2021b} that there exist examples of hypergraphs with unbounded instances of $\FOP_2$ and bounded $\VC_2$-dimension (see \cite[Example D.5]{Terry.2021b}). The quadratic Green-Sanders example $A_\QGS(p, n)$ is the first such example in the arithmetic setting.

\subsection{Relationships between $\FOP_2$ and other notions of tameness}\label{subsec:fopbasic}
In this section, we discuss in some detail the relationships between various notions of tameness appearing in this paper. For all the general implications presented, the proofs can be carried out with minor changes in the setting of formulas in first order structures, as we believe will be clear to the model theorist.

We begin by showing that sets with long functional order properties have large $\VC$-dimension (see Lemma \ref{lem:vcfop} in introduction).

\begin{lemma}[$\NIP \Rightarrow \NFOP_2$]\label{lem:fopvc}
For all $\ell\geq 2$, if $A\subseteq G$ has $\ell$-$\FOP_2$, then $A$ has $\ell$-$\IP$.
\end{lemma}
\begin{proof}
Suppose $A\subseteq G$ has $\ell$-$\FOP_2$.  Setting  $I=\{f:[\ell]^2\rightarrow [\ell]\}$, let 
$$
\{a_i^f, b_i,c_i: i\in [\ell], f\in I\}\subseteq G 
$$
witness that $A\subseteq G$ has $\ell$-$\FOP_2$.  For each $j\in [\ell]$, let $d_j=b_j\cdot c_2$, set $X=\{d_j: j\in [\ell]\}$, and for each $Y\subseteq [\ell]$, let $x_Y=a_2^{f_Y}$ where $f_Y:[\ell]\rightarrow [\ell]$ is defined by $f_Y(i,j)=2$ if $i=2$ and $j\in Y$ and $f_Y(i,j)=1$ otherwise.  Then $x_Y\cdot d_j\in A$ if and only if $2\leq f_Y(2,j)$ if and only if $j\in Y$.  
\end{proof}

We observe that this implication is strict: in the next subsection (see Corollary \ref{cor:quadatomfop2}), we shall show that a quadratic atom is $2$-$\NFOP_2$, while \cite[Lemma 3.9]{Terry.2021d} shows that an atom of a high-rank quadratic factor has unbounded $\VC$-dimension (in the sense that the length of the independence property grows as the dimension $n$ of the ambient vector space tends to infinity).\footnote{This fact was first proved in \cite[Example 5.1]{Terry.2021av2}.} 

We next verify that a set with large $\VC_2$-dimension has a long functional order property (see Lemma \ref{lem:nip2fop} in the introduction).

\begin{lemma}[$\NFOP_2 \Rightarrow \NIP_2$]\label{lem:vc2fop}
For all $\ell\geq 1$, if $A\subseteq G$ has $\ell$-$\IP_2$, then $A$ has $\ell$-$\FOP_2$. 
\end{lemma}
\begin{proof}
Suppose $\VC_2(A)\geq \ell$.  Then there are $c_1,\ldots, c_{\ell}, b_1,\ldots, b_{\ell}\in G$ and for each $X\subseteq [\ell]\times [\ell]$, an element $a_X\in G$ so that $a_X\cdot c_j\cdot b_k \in A$ if and only if $(j,k)\in X$.  Given $f:[\ell]\rightarrow [\ell]$, and $i\in [\ell]$, let $a_i^f$ be $a_{X_i^f}$ where $X_i^f=\{(j,k)\in [\ell]^2: k\leq f(i,j)\}$.  Then $ a_i^f\cdot c_j \cdot b_k\in A$ if and only if $k\leq f(i,j)$.  Thus $A$ has $\ell$-$\FOP_2$.
\end{proof}

Again, this implication is strict, with the quadratic Green-Sanders example discussed in Section \ref{subsec:qgs} being an example. 

It was proved by Chernikov, Palac\'{i}n and Takeuchi \cite{Chernikov.2019} that in any first order structure, the set of formulas of finite $\VC_2$-dimension is closed under finite boolean combinations.  In \cite[Corollary D.4]{Terry.2021b}, we showed that the analogous fact holds for $\NFOP_2$ formulas.  In the setting of subsets of groups, the result of Chernikov, Palac\'{i}n and Takeuchi implies that for all $k,\ell\geq 1$, there is $m\geq 1$ such that for any group $G$ and subsets $A_1,\ldots, A_k\subseteq G$, each of $\VC_2$-dimension at most $\ell$, any boolean combination of the sets $A_1,\ldots, A_k$ has $\VC_2$-dimension at most $m$.  Similarly, \cite[Corollary D.5]{Terry.2021b} implies the following.

\begin{theorem}[$\NFOP_2$ is closed under boolean combinations]\label{thm:bool}
For all $k,\ell\geq 1$, there is $m\geq 1$ such that for any subsets $A_1,\ldots, A_k\subseteq G$, each of which is $\ell$-$\NFOP_2$, any boolean combination of the sets $A_1,\ldots, A_k$ is $m$-$\NFOP_2$. 
\end{theorem}

Our proof in \cite{Terry.2021b} that $\NFOP_2$ formulas are closed under boolean combinations is finitary but extremely indirect, passing through hypergraph regularity lemmas.  An equally indirect  model theoretic proof that $\NFOP_k$ formulas are closed under boolean combinations (for any $k\geq 1$) later appeared later in \cite{Abd.2025}. It would be of interest to produce a direct combinatorial argument of this fact. In Appendix \ref{app:generalprops}, we provide such an argument for the special case of Theorem \ref{thm:bool} where $\ell=2$ and the boolean combination is simply a union (see Section \ref{subsec:atomsprops} for further discussion). 
We remark that it is very easy to show directly that $\NFOP_2$ formulas are closed under negations (see \cite{Terry.2021b} and earlier versions of the present paper \cite{Terry.2021av2}), so the difficulty really lies in showing closure under conjunctions. 

We next discuss the sense in which $\FOP_2$ can be viewed as a ternary analogue of the order property. This will lead us to investigating the relationship between $\FOP_2$ and a closely related definition which appears in unpublished work of Takeuchi \cite{Takeuchi.2018}.

First, we show that the order property is equivalent to a natural definition of ``$\FOP_1$".

\begin{lemma}\label{lem:opfop1}
Let $k\geq 1$, let $G$ be a group and let $A\subseteq G$.  The following are equivalent.
\begin{enumerate}[label=\normalfont(\roman*)]
\item $A$ has the $k$-order property.
\item There are $b_1,\ldots, b_k\in G$ and  for all $f:[k]\rightarrow [k]$, $a_1^f,\ldots, a_k^f$ such that $a_i^f\cdot b_j\in A$ if and only if $j\leq f(i)$.
\end{enumerate}
\end{lemma}
\begin{proof} Suppose first that there are $b_1,\ldots, b_k\in G$ and  for all $f:[k]\rightarrow [k]$, $a_1^f,\ldots, a_k^f$ such that $a_i^f\cdot b_j\in A$ if and only if $j\leq f(i)$.  Let $f_1(x)=x$.  Then $a_i^{f_1}\cdot b_j\in A$ if and only if $j\leq i $, so $A$ has the $k$-order property. Conversely, suppose $a_1,\ldots, a_k$, $b_1,\ldots, b_k\in G$ are such that $a_i\cdot b_j\in A$ if and only if $j\leq i$.  For each  $f:[k]\rightarrow [k]$, and $i\in [k]$, set $a^f_i=a_{f(i)}$.  Then given $1\leq i,j\leq k$, we have that $a_i^f\cdot b_j\in A$ holds if and only if $j\leq f(i)$.
\end{proof}

We next define a notion which turns out to be essentially equivalent to the functional order property.\footnote{Definition \ref{def:fopstar} is taken as the definition of $\FOP_2$ in the subsequent paper \cite{Abd.2025}.}

\begin{definition}[$\FOP_2^*$]\label{def:fopstar}
Given $m\geq 1$, a group $G$ and a subset $A\subseteq G$, we say that $A$ \emph{has $m$-$\FOP_2^*$} if there are $b_1,\ldots,b_m,c_1,\ldots,c_m\in G$ and for every $f:[m]\rightarrow [m]$, there is $a_f\in G$ such that  $a_f\cdot b_i\cdot c_j\in A$ if and only if $j\leq f(i)$.
\end{definition}

The properties $\FOP_2$ and $\FOP_2^*$ are qualitatively equivalent, as the following proposition shows.

\begin{proposition}[$\NFOP_2\Leftrightarrow \mathrm{NFOP}_2^*$]\label{prop:fopstar}
Let $m\geq 1$, let $G$ be a group and let $A\subseteq G$. 

Then the following statements hold.
\begin{enumerate}[label=\normalfont(\roman*)]
\item If $A$ has $m^m$-$\FOP_2$, then $A$ has $m$-$\FOP_2^*$.
\item If $A$ has $m$-$\FOP^*_2$, then $A$ has $m$-$\FOP_2$.
\end{enumerate}
\end{proposition}

\begin{proof}
(i)$\Rightarrow$(ii): Suppose that $A$ has $m^m$-$\FOP_2$, and let $s=m^m$.  Enumerate $[m]^{[m]}=\{f_1,\ldots, f_s\}$. Define $f\colon [m]^2\to [m]$ such that $f(i,j)=f_i(j)$ if $j\leq m$ and $f(i,j)=m$ if $j>m$. By assumption, there are $b_1,\ldots,b_m,c_1,\ldots,c_m\in G$ and for each $f:[m]^2\rightarrow [m]$, $a_1^f,\ldots, a_m^f\in G$ so that $a_i^f\cdot b_j\cdot c_k\in A$ if and only if $k\leq f(i,j)$. Then for all $j\leq \ell$, we have $a^f_i\cdot b_j\cdot c_k\in A$ holds if and only if $k\leq f_i(j)$. For each $f_u\in [m]^{[m]}$, set $a_{f_u}=a^f_u$. Given given $i,j\leq m$, we have $a_{f_u}\cdot b_i\cdot c_j\in A$  if and only if $a^f_u\cdot b_i\cdot c_j$ if and only if $j\leq f(u,i)$ if and only if $j\leq f_u(i)$. This shows that $A$ has $m$-$\FOP^*_2$.

(ii)$\Rightarrow$(i): Suppose that $A$ has $m$-$\FOP_2^*$.  Then there are $b_1,\ldots,b_m,c_1,\ldots,c_m\in G$, and for each $f:[m]\rightarrow [m]$ and $a_f\in G$ are such that $a_f\cdot b_j\cdot c_k\in A$ if and only if $k\leq f(j)$.  Fix $f\colon[m]^2\to[m]$, and for $i\leq m$, define $f_i\colon[m]\to[m]$ so that $f_i(j)=f(i,j)$. Set $a^f_i=a_{f_i}$. Then $a^f_i\cdot b_j\cdot c_k$ holds if and only if $k\leq f_i(j)$, if and only if $k\leq f(i,j)$. This shows that $A$ has $m$-$\FOP_2$.
\end{proof}

This reformulation makes clear how $\FOP_2$ is related to a notion called ``$\OP_2$" defined by Takeuchi in \cite{Takeuchi.2018}.  In particular, Takeuchi's $\OP_2$ corresponds to a variant of $\FOP_2^*$ in which one only insists on the existence of $a_f$ for \emph{increasing} functions $f:[n]\rightarrow [n]$.  We make a couple more definitions to clarify this relationship. 

\begin{definition}[$\IFOP_2$ and $\IFOP_2^*$]\label{def:ifop}
Let $m\geq 1$, let $G$ be a group and let $A\subseteq G$. 
\begin{enumerate}
\item We say that $A$ \emph{has the $m$-increasing functional order property} ($m$-$\text{IFOP}_2$) if there are $b_1,\ldots,b_m,c_1,\ldots,c_m\in G$ and for every $f:[m]^2\rightarrow [m]$ increasing in the second coordinate, there are  $a_1^f,\ldots, a_m^f\in G$ such that  $a_i^f\cdot b_j\cdot c_k\in A$ if and only if $k\leq f(i,j)$.
\item We say that $A$ \emph{has $m$-$\mathrm{IFOP}^*_2$} if there are $b_1,\ldots,b_m$, $c_1,\ldots,c_m\in G$ and for every increasing $f:[m]\rightarrow [m]$, there is $a_f\in G$ such that  $a_f\cdot b_i\cdot c_j\in A$ if and only if $j\leq f(i)$.
\end{enumerate}
\end{definition}

It is immediate from the definition that $m$-$\FOP_2$ implies $m$-$\IFOP_2$. Using an argument similar to the one used to prove Proposition \ref{prop:fopstar}, it is not difficult to see that for any $m$, $m$-$\text{IFOP}^*_2$ implies $m$-$\IFOP_2$ and $m^m$-$\IFOP_2$ implies $m$-$\text{IFOP}_2^*$, and consequently $\IFOP_2$ and $\IFOP_2^*$ are qualitatively equivalent, i.e. $\mathrm{NIFOP}_2 \Leftrightarrow \mathrm{NIFOP}_2^*$. 

Takeuchi's $\OP_2$ is exactly $\IFOP^*_2$ as defined above. Takeuchi showed that an unbounded instance of $\IFOP^*_2$ implies an unbounded instance of $\IP$. Combining all these facts yields the general implications
\[\NIP\Rightarrow \mathrm{NIFOP}_2\; (\mbox{equivalently, } \mathrm{NIFOP}^*_2) \Rightarrow \NFOP_2 \; (\mbox{equivalently, }  \NFOP_2^*)\Rightarrow \NIP_2.\]

We already discussed earlier in this section that the third implication is strict. This answers a question about the relationship between $\IP_2$ and $\OP_2$ asked by Takeuchi in \cite{Takeuchi.2018}. We shall see in Section \ref{subsec:atomsprops} that a single quadratic atom of a high-rank factor (which has unbounded $\VC$-dimension but $\VC_2$-dimension equal to 1) demonstrates that the first implication is strict (see remark following Corollary \ref{cor:quadatomfop2}). It is unclear to the authors at the time of writing whether the second implication is strict.

\begin{question}\label{que:op2vsfop2}
Is the implication $\mathrm{NIFOP}_2 \Rightarrow \mathrm{NFOP}_2$ strict? In particular, does there exist a subset of a group  with a large instance of $\mathrm{IFOP}_2$ but no large instances of $\FOP_2$?
\end{question}

\subsection{Properties of quadratic atoms}\label{subsec:atomsprops}

We begin with a discussion of the combinatorial properties of linear atoms. Given a subset of $A$ of a group $G$, we say that $A$ has the \emph{square auto-completion property} ($\mathrm{SAP}$) if  for all $x_1,x_2,y_1,y_2\in G$ satisfying $x_i\cdot y_j\in A$ for all $(i,j)\in [2]^2\setminus \{(2,1)\}$, one also has $x_2\cdot y_1\in A$. It is not difficult to check that $A\subseteq G$ has  $\mathrm{SAP}$ if and only if $A$ is a coset of a subgroup, or equivalently, if and only if $A$ is an atom of a purely linear factor.  It is also easy to see, after some reindexing, that having square auto-completion is equivalent to having no $2$-$\OP$. Thus a single atom of a linear factor has no $2$-OP, and thus also VC-dimension at most $1$.  On the other hand, an atom of a \emph{quadratic} factor may have large VC-dimension (see \cite[Lemma 3.9]{Terry.2021d}).\footnote{Example \ref{exa:quadvc} was originally proved as \cite[Example 5.1]{Terry.2021av2} in earlier versions of this paper, but has now been moved to \cite{Terry.2021d} with a new proof.}

\begin{example}[Quadratic atoms of high-rank factors have large $\VC$-dimension]\label{exa:quadvc}
For all $m\geq 1$, there exists a growth function $\rho_0=\rho_0(m)$ such that  for all growth functions $\rho\geq \rho_0$ the following holds.

Let $n\geq \ell\geq 1$ be integers and let $(\calL,\calQ)$ be a quadratic factor on $\F_p^n$ of complexity $(\ell,q)$ and rank at least $\rho(\ell+q)$. Then any atom of $(\calL,\calQ)$ has $\VC$-dimension at least $m$. 
\end{example}

Despite having large $\VC$-dimension, quadratic atoms exhibit rather strong restrictions on the kinds of ternary configurations they can encode.  Given a set $A\subseteq G$, we say that $A$ has the \emph{cube auto-completion property} ($\mathrm{CAP}$) if for all $x_1,x_2,y_1,y_2,z_1,z_2\in G$ satisfying $x_i\cdot y_j \cdot z_k\in A$ for each $(i,j,k)\in [2]^3\setminus \{(2,1,1)\}$, one also has $x_2 \cdot y_1 \cdot z_1\in A$.  We will show below that if $A\subseteq \F_p^n$ is an atom of a quadratic factor, then it has $\mathrm{CAP}$.  We begin by showing that this is the case for the simplest possible quadratic atom. 

\begin{example}\label{exa:quadtame}
Given $a\in \F_p$, the set $Q=\{x\in \F_p^n:x^T x=a\}\subseteq \F_p^n$ has $\mathrm{CAP}$.
\end{example}

\begin{proof}
For convenience, given $x\in \F_p^n$, we write $x^2$ for $x^Tx$.  Suppose that we have $x_1,x_2$, $y_1,y_2$ and $z_1,z_2$ such that  $x_1+y_1+z_2, x_2+y_1+z_2$, $x_1+y_2+z_1$, $x_2+y_2+z_1$, $x_1+y_2+z_2$,   $x_2+y_2+z_2$ and $x_1+y_1+z_1$ are all in $Q$. Since
\begin{align*}
(x_2+y_1+z_1)^2=&(x_1+y_1+z_1)^2 -(x_1+y_1+z_2)^2-(x_1+y_2+z_1)^2\\&+(x_2+y_2+z_1)^2+(x_1+y_2+z_2)^2+(x_2+y_1+z_2)^2 - (x_2+y_2+z_2)^2,
\end{align*}
we have that $(x_2+y_1+z_1)^T (x_2+y_1+z_1)=a$, and therefore $x_2+y_1+z_1\in Q$.  
\end{proof}

It is straightforward to extend this to any atom of a \emph{purely quadratic} factor.

\begin{remark}\label{rem:highrankat} 
Suppose $M_1,\ldots, M_q$ are symmetric $n\times n$ matrices with entries in $\F_p$, and $a_1,\ldots, a_q\in \F_p$.  Then the set 
$$
Q'=\{x\in \F_p^n:x^T M_1 x=a_1,\dots, x^T M_q x=a_q\}\subseteq \F_p^n
$$ 
has $\mathrm{CAP}$. This is simply a consequence of the fact that the quadratic identity used in Example \ref{exa:quadtame} above applies to each quadratic constraint individually.
\end{remark}

To deal with atoms of an arbitrary quadratic factor as defined in Definition \ref{def:quadfac} (which may include a linear part), the following lemma is helpful.

\begin{lemma}\label{lem:quadint}
If $A,B\subseteq G$ both have $\mathrm{CAP}$, then so does $A\cap B$. 
\end{lemma}

\begin{proof}
Suppose to the contrary that there are $x_1,x_2$, $y_1,y_2$, $z_1,z_2$ such that $x_i\cdot y_j\cdot z_k\in A\cap B$ if and only if $(i,j,k)\neq (2,1,1)$. Since $x_2\cdot y_1\cdot z_1\notin A\cap B$, we either have $x_2\cdot y_1\cdot z_1\notin A$ or $x_2\cdot y_1\cdot z_1\notin B$.  Without loss of generality, suppose $x_2\cdot y_1\cdot z_1\notin A$.  But now $x_i\cdot y_j\cdot z_k\in A$ if and only if $(i,j,k)\neq (2,1,1)$, a contradiction.  
\end{proof}

It is not difficult to show that SAP implies CAP.  Indeed, suppose a set $A$ fails CAP.  Then there are  $x_1,x_2$, $y_1,y_2$, $z_1,z_2$ such that $x_i\cdot y_j\cdot z_k\in A$ if and only if $(i,j,k)\neq (2,1,1)$.  Then letting $a_1=y_1+z_1$ and $a_2=y_1+z_2$, we have $x_i+a_j\in A$ if and only if $(i,j)\neq (2,1)$, i.e. $A$ fails SAP.  A consequence of this is that any atom of a linear factor has CAP, since, as we saw above, any atom of a linear factor has SAP.  We can now deduce that every atom of a quadratic factor has CAP.

 \begin{lemma}\label{lem:quadatomcap}
 For any quadratic factor $\calB=(\calL,\calQ)$ on $\F_p^n$, every atom of $\calB$  has $\mathrm{CAP}$.
 \end{lemma}
 \begin{proof}
 An atom of $\calB$ has the form $B(a,b)=L(a)\cap Q(b)$.  By the remarks at the start of this section, $L(a)$ has SAP, and thus CAP by the remarks following the proof of Lemma \ref{lem:quadint}.  By Remark \ref{rem:highrankat}, $Q(b)$ has CAP.  Consequently, $L(a)\cap Q(b)$ has CAP by Lemma \ref{lem:quadint}.
 \end{proof}
 
 The following corollary is immediate from the observation that a subset of $\F_p^n$ with CAP has no $2$-$\FOP_2$, by definition.
 
 \begin{corollary}[Quadratic atoms are $2$-$\NFOP_2$]\label{cor:quadatomfop2}
 For any quadratic factor $\calB=(\calL,\calQ)$ on $\F_p^n$, every atom of $\calB$  is $2$-$\NFOP_2$.
 \end{corollary}
 
In fact, any quadratic atom will be $2$-$\mathrm{NIFOP}_2$ (see Definition \ref{def:ifop}) since the function $f:[m]^2\ra [m]$ that takes $(i,j)$ to $\max\{i+j,m\}$ is increasing in the second coordinate. 
 
One might be tempted to think that a converse of Lemma \ref{lem:quadatomcap} should hold. However, atoms of quadratic factors as defined in Definition \ref{def:quadfac} are not the only sets that satisfy CAP. For instance, atoms of so-called \emph{generalised quadratic factors} \cite[Definition 4.3]{Green.2012}, which are of the form 
\[\{x\in \F_p^n: x^TMx +r^Tx=0\}\]
for some $n\times n$ symmetric matrix $M$ and an element $r\in \F_p^n$, are also easily seen to have CAP.

We conjecture that these are, in fact, the only examples.

\begin{conjecture}
Any subset $A\subseteq \F_p^n$ that satisfies $\mathrm{CAP}$ is an atom of a generalized quadratic factor.
\end{conjecture}

We end this section by discussing unions of quadratic atoms. By Corollary \ref{cor:quadatomfop2} and Theorem \ref{thm:bool}, it follows immediately that bounded unions of quadratic atoms are $\NFOP_2$ in the following sense.

\begin{theorem}[Unions of quadratic atoms are $\NFOP_2$]\label{thm:unionsatoms}
For all $k\geq 1$, there is $\ell\geq 2$ such that a union of $k$ quadratic atoms is $\ell$-$\NFOP_2$. 
\end{theorem}

As previously discussed, the proof of Theorem \ref{thm:bool} given in \cite{Terry.2021b} is extremely indirect. We include a more direct, combinatorial proof of a stronger version of Theorem \ref{thm:unionsatoms} in Appendix \ref{app:boolean}

To formulate this stronger version, we introduce a further ternary tameness property below. Extremely natural in the context of quadratic atoms, it played an important role in the authors' early attempts to frame the results in this paper, as well as in the companion paper \cite[Chapter 7]{Terry.2021b}. While it does not feature in the main results of this paper, it is combinatorially of independent interest.  

 \begin{definition}[Hyperplane order property ($\HOP_2$)]\label{def:hop}
Given $\ell\geq 1$, a group $G$ and a subset $A\subseteq G$, we say that $A$ \emph{has \emph{$\ell$-$\HOP_2$}} if there are  
$$
\{a_i: i\in [\ell]\}\cup \{b_i: i\in [\ell]\}\cup \{c_i: i\in [\ell]\}\subseteq G
$$
 such that $a_u \cdot b_v\cdot c_w\in A$ if and only if $u<v+w$.  

We say that $A$ \emph{is $\ell$-$\NHOP_2$} if it does not have $\ell$-$\HOP_2$.
\end{definition}

We observe that a set has CAP if and only if it has no $2$-$\HOP_2$, and thus by Lemma \ref{lem:quadatomcap}, all quadratic atoms are $2$-$\NHOP_2$.  It also follows straight from the definitions that a set which has $\ell$-$\FOP_2$ must also have $\ell$-$\HOP_2$.

\begin{lemma}[$\NHOP_2\Rightarrow \NFOP_2$]\label{lem:fophop}
For all $\ell\geq 2$, if $A\subseteq G$ has $\ell$-$\FOP_2$, then it has $\ell$-$\HOP_2$.
\end{lemma}

This implication is strict: there exist examples with unbounded instances of $\HOP_2$ but no unbounded instances $\FOP_2$.  In fact, there are sets with very small $\VC$-dimension and unbounded instances of $\HOP_2$ (which is a stronger statement by Lemma \ref{lem:fopvc}).

\begin{example}\label{exa:vchop}
For every $\ell\geq 2$, there is an abelian group $G$ and a subset $A\subseteq G$ such that $\VC(A)\leq 4$ and such that $A$ has $\ell$-$\HOP_2$. 
\end{example}
\begin{proof}
Fix $\ell\geq 2$.  Choose a prime $p>\ell^4$, and consider $G=\mathbb{Z}/p\mathbb{Z}$ and the subset $A=\{(p-1)/2,\ldots, p\}$.  It is shown in \cite{Alon.2018is}, for example, that this set has $\VC$-dimension at most $4$.  To see that $A$ has $\ell$-$\HOP_2$, let $x_i=(p-1)/2-i-1$ and $y_i=z_i=i$ for each $i\in [\ell]$.  Then $x_i+y_j+z_k\in A$ if and only if $(p-1)/2-i-1+j+k\in \{(p-1)/2,\ldots, p\}$.  Since $p>\ell^4$, this holds if and only if $(p-1)/2-i-1+j+k\geq  (p-1)/2$, i.e. if and only if $j+k>i$.
\end{proof}

Returning to unions of atoms, Theorem \ref{thm:unionsatoms} is implied by Corollary \ref{cor:quadatomfop2}, Lemma \ref{lem:fophop}, and the following.

\begin{proposition}\label{prop:twotameunion}
For all $k\geq 1$ there exists $\ell\geq 2$ such that whenever $A_1,\ldots, A_k$ are $2$-$\NHOP_2$, then $A_1\cup \ldots \cup A_k$ is $\ell$-$\NHOP_2$. 
\end{proposition}

The proof of Proposition \ref{prop:twotameunion} is finitary and explicit, and appears in Appendix \ref{app:boolean}.  It proceeds via a Ramsey theoretic argument employing the multi-dimensional version of Szemer\'{e}di's Theorem.

It is natural to ask whether $\NHOP_2$ is closed under boolean combinations, as we have seen $\NIP_2$ is by \cite{Chernikov.2019} and $\NFOP_2$ is by \cite{Terry.2021b}. It is easy to show that $\NHOP_2$ is closed under taking complements in the following sense.

\begin{lemma}[$\HOP_2$ is closed under complements]\label{lem:negations}
For all $\ell\geq 2$, if $A\subseteq G$ has $\ell$-$\HOP_2$, then $\neg A$ has $\lfloor \ell /2\rfloor$-$\HOP_2$. 
\end{lemma}
\begin{proof}
Suppose that $A$ has $\ell$-$\HOP_2$.  Then there are $x_1, \dots,x_\ell$, $y_1,\dots,y_\ell$, $z_1,\dots,z_\ell$ in $G$ such that $x_i\cdot y_j\cdot z_k\in A$ holds if and only if $i< j+k$.  For each $i,j,k\in [1,\lfloor\ell/2\rfloor]$, let $a_i=x_{\lfloor\ell/2\rfloor+i}$, $b_j=y_{\lfloor\ell/2\rfloor-j+1}$, and $c_k=z_k$.  Then $a_i\cdot b_j\cdot c_k\in \neg A$ holds if and only if 
\[\lfloor\ell/2\rfloor+i\geq \lfloor\ell/2\rfloor-j+1+k\Leftrightarrow i+j\geq k+1 \Leftrightarrow i+j>k.\]
This shows that $\neg A$ has $\lfloor \ell/2\rfloor$-$\HOP_2$.
\end{proof}

On the other hand, it remains an open question whether the intersection of two $\NHOP_2$ sets is $\NHOP_2$.  

\begin{question}\label{ques:inthop}
Does there exist a function $f:\mathbb{N}\rightarrow \mathbb{N}$ such that whenever $A,B\subseteq G$ are both $\ell$-$\NHOP_2$, then $A\cap B$ is $f(\ell)$-$\NHOP_2$?
\end{question}

We believe this to be an interesting question from a purely Ramsey theoretic point of view.  

Just like $\FOP_2$, $\HOP_2$ can be viewed as a natural ternary analogue of the order property. Indeed, with some simple reindexing it is possible to show that $A$ has $\ell$-$\OP$ if and only if there are $a_1,\ldots, a_\ell,b_1,\ldots, b_\ell$ such that $a_i+b_j\in A$ if and only if $i+j\leq \ell$.  On the other hand, one can show that $A$ has $\ell$-$\HOP_2$ if and only if there are $a_1,\ldots, a_\ell,b_1,\ldots, b_\ell,c_1,\ldots, c_\ell$ such that $a_i+b_j+c_k\in A$ if and only if $i+j+k\leq \ell$. 

In light of this, the reader may wonder whether there exists a yet stronger quadratic regularity lemma for sets without the hyperplane order property.  For instance, the reader may wonder if $\NHOP_2$ subsets of $\F_p^n$ admit quadratic decompositions with \emph{no irregular atoms}.  In the next section we shall show that this is false (see Proposition \ref{prop:GSnoHOP}). 

\subsection{The linear Green-Sanders example}\label{subsec:gs}

In this section we discuss several important properties of the linear Green-Sanders example, $A_{\GS}(p,n)$ of Definition \ref{def:GS}.

We begin by recalling its definition.  Given $i\in [n]$, we set 
$$
H_i:=\{x\in \mathbb{F}_p^n: x_1=\ldots=x_i=0\},
$$
 and  let $e_i$ denote the $i$th standard basis vector in $\F_p^n$.  

\begin{definition}[Green-Sanders example, recalled]\label{def:GSrecall}
For $p>2$ and $n\geq 1$, define the \emph{linear Green-Sanders $p$-example of dimension $n$} to be the subset $A_\GS(p,n)\subseteq\F_p^n$ defined by
$$
A_{\GS}(p,n)=\bigcup_{i=1}^n (H_i+e_i).
$$
\end{definition}

Noting that the $H_i+e_i$ are disjoint for distinct values of $i\in [n]$, it is easy to see that $A_\GS(p,n)$ has density $(1-p^{-n})/(p-1)=1/(p-1)+O(p^{-n})$. It is also not difficult to see that $A_{\GS}(p,n)$ has an order property of length increasing with $n$ (see \cite[Example 4]{Terry.2019}). On the other hand, we shall see below, in Theorem \ref{thm:gs3nip} (which appeared as Proposition \ref{prop:gs3vc} in the introduction), that the sets $A_{\GS}(p,n)$ have bounded $\VC$-dimension.  We provide a proof in the case $p=3$ only. We have chosen to omit the proof for $p>3$, which appeared in earlier versions of this paper (see \cite[Appendix A.2]{Terry.2021av2}), in this updated manuscript as a stronger and more elegant argument for these cases has since been found by Gladkova \cite{Gladkova.2024}.  

Throughout the discussion that follows, given $x\in \F_p^n$, let
$$
f(x)=\max\Big( \{1\}\cup \{2\leq j\leq n: \text{ for all }1\leq i< j, x_i=0\}\Big)
$$
denote the first non-zero coordinate of $x$. 

\begin{theorem}[Green-Sanders example has bounded $\VC$-dimension]\label{thm:gs3nip}
For all $n\geq 1$, $A_{\GS}(3,n)$ has $\VC$-dimension $3$.
\end{theorem}
\begin{proof}
We begin by showing that $A=A_{\GS}(3,n)$ has $\VC$-dimension less than $4$.  Suppose towards a contradiction that there exists a set $X=\{a^1,a^2,a^3,a^4\}\subseteq \F_p^n$ of size $4$ which is shattered by $A$.  So for all $Y\subseteq X$, there is some $e_Y\in \F_p^n$ so that $(A-e_Y)\cap X=Y$.  Let 
$$
m=\max\Big( \{0\}\cup \{i\in [n]: \text{ for all }1\leq j\leq i, a^1_j=a^2_j=a^3_j=a^4_j\}\Big)
$$
For each $i\in [4]$, set $b^i=a^i-a^1$.  Then $X'=\{b^1,b^2,b^3,b^4\}$ is still shattered by $A$, and 
$$m+1=\min\{f(b^1), f(b^2),f(b^3),f(b^4)\}<\max\{f(b^1), f(b^2),f(b^3),f(b^4)\}=f(b^1)=n.$$  Since the $b^i$ are pairwise distinct, we have that $f(b^2),f(b^3),f(b^4)<f(b^1)$.  

Suppose first that $f(b^2)<f(b^3)<f(b^1)$.   Let $u\in \{1,2\}$ be such that 
$$
b^2=ue_{f(b^2)}+\sum_{i=f(b^2)+1}^nb^2_ie_i.
$$
Suppose first that $u=1$.  By assumption, there exists some $y\in \F_p^n$ such that $(A-y)\cap X'=\{b^1\}$.  Since $b^1=0$, $y+b^1\in A$ implies that $y$ has the form $e_{f(y)}+\sum_{i=f(y)+1}^ny_ie_i$.  Then since $y+b^2\in \neg A$, we must have that $f(y)= f(b^2)<f(b^3)$.  But this forces $y+b^3\in A$, a contradiction.

Suppose now that $u=2$.  By assumption, there exists some $y\in \F_p^n$ such that $(A-y)\cap X'=\{b^2,b^3\}$.  Since $b^1=0$, $y+b^1\notin A$ implies that $y$ has the form $2e_{f(y)}+\sum_{i=f(y)+1}^ny_ie_i$.  Then since $y+b^2\in A$, we must have that $f(y)= f(b^2)<f(b^3)$.  But this forces $y+b^3\notin A$, a contradiction.

A symmetric argument arrives at a contradiction if $f(b^3)<f(b^2)<f(b^1)$.  We must therefore have that $f(b^2)=f(b^3)<f(b^1)$.  Repeating the same argument with $b^3,b^4,b^1$ yields that $f(b^3)=f(b^4)<f(b^1)$. So we have that for $r=f(b^2)=f(b^3)=f(b^4)<n$, there are $x_2,x_3,x_4\in \{2,3\}$, such that  
\begin{align*}
b^2=x_2e_r+\sum_{i=r+1}^nb^2_i,\text{ }b^3&=x_3e_r+\sum_{i=r+1}^nb^3_i,\text{ and }b^4=x_4e_r+\sum_{i=r+1}^nb^4_i.
\end{align*}
By the pigeonhole principle, two of $\{x_2,x_3,x_4\}$ are equal.  

Assume first that $x_2=x_3=1$.  Let $y$ be such that $(A-y)\cap X'=\{b^2,b^1\}$.  Clearly we must have that $f(y)=r$ and $y_r=2$.  But then $y+b^1\notin A$, a contradiction.  Assume instead that $x_2=x_3=2$.  Let $y$ be such that $(A-y)\cap X'=\{b^2,b^1\}$.  Clearly we must have that $f(y)=r$ and $y_r=1$.  But then $y+b^1\notin A$, a contradiction. This shows we cannot have $x_2=x_3$.  A similar argument shows that we cannot have $x_1=x_2$ or $x_1=x_3$, so we have arrived at a contradiction, and thus, we must have that $\VC(A)<4$.  

We end the proof with a quick example of a set of size $3$ which is shattered by $A_{\GS}(3,3)$. From this it follows that $A_{\GS}(3,n)$ shatters a set of size $3$ for any $n\geq 3$. Indeed, let $Z=\{(000), (012),(021)\}$.  It is straightforward to check that $Z$ is shattered by $A$, using its translates by the elements $(011), (020), (000), (010), (001), (022), (100), (200)$.  
\end{proof}

When $p=3$, Green and Sanders showed that $A_{\GS}(p,n)$ always has a large non-trivial Fourier coefficient with respect to the zero coset of any non-trivial subspace of $\F_p^n$. This immediately implies that for any linear factor $\calL$ of $\F_3^n$, there is an atom of $\calL$ on which $A_{\GS}(3,n)$ has density bounded away from $0$ and $1$. We give a direct proof of the latter fact for all primes $p>2$, stated as Proposition \ref{prop:gsvc} in the introduction.  We recall that combining this with Theorem \ref{thm:gs3nip} shows that the conclusions of Theorem \ref{thm:stablegps} cannot hold in general for sets of bounded VC-dimension.

\begin{proposition}[Green-Sanders example requires a linear atom of intermediate density]\label{prop:gsexirr}
For all primes $p>2$, all  $0<\e<p^{-1}$, and all integers $n>\ell\geq 0$,   there is no linear factor on $\F_p^n$ of complexity $\ell$ which is $\e$-atomic with respect to $A_{\GS}(p,n)$.  In fact, for any linear factor $\calL$ on $\F_p^n$ of complexity $\ell$, the density of $A_{\GS}(p,n)$ in $L(0)$ is always in  $[1/p,1-1/p]$.  

Moreover, for any union $Y$ of atoms of $\calL$, $|A_{\GS}(p,n)\Delta Y|>\e |L(0)|$.
\end{proposition}

\begin{proof}
Fix $0<\e<p^{-1}$, and integers $n>\ell\geq0$.    Suppose that $\calL$ is a linear factor on $\F_p^n$ of complexity $\ell$.  Let 
$$
i_0=\max\Big(\{0\}\cup \{i\in [n]: \text{ $e_j$ is contained in the span of $\calL$ for all $j\leq i$}\}\Big).
$$ 
Observe that $i_0\leq \ell<n$.  By definition of $i_0$, $e_{i_0+1}$ is not in the span of $\calL$, and  thus $\calL\cup \{e_{i_0+1}\}$ is a linear factor of complexity $\ell+1$.  Also by definition of $i_0$, $(H_{i_0+1}+e_{i_0+1})\cap L(0)$ is a non-empty atom of $\calL\cup \{e_{i_0+1}\}$, which therefore has size $p^{n-\ell-1}=p^{-1}|L(0)|$.  By definition of $A_{\GS}(p,n)$, $(H_{i_0+1}+e_{i_0+1})\cap L(0)\subseteq A_{\GS}(p,n)\cap L(0)$. Combining these facts, we have that
$$
|A_{\GS}(p,n)\cap L(0)|\geq | (H_{i_0+1}+e_{i_0+1})\cap L(0)|\geq p^{-1}|L(0)|.
$$
Similarly, $(H_{i_0+1}+2e_{i_0+1})\cap L(0)$ is a non-empty atom of $\calL\cup \{e_{i_0+1}\}$ and we have by definition of $A_{\GS}(p,n)$ that $ (H_{i_0+1}+2e_{i_0+1})\cap L(0)\subseteq L(0)\setminus A_{\GS}(p,n)$, so
$$
|L(0)\setminus A_{\GS}(p,n)|\geq | (H_{i_0+1}+2e_{i_0+1})\cap L(0)|\geq p^{-1}|L(0)|.
$$
This shows that the density of $A_{\GS}(p,n)$ in $L(0)$ is always in $[1/p,1-1/p]$.  Suppose now $Y$ is a union of atoms of $\calL$. If $L(0)\in Y$, then $|A_{\GS}(p,n)\Delta L(0)|\geq |L(0)\setminus A_{\GS}(p,n)|\geq p^{-1}|L(0)|$.  If on the other hand, $L(0)\notin Y$, then $|A_{\GS}(p,n)\Delta L(0)|\geq |L(0)\cap A_{\GS}(p,n)|\geq p^{-1}|L(0)|$.  In either case, we have shown $|A_{\GS}(p,n)\Delta Y|>\e |L(0)|$, which finishes the proof. \end{proof}

Next, we turn to Proposition \ref{prop:gslinear} from the introduction (restated below as Proposition \ref{prop:gslinearrecall}), which says that any quadratic factor necessarily has an irregular atom with respect to $A_{\GS}(p,n)$. We recall that this is of interest because, combining with Lemma \ref{lem:vcfop}, it shows that any version of Theorem \ref{thm:vc2finite} for $\NFOP_2$ sets must still allow for irregular atoms.

\begin{proposition}[Green-Sanders example requires a quadratic atom of intermediate density]\label{prop:gslinearrecall}
For all $\e\in (0,1)$, there exists a growth function $\omega=\omega(\e)$ such that for all $\ell,q\geq 0$ and all sufficiently large $n$, the following holds. 

For any quadratic factor $\calB=(\calL,\calQ)$ on $\mathbb{F}_p^n$ of complexity $(\ell,q)$ and rank at least $\omega(\ell+q)$, there is an atom $B\in \At(\calB)$ such that 
$$
\frac{|B\cap A_{\GS}(p,n)|}{|B|}\in \Big(\frac{1-\e}{p},\frac{1+\e}{p}\Big).  
$$
In fact, the above holds for all atoms of the form $B=B(0;b)$ with $b\in \F_p^q$.\end{proposition}

\begin{proof}
Given $\e\in (0,1)$, let $\omega_1=\omega_1(\e^2)$ be from Corollary \ref{cor:sizeofatoms}.  Then define $\omega:\mathbb{R}^+\rightarrow \mathbb{R}^+$ by setting $\omega(x)=\omega_1(x+1)$ for all $x\in \mathbb{R}^+$.   Let $\ell,q\geq 0$, and let $n$ be sufficiently large. Suppose $\calB=(\calL,\calQ)$ is a quadratic factor on $\mathbb{F}_p^n$ of complexity $(\ell,q)$ and rank at least $\omega(\ell+q)$.  Let
$$
i_0=\max\Big(\{0\}\cup \{i\in [n]: \text{ $e_j$ is contained in the span of $\calL$ for all $j\leq i$}\}\Big).
$$ 
Note $i_0\leq \ell<n$.  By definition of $i_0$, we have that $\calL'=\calL\cup \{e_{i_0+1}\}$ is a linear factor of complexity 
$\ell+1$, and thus $\calB'=(\calL',\calQ)$ is a quadratic factor of complexity $(\ell+1,q)$ and rank at least $\omega(\ell+q)\geq \omega_1(\ell+1+q)$.  

Fix any $b\in \F_p^q$.  We show $| A_{\GS}(p,n)\cap B(0;b)|/|B(0;b)|\in ((1-\e)p^{-1},1-(1-\e)p^{-1})$. By Corollary \ref{cor:sizeofatoms} and our choice of $\omega$,
\begin{align}\label{b0bsize}
|B(0;b)|=(1\pm \e^2)p^{n-\ell-q}.
\end{align}
By definition of $A_\GS(p,n)$, $B(0;b)\cap (H_{i_0+1}+e_{i_0+1})\subseteq A_{\GS}(p,n)$. We also have, by definition of $i_0$ and $\calB'$,  that $B(0;b)\cap (H_{i_0+1}+e_{i_0+1})$ is an atom of the factor $\calB'$, so by our choice of $\omega_1$, Corollary \ref{cor:sizeofatoms}, and (\ref{b0bsize}),
$$
|B(0;b)\cap (H_{i_0+1}+e_{i_0+1})|\geq (1- \e^2)p^{n-\ell-1-q}\geq \frac{1- \e^2}{1+ \e^2}p^{-1}|B(0;b)|\geq (1-\e)p^{-1}|B(0;b)|.
$$
Combining the above yields
\[|A_\GS(p,n)\cap B(0;b)|\geq |B(0;b)\cap (H_{i_0+1}+e_{i_0+1})|\geq (1-\e)p^{-1}|B(0;b)|.\]
Similarly, $B(0;b)\cap (H_{i_0+1}+2e_{i_0+1})$ is an  atom  of the factor $\calB'$, and  by definition of $A_\GS(p,n)$, $B(0;b)\cap (H_{i_0+1}+2e_{i_0+1})\cap A_{\GS}(p,n)=\emptyset$. Thus, as above, we have that
\[|B(0;b)\setminus A_\GS(p,n)|\geq |B(0;b)\cap (H_{i_0+1}+2e_{i_0+1})|\geq (1-\e)p^{-1}|B(0;b)|,\]
which concludes the proof.
\end{proof}

We next prove Proposition \ref{prop:gslinear1} (restated as Proposition \ref{prop:gslinear1res} below), which shows that while $A_{\GS}(p,n)$ forces the existence of atoms of intermediate density with respect to any high rank quadratic factor, they can be constrained to a single linear atom if the factor is chosen carefully.  We recall that this result hints at the correct form of decomposition theorem for $\NFOP_2$ sets.

\begin{proposition}\label{prop:gslinear1res}
For all $\e\in (0,1)$, there exists a growth function $\omega=\omega(\e)$ such that for all integers $\ell\geq 1$ and $q\geq 0$ the following holds. 

For all $n\geq \omega(\ell+q)+\ell+1$, there exists a quadratic factor $\calB=(\calL, \calQ)$ on $\mathbb{F}_p^n$ of complexity $(\ell,q)$ and rank at least $\omega(\ell+q)$ with the property that for every $B\in \At(\calB)$, either $B\subseteq L(0)$, in which case 
\[\frac{|B\cap A_{\GS}(p,n)|}{|B|}\in \left(\frac{1-\e}{p},1-\frac{1-\e}{p}\right),\] 
or $B\cap L(0)=\emptyset$, in which case $|B\cap A_{\GS}(p,n)|/|B|\in \{0,1\}$.  
\end{proposition}
\begin{proof}
Fix $\e\in (0,1)$ and let $\omega=\omega(\e)$ be as in Proposition \ref{prop:gslinearrecall}.  Fix $\ell\geq 1$, $q\geq 0$,  and $n\geq \omega(\ell,q)+\ell+1$.  Define $\calL=\{e_1,\ldots, e_{\ell}\}$ and let $\calQ$ be any purely quadratic factor of complexity $q$ so that $\calB=(\calL,\calQ)$ has rank at least $\omega(\ell+q)$ (if $q=0$, take $\calQ=\emptyset$, and if $q>0$, see Fact \ref{fact:basis} for why such a $\calQ$ exists).  Let $L(0)$ be the zero atom of $\calL$.  By definition of $\At(\calB)$, we know that any atom $B\in \At(\calB)$ is either contained in $L(0)$ or disjoint from $L(0)$.  When $B\subseteq L(0)$,  Proposition \ref{prop:gslinearrecall} implies $|B\cap A_{\GS}(p,n)|/|B|\in ((1-\e)p^{-1},1-(1-\e)p^{-1})$. 

Suppose now $B\in \At(\calL)$ and $B\cap L(0)=\emptyset$.  Then there is $a=(a_1,\ldots, a_{\ell})\in \F_p^{\ell}\setminus \{0\}$ such that $B\subseteq L(a)$.  Let 
$$
i_0=\max\Big( \{0\}\cup \{i\in [\ell]: \text{ for all }j\leq i, a_j=0\}\Big).
$$
Since $a\neq 0$, $i_0<\ell$.  Thus, $a_{i_0+1}\in  \{1,\ldots, p-1\}$ and
$$
B\subseteq L(a)\subseteq H_{i_0+1}+a_{i_0+1}e_{i_0+1}.
$$
Thus, either $a_{i_0+1}=1$ in which case  $B\subseteq A_{\GS}(p,n)$, or $a_{i_0+1}\in \{2,\ldots,p-1\}$, in which case $B\cap A_{\GS}(p,n)=\emptyset$.
\end{proof}

To conclude this section, we state the follow proposition, which says that the Green-Sanders  is $\NHOP_2$ in a strong sense.  

\begin{proposition}[Green-Sanders example is $4$-$\NHOP_2$]\label{prop:GSnoHOP}
For all $n\geq 4$, $A_{\GS}(3,n)$ has $3$-$\HOP_2$ but no $4$-$\HOP_2$.
\end{proposition}

The proof of this proposition proceeds by a rather tedious case analysis, and has therefore been relegated to Appendix \ref{app:no4hop}.  We conjecture that a similar result holds for sets of the form $A_{\GS}(p,n)$ for $p>3$.  

Combining Proposition \ref{prop:GSnoHOP} with Proposition \ref{prop:gslinear} demonstrates that any possible strengthening of Theorem \ref{thm:vc2finite} for $\NHOP_2$ sets must still allow for irregular atoms. It is unclear what, if any, structural conclusions can be deduced from the assumption that a set is $\ell$-$\NHOP_2$, beyond those which this paper shows hold because such a set is $\ell$-$\NFOP_2$.  

Both $\HOP_2$ and $A_{\GS}(3,n)$ provide important examples concerning certain questions about regular decompositions of 3-uniform hypergraphs. For an in-depth consideration of these examples from that perspective, we refer the reader to \cite[Section 7]{Terry.2021b}.

\subsection{The quadratic Green-Sanders example}\label{subsec:qgs}

In this section we prove Theorem \ref{thm:ex1}, which says that the quadratic Green-Sanders example of Definition \ref{def:QGS} does not admit the kind of decomposition appearing in the conclusion of Theorem \ref{thm:fop}.   This is in analogy to the linear Green-Sanders example, which does not admit the kind of decomposition appearing in the conclusion of Theorem \ref{thm:stablegps}.

Recall that for odd $n$, there is a basis $(M_i)_{i=1}^n$ of a vector space of $n\times n$ symmetric matrices of rank $n$ with entries in $\F_p$ (see e.g. \cite[Lemma 3]{Dumas.2011}). For every $i\in [n]$, let
\[Q_i(e_i)=\{x\in \F_p^n:x^TM_1x=x^TM_2x=\dots = x^TM_{i-1}x=0, x^TM_ix=1\}.\]

\begin{definition}[Quadratic Green-Sanders example, recalled]\label{def:QGSrecall}
For $p>2$ and odd $n\geq 1$, define the \emph{quadratic Green-Sanders $p$-example of dimension $n$} to be the subset $A_\QGS(p,n)\subseteq \F_p^n$ given by
\[A_{\QGS}(p,n)=\bigcup_{i=1}^n Q_i(e_i).\]
\end{definition}

Note that for every $i\in [n]$, the quadratic factor generated by $M_1,\dots, M_i$ has rank $n$, so Lemma \ref{lem:sizeofatoms} implies that for each $i\in [n]$, $Q_i(e_i)$ has density $p^{-i}+O(p^{-n/2})$ in $\F_p^n$. Since these atoms are all disjoint, the density of $A_{\QGS}(p,n)$ in $\F_p^n$ is 
\begin{align*}
\sum_{i=1}^{(n-1)/2}\left(p^{-i}+O(p^{-n/2})\right)+O\big(\frac{n+1}{2}\cdot p^{-n/2}\big)&=\frac{p^{-1}-p^{-((n-1)/2+1)}}{1-p^{-1}}+O\big(\frac{n+1}{2}\cdot p^{-n/2}\big)\\
&=\frac{1}{p-1}+O(p^{-n/3}).
\end{align*}

The remainder of this section is dedicated to proving Theorem \ref{thm:ex1} (restated below as Proposition \ref{prop:qgsbad}), which shows that a large quadratic Green-Sanders example does not admit a quadratic decomposition satisfying the conclusions of Theorem \ref{thm:fop}.\footnote{The authors are grateful to Val Gladkova for pointing out an error in an earlier version of the proof.}

\begin{proposition}\label{prop:qgsbad}
There exists a growth function $\psi$ so that for all primes $p>2$, there is a growth functions $\omega=\omega(p)$ such that for all integers $\ell\geq 1$ and $r\geq 0$ and all sufficiently large $n$, the following holds. 

Let $A=A_{\QGS}(p,n)$ be a quadratic Green-Sanders $p$-example of dimension $n$, and let $\calB=(\calL,\calR)$ be a quadratic factor on $\F_p^n$ of complexity $(\ell,r)$ and rank at least $\omega(\ell+r)$. Then for every $L\in \At(\calL)$, there is some $R\in \At(\calR)$ such that 
$$
\frac{|A\cap L\cap R|}{|L\cap R|}\in (p^{-\psi(r)}, 1-p^{-\psi(r)}).
$$ 
Moreover, for any union $Y$ of elements of $\At(\calB)$, $|A\Delta Y|>p^{-4\psi(r)}|G|$. 
\end{proposition}

\begin{proof}
Let $\psi$ be defined by setting $\psi(x)=x+2$.  Choose $\mu\in (0,1)$ sufficiently small compared to $p^{-1}$ and $p^{-\psi(0)}$, and let  $\omega(\mu)=\rho(\mu)$ be the function obtained in Corollary \ref{cor:sizeofatoms}.  Let $\tau=\tau(\omega)$ from Lemma \ref{lem:rank}. Fix integers $\ell\geq 1$ and $r\geq 0$, and choose an odd integer $n$ sufficiently large compared to $ \ell$, $r$, $\mu^{-1}$, and $f^{(r)}(\ell+r)$, where $f(x)=\tau(x)+x+1$. 

 Let $A=A_{\QGS}(p,n)$ be a quadratic Green-Sanders $p$-example of dimension $n$, and let $\calB=(\calL,\calR)$ be a quadratic factor on $\F_p^n$ of complexity $(\ell,r)$ and rank at least $\omega(\ell+r)$. Recall that $A=\bigcup_{i=1}^nQ_i(e_i)$, where for each $i$, $e_i$ is the $i$th standard basis vector in $\F_p^i$ and $\calQ_i=\{M_1,\ldots, M_i\}$, where $\{M_1,\ldots, M_n\}$ has rank $n$.  
 
Suppose first $r=0$. In this case, $\At(\calB)=\At(\calL)$, and $\calR$ has a unique atom, $\F_p^n$.  Fix any $L\in \At(\calL)$.  Since $(\calL,\{M_1\})$ has rank $n$, we have by definition of $A$ and Corollary \ref{cor:sizeofatoms} that
$$
|A\cap L|\geq |Q_1(1)\cap L|\geq (1-\mu)|L|>p^{-\psi(0)}|L|,
$$
where the last inequality is by our choice of $\mu$.  Similarly, 
$$
|L\setminus A|\geq |Q_1(2)\cap L|\geq (1-\mu)|L|>p^{-\psi(0)}|L|.
$$
For the ``moreover" part of the statement, let $Y$ be any union of atoms of $\calL$.  Then 
$$
|A\Delta Y|\geq \sum_{L\in Y}|A\Delta L|+\sum_{L\in \At(\calL)\setminus Y}|A\cap L|\geq p^{-\psi(0)}\sum_{L\in \At(\calL)}|L|=p^{-\psi(0)}|G|.
$$
This finishes the case $r=0$.
 
 Suppose now $r\geq 1$. This case will require much more care.   Enumerate $\calR$ as $\{R_1,\ldots, R_r\}$.  We begin with the following technical claim.

\begin{claim}\label{cl:gsbad}
There are integers $0\leq \ell'\leq f^{(r)}(\ell+r)$ and $0\leq q',r'\leq r$, and a quadratic factor $(\calL',\calQ')$ of complexity $(\ell',r')$ such that
\begin{enumerate}[label=\normalfont(\roman*)]
\item $\calL\subseteq \calL'$, and $(\calL',\calQ')\preceq (\calL,\calR)$;
\item  $(\calL',\calQ'\cup \{M_{q'+1}\})$ has rank at least $\omega(\ell'+r'+1)$;
\item one of the following holds.\label{cases1}
\begin{enumerate}
\item $q'=0$ and $(\calL',\calQ')=(\calL,\calR)$. \label{cases1a}
\item $q'=r=r'$ and $\calQ'=\{M_1,\ldots, M_r\}$.\label{cases1b}
\item $0<q'<r$ and,  after possibly reindexing $\calR$, we have $\calQ'=\{M_1,\ldots, M_{q'}\}\cup \calR'$, for some $\calR'\subseteq \{R_{q'+1},\ldots, R_r\}$.\label{cases1c}
\end{enumerate}
\end{enumerate}
\end{claim}
\begin{proof}
We proceed inductively. 

Step 1: If $\calR\cup \{M_1\}$ has rank at least $\omega(\ell+r+1)$, let $q'=0$, $\ell'=\ell$, $r'=r$, $\calL'=\calL$, and $\calQ'=\calR$.  Otherwise, apply Lemma \ref{lem:rank} to $(\calL,\calR\cup \{M_1\})$ using an enumeration listing $M_1$ first, to obtain $(\calL_1,\calS_1)\preceq (\calL,\calR\cup\{M_1\})$ of complexity $(\ell_1,s_1)$ and rank at least $\omega(\ell_1+s_1)$, such that $\calL_1\supseteq \calL$, $\calS_1 \subseteq \calR\cup \{M_1\}$, and $\ell_1\leq \tau(\ell+r+1)$.  Since $(\calL,\calR\cup \{M_1\})$ had rank less than $\omega(\ell+r+1)$, part \ref{start} of Lemma \ref{lem:rank} ensures $\calS_1\subsetneq \calR\cup \{M_1\}$.  Moreover, since $(\calL,\{M_1\})$ has rank $n$, part \ref{keep} of Lemma \ref{lem:rank} ensures that $M_1\in \calS_1$.  Possibly after reindexing $\calR$, we may thus assume that $\calS_1=\{M_1\}\cup \calR_1$ for some $\calR_1\subseteq \{R_2,\ldots, R_r\}$.   

Step $i+1$:  Suppose $i\geq 1$, $\calR_i\subseteq \{R_{i+1},\ldots, R_{r}\}$,  $\calL_i\supseteq \calL$, and $\calS_i=\{M_1,\ldots, M_{i}\}\cup \calR_i$  are such that $s_i:=|\calS_i|\leq r$, $\ell_i:=|\calL_i|\leq \tau(\ell_{i-1}+s_{i-1}+1)$ and $(\calL_i,\calS_i)\prec (\calL,\calR)$. If $\calS_i\cup \{M_{i+1}\}$ has rank at least $\omega(\ell_i+s_i+1)$, let $q'=i$, $\ell'=\ell_i$, $\calQ'=\calS_i$, $\calR'=\calR_i$, $\calL'=\calL_i$, and end the construction. 

Otherwise, apply Lemma \ref{lem:rank} to $(\calL_i,\calS_i\cup \{M_{i+1}\})$, using an enumeration listing the matrices $M_1,\ldots, M_{i+1}$ first.  From this we obtain $(\calL_{i+1},\calS_{i+1})$ of complexity $(\ell_{i+1},s_{i+1})$ and rank at least $\omega(\ell_{i+1},s_{i+1})$ such that $\calL_{i+1}\supseteq \calL_i$, $\calS_{i+1}\subseteq  \calS_i\cup \{M_{i+1}\}$, and $\ell_{i+1} \leq \tau(\ell_i + s_i+1)$.  Since $(\calL_i,\calS_i\cup \{M_{i+1}\})$ had rank less than $\omega(\ell_i+s_i+1)$, part \ref{start} of Lemma \ref{lem:rank} implies that $\calS_{i+1}\subsetneq \calS_i\cup \{M_{i+1}\}$.  Since $(\calL_i,\{M_1,\ldots, M_{i+1}\})$ has rank $n$, part \ref{keep} of Lemma \ref{lem:rank} ensures that $\{M_1,\ldots, M_{i+1}\}\subseteq \calS_{i+1}$.  After possibly reindexing $\calR$, we may therefore assume that $\calS_{i+1}=\{M_1,\ldots, M_{i+1}\}\cup \calR_{i+1}$, for some $\calR_{i+1}\subseteq \{ R_{i+2},\ldots, R_r\}$.

This process will either halt after at most $r$ steps, yielding the desired factor, or proceed all the way to step $r+1$. At step $r+1$, we will have deleted all the elements of $\calR$, yielding $\calS_r=\{M_1,\ldots, M_r\}$. Thus since $(\calL_r,\calS_r\cup \{M_{r+1}\})$ has rank $n$, the algorithm will halt and output the desired factor.   
 \end{proof}

Now let $\ell',q', r', \calL',\calQ',\calR'$ be as in the conclusion of Claim \ref{cl:gsbad}.  We proceed separately for each case appearing in Claim \ref{cl:gsbad} \ref{cases1}.  

Suppose first we are in case \ref{cases1a}, so  $(\calL',\calQ')=(\calL,\calR)$, and  $(\calL,\calR\cup \{M_1\})$ has complexity $(\ell, r+1)$ and rank at least $\omega(\ell+r+1)$.  Fix any $d\in \F_p^{\ell}$ and $a\in \F_p^r$. By Corollary \ref{cor:sizeofatoms} and our choice of $\omega$, for any $u\in \F_p$,
$$
|L(d)\cap R(a)\cap Q_1(u)|\geq (1-\mu)p^{n-\ell-r-1}\geq (1-\mu)(1+\mu)^{-1}p^{-1}|L(d)\cap R(a)|.
$$
By definition of $A$,  $L(d)\cap R(a)\cap Q_1(1)\cap A=\emptyset$ and $L(d)\cap R(a)\cap Q_1(2)\subseteq A$.  Combining these facts, we have that
$$
\min\Big\{\frac{|A\cap L(d)\cap R(a)|}{|L(d)\cap R(a)|},\frac{ | (L(d)\cap R(a))\setminus A|}{|L(d)\cap R(a)|}\Big\}\geq (1-\mu)(1+\mu)^{-1}p^{-1} \geq p^{-r-2},
$$
where the last inequality is because $\mu$ is sufficiently small compared to $p^{-1}$.  This finishes case \ref{cases1a}.

Assume now we are in case \ref{cases1b}, so $q'=r=r'$ and $\calQ'=\{M_1,\ldots, M_r\}$.  Fix $d\in \F_p^{\ell}$, and define 
$$
J_{d}=\{a\in \mathbb{F}_p^r: L(d)\cap Q_{r}(0)\cap R(a)\neq \emptyset\}.
$$
Since $n$ is large and $\calQ_r$ has rank $n$, Lemma \ref{lem:sizeofatoms}, implies  
$$
|L(d)\cap Q_{r}(0)|=(1\pm \mu)p^{n-\ell-r}.
$$
Consequently, since $|J_{d}|\leq p^r$ and $L(d)\cap Q_{r}(0)=\bigcup_{a\in J_{d}}L(d)\cap Q_{r}(0)\cap R(a)$, there is some $a\in J_{d}$ with $|L(d)\cap Q_{r}(0)\cap R(a)|\geq (1-\mu)p^{n-\ell-2r}$.    Fix such an $a\in \F_p^r$ and define 
\[J_{d,a}=\{b\in \mathbb{F}_p^{\ell'-\ell}: L'(db)\cap Q_{r}(0)\cap R(a)\neq \emptyset\}.\]
Observe that for any $b\in \F_p^{\ell'-\ell}$, if $L'(db)\cap Q_{r}(0)\cap R(a)\neq \emptyset$ then since $(\calL',\calQ')\preceq (\calL,\calR)$, 
$$
L'(db)\cap Q_{r}(0)\cap R(a)=L'(db)\cap Q_{r}(0).
$$
Because $n$ is large and $(\calL',\calQ'\cup \{M_{r+1}\})$ has rank $n$,  Lemma \ref{lem:sizeofatoms} implies that for all $b\in J_{d,a}$ and $u\in \F_p$,
\begin{align*}
 |L'(db)\cap Q_{r}(0)\cap \{x\in \mathbb{F}_p^n: x^TM_{r+1}x=u\}|\geq (1-\mu)p^{n-\ell'-r-1}.
\end{align*}
Note that by definition of $A$, $L'(db)\cap Q_{r}(0) \cap \{x\in \mathbb{F}_p^n: x^TM_{r+1}x=u\}$ is contained in $A$ or disjoint from $A$, when $u=1$ or $u=2$, respectively. Hence for any $b\in J_{d,a}$,
\begin{align}\label{alqLB}
\min\{|A\cap L'(db)\cap Q_{r}(0)|, |(\neg A)\cap L'(db)\cap Q_{r}(0)|\}&\geq (1-\mu)p^{n-\ell'-r-1}.
\end{align}
Note that 
$$
L(d)\cap R(a)\cap Q_{r}(0)=\bigsqcup_{b\in J_{d,a}}L'(db)\cap R(a)\cap Q_{r}(0)=\bigsqcup_{b\in J_{d,a}}L'(db)\cap Q_{r}(0).
$$
 Consequently, 
 \begin{align}\label{lrq}
|L(d)\cap R(a)\cap Q_{r}(0)|=\Big|\bigcup_{b\in J_{d,a}}L'(db)\cap Q_{r}(0)\Big|=|J_{d,a}|(1\pm \mu )p^{n-\ell'-r},
\end{align}
where the second equality uses Lemma \ref{lem:sizeofatoms}, the fact $\calQ_r$ has rank $n$, and the fact $n$ is large.  This implies 
\begin{align}\label{jadbound}
|J_{d,a}|\geq p^{-n+\ell'+r}|L(d)\cap R(a)\cap Q_{r}(0)|(1+ \mu)^{-1}\geq \frac{(1-\mu)}{(1+\mu)}p^{\ell'-\ell-r}|L(d)\cap R(a)\cap Q_{r}(0)|,
\end{align}
where the second inequality uses that $|L(d)\cap R(a)\cap Q_{r}(0)|\geq (1-\mu)p^{n-\ell-2r}$ holds by our choice of $a$.  We also have 
 \begin{align*}
 |A\cap L(d)\cap R(a)\cap Q_r(0)|=\sum_{b\in J_{d,a}}|A\cap L'(db)\cap R(a)\cap Q_r(0)|&=\sum_{b\in J_{d,a}}|A\cap L'(db)\cap Q_r(0)|\\
 &\geq (1-\mu)p^{n-\ell'-r-1}|J_{d,a}|,
 \end{align*}
where the last inequality is by (\ref{alqLB}). The same argument shows  
$$
 |(\neg A)\cap L(d)\cap R(a)\cap Q_r(0)|\geq (1-\mu)p^{n-\ell'-r-1} |J_{d,a}|.
 $$
   Consequently,
 \begin{align}\label{alra}
\min\{|A\cap L(d)\cap R(a)|,|(\neg A)\cap L(d)\cap R(a)|\}&\geq (1-\mu)p^{n-\ell'-r-1}|J_{d,a}|.
\end{align}
Using (\ref{jadbound}), we have
\begin{align*}
(1-\mu)p^{n-\ell'-r-1} |J_{d,a}|&\geq (1-\mu)p^{n-\ell'-r-1} \Big(\frac{(1-\mu)}{(1+\mu)}p^{\ell'-\ell-r}|L(d)\cap R(a)\cap Q_{r}(0)|\Big)\\ 
&=\frac{(1-\mu)^2}{(1+\mu)} p^{n-\ell-2r-1}\\
&\geq \frac{(1-\mu)^2}{(1+\mu)^2}p^{-r-1}|L(d)\cap R(a)|\\
&\geq p^{-r-2}|L(d)\cap R(a)|,
\end{align*}
where the second inequality  holds by Corollary \ref{cor:sizeofatoms},  our choice of $\omega$, and since $(\calL,\calR)$ has rank at least $\omega(\ell+r)$, and the last  inequality  holds because $\mu$ is sufficiently small compared to $p^{-1}$.  Combining with (\ref{alra}), we find that 
 \begin{align*}
\min \{|A\cap L(d)\cap R(a)|, |(\neg A)\cap L(d)\cap R(a)| \}&\geq  p^{-r-2}|L(d)\cap R(a)|,
\end{align*} 
which finishes the proof in case \ref{cases1b}.

Suppose now we are in case \ref{cases1c}, so $0<q' <r$, and $\calQ'=\{M_1,\ldots, M_{q'}\}\cup \calR'$, for some $\calR'\subseteq \{R_{q'+1},\ldots, R_r\}$. Given $d\in \F_p^\ell$, define
$$
J_{d}=\{a\in \mathbb{F}_p^r: L(d)\cap Q_{q'}(0)\cap R(a)\neq \emptyset\}.
$$
By Corollary \ref{cor:sizeofatoms} and our choice of $\omega$,
$$
|L(d)\cap Q_{q'}(0)|=(1\pm \mu)p^{n-\ell-q'}.
$$
Consequently, since $|J_{d}|\leq p^r$ and $L(d)\cap Q_{q'}(0)=\bigsqcup_{a\in J_{a}}L(d)\cap Q_{q'}(0)\cap R(a)$, there is some $a\in J_{d}$ with $|L(d)\cap Q_{q'}(0)\cap R(a)|\geq (1-\mu)p^{n-\ell-q'-r}$.   Fix such an $a\in \F_p^r$ and define 
\[J_{d,a}=\{b\in \mathbb{F}_p^{\ell'-\ell}: L'(db)\cap Q_{q'}(0)\cap R(a)\neq \emptyset\}.\]
Observe that for any $b\in \F_p^{\ell'-\ell}$, if $L'(db)\cap Q_{q'}(0)\cap R(a)\neq \emptyset$ then 
$$
L'(db)\cap Q_{q'}(0)\cap R(a)=L'(db)\cap Q_{q'}(0)\cap R'(a|_{\calR'}),
$$
where $a|_{\calR'}=(a_i)_{R_i\in \calR'}$.  Because $(\calL',\calQ'\cup \{M_{q'+1}\})$ has rank at least $\omega(\ell'+r'+1)$, Corollary \ref{cor:sizeofatoms} implies that for all $b\in J_{d,a}$ and $u\in \F_p$,
\begin{align*}
(1\pm \mu) p^{n-\ell'-r'-1}&= |L'(db)\cap Q_{q'}(0)\cap R'(a|_{\calR'})\cap \{x\in \mathbb{F}_p^n: x^TM_{q'+1}x=u\}|\\
&=p^{-1}(1\pm \mu)^2 |L'(db)\cap Q_{q'}(0)\cap R'(a|_{\calR'})|.
\end{align*}
Note that by definition of $A$, $L'(db)\cap Q_{q'}(0)\cap R'(a|_{\calR'})\cap \{x\in \mathbb{F}_p^n: x^TM_{q'+1}x=u\}$ is contained in $A$ or is disjoint from $A$, when $u=1$ or $u=2$, respectively. Hence the above shows that for any $b\in J_{d,a}$,
\begin{align*}
\min\Big\{\frac{|A\cap L'(db)\cap Q_{q'}(0)\cap R'(a|_{\calR'})|}{|L'(db)\cap Q_{q'}(0)\cap R'(a|_{\calR'})|}, \frac{|(\neg A)\cap L'(db)\cap Q_{q'}(0)\cap R'(a|_{\calR'})|}{|L'(db)\cap Q_{q'}(0)\cap R'(a|_{\calR'})|}\Big\}\geq \frac{(1-\mu)^2}{p}.
\end{align*}
Note that $L(d)\cap R(a)\cap Q_{q'}(0)=\bigsqcup_{b\in J_{d,a}}L'(db)\cap Q_{q'}(0)\cap R'(a|_{\calR'})$.    Consequently, 
 
 \begin{align*}
\min\Big\{\frac{|A\cap L(d)\cap R(a)|}{|L(d)\cap R(a)|},\frac{|(\neg A)\cap L(d)\cap R(a)|}{|L(d)\cap R(a)|}\Big\}&\geq \sum_{b\in J_{d,a}}\frac{(1-\mu)^2}{p}\frac{|L'(db)\cap Q_{q'}(0)\cap R'(a|_{\calR'})|}{|L(d)\cap R(a)|}\\ 
&\geq \frac{(1-\mu)^2}{p}\frac{|J_{d,a}|p^{n-\ell'-r'}}{|L(d)\cap R(a)|}\\
&\geq  \frac{(1-\mu)^2}{(1+\mu)} |J_{d,a}|p^{\ell-\ell'+r-r'-1},
\end{align*}
where both inequalities use Corollary \ref{cor:sizeofatoms} and our choices of $\omega$ (the first applied to $(\calL',\calQ')$ and the second applied to $(\calL,\calR)$).  Using Corollary \ref{cor:sizeofatoms} again, we have that 
$$
|L(d)\cap R(a)\cap Q_{q'}(0)|=\Big|\bigsqcup_{b\in J_{d,a}}L'(db)\cap Q_{q'}(0)\cap R'(a|_{\calR'})\Big|=|J_{d,a}|(1\pm \mu)p^{n-\ell'-r'},
$$
and thus,
\begin{align*}
|J_{d,a}|\geq p^{-n+\ell'+r'}|L(d)\cap R(a)\cap Q_{q'}(0)|(1+ \mu)^{-1}.
\end{align*}
Consequently, it follows that 
\begin{align*}
\frac{(1-\mu)^2}{(1+\mu)}  |J_{d,a}|p^{\ell-\ell'+r-r'-1}&\geq \frac{(1-\mu)^2}{(1+\mu)^2} p^{-n+\ell+r-1}|L(d)\cap R(a)\cap Q_{q'}(0)|\\
&\geq \frac{(1-\mu)^3}{(1+\mu)^2} p^{-q'-1}\\
&\geq p^{-r-2},
\end{align*}
where the second inequality uses that $|L(d)\cap Q_{q'}(0)\cap R(a)|\geq (1-\mu)p^{n-\ell-q'-r}$ holds by our choice of $a$,  and the last inequality uses that $q'\leq r$ and $\mu$ is sufficiently small.  Combining, we have that
 \begin{align*}
\min\Big\{\frac{|A\cap L(d)\cap R(a)|}{|L(d)\cap R(a)|},\frac{|(\neg A)\cap L(d)\cap R(a)|}{|L(d)\cap R(a)|}\Big\}&\geq  p^{-r-2},
\end{align*} 
which finishes case \ref{cases1c}.  

We conclude by proving the ``moreover" part of the desired statement.  To this end, let $Y$ be a union of atoms of $\calB$.  Given any $d\in \F_p^{\ell}$,  we have proved there exists some $f(d)\in \F_p^r$ so that 
$$
|A\cap L(d)\cap R(f(d))|/|L(d)\cap R(f(d))|\in (p^{-\psi(r)}, 1-p^{-\psi(r)}),
$$
where we recall $\psi(x)=x+2$.  If $L(d)\cap R(f(d))\in Y$, then this implies 
$$
|(A\Delta Y)\cap L(d)|\geq | (L(d)\cap R(f(d)))\setminus A|>p^{-\psi(r)}|L(d)\cap R(f(d))|.
$$
 On the other hand, if $L(d)\cap R(f(d))\notin Y$, then this implies 
$$
|(A\Delta Y)\cap L(d)|\geq | (L(d)\cap R(f(d)))\cap A|> p^{-\psi(r)}|L(d)\cap R(f(d))|.
$$
In either case, we have shown 
$$
|(A\Delta Y)\cap L(d)|>p^{-\psi(r)}|L(d)\cap R(f(d))|\geq p^{-r-\psi(r)}(1-\mu) |L(d)|=(1-\mu)p^{n-\ell-r-\psi(r)},
$$
 where the second inequality uses Corollary \ref{cor:sizeofatoms} and our assumption on the rank of $\calB$.  We can conclude that
\begin{align*}
|A\Delta Y|\geq \sum_{d\in \F_p^{\ell}}|(A\Delta Y)\cap L(d)|\geq (1-\mu)p^{n-r-\psi(r)}&=(1-\mu)p^{-r-\psi(r)}|G|\geq p^{-4\psi(r)}|G|,
\end{align*}
where the final inequalities use that $p>2$, $r<\psi(r)$, and $\mu$ is sufficiently small.  This finishes the case $r>0$, and thus the proof.
\end{proof}

Proposition \ref{prop:qgsbad} implies Theorem \ref{thm:ex1} as stated in the introduction.  As explained there, in conjunction with Theorem \ref{thm:fop}, this shows that for every $\ell \geq 1$, a sufficiently  large quadratic Green-Sanders example has $\ell$-$\FOP_2$.  We also give a quick explicit proof of this fact in the next section (see Corollary \ref{cor:qgscor}).

\section{Sufficient conditions for $k$-$\FOP_2$ and $k$-$\IP_2$}\label{sec:sufficient}

In this subsection, we give sufficient conditions for when a set $A\subseteq \F_p^n$ has $k$-$\IP_2$ and $k$-$\FOP_2$.  The sufficient condition for $k$-$\FOP_2$ will play a crucial role in the proof of our main theorem.  We include the results for $k$-$\IP_2$ to illustrate the analogy ``stability is to $\FOP_2$ as NIP is to $\NIP_2$."  The sufficient condition for $k$-$\IP_2$ will also be used in a subsequent paper of the authors addressing the bounds in Theorem \ref{thm:vc2finite} (see \cite{Terry.2025}).  The ingredients required for this section include a counting lemma from \cite{Terry.2021d} (covered in Section \ref{ss:counting}) along with notions of reduced structures associated to a set $A$ in a quadratic factor $\calB$ (Section \ref{ss:reduced}).  

\subsection{Counting lemmas}\label{ss:counting} In this subsection we discuss both linear and quadratic counting lemmas.   The linear counting lemma  shows up in several places in Section \ref{sec:FOP}, and also serves as a warm up for the quadratic version, which we will use to prove the sufficient conditions for $\FOP_2$ and $\IP_2$ in Section \ref{ss:reduced}.

We begin with the linear case.  To state this, we require the local versions of the $U^2$ inner product and the $U^2$ norm defined in \cite{Terry.2021d}.  For this and what follows, we will use the following notation. Given $k\geq 1$ and $\sigma=(\sigma_1,\ldots, \sigma_k)\in \{0,1\}^k$, define for each $i\in [k]$,  $\sigma(i)=\sigma_i$,  and let $|\sigma|=\sigma_1+\ldots+\sigma_k$.  We write $\calC$ to denote the complex conjugate operator.

\begin{definition}[Local $U^2$ inner product]\label{def:localinner}
Given a linear factor $\calL$ on $\F_p^n$ of complexity $\ell$, a tuple $d=(a_1,a_2)\in \mathbb{F}_p^{\ell}\times \mathbb{F}_p^{\ell}$ and functions $(f_{\omega})_{\omega\in \{0,1\}^2}:\F_p^n\ra \C$, we define
\[\langle (f_{\omega})_{\omega\in \{0,1\}^2}\rangle_{U^2(d)}=\E_{\begin{subarray}{l} x_0,x_1\in L(a_1)\end{subarray}}\E_{\begin{subarray}{l}y_0,y_1\in L(a_2)\end{subarray}}\prod_{\omega\in \{0,1\}^2}\calC^{|\omega|}f_{\omega}(x_{\omega(1)}+y_{\omega(2)}).\]
\end{definition}

\begin{definition}[Local $U^2$ semi-norm]\label{def:localu2}
Given a linear factor $\calL$ on $\F_p^n$ of complexity $\ell$, a tuple $d=(a_1,a_2)\in \mathbb{F}_p^{\ell}\times \mathbb{F}_p^{\ell}$ and $f:\F_p^n\ra\C$, we define
\[\|f\|_{U^2(d)}= \langle (f)_{\omega\in \{0,1\}^2}\rangle_{U^2(d)}^{1/4}.\]
\end{definition}

It was shown in \cite[Lemma 3.4]{Terry.2021d} that this expression does indeed define a semi-norm on the space of functions on $\F_p^n$.  

We will use this semi-norm to count  instances of bipartite graphs in sum-graphs associated to a set $A$, where the vertices involved fall in certain pre-determined atoms of fixed linear factor.  To describe this in more detail we fix some notation. Given a bipartite graph  $F=(U\cup V, E)$,   a linear factor $\calL$ on $\F_p^n$ of complexity $\ell$, and a tuple of labels $d=(d_u)_{u\in U}(d_v)_{v\in V}\in \F_p^{|U|\ell+|V|\ell}$, define
$$
\calI_F(d)=\prod_{u\in U}L(d_u)\times \prod_{v\in V}L(d_v).
$$
In other words, $\calI_F(d)$ consists of tuples from $\F_p^n$ indexed by the vertices of $F$, which also fall inside atoms of $\calL$ whose labels are determined by the tuple $d$.  Note that since all atoms of a linear factor have the same size, 
$$
|\calI_F(d)|=\prod_{u\in U}|L(d_u)||\prod_{v\in V}L(d_v)|=|L(0)|^{|U|+|V|}.
$$
  Given a set $A\subseteq \F_p^n$, we wish to understand the size of the set
$$
\{(a_u)_{u\in U}(b_v)_{v\in V}\in \calI_F(d): a_u+b_v\in A\text{ if and only if }uv\in E\}.
$$
The linear counting lemma below,  proved in \cite[Proposition D.4]{Terry.2021d}, allows us approximate the size of this set under appropriate uniformity assumptions, in terms of the density of $A$ on the atoms whose labels are sums from the tuple $d$.

\begin{proposition}[Counting lemma for induced binary sums across linear atoms]\label{prop:countinggplinear}
For all $\e\in (0,1)$ and bipartite graphs $F=(U\cup V,E)$ the following holds.  Let $A\subseteq \mathbb{F}_p^n$,  let $\calL$ be a linear factor of complexity $\ell$ on $\F_p^n$, and let $d=(d_u)_{u\in U}(d_v)_{v\in V}\in \F_p^{|U|\ell+|V|\ell}$ be a tuple of labels. For each $(u,v)\in U\times V$, define 
$$
\alpha_{uv}=\frac{|A\cap L(d_u+d_v)|}{|L(d_u+d_v))|}.
$$
Suppose that for each $(u,v)\in U\times V$, $\|1_A-\alpha_{uv}\|_{U^2((d_u,d_v))}\leq \e$.  Then the number of $(x_u)_{u\in U}(y_v)_{v\in V}\in \calI_F(d)$ such that $x_u+y_v\in A$ if and only if $uv\in E$ is 
$$
\Big(\prod_{u\in U, v\in V: uv\in E}\alpha_{uv}\prod_{u\in U, v\in V: uv\notin E} (1-\alpha_{uv}) \pm \mu |U||V|\Big)|\calI_F(d)|.
$$
\end{proposition}

 Before stating the quadratic  analogues of the above, we require some more notation. Given a set of pairs $\Gamma\subseteq \F_p^n\times \F_p^n$, define the \emph{characteristic measure} $\mu_{\Gamma}: \F_p^n\times \F_p^n\rightarrow [0,1]$ of $\Gamma$ by setting $\mu_{\Gamma}(x,y)=1_{\Gamma}(x,y)|\F_p^{n}|^2/|\Gamma|$. This function has the property that its normalised average over $(x,y)\in \F_p^n\times \F_p^n$ is 1. In practice, $\Gamma$ will typically be a bilinear level set associated with a quadratic factor.

\begin{notation}[Bilinear level sets]\label{not:beta}
Let $q\geq 1$ be an integer, and let $\calB=(\calL,\calQ)$ be a quadratic factor on $\F_p^n$ with $\calQ=\{M_1,\ldots, M_q\}$. Given $b=(b_1,\ldots, b_q)\in \F_p^q$, define 
$$
\beta_{\calQ}(b)=\{(x,y)\in \F_p^n\times \F_p^n: \text{ for each }i\in[q], x^TM_iy=b_i\}.
$$
When $\calQ$ is clear from the context, we will drop the subscript on $\beta$.
\end{notation}

It will be helpful to know that if the quadratic factor is high rank, then the associated bilinear level set has roughly the expected size (see \cite[Corollary 2.24]{Terry.2021d}). 

\begin{lemma}\label{lem:sizeofbeta}
Let $\ell \geq 0$ and $q\geq 1$ be integers, and let $\calB=(\calL,\calQ)$ a quadratic factor on $\F_p^n$ of complexity $(\ell,q)$ and rank $\tau$. Then for any $b\in \F_p^q$,
$$
\big| |\beta_{\calQ}(b)|-p^{2n-q}\big|\leq p^{2n+q-\tau}.
$$
\end{lemma}

We now define the quadratic analogue of  Definition \ref{def:localinner} (see \cite[Definition 3.11]{Terry.2021d}).

\begin{definition}[Local $U^3$ inner product]\label{def:localu3}
Given integers $\ell,q\geq 1$, a quadratic factor $\calB=(\calL,\calQ)$ on $\F_p^n$ of complexity $(\ell,q)$, a tuple 
$$
d=(a_1,a_2,a_3,b_{12},b_{13},b_{23})\in \F_p^{\ell+q}\times \F_p^{\ell+q}\times \F_p^{\ell+q}\times \F_p^q\times \F_p^q\times \F_p^q,
$$
 and functions $(f_{\omega})_{\omega\in \{0,1\}^3}:\F_p^n\ra \C$, define
\begin{align*}\langle (f_{\omega})_{\omega\in \{0,1\}^3}\rangle_{U^3(d)}=\E_{ x_0,x_1\in B(a_1)}&\E_{y_0,y_1\in B(a_2)}\E_{z_0,z_1\in B(a_3)}\\
\prod_{(i,j)\in \{0,1\}^2} &\mu_{\beta(b_{12})}(x_i,y_j)\prod_{(i,k)\in \{0,1\}^2} \mu_{\beta(b_{13})}(x_i,z_k)\prod_{(j,k)\in \{0,1\}^2} \mu_{\beta(b_{23})}(y_j,z_k)\\
&\prod_{\omega\in \{0,1\}^3}\calC^{|\omega|}f_{\omega}(x_{\omega(1)}+y_{\omega(2)}+z_{\omega(3)}).
\end{align*}
\end{definition}

Using this inner product, we now define the local analogue of the $U^3$ norm (see \cite[Definition 3.13]{Terry.2021d}).

\begin{definition}[Local $U^3$ semi-norm]\label{def:localu3}
Given integers $\ell,q\geq 1$, a quadratic factor $\calB=(\calL,\calQ)$ on $\F_p^n$ of complexity $(\ell,q)$, a tuple 
$$
d=(a_1,a_2,a_3,b_{12},b_{13},b_{23})\in \F_p^{\ell+q}\times \F_p^{\ell+q}\times \F_p^{\ell+q}\times \F_p^q\times \F_p^q\times \F_p^q,
$$
 and  a function $f:\F_p^n\ra \C$, we define
\[\|f\|_{U^3(d)}= \langle (f)_{\omega\in \{0,1\}^3}\rangle_{U^3(d)}^{1/8}.\]
\end{definition}

We showed in \cite[Lemma 3.16]{Terry.2021d} that this expression defines a semi-norm on the space of functions on $\F_p^n$.  It was also shown that when a set $A$ has density near $0$ or $1$ on a fixed atom of a quadratic factor, then it is uniform in the above sense (see \cite[Proposition 3.20]{Terry.2021d}).

\begin{proposition}[Density near $0$ or $1$ implies locally uniform]\label{prop:locsparseuni}
For all $\e\in (0,1)$, there exists a growth function $\rho_0=\rho_0(\e)$ such that for any growth function $\rho\geq \rho_0$   the following holds. 

Let $\ell,q\geq1$ be integers, let $\calB=(\calL,\calQ)$ be a quadratic factor on $\F_p^n$ of complexity $(\ell,q)$ and rank at least $\rho(\ell+q)$, and let $d=(a_1,a_2,a_3,b_{12},b_{13},b_{23})\in  \F_p^{\ell+q}\times \F_p^{\ell+q}\times \F_p^{\ell+q}\times\F_p^q\times \F_p^q\times \F_p^q$. Given any subset $A\subseteq \F_p^n$, denote by $\alpha_{B(\Sigma(d))}$ the density of $A$ on the atom $B(\Sigma(d))$.

If $\alpha_{B(\Sigma(d))}\in[0,\e)\cup(1-\e,1]$, then $\|1_A-\alpha_{B(\Sigma(d))}\|_{U^{3}(d)}\leq 2\e^{1/8}$. 
\end{proposition}

The goal of our local quadratic counting lemma is to  count copies of a fixed $3$-partite $3$-uniform hypergraph $F$ with vertices behaving in pre-specified ways with respect to a quadratic factor $\calB$.  We describe the relevant type of counting problem in more detail below.

\begin{definition}\label{def:if}
Let $F=(U\cup V\cup W,E)$ be a $3$-partite $3$-uniform hypergraph, let $\ell,q\geq 1$ be integers, let $\calB=(\calL,\calQ)$ be a quadratic factor on $\F_p^n$ of complexity $(\ell,q)$, and let $A\subseteq \F_p^n$.  Assume  that for each $(u,v,w)\in U\times V
\times W$,  we have a tuple 
$$
e_{uvw}=(a_u,b_v,c_w,d_{uv},d_{uw},d_{vw})\in \F_p^{\ell+q}\times \F_p^{\ell+q}\times \F_p^{\ell+q}\times \F_p^q\times \F_p^q\times \F_p^q.
$$
Setting $e=(e_{uvw})_{u\in U, v\in V, w\in W }$, define
\begin{align*}
\calI_F(e)=\{(x_u)_{u\in U}(y_v)_{v\in V}(z_w)_{w\in W}\in  \prod_{u\in U}B(a_u)&\times\prod_{v\in V}B(b_v)\times\prod_{w\in W}B(c_w):\\
 \text{ for each }&(u,v,w)\in U\times V\times W,
(x_u,y_v)\in \beta_{\calQ}(d_{uv}), \\
&(x_u,z_w)\in \beta_{\calQ}(d_{uw}), \text{ and }(y_v,z_w)\in \beta_{\calQ}(d_{vw})\}.
\end{align*}
\end{definition}

When the factor $\calB$ has sufficiently high rank, the sets of the form $\calI_F(e)$ are always non-empty (see \cite[Lemma D.8]{Terry.2021d}). 

In the notation of Definition \ref{def:if}, we are interested in understanding what proportion of tuples $(x_u)_{u\in U}(y_v)_{v\in V}(z_w)_{w\in W}\in \calI_F(e)$ satisfy $ x_u+y_v+z_w\in A$ if and only if $uvw\in E$, for a fixed subset $A\subseteq \F_p^n$.  The answer to this question is closely related to the density of the set $A$ on the atoms containing the sums in question.  To make this precise, we will use the following notation.

\begin{notation}[$\Sigma_{\calB}(d)$]\label{not:sigmad}
Let $\ell,q\geq 1$ be integers, and let $\calB$ be a quadratic factor on $\F_p^n$ of complexity $(\ell,q)$. Given $d=(a_1,a_2,a_3;b_{12},b_{13},b_{23})\in \mathbb{F}_p^{\ell+q}\times  \mathbb{F}_p^{\ell+q}\times  \mathbb{F}_p^{\ell+q}\times  \mathbb{F}_p^{q}\times  \mathbb{F}_p^{q}\times  \mathbb{F}_p^q$, define
$$ 
\Sigma_{\calB}(d)=a_1+a_2+a_3+0(2b_{12})+0(2b_{13})+0(2b_{23})\in \F_p^{\ell+q},
$$
 where each of the vectors $0(2b_{ij})$ is $2b_{ij}$ augmented by $\ell$ initial zeros. When the factor $\calB$ is clear from context, we will drop the subscript.
\end{notation}

The motivation behind this notation is to track the labels of sums in the following sense.  Given a quadratic factor $\calB=(\calL,\calQ)$ on $\F_p^n$ of complexity $(\ell,q)$, and a triple $(x_1,x_2,x_3)\in \F_p^{3n}$, if we know that for each $i\in [3]$, $x_i\in B(a_i)$, and for each $1\leq i<j\leq 3$, $(x_i^TM_u x_j)_{u=1}^q=b_{ij}$, then $x_1+x_2+x_3$ lies in $B(\Sigma(d))$.    We now state our counting lemma, proved in \cite[Proposition D.11]{Terry.2021d}.

\begin{proposition}[Counting lemma for induced ternary sums across quadratic atoms]\label{prop:countinggp}
For every positive integer $m$ and all $\e\in (0,1)$, there exists a growth function $\tau_0=\tau_0(m,\e)$ such that the following holds for any growth function $\tau\geq \tau_0$. 

Let $U,V,W\subseteq [m]$, and let $F=(U\cup V\cup W,E_F)$ be a $3$-partite $3$-uniform hypergraph. Let $\ell,q\geq 1$ be integers, let $\calB=(\calL,\calQ)$ be a quadratic factor on $\F_p^n$ of complexity $(\ell,q)$ and rank at least $\tau$, and let $A\subseteq \mathbb{F}_p^n$. For each $(u,v,w)\in U\times V
\times W$,  let 
$$
e_{uvw}=(a_u,b_v,c_w,d_{uv},d_{uw},d_{vw})\in \F_p^{\ell+q}\times \F_p^{\ell+q}\times \F_p^{\ell+q}\times \F_p^q\times \F_p^q\times \F_p^q 
$$
be a tuple of labels, and let $e=(e_{uvw})_{u\in U, v\in V, w\in W }$.

Suppose that for each $(u,v,w)\in U\times V\times W$, $\|1_A-\alpha_{uvw}\|_{U^3(e_{uvw})}\leq \e$, where 
$$
\alpha_{uvw}=\frac{|A\cap B(\Sigma(e_{uvw}))|}{|B(\Sigma(e_{uvw}))|}.
$$
Then
the number of tuples $(x_u)_{u\in U}(y_v)_{v\in V}(z_w)_{w\in W}\in \calI_F(e)$ such that $ x_u+y_v+z_w\in A$ if and only if $uvw\in E_F$ differs from 
\[\prod_{u\in U, v\in V, w\in W: uvw\in E_F}\alpha_{uvw}\prod_{u\in U, v\in V, w\in W: uvw\notin E_F} (1-\alpha_{uvw}) |\calI_F(e)|\]
by at most $O_m(\e)|\calI_F(e)|$.
\end{proposition}

We will use the following corollary (see \cite[Corollary D.12]{Terry.2021d}), which states that if $\tau$ grows sufficiently fast (depending on $\e$, $|U|$, $|V|$, and $|W|$), and every atom labeled by $e$ is uniform in the local $U^3$ sense and of density in $(\e,1-\e)$,  then there exists some $(x_u)_{u\in U}(y_v)_{v\in V}(z_w)_{w\in W}\in \calI_F(e)$ such that $ x_u+y_v+z_w\in A$ if and only if $uvw\in E_F$.

\begin{corollary}\label{cor:countinggp}
For every positive integer $m$, there exists a constant $C(m)$ and for all $\e\in (0,1)$, there exists a growth function $\tau_0=\tau_0(m,\e)$ such that the following holds for any growth function $\tau\geq \tau_0$. 

Let $U,V,W\subseteq [m]$, and let $F=(U\cup V\cup W,E_F)$ be a $3$-partite $3$-uniform hypergraph. Let $\ell,q\geq 1$ be integers, let $\calB=(\calL,\calQ)$ be a quadratic factor on $\F_p^n$ of complexity $(\ell,q)$ and rank at least $\tau$, and let $A\subseteq \mathbb{F}_p^n$. For each $(u,v,w)\in U\times V
\times W$,  let 
$$
e_{uvw}=(a_u,b_v,c_w,d_{uv},d_{uw},d_{vw})\in \F_p^{\ell+q}\times \F_p^{\ell+q}\times \F_p^{\ell+q}\times \F_p^q\times \F_p^q\times \F_p^q 
$$
 be a tuple of labels, and let $e=(e_{uvw})_{u\in U, v\in V, w\in W }$.

Suppose that for each $(u,v,w)\in U\times V\times W$, $\|1_A-\alpha_{uvw}\|_{U^3(e_{uvw})}\leq C \e^{m^2 2^{m^2}}$, where
$$
\alpha_{uvw}=\frac{|A\cap B(\Sigma(e_{uvw}))|}{|B(\Sigma(e_{uvw}))|}\in (\e,1-\e).
$$ 
Then there exists a tuple $(x_u)_{u\in U}(y_v)_{v\in V}(z_w)_{w\in W}\in \calI_F(e)$ such that $ x_u+y_v+z_w\in A$ if and only if $uvw\in E_F$.
\end{corollary}

\subsection{Reduced structure}\label{ss:reduced} We now turn to defining a notion of a reduced structure associated to a set $A$ and quadratic factor $\calB$. Throughout this subsection, we use some additional notational conventions. First, we will use bars to denote elements we consider as vectors in some $\F_p^m$, and we will denote elements of $\F_p$ without bars.  In keeping with this, we denote the standard basis vectors in $\F_p^m$ by $\ebar_1,\ldots, \ebar_m$.  Given a vector  $\dbar$, we write $\overline{0}\dbar$ to denote the vector obtained by adjoining a string of $0$'s to the left of $\dbar$.  Similarly, given $u\in \F_p$,  $\overline{0}u$ denotes a vector obtained by adjoining a string of $0$'s to the left of $u$. The length of the string $\overline{0}$ will be clear from context in each usage of the notation.  

Our notion of reduced structure will be phrased in terms of a $3$-coloring of a group.

\begin{definition}[3-coloring of a group]
A \emph{$3$-coloring of a group $G$} is a tuple $(G; P_0,P_1,P_2)$ where  $G=P_0\cup P_1\cup P_2$ is a partition.
\end{definition}

Given a positive parameter $\e$, a set $A\subseteq\F_p^n$, and a quadratic factor $\calB$, we naturally obtain the following $3$-coloring in the configuration space associated to a quadratic factor $\calB$. 

\begin{definition}[$\e$-reduced coloring]\label{def:reduced}
Given $A\subseteq \mathbb{F}_p^n$, integers $\ell,q\geq 1$, and a quadratic factor $\calB$ on $\mathbb{F}_p^n$ of complexity $(\ell,q)$, set 
\begin{align*}
\mathbf{A}_0&=\{\abar \in \mathbb{F}_p^{\ell+q}: |A\cap B(\abar)|/|B(\abar)|\leq \e\},\\
\mathbf{A}_1&=\{\abar\in \mathbb{F}_p^{\ell+q}: |A\cap B(\abar)|/|B(\abar)|\geq 1-\e\},\text{ and }\\
\mathbf{A}_2&=\F_p^{\ell+q}\setminus (\mathbf{A}_1\cup \mathbf{A}_0).
\end{align*}
The \emph{$\e$-reduced coloring associated to $A$ and $\calB$}  is the $3$-coloring 
$$
\Red_{\calB,\e}(\F_p^n,A):=(\F_p^{\ell+q}; \mathbf{A}_0,\mathbf{A}_1,\mathbf{A}_2).
$$
\end{definition}

Note that in Definition \ref{def:reduced}, $\mathbf{A}_0$, $\mathbf{A}_1$, and $\mathbf{A}_2$ depend on both $\calB$ and $\e$, but we suppress these from the notation as they will be clear from the context.  Observe further that in Definition \ref{def:reduced}, $\F_p^{\ell+q}$ need not be equal to $\mathbf{A}_0\cup \mathbf{A}_1$, as there could be atoms in $\At(\calB)$ on which the density of $A$ is not near $0$ or $1$.  We will informally refer to these as ``error atoms''. When $\calB$ is almost $\e$-atomic with respect to $A$ (as in the conclusion of Theorem \ref{thm:vc2finite}), there will be very few error atoms, in which case  $\Red_{\calB,\e}(\F_p^n,A)$ will contain a great deal of structural information about the set $A$.  

Our sufficient conditions for $k$-$\IP_2$ and $k$-$\FOP_2$ will deal with certain copies of bipartite graphs in reduced colorings.  

\begin{definition}[$(P_0|P_1)$-copy of a bipartite graph]\label{def:copy}
Given a $3$-coloring $P=(G; P_0,P_1,P_2)$ and a bipartite graph $F=(U\cup W, E)$, a \emph{$(P_0|P_1)$-copy of $F$ in $P$} is a tuple 
$$
(x_u)_{u\in U}(y_w)_{w\in W}\in G^{|U|}\times G^{|W|}
$$
 such that $x_u\cdot y_w\in P_1$ if $uw\in E$ and $x_u\cdot y_w\in P_0$ if $uw\notin E$. 
\end{definition}
Note that in Definition \ref{def:copy}, we impose no requirements on terms of the form $x_u\cdot x_{u'}$ or $y_w\cdot y_{w'}$.  In practice, we will be interested in such copies of the following bipartite graphs, which correspond to stability and NIP respectively.

\begin{definition}[$H(k)$ and $U(k)$]\label{def:ukhk}
Given $k\geq 1$, define the bipartite graphs
\begin{align*}
H(k)&=(\{a_i:i\in [k]\}\cup \{b_i: i\in [k]\}, \{a_ib_j: i\leq j\})\text{ and }\\
U(k)&= (\{a_i:i\in [k]\}\cup \{b_S:  S\subseteq [k]\},\{a_ib_S: i\in S\}).
\end{align*}
\end{definition}

We will use $(\mathbf{A}_0|\mathbf{A}_1)$-copies of $H(k)$ and $U(k)$ in $\Red_{\calB,\e}(\F_p^n,A)=(\F_p^{\ell+q};\mathbf{A}_0,\mathbf{A}_1,\mathbf{A}_2)$ to produce sufficient conditions for when the set $A$ has $\FOP_2$ and $\IP_2$, respectively. For this idea to work,  care must be taken when selecting the labels used to build the $(\mathbf{A}_0|\mathbf{A}_1)$-copy.  

For example, to show a set $A$ has $k$-$\FOP_2$, it is not enough to merely find a $(\mathbf{A}_0|\mathbf{A}_1)$-copy of $H(k)$ in some $\Red_{\calB,\e}(\F_p^n,A)$,  even if the quadratic factor $\calB$ has high rank. Indeed, we show in Proposition \ref{prop:gsred} below that a large linear Green-Sanders example will have many $(\mathbf{A}_0|\mathbf{A}_1)$-copies of $H(k)$ in any $3$-coloring of the form $\Red_{\calB,\e}(\F_p^n,A)$, yet has no $4$-$\FOP_2$ (by Section \ref{subsec:gs}).

\begin{proposition}\label{prop:gsred}
For all primes $p>2$ and integers $k\geq 1$, there exist $\e\in (0,1)$, $\ell_0\geq 1$, and a growth function $\omega$ with the property that for all integers $\ell_0\leq \ell$ and $1\leq q$, and all sufficiently large $n$, the following holds.  Suppose  $\calB$ is a quadratic factor on $\F_p^n$ of complexity $(\ell,q)$ and rank at least $\omega(\ell+q)$ that is almost $\e$-atomic with respect to $A_{GS}(p,n)$. 

Let $\Red_{\calB,\e}(\F_p^n,A_{GS}(p,n))=(\F_p^{\ell+q};\mathbf{A}_0,\mathbf{A}_1,\mathbf{A}_2)$. Then there are at least $\e (p^{\ell+q})^{2k}$ many $(\mathbf{A}_0|\mathbf{A}_1)$-copies of $H(k)$ in $\Red_{\calB,\e}(\mathbb{F}_p^n, A_{\GS}(p,n))$.
\end{proposition}
\begin{proof}
Fix $k\geq 1$ and a prime $p>2$.  Let $\e\in (0,1)$ be sufficiently small based on $p$ and $k$, and let $\omega_1=\omega_1(\e)$ be as in Corollary \ref{cor:sizeofatoms}.  Set $\ell_0=2k$ and $\omega(x)=\omega_1(x+2)$.  Suppose $\ell\geq \ell_0$ and $q\geq 1$.  Let $N=\ell+2k$.

Suppose $n\geq N$ and $\calB=(\calL,\calQ)$ is a quadratic factor  on $\F_p^n$  of complexity $(\ell,q)$ and rank at least $\omega(\ell+q)$, and assume $\calB$ is almost $\e$-atomic with respect to $A=A_{\GS}(p,n)$.  Let $Red_{\calB,\e}(\F_p^n,A)=(\F_p^n;\mathbf{A}_0,\mathbf{A}_1, \mathbf{A}_2)$.

Given $i\in [n]$, let $e_i$ denote the $i$th standard basis vector in $\F_p^n$, and for each $i\in [n]$, set $\calL_i=\{e_j: j\leq i\}$.  It will be convenient to also let $\calL_0=\emptyset$.  Define
$$
i_0=\min\Big\{\{0\}\cup \{i\in [n]: e_j\in \Span(\calL)\text{ for all }j\leq i\}\Big\}.
$$
Note $i_0\leq \ell$, since $\Span(\calL)$ has dimension $\ell$.  By definition of $i_0$, we can find some $\calL'=\{r_1,\ldots, r_{\ell}\}$ so that $\At(\calL')=\At(\calL)$ and such that for each $1\leq i\leq i_0$, $r_i=e_i$ (note this last condition is vacuous if $i_0=0$). Set $\calB'=(\calL',\calQ)$.  By construction, $\At(\calB)=\At(\calB')$, so $\calB'$ is almost $\e$-atomic with respect to $A$, and $\Red_{\calB',\e}(\F_p^n,A)\cong \Red_{\calB,\e}(\F_p^n,A)$.  Note that the factor $(\calL'\cup \calL_{i_0+1},\calQ)=(\calL'\cup \{e_{i_0+1}\},\calQ)$ has complexity $(\ell+1,q)$ and rank at least $\omega(\ell+q)\geq \omega_1(\ell+1+q)$.

Fix any $\bbar\in \F_p^q$ and $\xbar\in \F_p^{\ell-i_0}$.  By the observations above, we have that for any $u\in \F_p$ the set  $L'(\overline{0}\xbar)\cap Q(\bbar)\cap L_{i_0+1}(\overline{0}u)$ is an atom of $(\calL'\cup \{e_{i_0+1}\},\calQ)$ (note $\overline{0}\xbar$ denotes a string of $i_0$ many $0$'s followed by $\xbar$, and $\overline{0}u$ denotes $i_0$ many $0$'s followed by $u$).  Consequently, any such set has size at least $(1-\e)p^{n-\ell-1-q}$ by Corollary \ref{cor:sizeofatoms}.  Applying Corollary \ref{cor:sizeofatoms} again to estimate the size of $L'(\overline{0}\xbar)\cap Q(\bbar)$, this implies 
\begin{align*}
\min\{|L'(\overline{0}\xbar)\cap Q(\bbar)\cap L_{i_0+1}(\overline{0}1)|,|L'(\overline{0}\xbar)\cap Q(\bbar)\cap L_{i_0+1}(\overline{0}2)|\}&\geq \frac{(1-\e)}{(1+\e)}p^{-1}|L'(\overline{0}\xbar)\cap Q(\bbar)|\\
&>\e|L'(\overline{0}\xbar)\cap Q(\bbar)|,
\end{align*}
where the second inequality uses that $\e$ is sufficiently small.  Since $L_{i_0+1}(\overline{0}1)\subseteq A$ and $L_{i_0+1}(\overline{0}2)\subseteq \neg A$ by definition of $A$, we have shown that for all $\bbar\in \F_p^q$ and $\xbar\in \F_p^{\ell-i_0}$,  $(\overline{0}\xbar; \bbar) \in \mathbf{A}_2$, and thus $|\mathbf{A}_2|\geq p^{\ell-i_0+q}$.  By assumption, $|\mathbf{A}_2|\leq \e p^{\ell+q}$, so we must have that  $p^{\ell-i_0}\leq \e p^{\ell}$, which then implies $i_0\geq 2k$, as $\e$ is sufficiently small.  

Let $\ebar_1,\ldots, \ebar_{\ell+q}$ be the standard basis vectors in $\F_p^{\ell+q}$. Fix any tuple $(\xbar_i\ybar_i)_{i\in [k]}\in (\F_p^{\ell+q-k-1}\times \F_p^{\ell+q-k-1})^{k}$, and for each $i\in [k]$, set $\abar_i=\ebar_i+\ebar_{i+1}+\overline{0}\xbar_i$ and $\bbar_i=(p-1)(\ebar_i+\overline{0}\ybar_i)$, where $\overline{0}\xbar_i$ and $\overline{0}\ybar_i$ are obtained by adding a string of $k+1$ many $0$'s to the left of $\xbar_i$ and $\ybar_i$, respectively.  Then $B'(\abar_i+\bbar_j)$ is contained in $A$ if $i\leq j$ and disjoint from $A$ otherwise, so $(\abar_i\bbar_i)_{i\in [k]}$ is a $(\mathbf{A}_0|\mathbf{A}_1)$-copy of $H(k)$ in $\Red_{\calB',\e}(\F_p^n,A)$.  The number of such copies is the number of ways of choosing the tuple $(\xbar_i\ybar_i)_{i\in [k]}$, which is at least $(p^{\ell+q-k-1})^{2k}=p^{-2k(k+1)}(p^{\ell+q})^{2k}\geq \e (p^{\ell+q})^{2k}$, where the final inequality is because $\e$ is sufficiently small.
\end{proof}

We see from the proof of Proposition \ref{prop:gsred} that we can build many $(\mathbf{A}_0|\mathbf{A}_1)$-copies of $H(k)$ in the reduced coloring associated to $A_{\GS}(p,n)$ using the linear part of the factor, with the purely quadratic part not contributing much.  It turns out that in order for a  $(\mathbf{A}_0|\mathbf{A}_1)$-copy of $H(k)$ in a reduced coloring to give rise to an instance of $k$-$\FOP_2$, it must use the purely quadratic structure in a more fundamental way.  Formally, this can be accomplished by considering  $(\mathbf{A}_0|\mathbf{A}_1)$-copies of $H(k)$ where one side of the bipartite graph sits inside a distinguished subspace.  

\begin{definition}[Right side in a subspace]\label{def:oneside}
Let $F=(V\cup W,E)$ be a bipartite graph.  Suppose that $m\geq 1$, $\e\in (0,1)$, $(\F_p^m; P_0,P_1,P_2)$ is a $3$-coloring,  $H\leqslant \F_p^m$ is a subspace, and $(\abar_v)_{v\in V}(\bbar_w)_{w\in W}\in (\F_p^m)^{|V|+|W|}$ is a $(P_0|P_1)$-copy of $F$.  

We say that $(\abar_v)_{v\in V}(\bbar_w)_{w\in W}$  \emph{has the right side in $H$} if $\{\bbar_w: w\in W\}\subseteq H$.
\end{definition}

Where ambiguity might arise, we will always specify which side of a bipartite graph we consider to be the ``right hand'' side.  In particular, we will only use this definition with the bipartite graph $F$ being either  $H(k)$ or $U(k)$ (see Definition \ref{def:ukhk}).  The two sides of $H(k)$ are symmetric, so for these applications no ambiguity arises.  On the other hand, the two sides of $U(k)$ are not symmetric.  In this case, by the ``right side" of $U(k)$, we  will always mean the side indexed by \emph{subsets of $[k]$}.  

Our applications will consider $(\mathbf{A}_0|\mathbf{A}_1)$-copies of $H(k)$ and $U(k)$ in $\Red_{\calB,\e}(\F_p^n,A)$ with right side in the following special subspace defined using the linear part of the factor $\calB$. 

\begin{definition}[$\calL_\calB$ and $H_\calB$]\label{def:hb}
Suppose $n,\ell,q\geq 1$ are integers and  $\calB$ is a quadratic factor on $\F_p^n$ of complexity $(\ell,q)$. Define $H_{\calB}:=\{\overline{0}\bbar\in \F_p^{\ell+q}: \bbar\in \F_p^q\}$.  Let $\calL_{\calB}=\{\ebar_1,\ldots, \ebar_{\ell}\}$ denote the linear factor on $\F_p^{\ell+q}$ consisting of the first $\ell$ basis vectors.
\end{definition}

Observe that the set $H_{\calB}$ in Definition \ref{def:hb}  is a subspace of $\F_p^{\ell+q}$ of codimension $\ell$, and is in fact the atom of $\calL_{\calB}$ labelled by $0\in \F_p^\ell$.

  Our sufficient conditions for $\FOP_2$ and $\IP_2$ will say roughly the following. If there exists an $(\mathbf{A}_0|\mathbf{A}_1)$-copy of $H(k)$ (respectively $U(k)$) in $\Red_{\calB,\e}(\F_p^n,A)$ with right side in $H_{\calB}$, then $A$ has $k$-$\FOP_2$ (respectively $k$-$\IP_2$).   The proofs will be very short from Corollary \ref{cor:countinggp}.

\begin{theorem}[Sufficient condition for $\IP_2$]\label{thm:extract3}
For all primes $p>2$ and integers $k\geq 1$,  there are $\e\in (0,1)$ and a growth function $\omega$ such that for all $\ell, q\geq 1$ and all sufficiently large $n$, the following holds.  Let $A\subseteq \mathbb{F}_p^n$, let $\calB=(\calL,\calQ)$ be a quadratic factor on $\F_p^n$ of complexity $(\ell,q)$ and rank at least $\omega(\ell+q)$, and let $\Red_{\calB,\e}(\F_p^n,A)=(\F_p^{\ell+q};\mathbf{A}_0,\mathbf{A}_1,\mathbf{A}_2)$.  Suppose that there is an $(\mathbf{A}_0|\mathbf{A}_1)$-copy of $k$-$\IP$ in $\Red_{\calB,\e}(\F_p^n,A)$ with right side in $H_{\calB}$.  Then $A$ has $k$-$\IP_2$.
\end{theorem}

\begin{proof}
Fix $k\geq 1$ and a prime $p>2$, and let $t=2k+2^{k^2}$. Let $\e\in (0,1)$ be sufficiently small based on $p$ and $t$,  and let $\omega$ be a sufficiently fast growing growth function. Fix $\ell, q\geq 1$ and let $n$ be sufficiently large.  Suppose $A\subseteq \mathbb{F}_p^n$,  $\calB=(\calL,\calQ)$ is a quadratic factor on $\F_p^n$ of complexity $(\ell,q)$ and rank at least $\omega(\ell+q)$, and $\Red_{\calB,\e}(\F_p^n,A)=(\F_p^{\ell+q};\mathbf{A}_0,\mathbf{A}_1,\mathbf{A}_2)$. Suppose further that there is an $(\mathbf{A}_0|\mathbf{A}_1)$-copy of $k$-$\IP$ in $\Red_{\calB,\e}(\F_p^n,A)$ with right side in $H_{\calB}$.  

Then there exist $\{\fbar_S:S\subseteq [k]\}\cup \{\gbar_j: j\in [k]\}\subseteq \F_p^{\ell+q}$ with $\{\fbar_S:S\subseteq  [k]\}\subseteq H_{\calB}$, and such that for all $j\in [k]$ and $S\subseteq [k]$, if $j\in S$, then $|A\cap B(\gbar_j+\fbar_S)|/|B(\gbar_j+\fbar_S)|\geq 1-\e$, and if $j\notin S$, then $|A\cap B(\gbar_j+\fbar_S)|/|B(\gbar_j+\fbar_S)|\leq \e$.  For each $j\in [k]$, let $\abar_j\in \F_p^{\ell}$ and $\bbar_j\in \F_p^q$ be such that $\gbar_j=\abar_j \bbar_j$, and for each $S\subseteq [k]$, let $\dbar_S\in \F_p^q$ be such that $\fbar_S=\overline{0} \dbar_S$ (each $\fbar_S$ has this form because it is in $H_{\calB}$ by assumption).
 
For each $i\in [k]$ and $X\subseteq [k]^2$, set $\overline{h}_i=\kbar_X=\overline{0}\in \F_p^{\ell+q}$.   Given $X\subseteq [k]^2$, and $i\in [k]$, define $X(i)=\{j\in [k]: (i,j)\in X\}$.  Then for each $i,j\in [k]$ and $X\subseteq[k]^2$, define $\ubar_{ij}=\wbar_{jX}=\overline{0}\in \F_p^q$, and $\vbar_{iX}=\frac{1}{2}\dbar_{X(i)}\in \F_p^q$, and set $\ebar_{ijX}=\overline{h}_i \gbar_j \kbar_{X}\ubar_{ij}\vbar_{iX}\wbar_{jX}$.  Note that 
\[B(\Sigma(\ebar_{ijX}))=L(\abar_{j})\cap Q(\bbar_j+\dbar_{X(i)})=B(\gbar_j+\fbar_{X(i)}).\]
By assumption, the density of $A$ on $B(\Sigma(\ebar_{ijX}))$ is at least $1-\e$ if $j\in X(i)$, and at most $\e$ if $j\notin X(i)$.  By Proposition \ref{prop:locsparseuni} and Corollary \ref{cor:countinggp} (and because $\e$ is sufficiently small and $\omega$ grows sufficiently fast), we see there are $(\sbar_i)_{i\in [k]}(\tbar_i)_{j\in [k]}(\rbar_X)_{X\subseteq [k]^2}$  such that for each $i,j\in [k]$ and $X\subseteq [k]^2$, $\sbar_i+\tbar_j+\rbar_X\in A$ if and only if $j\in X(i)$, i.e. if and only if $(i,j)\in X$.  Thus $A$ has $k$-$\IP_2$.
\end{proof}

\begin{theorem}[Sufficient conditions for $\FOP_2$]\label{thm:extract2}
For all primes $p>2$ and integers $k\geq 1$,  there are $\e\in (0,1)$ and a growth function $\omega$ such that for all $\ell, q\geq 1$ and all sufficiently large $n$,  the following holds.  Let $A\subseteq \mathbb{F}_p^n$, let $\calB=(\calL,\calQ)$ be a quadratic factor on $\F_p^n$ of complexity $(\ell,q)$ and rank at least $\omega(\ell+q)$, and let $\Red_{\calB,\e}(\F_p^n,A)=(\F_p^{\ell+q};\mathbf{A}_0,\mathbf{A}_1,\mathbf{A}_2)$.  Suppose that there is an $(\mathbf{A}_0|\mathbf{A}_1)$-copy of $H(k)$ in $\Red_{\calB,\e}(\F_p^n,A)$ with right side in $H_{\calB}$. Then $A$ has $k$-$\FOP_2$.
\end{theorem}
\begin{proof}
Fix an integer $k\geq 1$ and a prime $p>2$, and set $t=2k+k^{k^2}$. Let $\e\in (0,1)$ be sufficiently small based on $t$ and $p$, and let $\omega$ be a growth function which grows at a sufficiently fast rate. Fix $\ell,q\geq 1$ and let $n$ be sufficiently large.  

Suppose $A\subseteq \mathbb{F}_p^n$, $\calB=(\calL,\calQ)$ is a factor of complexity $(\ell,q)$ and rank at least $\omega(\ell+q)$, and $\Red_{\calB,\e}(\F_p^n,A)=(\F_p^{\ell+q};\mathbf{A}_0,\mathbf{A}_1,\mathbf{A}_2)$.  Assume there is an $(\mathbf{A}_0|\mathbf{A}_1)$-copy of $H(k)$ in $\Red_{\calB,\e}(\F_p^n,A)$ with right side in $\calB$.  By definition, this means there exist $\{\ubar_i : i\in [k]\}\subseteq  \mathbb{F}_p^{\ell+q}$ and $\{\wbar_i: i\in [k]\}\subseteq H_{\calB}$ such that if $i\leq j$, then $|A\cap B(\ubar_i+\wbar_j)|/|B(\ubar_i+\wbar_j)|\geq 1-\e$, and if $i>j$, then $|A\cap B(\ubar_i+\wbar_j)|/|B(\ubar_i+\wbar_j)|\leq \e$.  For each $i\in [k]$, let $\abar_i\in \F_p^{\ell}$ and $\bbar_i,\cbar_i\in \F_p^q$ be such that $\ubar_i=\abar_i\bbar_i$ and $\wbar_i=\overline{0}\cbar_i$ (each $\wbar_i$ has this form since $\wbar_i\in H_{\calB}$).  To ease notation below, we let $I=\{f: [k]^2\rightarrow [k]\}$.

For each $i\in [k]$ and $(i,f)\in [k]\times I$, define $ \xbar_i=\ybar_i^f=\overline{0}\in \F_p^{\ell+q}$.  Then for each $i,j\in [k]$ and $f\in I$, set $\dbar_{ij}^{xu}=\dbar_{ij}^{y^fu}=\overline{0}\in \mathbb{F}_p^q$, and $\dbar_{ij}^{xy^f}=\frac{1}{2}\cbar_{f(i,j)}\in \mathbb{F}_p^q$.  Given $1\leq i,j,s\leq k$, define $\ebar_{ijs}^{f}=\xbar_i\ybar^f_j\ubar_s\dbar_{ij}^{xy^f}\dbar_{is}^{xu}\dbar_{js}^{y^fu}$.  Observe that $B(\Sigma(\ebar^f_{ijs}))=B(\abar_s; \bbar_s+\cbar_{f(i,j)})=B(\ubar_s+\wbar_{f(i,j)})$.  Thus,
\begin{align*}
&\frac{|A\cap B(\Sigma(\ebar^f_{ijs}))|}{|B(\Sigma(\ebar^f_{ijs})|}\geq 1-\e\text{ if }s\leq  f(i,j)\text{ and }\frac{|A\cap B(\Sigma(\ebar^f_{ijs}))|}{|B(\Sigma(\ebar^f_{ijs})|}\leq \e\text{ if }s>f(i,j).
\end{align*}
By Proposition \ref{prop:locsparseuni} and Corollary \ref{cor:countinggp} (and because $\e$ is sufficiently small and $\omega$ grows sufficiently fast),   we see there are $(s_i)_{i\in [k]}(t_i^f)_{(i,f)\in [k]\times I}(r_i)_{i\in [k]}$  such that $s_i+t_j^f+r_m\in A$ if and only if $m\leq f(i,j)$.  Thus $A$ has $k$-$\FOP_2$. 
\end{proof}

Theorem \ref{thm:extract2} is crucial for our main results about $\FOP_2$.  We can also use these ideas to show explicitly that the quadratic Green-Sanders example has unbounded instances of the functional order property.

\begin{corollary}\label{cor:qgscor}
For all primes $p>2$ and integers $k\geq 1$, and and all sufficiently large odd integers $n$, $A_{\QGS}(p,n)$ has $k$-$\FOP_2$.
\end{corollary}
\begin{proof}
Fix an integer $k\geq 1$ and a prime $p>2$. Let $\e$ and $\omega$ be as in Theorem \ref{thm:extract2}, and suppose $n$ is a sufficiently large odd integer.  Set $A=A_{\QGS}(p,n)$, so $A=\bigcup_{i=1}^nQ_i(e_i)$, where for each $i\in [n]$, $\calQ_i=\{M_1,\ldots, M_i\}$, for some set $\{M_1,\ldots, M_n\}$ of rank $n$. Note that $\calQ_k$ is a quadratic factor of complexity $(0,k)$ and rank $n\geq \omega(k)$. Let $\Red_{\calB,\e}(\F_p^n,A)=(\F_p^{\ell+q};\mathbf{A}_0,\mathbf{A}_1,\mathbf{A}_2)$.

For each $1\leq i\leq k$, let $\bbar_i=(p-1)\ebar_{i}$ and $\abar_i=\ebar_{i}+\ebar_{i+1}$.  Note that for any $1\leq i,j\leq k$, $Q_k(\abar_i+\bbar_j)=Q_k(\ebar_{i}+\ebar_{i+1}+(p-1)\ebar_{j})$.  

If $i< j$, then $Q_k(\abar_i+\bbar_j)\subseteq Q_i(\ebar_i)\subseteq A$.  If $i=j$, then $Q_k(\abar_i+\bbar_j)\subseteq Q_{i+1}(\ebar_{i+1})\subseteq \neg A$.  If $j<i$, then $Q_k(\abar_i+\bbar_j)\subseteq Q_{j}((p-1)\ebar_{j})\subseteq\neg A$.  Therefore, $(\abar_i)_{i\in [k]}(\bbar_i)_{i\in [k]}$ is an $(\mathbf{A}_0|\mathbf{A}_1)$-copy of $H(k)$ in $\Red_{\calQ_k, \e}(\F_p^n,A)$ with right side in $H_{\calB}$. By Theorem \ref{thm:extract2}, $A$ has $k$-$\FOP_2$. 
\end{proof}

\section{Proof of Theorem \ref{thm:fop}}\label{sec:FOP}

In this section we prove our main theorem about $\NFOP_2$ sets, Theorem \ref{thm:fop}.  Our strategy will be to work with reduced colorings associated to a set $A$ and a quadratic factor $\calB$ (see Definition \ref{def:reduced}).  By Section \ref{sec:sufficient}, we know that when the set $A$ has no $k$-$\FOP_2$,  the reduced coloring exhibits an approximate version of stability, localized to a distinguished subspace (Theorem \ref{thm:extract2}).  In this section, we will prove a structure theorem for sets satisfying this localized, approximate  form of stability. This result, Theorem \ref{thm:almoststableinhstrong}, is a generalization of the stable regularity lemma of \cite{Terry.2019}, and its proof will use a standard translation between the order property and  certain model theoretic trees. This translation between order properties and trees is well known in the model theory community, having first been proved by Shelah \cite{Shelah.1990o5n}. It  also played a key role in the authors' proof of a linear arithmetic regularity lemma for stable subsets of $\F_p^n$ \cite{Terry.2019} and subsets of abelian groups \cite{Terry.2020}. 

The paper's main result, Theorem \ref{thm:fop}, will then be proved in Section \ref{subsec:foppfs} by applying Theorem \ref{thm:almoststableinhstrong} to a reduced coloring associated to an $\NFOP_2$ set $A$ and a quadratic factor $\calB$ with respect to which $A$ is highly uniform.  Sections \ref{sec:stablereg} and  \ref{subsec:iterate} are then devoted to proving Theorem \ref{thm:almoststableinhstrong}.

\subsection{Translating orders to trees}\label{ss:trees}

In this subsection, we prove a key corollary of Theorem \ref{thm:extract2} using an equivalent formulation of the order property involving trees.  A foundational result in stability theory, due to Shelah \cite{Shelah.1990o5n}, says that there is a correspondence between instances of a certain binary ``tree" and the order property.  In order to describe this correspondence, we will first set some notation. We refer the reader to the literature for further details (see \cite{Hodges.1981, Malliaris.2014, Terry.2019} for treatments in a finitary setting). 

Given an integer $d\geq 1$, let $2^d=\{0,1\}^d$.  In other words, $2^d$ consists of sequences of $0$s and $1$s of length $d$.  By convention, $2^0=\{<>\}$, where $<>$ denotes the so-called empty string.  Then for $d\geq 1$, $2^{<d}:=\bigcup_{i=0}^{d-1}2^i$.  Given $0\leq i\leq j$, $\sigma\in 2^i$ and $\tau\in 2^j$, we write $\sigma \triangleleft \tau$ to denote that $\sigma$ is a proper initial segment of $\tau$, and $\sigma\trianglelefteq \tau$ denotes that either $\sigma=\tau$ or $\sigma\triangleleft \tau$.  By convention, $<>\trianglelefteq \sigma$ holds for any $\sigma\in 2^d$ and $d\geq 0$.  For any $\sigma \in 2^i$ and $\tau\in 2^j$, we let $\sigma \wedge \tau$ denote the element of $2^{i+j}$ obtained by adjoining $\tau$ to the end of $\sigma$ (for example, $(01)\wedge (001)=(01001)$). 

\begin{definition}[$T(d)$]
Given $d\geq 1$, define $T(d)$ to be the bipartite graph with vertex set $\{a_\eta: \eta\in 2^d\}\cup \{b_{\sigma}: \sigma\in 2^{<d}\}$ and edge set $\{a_{\eta}b_{\sigma}: \sigma \wedge 1\trianglelefteq \eta\}$. 
\end{definition}

We think of $T(d)$ as a binary branching tree, and refer to the elements $\{a_\eta: \eta\in 2^d\}$ as the \emph{leaves}, and the elements $\{b_{\sigma}: \sigma\in 2^{<d}\}$ as the \emph{nodes}.  We will use throughout that there are $2^d$ leaves and $2^d-1$ nodes in $T(d)$.  For our purposes, it will only be important for us to know what happens  \emph{along branches} of $T(d)$, i.e. that $\sigma\wedge 1\trianglelefteq \eta$ implies $b_{\sigma}a_{\eta}\in E(T(d))$ and $\sigma\wedge 0\trianglelefteq \eta$ implies $b_{\sigma}a_{\eta}\notin E(T(d))$.  For this reason, we will not work with  copies of $T(d)$, but rather with what we will call \emph{encodings} of $T(d)$.

\begin{definition}[Encoding of $T(d)$]\label{def:encoding}
Suppose $n\geq 1$ and $A\subseteq \F_p^n$. An \emph{encoding of $T(d)$ in $(\F_p^n,A)$} is a tuple $(g_{\eta})_{\eta\in 2^d}(h_{\sigma})_{\sigma\in 2^{<d}}\in (\F_p^n)^{2^d+2^d-1}$ such that if $\sigma \wedge 1\trianglelefteq \eta$, then $g_{\eta}+h_{\sigma}\in A$, and if $\sigma \wedge 0\trianglelefteq \eta$, then $g_{\eta}+h_{\sigma}\notin A$.
\end{definition}

In the notation of Definition \ref{def:encoding}, we will refer to the group elements $(h_{\sigma})_{\sigma\in 2^{<d}}$ as the \emph{nodes} of the encoding, and the group elements $(g_{\eta})_{\eta\in 2^{d}}$ as the \emph{leaves} of the encoding.   It will be important to keep track of nodes and leaves, which motivates the following definition.

\begin{definition}\label{def:nodesleaves1}
Let $X,Y,A\subseteq \F_p^n$,  and let $(a_{\eta})_{\eta\in 2^d}(b_{\sigma})_{\sigma\in 2^{<d}} \in (\F_p^n)^{2^{d}+2^d-1}$ be an encoding of $T(d)$ in $(\F_p^n,A)$.  We say that the encoding has \emph{leaves in $X$ if $\{a_\eta: \eta\in 2^d\}\subseteq X$} and \emph{nodes in $Y$ if $\{b_{\sigma}:\sigma\in 2^{<d}\}\subseteq Y$}.
\end{definition}

Note the definition of an encoding specifies nothing about the sums the form $g_{\eta}+h_{\sigma}$ with $\sigma \ntrianglelefteq \eta$ (which is why we use the word ``encoding" rather than ``copy").  We state below the required correspondence between encodings of trees and the order property. With some minor adjustments to the setting of this paper, this version can be deduced directly from \cite{Hodges.1981}.

\begin{theorem}[Correspondence between trees and the order property]\label{thm:hodges}
For all $k\geq 1$, there exists an integer $d=d(k)\geq 1$ so that the following holds.  Suppose $G=(V,E)$ is a graph containing vertices $(g_{\eta})_{\eta\in 2^d}(h_{\sigma})_{\sigma\in 2^{<d}}$ so that  if $\sigma \wedge 1\trianglelefteq \eta$, then $g_{\eta} h_{\sigma}\in E$, and if $\sigma \wedge 0\trianglelefteq \eta$, then $g_{\eta} h_{\sigma}\notin E$.   

 Then there are $\{b_1,\ldots, b_k\}\subseteq \{g_{\eta}: \eta\in 2^d\}$ and $\{c_1,\ldots, c_k\}\subseteq \{h_{\sigma}: \sigma\in 2^{<d}\}$ such that $c_i b_j\in E$ if and only if $i\leq j$.
\end{theorem}

In other words, this tells us that for $d$ sufficiently large, an encoding of $T(d)$ in $(\F_p^n,A)$ always produces a copy of $H(k)$ in $(\F_p^n,A)$, with ``left side'' contained in the nodes of the encoding and ``right side'' contained in the leaves of the encoding.\footnote{For the curious, Hodges showed one can take $d(k)=2^{k+2}-2$ in Fact \ref{thm:hodges}.  We will not use it here, but there is a similar result going in the reverse direction (see \cite{Hodges.1981}), i.e. if one has a sufficiently large copy of $H(k)$, then one has an encoding of $T(d)$.}  

In light of Theorem \ref{thm:hodges}, we can restate the main result of \cite{Terry.2019} as follows.

\begin{theorem}[Stable linear regularity lemma]\label{thm:linregtree}
For all $d\geq 1$ and $\e\in (0,1)$, there is $M=M(\e,d)$ so that the following holds. Suppose that $A\subseteq \F_p^n$ has the property that there exists no encoding of $T(d)$ in $(\F_p^n,A)$.  Then there exists a subspace $H\leqslant \F_p^n$ of codimension at most $M$ which is $\e$-atomic with respect to $A$.
\end{theorem}

\subsection{Encoding trees in $3$-colorings}

In this section we define encodings of trees in $3$-colorings, and prove a reformulation of Theorem \ref{thm:extract2} in terms of trees.  We will use the following analogues of encodings of $T(d)$ for $3$-colorings. 

\begin{definition}[$(P_0|P_1)$-encoding of $T(d)$]\label{def:goodenc}
Let $n\geq 1$, $\e\in (0,1)$, and suppose $(\F_p^n; P_0,P_1,P_2)$ is a $3$-coloring.  A \emph{$(P_0|P_1)$-encoding of $T(d)$ in $(\F_p^n; P_0,P_1,P_2)$} is a tuple $(a_{\eta})_{\eta\in 2^d}(b_{\sigma})_{\sigma\in 2^{<d}}$ in $(\F_p^n)^{2^d+2^d-1}$ such that 
\begin{itemize}
\item for all $\eta\in 2^d$ and $\sigma\in 2^{<d}$, if $\sigma\wedge 1\triangleleft \eta$, then $a_{\eta}+b_{\sigma}\in P_1$, and if $\sigma\wedge 0\triangleleft \eta$, then $a_{\eta}+b_{\sigma}\in P_0$, and
\item for all $\sigma\in 2^{<d}$ and all $\eta\in 2^d$, $a_{\eta}+b_{\sigma}\in P_1\cup P_0$.
\end{itemize}
\end{definition}

In the notation above, we observe that a $(P_0|P_1)$-encoding of $T(d)$ in $(\F_p^n; P_0,P_1,P_2)$ need not be a $(P_0|P_1)$-copy of $T(d)$ (in the sense of Definition \ref{def:copy}).  Indeed, a $(P_0|P_1)$-copy would insist that $a_{\eta}+b_{\sigma}\in  P_0$ for all $\sigma\ntrianglelefteq \eta$, while a $(P_0|P_1)$-encoding merely requires that $a_{\eta}+b_{\sigma}\in  P_1\cup P_0$ for all $\sigma\ntrianglelefteq \eta$.  We will need to keep track of where the leaves and nodes of our encodings end up, for which we will use the following terminology.

\begin{definition}\label{def:nodesleaves}
Let $X,Y\subseteq \F_p^n$, let $(\F_p^n; P_0,P_1,P_2)$ be a $3$-coloring, and let $(a_{\eta})_{\eta\in 2^d}(b_{\sigma})_{\sigma\in 2^{<d}}$ be a $(P_0|P_1)$-encoding of $T(d)$ in $(\F_p^n; P_0,P_1,P_2)$.  We say that the encoding has \emph{leaves in $X$ if $\{a_\eta: \eta\in 2^d\}\subseteq X$} and \emph{nodes in $Y$ if $\{b_{\sigma}:\sigma\in 2^{<d}\}\subseteq Y$}.
\end{definition}

We next deduce Fact \ref{thm:hodges1}, which is a version of Theorem \ref{thm:hodges} for encodings of $T(d)$ in $3$-colorings.  

\begin{fact}\label{thm:hodges1}
For all $k\geq 1$,  the following holds, where $d=d(k)$ is from Theorem \ref{thm:hodges}.  Suppose $(\F_p^n;P_0,P_1,P_2)$ is a $3$-coloring, and $(g_{\eta})_{\eta\in 2^d}(h_{\sigma})_{\sigma\in 2^{<d}}$ is a $(P_0|P_1)$-encoding of $T(d)$ in $(\F_p^n;P_0,P_1,P_2)$.  

Then there are $\{b_1,\ldots, b_k\}\subseteq \{g_{\eta}: \eta\in 2^d\}$ and $\{c_1,\ldots, c_k\}\subseteq \{h_{\sigma}: \sigma\in 2^{<d}\}$ such that $c_i+b_j\in P_1$ if  $i\leq j$ and $c_i+b_j\in P_0$ if $i>j$.
\end{fact}
\begin{proof}
Apply Theorem \ref{thm:hodges} to the graph $(V,E)$ where 
$$
V=\{g_{\eta}: \eta\in 2^d\}\cup \{h_{\sigma}: \sigma\in 2^{<d}\},
$$
 and $E=\{g_{\eta}h_{\sigma}: g_{\eta}+h_{\sigma}\in P_1\}$.  From this we obtain $\{b_1,\ldots, b_k\}\subseteq \{g_{\eta}: \eta\in 2^d\}$ and $\{c_1,\ldots, c_k\}\subseteq \{h_{\sigma}: \sigma\in 2^{<d}\}$ such that $c_i+b_j\in P_1$ if  $i\leq j$ and $c_i+b_j\notin P_1$ if $i>j$.  Since  $(g_{\eta})_{\eta\in 2^d}(h_{\sigma})_{\sigma\in 2^{<d}}$ is a $(P_0|P_1)$-encoding of $T(d)$, we must have that for all $i>j$ that $c_i+b_j\notin P_1 \Rightarrow c_i+b_j\in P_0$.  This finishes the proof.
\end{proof}

As an immediate corollary, we obtain a repackaging of Theorem \ref{thm:extract2} in terms of trees rather than half-graphs.  

\begin{corollary}[Sufficient condition for $\FOP_2$ in terms of tree encodings]\label{cor:extract2}
For all integers $k\geq 1$ and primes $p>2$, there exists  $\e\in (0,1)$ and a growth function $\omega$, such that for all integers $\ell, q\geq 1$ and all sufficiently large $n$, the following holds, where $d=d(k)$ is from Theorem \ref{thm:hodges}.

Let $A\subseteq \mathbb{F}_p^n$,  let $\calB=(\calL,\calQ)$ be a factor on $\F_p^n$ of complexity $(\ell,q)$ and rank at least $\omega(\ell+q)$, and let $\Red_{\calB,\e}(\F_p^n,A)=(\F_p^{\ell+q}; \mathbf{A}_0,\mathbf{A}_1,\mathbf{A}_2)$.  Suppose that there is an $(\mathbf{A}_0|\mathbf{A}_1)$-encoding of $T(d)$ in $\Red_{\calB,\e}(\F_p^n,A)$ with leaves in $H_{\calB}$. Then $A$ has $k$-$\FOP_2$.

\end{corollary}
\begin{proof}
Fix an integer $k\geq 1$, a prime $p>2$, and let $d=d(k)$ be as in Theorem \ref{thm:hodges}. Choose $\e$ and $\omega$ as in Theorem \ref{thm:extract2}.  Fix $\ell,q\geq 1$ and let $n$ be sufficiently large.  

Suppose $A\subseteq \F_p^n$  and $\calB$ is a quadratic factor on $\F_p^n$ of complexity $(\ell,q)$ and rank at least $\omega(\ell+q)$.  Suppose that there is an $(\mathbf{A}_0|\mathbf{A}_1)$-encoding of $T(d)$ in $\Red_{\calB,\e}(\F_p^n,A)$ with leaves in $H_{\calB}$.  By Fact \ref{thm:hodges1}, there are $\{b_1,\ldots, b_k\}\subseteq \{w_{\eta}: \eta\in 2^d\}$ and $\{c_1,\ldots, c_k\}\subseteq \{v_{\sigma}:\sigma\in 2^{<d}\}$ such that $c_i+b_j\in \mathbf{A}_1$ if $i\leq j$, and $c_i+b_j\in \mathbf{A}_0$ if $i>j$. Consequently, $(c_i)_{i\in [k]}(b_i)_{i\in [k]}$ is an $(\mathbf{A}_0|\mathbf{A}_1)$-copy of $H(k)$ in $\Red_{\calB,\e}(\F_p^n,A)$ (see Definition \ref{def:copy}).  It has right side in $H_{\calB}$ because $\{b_1,\ldots, b_k\}\subseteq \{w_{\eta}: \eta\in 2^d\}\subseteq H_{\calB}$ (since the original encoding of $T(d)$ had leaves in $H_{\calB}$). By Theorem \ref{thm:extract2}, $A$ has $k$-$\FOP_2$.
\end{proof}

\subsection{Distribution of tree encodings in pairs arising from $3$-colorings}

In this section we state and prove Corollary \ref{cor:fewbad}, which gives strong restrictions on the distribution of encodings of trees in pairs of the form $(\F_p^{\ell+q}, \mathbf{A}_1)$, where $\mathbf{A}_1$ arises from the reduced $3$-coloring associated to a set $A$ with no $k$-$\FOP_2$.  Informally, we will see that if $A$ has no $k$-$\FOP_2$, and $\calB$ is a high rank quadratic factor of complexity $(\ell,q)$ which is almost $\e$-atomic with respect to $A$, then there are few encodings of $T(d)$ in $(\F_p^{\ell+q}, \mathbf{A}_1)$ with leaves in $H_{\calB}$, and further, for almost all cosets $C$ of $H_{\calB}$, there are very few encodings of $T(d)$ in $(\F_p^{\ell+q}, \mathbf{A}_1)$ with leaves in $H_{\calB}$ and nodes in $C$.  The idea of the proof is as follows. Recall that Corollary \ref{cor:extract2} says that there can be no $(\mathbf{A}_0|\mathbf{A}_1)$-encodings of $T(d)$ in $\Red_{\calB,\e}(\F_p^n,A)$ with leaves in $H_{\calB}$. Consequently, any encoding  of $T(d)$ in $(\F_p^{\ell+q}, \mathbf{A}_1)$  must use the label of an error atom, i.e. an element of $\mathbf{A}_2$.  Since $\calB$ is almost $\e$-atomic, the set $\mathbf{A}_2$ is very small, so there can be only a small number of such encodings.

\begin{corollary}\label{cor:fewbad}
For all integers $k\geq 1$ and primes $p>2$, there exists $\e^*\in (0,1)$ and a growth function $\omega^*$ such that for all $0<\e\leq \e^*$, there is $\mu=\mu(k,\e)\in (0,1)$ so that for all growth functions $\omega\geq \omega^*$, all integers $\ell,q\geq 1$, and all sufficiently large $n$, the following holds.

Assume $A\subseteq \F_p^n$ has no $k$-$\FOP_2$ and $\calB=(\calL,\calQ)$ is a quadratic factor on $\F_p^n$ of complexity $(\ell,q)$ and rank at least $\omega(\ell+q)$ which is almost $\mu$-atomic with respect to $A$.  Then the following hold, where $(\F_p^{\ell+q}; \mathbf{A}_0,\mathbf{A}_1,\mathbf{A}_2)=\Red_{\e,\calB}(\F_p^n, A)$ and  $d=d(k)$ is from Theorem \ref{thm:hodges}.
\begin{enumerate}[label=\normalfont(\arabic*)]
\item There are at most $\e |\calH_{\calB}|^{2^d}|\F_p^{\ell+q}|^{2^d-1}$ encodings of $T(d)$ in $(\F_p^{\ell+q},\mathbf{A}_1)$ with leaves in $H_{\calB}$.
\item There is a set $\Omega\subseteq \F_p^{\ell+q}/H_{\calB}$ satisfying $|\Omega|\leq \e p^{\ell}$ such that for all $g\notin \Omega$, the number of encodings of $T(d)$ in $(\F_p^{\ell+q},\mathbf{A}_1)$ with leaves in $H_{\calB}$ and nodes in $H_{\calB}+g$ is at most $\e |H_{\calB}|^{2^d+2^d-1}$.
\end{enumerate}
\end{corollary}
\begin{proof}
Fix an integer $k\geq 1$ and a prime $p>2$, and let $d=d(k)$ be as in Theorem \ref{thm:hodges}.   Let $\e^*$ and $\omega^*$ be as in Corollary \ref{cor:extract2} for $k$ and $p$.  Given $0<\e\leq \e^*$, set $\mu=\e^{4(2^d+2^d-1)}$.  Fix a growth function $\omega\geq \omega^*$, integers  $\ell,q\geq 1$, and a sufficiently large integer $n$.

Suppose that $A\subseteq \F_p^n$ does not have $k$-$\FOP_2$, and that $\calB=(\calL,\calQ)$ is a quadratic factor on $\F_p^n$ of complexity $(\ell,q)$ and rank at least $\omega(\ell+q)$ which is almost $\mu$-atomic with respect to $A$.  Let  $\Red_{\calB,\mu}(\F_p^n,A)=(\mathbb{F}_p^{\ell+q};\mathbf{A}_1,\mathbf{A}_0,\mathbf{A}_2)$.  Since $\calB$ is almost $\mu$-atomic with respect to $A$, we know $|\mathbf{A}_2|\leq \mu p^{\ell+q}$.   Let $\calL_1$ be a linear factor in $\F_p^{\ell+q}$ such that the zero-atom is $H_{\calB}$, i.e. such that $L_1(0)=H_{\calB}$.  Note $\calL_1$ has complexity $\ell$ as $H_{\calB}$ has codimension $\ell$ in $\F_p^{\ell+q}$.  For ease of notation, set $G'=\F_p^{\ell+q}$, and define 
$$
\Omega=\{L\in \At(\calL_1): |L\cap \mathbf{A}_2|\geq \mu^{1/2}p^{q}\}.
$$
Note that every element of $\Omega$ is a coset of $H_{\calB}$, and since $|\mathbf{A}_2|\leq \mu p^{\ell+q}$, we must have that $|\Omega|\leq \mu^{1/2}p^{\ell}$.  We will show (1) and (2) hold with this choice of $\Omega$.

Since $A$ does not have $k$-$\FOP_2$, Corollary \ref{cor:extract2} implies that if $(a_{\eta})_{\eta\in 2^d}(b_{\sigma})_{\sigma\in 2^{<d}}$ is an encoding of $T(d)$ in $(\F_p^n,\mathbf{A}_1)$   with leaves in $H_{\calB}$, then there must be some $\sigma\in 2^{<d}$ and $\eta\in 2^d$, such that $b_{\sigma}+a_{\eta}\in \mathbf{A}_2$.  Thus each encodings of $T(d)$ in $(\F_p^{\ell+q},\mathbf{A}_1)$ with leaves in $H_{\calB}$ will be constructed by the following procedure.
\begin{itemize}
\item Choose $\sigma\in 2^{<d}$ and $\eta\in 2^d$.  There are $(2^d)(2^d-1)$ choices.
\item Choose $g\in \mathbf{A}_2$.  There are at most $\mu p^{\ell+q}$ choices.
\item Choose $a_{\eta}\in H_{\calB}$.  There are $|H_{\calB}|$ choices.
\item Set $b_{\sigma}=g-a_{\eta}$.  There is one choice.
\item Choose any $b_{\sigma'}\in G'$ for each $\sigma'\in 2^{<d}\setminus\{\sigma\}$ and any $a_{\eta'}\in H_{\calB}$ for each $\eta\in 2^d\setminus\{\eta\}$.  There are $|G'|^{2^d-2}|H_{\calB}|^{2^d-1}$ choices.
\end{itemize}
This yields a bound of 
$$
(2^d)(2^d-1)\mu p^{\ell+q}|H_{\calB}|^{2^d}|G'|^{2^d-2}=(2^d)(2^d-1)\mu |H_{\calB}|^{2^d}|G'|^{2^d-1}<\e |H_{\calB}|^{2^d}|G'|^{2^d-1}
$$
for the number of encodings of $T(d)$ in $(\F_p^{\ell+q},\mathbf{A}_1)$ with leaves in  $H_{\calB}$, where the last inequality is by our choice of $\mu$.  This completes the verification of (1).  

We now prove (2).  We have already shown that $|\Omega|\leq \mu^{1/2}p^q<\e p^q$, where the second inequality is by our choice of $\mu$.   Assume now that $H_{\calB}+g\notin \Omega$ and $(a_{\eta})_{\eta\in 2^d}(b_{\sigma})_{\sigma\in 2^{<d}}$ is an encoding of $T(d)$ in $(\F_p^{\ell+q}, \mathbf{A}_1)$  with leaves in $H_{\calB}$ and nodes in $H_{\calB}+g$.  Since $A$ has no $k$-$\FOP_2$, Corollary \ref{cor:extract2} implies there must be some $\sigma\in 2^{<d}$ and $\eta\in 2^d$ such that $b_{\sigma}+a_{\eta}\in \mathbf{A}_2$.  Note that since $b_{\sigma}\in H_{\calB}+g$ and $a_{\eta}\in H_{\calB}$, we have $b_{\sigma}+a_{\eta}\in \mathbf{A}_2 \cap (H_{\calB}+g)$.  This shows that for every $H_{\calB}+g\notin \Omega$, every encoding of $T(d)$ in $(\F_p^{\ell+q},\mathbf{A}_1)$ with leaves in $H_{\calB}$ and nodes in $H_{\calB}+g$ will be constructed by the following procedure.
\begin{itemize}
\item Choose $\sigma\in 2^{<d}$ and $\eta\in 2^d$.  There are $(2^d)(2^d-1)$ choices.
\item Choose $f\in \mathbf{A}_2\cap (H_{\calB}+g)$.  Since $H_{\calB}+g\notin \Omega$, the number of choices is at most $\mu^{1/2} p^{q}$. 
\item Choose $a_{\eta}\in H_{\calB}$.  There are $|H_{\calB}|$ choices.
\item Set $b_{\sigma}=f-a_{\eta}$.  There is one choice.
\item Choose any $b_{\sigma'}\in H_{\calB}+g$ for each $\sigma'\in 2^{<d}\setminus \{\sigma\}$ and any $a_{\eta'}\in H_{\calB}$ for each $\eta'\in 2^d\setminus\{\eta\}$.  There are $|H_{\calB}|^{2^d+2^d-3}$ choices.
\end{itemize}

Thus, for every $g\notin \Omega$, the number of encodings of $T(d)$ in $(\F_p^{\ell+q},\mathbf{A}_1)$ with leaves in $H_{\calB}$ and nodes in $H_{\calB}+g$ is at most
\begin{align*}
(2^d)(2^d-1)\mu^{1/2} |H_{\calB}|p^{q}|H_{\calB}|^{2^d+2^d-3}=(2^d)(2^d-1)\mu^{1/2} |H_{\calB}|^{2^d+2^d-1}<\e |H_{\calB}|^{2^d+2^d-1},
\end{align*}
where the equality uses that $|H_{\calB}|=p^q$, and the final inequality is by definition of $\mu$.  This completes the verification of (2). 
\end{proof}

In light of Corollary \ref{cor:fewbad}, it is a good time to state the generalized stable regularity lemma we will use to prove Theorem \ref{thm:fop}.  Its hypotheses are designed to reflect the conclusions of Corollary \ref{cor:fewbad}, and it will later be applied to structures of the form $(\F_p^{\ell+q},\mathbf{A}_1)$, where $\mathbf{A}_1$ comes from $\Red_{\calB,\e}(\F_p^n,A)$ for some $A\subseteq \F_p^n$ with no $k$-$\FOP_2$.

\begin{theorem}[Generalised stable arithmetic regularity lemma]\label{thm:almoststableinhstrong}
For all primes $p>2$, all increasing functions $\psi:\mathbb{N}\rightarrow \mathbb{R}^+$, and all integers $d\geq 1$, there exist positive integers $D=D(p,\psi,d)$, $N=N(p,\psi,d)$, and $C=C(p,\psi,d)$, and a real number $\mu_0=\mu_0(p,\psi,d)\in (0,1)$   such that for all $0<\mu\leq \mu_0$, the following holds.

Suppose $\ell\geq 0$, $n-\ell \geq N$, $A\subseteq \F_p^n$,  $H\leqslant \F_p^n$, and $\Omega\subseteq G/H$ satisfy the following.
\begin{enumerate}[label=\normalfont(\roman*)]
\item  $H$ has codimension $\ell$ in $\F_p^n$ and $|\Omega|\leq \mu p^{\ell}$;
\item the number of encodings of $T(d)$ in $(\F_p^n, A)$ with leaves in $H$ is at most $\mu^C |H|^{2^d}|G|^{2^d-1}$;
\item for all $H+g\notin \Omega$, the number of encodings of $T(d)$ in $(\F_p^n, A)$ with leaves in $H$ and nodes in $H+g$ is at most $\mu |H|^{2^d+2^d-1}$. 
\end{enumerate}
Then there exists a subspace $H'\leqslant H$ of codimension $m\leq D$ in $H$, and a set $\Omega'\subseteq G/H'$ satisfying $|\Omega'|\leq (\mu +p^{-\psi(m)})p^{m+\ell}$ such that for every $K\in (G/H')\setminus \Omega'$, 
$$
\frac{ |A\cap K|}{|K|}\in [0,p^{-\psi(m)})\cup (1-p^{-\psi(m)},1].
$$
\end{theorem}

The proof of this theorem appears in Section \ref{sec:stablereg}, and is based on the original proof of Theorem \ref{thm:stablegps} from \cite{Terry.2019}.  We end this subsection with a brief discussion of the key features of Theorem \ref{thm:almoststableinhstrong}.  First, it is extremely important that the bound $D$ on the codimension of $H'$ in $H$ depends only on $p$, $\psi$, and $d$, and not on $\ell$ or on $\mu$.  For instance, Theorem \ref{thm:almoststableinhstrong} provides another proof of Theorem \ref{thm:stablestructure}, a fact which relies on this dependence among the parameters.  The reader may wish to compare the following with Theorem \ref{thm:stablestructure}. 

\begin{theorem}\label{thm:stablegpsstrong}
For all integers $k\geq 1$, primes $p>2$, and increasing functions $\theta:\mathbb{N}\rightarrow \mathbb{R}^+$, there  $D\geq 1$ such that the following holds for all sufficiently large $n$.  Let $A\subseteq \F_p^n$ be $k$-stable.  Then there exists a linear factor $\calL$ of complexity $m\leq D$ which is $p^{-\theta(m)}$-atomic with respect to $A$.
\end{theorem}
\begin{proof}
Fix an  integer  $k\geq 1$, prime  $p>2$, and increasing function $\theta:\mathbb{N}\rightarrow \mathbb{R}^+$.  Let $D$, $N$, $\mu_0$, and $C$ be as in Theorem \ref{thm:almoststableinhstrong} for $\psi(x)=2\theta(x)+2x$ and $d=d(k)$ from Theorem \ref{thm:hodges}.  Choose any $0<\mu<\mu_0p^{-2D}$.  

Let $n\geq N$ and suppose that $A\subseteq G=\F_p^n$ is $k$-stable. Observe the hypotheses of Theorem \ref{thm:almoststableinhstrong} are satisfied with  $H=G$, $\ell=0$, and $\mu$.  Indeed, because $A$ is $k$-stable, there are no encodings at all of $T(d)$ in $(\F_p^n,A)$.  Theorem \ref{thm:almoststableinhstrong} thus yields a subspace $H'\leqslant G$ of codimension $m\leq D$ and a set $\Omega'\subseteq G/H'$ of size at most $(\mu +p^{-\psi(m)})|G/H'|$, so that $A$ has density within $p^{-\psi(m)}$ of $0$ or $1$ on all cosets of $H'$ outside $\Omega'$.  By our choice of $\mu$ and $\psi$, $|\Omega'|\leq (\mu+p^{-\psi(m)})p^m<1$, so $\Omega'=\emptyset$.  Thus the partition generated by the cosets of $H'$ is in fact $p^{-\psi(m)}$-atomic with respect to $A$. 
\end{proof}

Note that this proof relies on the fact that we could choose $\mu$ last among the parameters.  A more sophisticated use of this dependence will appear in the proof of our main theorem, Theorem \ref{thm:fop}.  It is important for our applications that Theorem \ref{thm:almoststableinhstrong} can be applied  in a wider range of circumstances than Theorem \ref{thm:stablegpsstrong}, namely in cases where the set $A$ is not necessarily stable.

\subsection{Reduction to Theorem \ref{thm:almoststableinhstrong}}\label{subsec:foppfs} 

In this section we prove our main result, assuming Theorem \ref{thm:almoststableinhstrong}. We begin with an informal outline of the strategy.   Given a set $A\subseteq \F_p^n$ with no $k$-$\FOP_2$, we will first choose a high rank quadratic factor $\calB$ of bounded complexity $(\ell,q)$, such that $\calB$ is almost $\e$-atomic with respect to $A$ (such a $\calB$ will exist by Theorem \ref{thm:vc2finite}).  We will then consider the reduced $3$-coloring $\Red_{\calB,\e}(\F_p^n,A)= (\F_p^{\ell+q};\mathbf{A}_0,\mathbf{A}_1,\mathbf{A}_2)$.  By Corollary \ref{cor:extract2}, the hypotheses of Theorem \ref{thm:almoststableinhstrong} will be satisfied by the pair $(\F_p^{\ell+q},\mathbf{A}_1)$ and subgroup $H_{\calB}$.  We will then apply Theorem \ref{thm:almoststableinhstrong} to $(\F_p^{\ell+q},\mathbf{A}_1)$ and  $H_{\calB}$, from which we will obtain a subspace $H\leqslant H_{\calB}$ of bounded codimension in $H_{\calB}$, with the property that almost all cosets of $H$ are $\e'$-atomic with respect to $\mathbf{A}_1$ for a very small $\e'$.  We will then ``pull back'' the partition of $\F_p^{\ell+q}$ generated by the cosets of $H$ to find a partition of the original group $\F_p^n$ satisfying the conclusions of Theorem \ref{thm:fop}.

 There are two remaining technicalities to be addressed before we can carry out this plan.  The first technicality is how to ``pull  back'' a linear factor in configuration space to obtain a quadratic factor in the original vector space.   This will be dealt with in Proposition \ref{prop:pullbackfactor}.  Second, we will be applying Theorem \ref{thm:almoststableinhstrong} to $(\F_p^{\ell+q},\mathbf{A}_1)$, where $\mathbf{A}_1$ comes from $\Red_{\calB,\e}(\F_p^n,A)$, and $(\ell,q)$ is the complexity of $\calB$.  Theorem \ref{thm:almoststableinhstrong} contains a lower bound on the size of the vector space involved, so we will need to be able to ensure that $\ell$ and $q$ are sufficiently large.  We will deal with this  in Theorem \ref{thm:quadlb}.   
 
Below, and throughout this subsection, we will use bars to denote tuples of elements from $\F_p^n$,  while we will   denote single elements from $\F_p$ as letters without bars. We will frequently be attaching subscripts to vectors.  To avoid confusion in this case, we will write $(\xbar_{\alpha})_i$ to denote the $i$th coordinate of a  vector $\xbar_{\alpha}$.
 
  We turn now to proving Proposition \ref{prop:pullbackfactor}. Given a quadratic factor $\calB$ on $\F_p^n$, any linear factor in the configuration space of $\calB$ naturally gives rise to a partition of $\F_p^n$.  
 
\begin{definition}\label{def:pullback}
Let $\ell,q\geq 1$ be integers, let  $\calB=(\calL,\calQ)$ be a quadratic factor on $\F_p^n$ of complexity $(\ell,q)$, and let $\calR$ be a linear factor on $\F_p^{\ell+q}$.  For each $R\in \At(\calR)$, set $X_{R}=\bigcup_{\cbar\in R}B(\cbar)$, where each $B(\cbar)$ is an atom of $(\calL,\calQ)$.  Define
 $$
 \calX_{\calR}=\{X_{R}: R\in \At(\calR)\}.
 $$
\end{definition}
Note that for each $R\in \At(\calR)$, $X_{R}\subseteq \F_p^n$.  Clearly $\calX_{\calR}$ is a partition of $\F_p^n$, since $\At(\calR)$ is a partition of $\At(\calB)$ and $\At(\calB)$ is a partition of $\F_p^n$.  Thus Definition \ref{def:pullback} gives us a way of building a partition $\calX_{\calR}$ of $\F_p^n$ from the atoms of the linear factor $\calR$ in the configuration space associated to a quadratic factor in $\F_p^n$.   It will be important to note that, in the notation of Definition \ref{def:pullback}, if two linear factors $\calR,\calR'$ on $\F_p^{\ell+q}$ satisfy $\At(\calR)=\At(\calR')$, then $\calX_{\calR}=\calX_{\calR'}$.

Our next goal is to show that under certain hypotheses, the partition $\calX_{\calR}$ in Definition \ref{def:pullback} is equal to the set of atoms of a quadratic factor.  The relevant hypotheses will always be that the factor $\calR$ refines the linear factor $\calL_{\calB}$ from Definition \ref{def:hb}.  Our first lemma in this direction shows that if $\Span(\calR)=\Span(\calL_{\calB})$, then $\calX_{\calR}$ is simply $\At(\calL)$.

\begin{lemma}\label{lem:pullback1}
Let $\ell,q\geq 1$ be integers. Suppose $\calB=(\calL,\calQ)$ is a quadratic factor on $\F_p^n$ of complexity $(\ell,q)$, and $\calR $ is a linear factor on $\F_p^{\ell+q}$ satisfying $\calR \preceq\calL_{\calB}$. If $\Span(\calR)=\Span(\calL_{\calB})$, then $\calX_{\calR}=\At(\calL)$.
\end{lemma}
\begin{proof}
Let $\calL=\{\sbar_1,\ldots, \sbar_{\ell}\}$, $\calQ=\{M_1,\ldots, M_q\}$ and $\calR=\{\rbar_1,\ldots, \rbar_k\}$.  Recall that, by definition, $\calL_{\calB}=\{\ebar_1,\ldots, \ebar_{\ell}\}$ is the linear factor on $\F_p^{\ell+q}$ consisting of the first $\ell$ basis vectors in $\F_p^{\ell+q}$.  Since  $\Span(\calR)=\Span(\calL_{\calB})$,  $\At(\calR)=\At(\calL_{\calB})$, so by the remarks following Definition \ref{def:pullback},  $\calX_{\calR}=\calX_{\calL_{\calB}}$.  Therefore, it suffices to show  $\calX_{\calL_{\calB}}=\At(\calL)$. We accomplish this by showing that for all $\abar\in \F_p^{\ell}$, $X_{L_{\calB}(\abar)}=L(\abar)$ (note here $L_{\calB}(\abar)\subseteq \F_p^{\ell+q}$ is an atom of $\calL_{\calB}$, while $L(\abar)\subseteq\F_p^n$ is an atom of $\calL$).

To this end, fix $\abar=(a_1,\ldots, a_{\ell})\in \mathbb{F}_p^{\ell}$. By definition, $\xbar\in X_{L_{\calB}(\abar)}$ if and only if there is some $\dbar\cbar\in L_{\calB}(\abar)$ such that $\xbar\in B(\dbar;\cbar)$.  Note that for any $\dbar\cbar\in \F_p^{\ell}\times \F_p^q$, 
\begin{align*}
\xbar\in B(\dbar;\cbar)\Leftrightarrow(\xbar^T \sbar_1,\ldots, \xbar^T\sbar_{\ell};\xbar^TM_1\xbar, \ldots, \xbar^TM_q\xbar)=(\dbar;\cbar).
\end{align*}
Consequently, $\xbar\in X_{L_{\calB}(\abar)}$ if and only if 
$$
(\xbar^T \sbar_1,\ldots, \xbar^T\sbar_{\ell},\xbar^TM_1\xbar, \ldots, \xbar^TM_q\xbar)\in L_{\calB}(\abar).
$$
Since $\calL_{\calB}$ consists of the first $\ell$ basis vectors in $\F_p^{\ell}$, this holds if and only $(\xbar^T \sbar_1,\ldots, \xbar^T\sbar_{\ell})=\abar$, i.e. if and only if $\xbar\in L(\abar)$.      Thus we have shown that $\xbar\in X_{ L_{\calB}(\abar)}$ if and only if  $\xbar\in L(\abar)$, as desired.  
\end{proof}

Our next  lemma gives us a description of the elements in $\calX_{\calR}$  in terms of linear and bilinear forms, in the case where $\Span(\calR)\neq \Span(\calL_{\calB})$.

\begin{lemma}\label{lem:pullback2}
Let $\ell,q\geq 1$ be integers, and let $\calB=(\calL,\calQ)$ be a quadratic factor on $\F_p^n$ of complexity $(\ell,q)$, where $\calL=\{\sbar_1,\ldots, \sbar_{\ell}\}$ and $\calQ=\{M_1,\ldots, M_q\}$.  

Suppose $\calR=\{\rbar_1,\ldots, \rbar_k\}\preceq\calL_{\calB}$ is a linear factor on $\F_p^{\ell+q}$ satisfying $\Span(\calR)\neq \Span(\calL_{\calB})$, and for each $\ell<\alpha\leq k$, let $M^{\alpha}:=\sum_{j=1}^q(\rbar_{\alpha})_{\ell+j}M_j$. For all $\abar=(a_1,\ldots, a_k)\in \F_p^k$ and  $\xbar\in \F_p^{n}$, the following are equivalent.
\begin{enumerate}
\item $\xbar\in X_{R(\abar)}$. 
\item $\xbar\in L(a_1,\ldots, a_{\ell})$ and for all $\ell < \alpha\leq k$, $a_{\alpha}=\sum_{i=1}^{\ell}(\rbar_{\alpha})_i(\xbar^T\sbar_i)+\xbar^TM^{\alpha}\xbar$. 
\end{enumerate}
\end{lemma}
\begin{proof}
Recall $\calL_{\calB}=\{\ebar_1,\ldots, \ebar_{\ell}\}$ consists of the first $\ell$ basis vectors in $\F_p^{\ell+q}$.  Fix $\xbar\in \F_p^n$ and $\abar=(a_1,\ldots, a_{k})\in \mathbb{F}_p^k$.  By definition, $\xbar\in X_{R(\abar)}$ if and only if there is $\dbar\cbar\in R(\abar)$ such that $\xbar\in B(\dbar;\cbar)$.  By definition, for any $(\dbar,\cbar)\in \F_p^{\ell}\times \F_p^q$,
\begin{align*}
\xbar\in B(\dbar;\cbar)\Leftrightarrow(\xbar^T \sbar_1,\ldots, \xbar^T\sbar_{\ell},\xbar^TM_1\xbar, \ldots, \xbar^TM_q\xbar)=\dbar\cbar.
\end{align*}
Thus  
\begin{align}\label{al:3}
\xbar\in X_{R(\abar)}&\Leftrightarrow (\xbar^T \sbar_1,\ldots, \xbar^T\sbar_{\ell},\xbar^TM_1\xbar, \ldots, \xbar^TM_q\xbar)\in R(\abar)\nonumber \\
&\Leftrightarrow \text{  for all $1\leq \alpha\leq k$}, (\xbar^T \sbar_1,\ldots, \xbar^T\sbar_{\ell},\xbar^TM_1\xbar, \ldots, \xbar^TM_q\xbar) \cdot \rbar_\alpha=a_\alpha,
\end{align}
where the second ``$\Leftrightarrow$" is by definition of $R(\abar)$.    When $1\leq \alpha\leq \ell$, we know that $\rbar_{\alpha}=\ebar_{\alpha}$, in which case  (\ref{al:3}) simplifies to the equation $\xbar^T\sbar_{\alpha}=a_{\alpha}$.  For $\ell<\alpha\leq k$, expanding (\ref{al:3}) yields 
\begin{align}\label{al:4}
a_\alpha=\sum_{i=1}^{\ell} (\rbar_\alpha)_i (\xbar^T\sbar_i) + \sum_{j=1}^q (\rbar_\alpha)_{\ell+j}(\xbar^TM_j\xbar)=\sum_{i=1}^{\ell} (\rbar_\alpha)_i (\xbar^T \sbar_i) +\xbar^TM^\alpha \xbar,
\end{align}
where the last equality is by definition of $M^{\alpha}$.  This finishes the proof.
\end{proof}

 We now prove our main result about partitions of the form $\calX_{\calR}$, when $\calR\preceq \calL_{\calB}$.\footnote{A more general version of Proposition \ref{prop:pullbackfactor} appears in earlier versions of this paper (see e.g. Lemma 4.24 in \cite{Terry.2021av2}). While this more general version may be of interest, we have updated our proof so that the above  simpler version suffices.  }

 \begin{proposition}\label{prop:pullbackfactor}
Let $\ell,q\geq 1$ and $r\geq 0$ be integers.  Suppose $\calB=(\calL,\calQ)$ is a quadratic factor on $\F_p^n$ of complexity $(\ell,q)$ and rank $\rho$, and $\calR\preceq\calL_{\calB}$ is a linear factor of complexity $\ell+r$ on $\F_p^{\ell+q}$.  

Then there exists a quadratic factor $\tilde{\calB}=(\calL,\tilde{\calQ})$ on $\F_p^n$ of complexity $(\ell,r)$ and rank at least $\rho$ such that $\calX_{\calR}=\At(\tilde{\calB})$.   
\end{proposition}
\begin{proof}
Fix integers $\ell,q\geq 1$ and a quadratic factor $\calB=(\calL,\calQ)$ on $\F_p^n$ of complexity $(\ell,q)$ and rank $\rho$, and suppose that $\calR\preceq\calL_{\calB}$ is a linear factor in $\F_p^{\ell+q}$ of complexity $k=\ell+r$.  Recall $\calL_{\calB}=\{\ebar_1,\ldots, \ebar_{\ell}\}$, where $\ebar_i$ denotes the $i$-th standard basis vector in $\F_p^{\ell+q}$.  Recall that $\calX_{\calR}$ depends only on the partition $\At(\calR)$ of $\F_p^{\ell+q}$, which in turn depends only on $\Span(\calR)$. Thus, after possibly replacing $\calR$ with a different basis for $\Span(\calR)$, we may assume  that $\calR\supseteq \calL_{\calB}$, and that for all $\rbar\in \calR\setminus \calL_{\calB}$, $(\rbar)_j=0$ for all $1\leq j\leq \ell$. 

Suppose first $\calR=\calL_{\calB}$. Then $r=0$, and Lemma \ref{lem:pullback1} implies $\calX_{\calR}=\At(\calL)$.  In this case, set $\tilde{\calB}=(\calL,\emptyset)$. Then $\At(\tilde{\calB})=\At(\calL)=\calX_{\calR}$ and $\tilde{\calB}$ has  rank $n\geq \rho$ by convention (since $\tilde{\calB}$ is purely linear).   

 Assume now $\calR\neq \calL_{\calB}$. Then we must have $r>0$ and $k>\ell$.  Fix an enumeration  $\calR=\{\rbar_1,\ldots, \rbar_k\}$ where $\rbar_1=\ebar_1,\ldots, \rbar_{\ell}=\ebar_{\ell}$, and let  $\calQ=\{M_1,\ldots, M_q\}$.  For each $\ell+1\leq \alpha\leq k$, set $M^\alpha=\sum_{j=1}^q (\rbar_\alpha)_{\ell+j}M_j$, and define $\tilde{\calB}=(\calL,\tilde{\calQ})$, where $\tilde{\calQ}=\{M^{\alpha}: \ell+1\leq \alpha \leq k\}$.   We will first show  that $\calX_{\calR}=\At(\tilde{\calB})$.
 
 \underline{$\calX_{\calR}\subseteq \At(\tilde{\calB})$:} Fix  $\abar=(a_1,\ldots, a_k)\in \F_p^k$. We show $X_{R(\abar)}=\tilde{B}(a_1,\ldots, a_{\ell}; a_{\ell+1}',\ldots, a_k')$, where for each $\ell+1\leq \alpha\leq k$, $a_{\alpha}'=a_{\alpha}-\sum_{i=1}^{\ell}(\rbar_{\alpha})_ia_i$.  To this end, fix $\xbar\in X_{R(\abar)}$.   By Lemma \ref{lem:pullback2},  $\xbar\in L(a_1,\ldots, a_{\ell})$ and for all $\ell+1\leq \alpha\leq k$,  
$$
a_{\alpha}=\sum_{i=1}^{\ell} (\rbar_\alpha)_i (\xbar^T \sbar_i) +\xbar^TM^\alpha \xbar=\sum_{i=1}^{\ell}(\rbar_{\alpha})_ia_i+\xbar^TM^\alpha \xbar,
$$
where the last equality is because $\xbar\in L(a_1,\ldots,a_{\ell})$.  Rearranging yields that $\xbar^TM^{\alpha}\xbar=a_{\alpha}'$.  Thus $\xbar\in \tilde{B}(a_1,\ldots, a_{\ell}; a_{\ell+1}',\ldots, a_k')$.

Conversely, suppose $\xbar\in \tilde{B}(a_1,\ldots, a_{\ell}; a_{\ell+1}',\ldots, a_k')$.    By definition of $\tilde{\calB}$, this implies  $\xbar\in L(a_1,\ldots, a_{\ell})$.   For $\ell+1\leq \alpha\leq k$, by definition of $a_{\alpha}'$, and since $\xbar^TM^{\alpha}\xbar=a_{\alpha}'$,  we have
$$
a_{\alpha}=a_{\alpha}'+\sum_{i=1}^{\ell}(\rbar_{\alpha})_ia_i=\xbar^TM^\alpha \xbar+\sum_{i=1}^{\ell}(\rbar_{\alpha})_ia_i=\xbar^TM^\alpha \xbar+\sum_{i=1}^{\ell} (\rbar_\alpha)_i (\xbar^T \sbar_i),
$$
where the last equality uses that $\xbar\in L(a_1,\ldots, a_{\ell})$.  By Lemma \ref{lem:pullback2}, $\xbar\in X_{R(\abar)}$. This finishes the proof of $\calX_{\calR}\subseteq \At(\tilde{\calB})$.
 
 \underline{$\calX_{\calR}\subseteq \At(\tilde{\calB})$:} Fix   $\bbar=(b_1,\ldots, b_k)\in \F_p^k$. We show  $\tilde{B}(b_1,\ldots, b_k)=X_{R(b_1,\ldots, b_{\ell},b_{\ell+1}',\ldots, b_k')}$,   where for $\ell+1\leq \alpha\leq k$, $b_{\alpha}'=b_{\alpha}+\sum_{i=1}^{\ell}(\rbar_{\alpha})_ib_i$.  Our proof of the ``$\calX_{\calR}\subseteq \At(\tilde{\calB})$" direction shows  
 $$
 X_{R(b_1,\ldots, b_{\ell},b_{\ell+1}',\ldots, b_k')}=\tilde{B}(b_1,\ldots, b_{\ell},b_{\ell+1}'',\ldots, b_k''),
 $$
  where for each $\ell+1\leq \alpha\leq k$, $b_{\alpha}''=b'_{\alpha}-\sum_{i=1}^{\ell}(\rbar_{\alpha})_ib_i$.   However, for each $\ell+1\leq \alpha\leq k$, $b'_{\alpha}-\sum_{i=1}^{\ell}(\rbar_{\alpha})_ib_i=b_{\alpha}$ by definition of $b_{\alpha}'$.  Thus  $\tilde{B}(b_1,\ldots, b_k)=X_{R(b_1,\ldots, b_{\ell},b_{\ell+1}',\ldots, b_k')}$. This shows $\At(\tilde{\calB})\subseteq \calX_{\calR}$, which finishes our verification that $\calX_{\calR}=\At(\tilde{\calB})$.
  
We just have left to show $\tilde{\calB}$ has rank at least $\rho$. Suppose towards a contradiction this is false.  We now consider two cases.

Suppose first $r=1$. Then $k=\ell+1$ and $\tilde{\calQ}=\{M^{\ell+1}\}$. Our assumption that $\rk(\tilde{\calB})<\rho$ implies there exists some $\lambda\in \F_p\setminus \{ 0\}$ so that  $\lambda M^{\ell+1}=\lambda \sum_{j=1}^q (\rbar_{\ell+1})_{\ell+j}M_j$  has rank less than $\rho$. Since $\calB$ has rank $\rho$ and $\lambda\neq 0$, we must have that for each $1\leq j\leq q$, $(\rbar_{\ell+1})_{\ell+j}=0$. However, this implies $\calR\subseteq \Span(\calL_{\calB})$ and thus $\calR=\calL_{\calB}$, contradicting our assumption.   

Suppose now $r>1$. Then $k>\ell+1$, and our assumption that $\rk(\tilde{\calB})<\rho$ implies there exist $\lambda_{\ell+1},\ldots, \lambda_k\in \F_p$, not all $0$, so that $ \lambda_{\ell+1}M^{\ell+1}+\ldots+\lambda_kM^k$ has rank less than $\rho$.  By reindexing if necessary, we may assume $\lambda_{\ell+1}\neq 0$.  Using the definition of the $M^{\alpha}$ and rearranging, we have
\begin{align}\label{factorM}
\lambda_{\ell+1}M^{\ell+1}+\ldots+\lambda_kM^k=\sum_{j=1}^q\Big(\sum_{\alpha=\ell+1}^k\lambda_{\alpha}(\rbar_{\alpha})_{\ell+j} \Big)M_j.
\end{align}
Since $\calB$ has rank $\rho$, it must be the case that for all $j\in [q]$, the coefficient of $M_j$ in (\ref{factorM})  is zero, i.e. $\sum_{\alpha=\ell+1}^k\lambda_{\alpha} (\rbar_{\alpha})_{\ell+j}=0$ for all $j\in [q]$.   Consequently, for each $j\in [q]$, $(\rbar_{\ell+1})_{\ell+j}=-\lambda_{\ell+1}^{-1}\sum_{\alpha=\ell+2}^k \lambda_{\alpha}(\rbar_{\alpha})_{\ell+j}$. Therefore, we can write
\begin{align*}
\rbar_{\ell+1}&=\sum_{j=1}^{\ell}(\rbar_{\ell+1})_j\ebar_j+\sum_{j=1}^q\Big(-\lambda_{\ell+1}^{-1}\sum_{\alpha=\ell+2}^k \lambda_{\alpha}(\rbar_{\alpha})_{\ell+j}\Big)\ebar_{\ell+j}\\
&=\sum_{j=1}^{\ell}(\rbar_{\ell+1})_j\ebar_j-\lambda_{\ell+1}^{-1}\sum_{j=1}^q\sum_{\alpha=\ell+2}^k \lambda_{\alpha}(\rbar_{\alpha})_{\ell+j}\ebar_{\ell+j}\\
&=\sum_{j=1}^{\ell}(\rbar_{\ell+1})_j\ebar_j-\lambda_{\ell+1}^{-1} \sum_{\alpha=\ell+2}^k \lambda_{\alpha}\sum_{j=1}^q(\rbar_{\alpha})_{\ell+j}\ebar_{\ell+j}\\
&=\sum_{j=1}^{\ell}(\rbar_{\ell+1})_j\ebar_j-\lambda_{\ell+1}^{-1} \sum_{\alpha=\ell+2}^k \lambda_{\alpha}\sum_{j=1}^q \rbar_{\alpha}, 
\end{align*}
where the last equality uses our assumption that for all $\ell+1\leq \alpha\leq k$ and all $1\leq j\leq q$, $(\rbar_{\alpha})_j=0$.   We have now expressed $\rbar_{\ell+1}$ as a linear combination of $\{\ebar_1,\ldots, \ebar_{\ell},\rbar_{\ell+2},\ldots, \rbar_k\}$, contradicting that $\calR$ has complexity $k$.  This finishes the proof.
\end{proof}

We now turn to proving Theorem \ref{thm:quadlb}, which will show that we can ensure the existence of  almost $\e$-atomic quadratic factors of sufficiently large complexity for sets of bounded $\VC_2$-dimension.  These results rely on the following simple fact about refinements of almost $\e$-atomic partitions. The proof is elementary, and included in Appendix \ref{app:easy} for completeness.

\begin{fact}\label{fact:aveofavecor}
Suppose $A\subseteq \F_p^n$ and $\calY,\calY'$ are partitions of $\F_p^n$ with $\calY'$ refining $\calY$.  If $\calY$ is almost $\e$-atomic with respect to $A$, then $\calY'$ is almost $\sqrt{2\e}$-atomic with respect to $A$. 
\end{fact}

As a warmup, we show how to use Fact \ref{fact:aveofavecor} to prove a version of Theorem \ref{thm:vc} with a lower bound on the complexity of the resulting linear factor.

\begin{theorem}\label{thm:linearlower}
For all primes $p$, integers $k,f_1\geq 1$, and reals $\mu\in (0,1)$, there is an integer $D\geq 1$ such that the following holds. Suppose $n$ is sufficiently large and $A\subseteq \F_p^n$ satisfies $\VC(A)<k$.  Then there is $f_1\leq \ell\leq D$ and a linear factor $\calL$ of complexity $\ell$ which is almost $\mu$-atomic with respect to $A$.
\end{theorem}
\begin{proof} Choose $D_1=M(p,k,\mu^2/4)$ as in Theorem \ref{thm:vc} and set $D_2=\max\{D_1,f_1\}$. Let $n$ be sufficiently large, and suppose $A\subseteq \F_p^n$ satisfies $\VC(A)<k$. By Theorem \ref{thm:vc}, there exists a linear factor $\calL$ on $\F_p^n$ of complexity $\ell\leq D_1$   which is almost $\mu^2/4$-atomic with respect to $A$.  If $\ell\geq f_1$, then we are done, so assume $\ell<f_1$.  Since $n$ is large, there exists a set $\calX\subseteq\F_p^n$ of size $f_1-\ell$ such that $\calL'=\calX\cup \calL$ is a set of linearly independent vectors.  Since $\At(\calL')$ refines $\At(\calL)$, Fact \ref{fact:aveofavecor} implies $\calL'$ is almost $\mu$-atomic with respect to $A$.
\end{proof}

We want to prove the following analogous result for sets of bounded $\VC_2$-dimension. In particular, this is a version of Theorem \ref{thm:vc2finite} with a lower bound on the complexity of the resulting quadratic factor.

\begin{theorem}\label{thm:quadlb}
For all primes $p>2$,  integers $k,f_1,f_2\geq 1$, reals $\mu\in (0,1)$, and growth functions $\omega$, there is an integer $D\geq 1$ so that the following holds. Let $n$ be sufficiently large, and suppose  $A\subseteq \F_p^n$ satisfies $\VC_2(A)<k$.  Then there is are integers $f_1\leq \ell\leq D$ and $f_2\leq q\leq D$, and a quadratic factor $\calB=(\calL,\calQ)$ on $\F_p^n$ of complexity $(\ell,q)$ and rank at least $\omega(\ell+q)$, such that $\calB$ is almost $\mu$-atomic with respect to $A$.
\end{theorem}

The proof of this will be more complicated than the proof of the linear analogue, Theorem \ref{thm:linearlower}, due to rank considerations.  In particular, we require the following lemma about adding high rank matrices to an existing factor (see Appendix \ref{app:easy} for a proof).

\begin{lemma}\label{lem:padding}
Let $q,q'\geq 1$ be integers and $r>0$ a real.  Let $n$ be sufficiently large compared to $q,q'$ and $r$, let $\calS=\{M_1,\ldots, M_{n-1}\}$ be a purely quadratic factor on $\F_p^n$ of rank at least $n-1$, and let $\calQ$ be a purely quadratic factor on $\F_p^n$ of complexity $q$ and rank at least $r$.  Then there is a set $\calS'\subseteq \calS$ so that $\calQ\cup \calS'$ has complexity $q+q'$ and rank at least $r$.
\end{lemma}

\vspace{2mm}

The last thing we need for  Theorem \ref{thm:quadlb} is the following fact about existence of high rank factors. This will be necessary for us to apply Lemma \ref{lem:padding} above. It is an immediate corollary of  \cite[Lemma 3]{Dumas.2011}. For completeness, we include a proof in Appendix \ref{app:easy}.

\begin{fact}\label{fact:basis}
For all $n\geq 1$, there exists a purely quadratic factor $\{M_1,\ldots, M_{n-1}\}$ on $\F_p^n$ of complexity $n-1$ and rank at least $n-1$.
\end{fact}

We now prove  Theorem \ref{thm:quadlb}.

\vspace{2mm}

\begin{proofof}{Theorem \ref{thm:quadlb}} Fix a prime $p>2$, a real $\mu\in (0,1)$, integers $k,f_1,f_2\geq 1$, and a growth function $\omega$.  Let $\omega_1$ be the growth function defined by  setting $\omega_1(x)= \omega(x+f_1+f_2)$.  Choose  $C=C(p,k, \mu^2/2,\omega_1)$ as in Theorem \ref{thm:vc2finite}, and set $D=C+f_1+f_2$.  Let $n$ be sufficiently large compared to $D$, and $\omega_1(2D)$.

Suppose $A\subseteq \F_p^n$ satisfies $\VC_2(A)<k$.  By Theorem \ref{thm:vc2finite}, there exist integers $\ell',q'\geq 0$ and a quadratic factor $\calB'=(\calL',\calQ')$ of complexity $(\ell',q')$ and rank at least $\omega_1(\ell'+q')$, such that $\ell'+q'\leq C$, and such that $\calB'$ is which is almost $\mu^2/2$-atomic with respect to $A$.  Because $n$ is sufficiently large, we can choose   a  (possibly empty) set $X\subseteq \F_p^n$ so that $\calL:=\calL'\cup X$ is linearly independent and $ |\calL|= \max\{\ell',f_1\}$.    By Fact \ref{fact:basis}, there exists a purely quadratic factor $\calS=\{M_1,\ldots, M_{n-1}\}$ on $\F_p^n$ of rank at least $n-1$.   Since $n$ is sufficiently large,  Lemma \ref{lem:padding} implies the existence of a subset $\calS'\subseteq \calS$ such that  $\calQ:=\calQ'\cup \calS'$ has size $q'+f_2$ and rank at least $\omega_1(\ell'+q')$.   

Define $\calB=(\calL,\calQ)$. By construction, this factor has complexity $(\ell,q)$ for $\ell=\max\{\ell',f_1\}$ and $q=q'+f_2$ (both bounded by $D$), and rank at least 
$$
\omega_1(\ell'+q')=\omega(\ell'+q'+f_1+f_2)\geq \omega(\ell+q),
$$
where the last inequality is because $\omega$ is increasing.  By Fact \ref{fact:aveofavecor}, $\calB$ is almost $\mu$-atomic with respect to $A$ (since it refines $\calB'$). Thus $\calB$ has all the desired properties.
\end{proofof}

\vspace{2mm}

We can now carry out our strategy to prove Theorem \ref{thm:fop}.
\vspace{2mm}

\begin{proofof}{Theorem \ref{thm:fop}} Fix a prime $p>2$, an integer $k\geq 1$, an increasing function $\psi:\mathbb{N}\rightarrow \mathbb{R}^+$, and a growth function $\omega$.  Let $d=d(k)$ be as in Theorem \ref{thm:hodges}.  Let $\e^*=\e^*(k,p)$ and $\omega^*=\omega^*(k,p)$ be as in Corollary \ref{cor:fewbad}.  Let $\psi':\mathbb{N}\rightarrow \mathbb{R}^+$ be an increasing function satisfying the following:
\begin{itemize}
\item $p^{-\psi'(0)}<\e^*$.
\item $\psi'(x)>\psi(x)+2x+4$ for all $x\in \mathbb{N}$.
\end{itemize}
Let $D_0=D_0(p,\psi',d)$, $N_0=N_0(p,\psi',d)$, $C=C(p,\psi',d)$, and $\mu_0=\mu_0(p,\psi',d)$,  be as in Theorem \ref{thm:almoststableinhstrong}.  Choose any $0<\tau<1/7$ and let $\rho=\rho(\tau)$ be the growth function from Corollary \ref{cor:sizeofatoms}.  Define a new growth function $\omega_1$ by setting $\omega_1(x)=\max\{\omega^*(x),\rho(x),\omega(x)\}$.

Choose any $0<\mu_1<\min\{\e^*, p^{-\psi'(D_0)}\}$.  Since it will be used implicitly in our proof, we point out here that $\mu_1^C\leq \mu_1$ holds because $C\geq 1$ and $\mu_1\in (0,1)$.   Define now $\mu_2=\mu_2(k,\mu_1^C)$  as in Corollary \ref{cor:fewbad}, and set $\mu=\min\{\mu_0, \mu_1,\mu_2\}$.  Set $f_1=1$ and $f_2=\max\{D_0,N_0\}$, and choose $D_1=D_1(p,k, \mu,f_1,f_2, \omega_1)$   as in Theorem \ref{thm:quadlb}.  Finally, assume $n$ is sufficiently large, based on all the previously chosen parameters.

Suppose $A\subseteq \F_p^n$ does not have $k$-$\FOP_2$.  By Lemma \ref{lem:vc2fop}, $\VC_2(A)<k$.  By Theorem \ref{thm:quadlb}, there exist integers $\ell,q$ satisfying $f_1=1\leq \ell\leq D_1$, $f_2\leq q\leq D_1$, and a quadratic factor $\calB=(\calL,\calQ)$ on $\F_p^n$ of complexity $(\ell,q)$ and rank at least $\omega_1(\ell+q)$, such that $\calB$ is almost $\mu$-atomic with respect to $A$.  Let $\Red_{\calB,\mu}(\F_p^n, A)=(\F_p^{\ell+q};\mathbf{A}_0,\mathbf{A}_1,\mathbf{A}_2)$. Recall $\calL_{\calB}=\{\ebar_1,\ldots, \ebar_{\ell}\}$ is the linear factor on $\F_p^{\ell+q}$ consisting of the first $\ell$ standard basis vectors in $\F_p^{\ell+q}$, and $H_{\calB}=L_{\calB}(0)$ is the associated codimension $\ell$ subspace of $\F_p^{\ell+q}$.  

In many of the arguments below, we think of $\F_p^{\ell+q}$ not just as the collection of labels associated to atoms of $\calB$, but as itself a group and an $\F_p$-vector space.   For this reason, we also let $G'=\F_p^{\ell+q}$.  Since $\mu_1^C<\e^*$, since $\mu\leq \mu_2$, and since $\calB$ is $\mu$-atomic with respect to $A$,  Corollary \ref{cor:fewbad} implies the following.
\begin{itemize}
\item There are at most $\mu_1^C|H_{\calB}|^{2^d}|G'|^{2^d-1}$ encodings of $T(d)$ in $(G',\mathbf{A}_1)$ with leaves in $H_{\calB}$, and
\item There is a set $\Omega\subseteq G'/H_{\calB}$ of size at most $\mu_1^Cp^{\ell}$, such that for all $\gbar+H_{\calB}\notin \Omega$, there are at most $\mu_1^C|H_{\calB}|^{2^d+2^d-1}$ many encodings of $T(d)$ in $(G',\mathbf{A}_1)$ with leaves in $H_{\calB}$ and nodes in $H_{\calB}+\gbar$.  
\end{itemize}

Since $\mu_1<\mu_0$ and $\dim(G')-\codim(H_{\calB})=(\ell+q)-\ell=q\geq f_2\geq N_0$, we can now apply Theorem \ref{thm:almoststableinhstrong} to the group $G'$, the subspace $H_{\calB}\leq G'$, and the subset $\mathbf{A}_1\subseteq G'$.  From this we obtain an integer $0\leq m\leq D_0$, a subspace $H'\leq H_{\calB}$ of codimension $m$ in $H_{\calB}$, and a set $\Omega'\subseteq G'/H'$ satisfying
\begin{enumerate}[label=\normalfont(\arabic*)]
\item  $|\Omega'|\leq (\mu_1+p^{-\psi'(m)})p^{m+\ell}$;
\item for every coset $K\in G'/H'$ with $K\notin \Omega'$, $\frac{|\mathbf{A}_1\cap K|}{|K|}\in [0,p^{-\psi'(m)})\cup (1-p^{-\psi'(m)},1]$.   
\end{enumerate}
Since $H'\leq H_{\calB}$, we may fix a linear factor $\calR$ on $\F_p^{\ell+q}$ satisfying $\calR\preceq_{syn}\calL_{\calB}$, such that the cosets of $H'$ in $G'$ are exactly the atoms of $\calR$.   Letting $r$ denote the complexity of $\calR$, we have  $r=\ell+m$.   Observe the atoms of $\calR$ have size $p^{-r}|G'|=p^{\ell+q-r}=p^{q-m}$.  We now define a partition $\At(\calR)=\mathbf{\Omega}_{err}\cup \mathbf{\Omega}_0\cup \mathbf{\Omega}_1$ by setting 
\begin{align*}
\mathbf{\Omega}_{err}&=\{R\in \At(\calR): R\in \Omega'\text{ or }|R\cap \mathbf{A}_2|\geq \sqrt{\mu}|R|\},\\   
\mathbf{\Omega}_1&=\{R\in \At(\calR)\setminus \mathbf{\Omega}_{err}: |\mathbf{A}_1\cap R| \geq (1- \e p^{-\psi'(m)})|R|\},\text{ and }\\
\mathbf{\Omega}_0&=\{R\in \At(\calR)\setminus \mathbf{\Omega}_{err}: |\mathbf{A}_1\cap R| \leq  \e p^{-\psi'(m)}|R|\}. 
\end{align*}
By (2) above, this is indeed a partition of $\At(\calR)$. It will be useful later to know the set $\mathbf{\Omega}_{err}$ is small. To this end, observe that
$$
\mu |\F_p^{\ell+q}|\geq |\mathbf{A}_2|\geq \sqrt{\mu}p^{\ell+q-r}|\{R\in \At(\calR): |\mathbf{A}_2\cap R|\geq \sqrt{\mu}|R|\}|, 
$$
which implies $|\{R\in \At(\calR): |\mathbf{A}_2\cap R|\geq \sqrt{\mu}|R|\}|\leq \sqrt{\mu}p^r$.   Combining this with (1) above, we have 
\begin{align}\label{omegabd}
 |\mathbf{\Omega}_{err}|\leq (\mu_1+p^{-\psi'(m)}+\sqrt{\mu})p^r<3p^{-\psi'(m)}p^r,
\end{align}
where the last inequality is by our choices of $\mu$ and $\mu_1$, and since $m\leq D_0$ and $\psi'$ is increasing.

 By Proposition \ref{prop:pullbackfactor}, there is a quadratic factor $\calB'=(\calL,\calQ')$ of complexity $(\ell,m)$ and rank at least $\rk(\calB)$ so that $\calX_{\calR}=\At(\calB')$.  We now define a distinguished subset of $\At(\calL)$ by letting  
\begin{align*}
\Sigma=\{L\in \At(\calL): \text{ there exists }R\in \mathbf{\Omega}_{err}\text{ such that }X_R\subseteq L\}.
\end{align*}
We show the conclusions of Theorem \ref{thm:fop} hold for this $\calB'$ and $\Sigma$.     We throughout be conflating $\calX_{\calR}$ and $\At(\calB')$ (since they are equal),  using both the notations $X_R$ and $B'(\abar;\bbar)$ for atoms of $\calB'$.
 
 First, we note that by construction, $\ell\leq D_1$ and $m\leq D_0$, so (i) holds.  For (ii), since $\rk(\calB')\geq \rk(\calB)$,  since $m\leq D_0$, and since $q\geq f_2 \geq D_0$, we have
 $$
 \rk(\calB')\geq \rk(\calB)\geq \omega_1(\ell+q)\geq \omega_1(\ell+D_0)\geq \omega_1(\ell+m)\geq \omega(\ell+m),
 $$
as desired.  For (iii), note that by definition,  $|\Sigma|\leq |\mathbf{\Omega}_{err}|$, so by (\ref{omegabd}),
 \begin{align}\label{sigmabound}
|\Sigma|<3p^{-\psi'(m)}p^r=3p^{\ell+m-\psi'(m)}\leq p^{\ell-\psi(m)-1}<p^{-\psi(m)}|\At(\calL)|,
\end{align}
where the last inequality is by our choice of $\psi'$.

We next turn to proving (iv).  Define 
\begin{align*}
\Sigma_0 &=\{X_R: R\in \mathbf{\Omega}_0\text{ and }X_R\cap L=\emptyset \text{ for all }L\in \Sigma\},\text{ and }\\\Sigma_1 &=\{X_R: R\in \mathbf{\Omega}_1\text{ and }X_R\cap L=\emptyset \text{ for all }L\in \Sigma\}.
\end{align*}
Since $\mathbf{\Omega}_{err},\mathbf{\Omega}_0,\mathbf{\Omega}_1$ partition $\At(\calR)$ and $\At(\calB')=\calX_{\calR}$, we have that for  any $\calB'$-atom $B'(\abar;\bbar)=L(\abar)\cap Q'(\bbar)$, we have  $B'(\abar;\bbar)\in \Sigma_0\cup \Sigma_1$ if and only if $L(\abar)\notin \Sigma$.  Therefore, to show (iv), it suffices to show the following.
\begin{align}\label{suffxr}
X_R\in \Sigma_0\Rightarrow |A\cap X_R|\leq p^{-\psi(m)-1}|X_R|\text{ and }X_R\in \Sigma_1\Rightarrow |A\cap X_R|\geq (1-p^{-\psi(m)-1})|X_R|.
\end{align}
  Our verification of this will frequently use the fact that for every  $X_R\in \At(\calB')$,  Corollary \ref{cor:sizeofatoms} implies
\begin{align}\label{align:XL}
|X_{R}|=(1\pm \tau)p^{n-\ell-m}=(1\pm \tau)p^{n-r}.
\end{align}
Fix $X_R\in \Sigma_0$. Using the definition of $X_R$,  we have
  \begin{align*} 
\nonumber|A\cap X_{R}|=\sum_{\dbar\cbar\in R }|A\cap B(\dbar;\cbar)|&=\sum_{\dbar\cbar\in R \cap \mathbf{A}_0}|A\cap B(\dbar;\cbar)|+\sum_{\dbar\cbar\in R\cap (\mathbf{A}_1\cup \mathbf{A}_2)}|B(\dbar;\cbar)|\\
\nonumber&\leq \sum_{\dbar\cbar\in R \cap \mathbf{A}_0}\mu |B(\dbar;\cbar)|+\sum_{\bar\cbar\in R\cap (\mathbf{A}_1\cup \mathbf{A}_2)}|B(\dbar;\cbar)|.
\end{align*}
where second inequality uses the definition of $\mathbf{A}_0$. Combining with Corollary \ref{cor:sizeofatoms}, this implies
\begin{align}\label{112sum}
\nonumber|A\cap X_{R}|&\leq \sum_{\dbar\cbar\in R \cap \mathbf{A}_0}\mu (1+\tau)p^{n-\ell-q} +\sum_{\dbar\cbar\in R\cap (\mathbf{A}_1\cup \mathbf{A}_2)}(1+\tau)p^{n-\ell-q}\\
&\leq \mu |R|(1+\tau)p^{n-\ell-q}+|R\cap (\mathbf{A}_1\cup \mathbf{A}_2)|(1+\tau)p^{n-\ell-q}.
\end{align}
Since $R\in \mathbf{\Omega}_0$,  $|\mathbf{A}_1\cap R|\leq p^{-\psi'(m)}|R|$, and since $R\notin \mathbf{\Omega}_{err}$,  $|\mathbf{A}_{2}\cap R|\leq \sqrt{\mu}|R|$. Therefore,
\begin{align}\label{ra1a2}
|R\cap (\mathbf{A}_1\cup \mathbf{A}_2)|\leq (p^{-\psi'(m)}+\sqrt{\mu})|R|.
\end{align}
Combining (\ref{ra1a2}) with (\ref{112sum}) and the fact that $|R|=p^{\ell+q-r}$, we have 
\begin{align*}
|A\cap X_{R}|\leq (\mu+p^{-\psi'(m)}+\sqrt{\mu})|R|(1+\tau)p^{n-\ell-q}&= (\mu+p^{-\psi'(m)}+\sqrt{\mu})(1+\tau)p^{n-r}\\
&\leq \frac{(\mu+p^{-\psi'(m)}+\sqrt{\mu})(1+\tau)}{(1-\tau)}|X_R|\\
&<p^{-\psi(m)-1}|X_R|,
\end{align*}
where the second inequality uses (\ref{align:XL}), and the final inequality uses our choices of $\mu$, $\tau$, and $\psi'$.

Now fix $X_R\in \Sigma_1$.   By definition of $X_R$ and $\mathbf{A}_1$, we have
\begin{align}\label{ra1}
\nonumber|A\cap X_{R}|\geq \sum_{\cbar\dbar\in R\cap \mathbf{A}_1}|A\cap B(\cbar;\dbar)|&\geq  \sum_{\cbar\dbar\in R\cap \mathbf{A}_1}(1- \mu)|B(\abar;\bbar)|\\
\nonumber&\geq \sum_{\cbar\dbar\in R\cap \mathbf{A}_1}(1- \mu)(1-\tau)p^{n-\ell-q}\\
&=|R\cap \mathbf{A}_1|(1- \mu)(1-\tau)p^{n-\ell-q},
\end{align}
where the third inequality is by Corollary \ref{cor:sizeofatoms}.   Since $R\in \mathbf{\Omega}_1$, $|\mathbf{A}_1\cap R|\geq (1-p^{-\psi'(m)})|R|$.  Combining this with (\ref{ra1}), we have
\begin{align*}
|A\cap X_R|\geq (1- p^{-\psi'(m)})(1-\mu)(1-\tau)p^{n-\ell-q}|R|=(1- p^{-\psi'(m)})(1-\mu)(1-\tau)p^{n-r},
\end{align*}
where the last inequality follows from $|R|=p^{\ell+q-r}$.  Combining with (\ref{align:XL}), we conclude
$$
|A\cap X_R|\geq \frac{(1- p^{-\psi'(m)})(1-\mu)(1-\tau)}{(1+\tau)}|X_R|\geq (1-p^{-\psi(m)-1})|X_R|,
$$
where the final inequality is by and our choices of $\mu$, $\psi'$, and $\tau$. This finishes the verification of (\ref{suffxr}), and thus of part (iv).    

Finally, we verify (v). Set $Y=\bigcup_{B\in \Sigma_1}B$.  Using the definition of $Y$ and (\ref{suffxr}), we have
\begin{align*}
|A\Delta Y|&\leq \sum_{B\in \Sigma_1}|B\setminus A|+\sum_{B\in \Sigma_0}|B\cap A|+\sum_{L\in \Sigma}|L|\\
&\leq \sum_{B\in \Sigma_1}p^{-\psi(m)-1}|B|+ \sum_{B\in \Sigma_0}p^{-\psi(m)-1}|B|+|\Sigma|p^{n-\ell}\\
&\leq  \sum_{B\in \Sigma_1\cup \Sigma_0}p^{-\psi(m)-1}p^{n-\ell-m}(1+\tau)+|\Sigma|p^{n-\ell}\\
& \leq 2|\Sigma_1\cup \Sigma_0|p^{-\psi(m)-1}p^{n-\ell-m}+|\Sigma|p^{n-\ell},
\end{align*}
where the third inequality uses Corollary \ref{cor:sizeofatoms} and the last uses $\tau<1$.  Using (\ref{sigmabound}) and the fact $|\Sigma_0\cup \Sigma_1|\leq p^{\ell+m}$, we can conclude that
$$
|A\Delta Y|\leq 2p^{-\psi(m)-1}p^{n}+p^{n-\psi(m)-1}=3p^{-\psi(m)-1}|\F_p^n|\leq p^{-\psi(m)}|\F_p^n|,
$$
where the last inequality uses that $p>2$.  This finishes the proof.
\end{proofof}

\subsection{A preliminary generalized stable regularity lemma}\label{sec:stablereg}

The remainder of this section is devoted to proving Theorem \ref{thm:almoststableinhstrong}.  We will first set out some preliminaries, after which we will prove a weak version of Theorem \ref{thm:almoststableinhstrong} (see Theorem \ref{thm:stableinh} below).  We subsequently iterate the weak version (Theorem \ref{thm:stableinh}) to obtain the strong version, Theorem \ref{thm:almoststableinhstrong}. 

Our first goal is to prove Lemma \ref{fact:densesubspace}, which says that if a set $A\subseteq \F_p^n$ has non-trivial density and admits few encodings of $T(d)$ with one side in a subspace $H$, then there is a subspace $H'\leqslant H$ of bounded codimension in $H$ so that $A$ has density near $1$  on some coset of $H'$ in $H$.  We will  use the following fact, which shows we can efficiently find uniform subspaces of small codimension (here ``uniform" is in the Fourier sense, first defined in \cite{Green.2005} and used in \cite{Terry.2019}). This can be proved via a standard density increment argument (see e.g. Proposition 2 in \cite{Terry.2019}). 

\begin{fact}\label{lem:lemma7}
For all $\e\in (0,1)$ there is $N=N(\e)$ so that the following holds.  Suppose $n\geq N$ and $A\subseteq G=\F_p^n$.  Then there are $H\leqslant G$ of codimension at most $\lfloor 2/\e\rfloor$ in $G$ and $y\in G$ such that $|A\cap (H+y)|/|H|\geq |A|/|G|$, and $H+y$ is $\e$-uniform with respect to $A$.
\end{fact}

We will not define uniformity here, but simply note that it can be rephrased in terms of the local $U^2$-norm in the following sense (see \cite[Lemma 3.4]{Terry.2021d}).

\begin{observation}\label{ob:uniformu2}
Suppose $A\subseteq \F_p^n$,  $H \leqslant \F_p^n$, $\calL$ is a linear factor on $\F_p^n$ of complexity $\ell$ with zero atom $H=L(0)$, and $H+y$ is a coset of $H$, equal to the linear atom $L(a)$.  Then $H+y$ is $\e$-uniform with respect to $A$ if and only if $\|1_A-\alpha_y\|_{U^2(d)}<\e$, where $\alpha_y$ is the density of $A$ on $H+y$, and $d$ is the label $(0,a)\in \F_p^{2\ell}$.
\end{observation}

We can now prove Lemma \ref{fact:densesubspace}.  We will use in the proof below, and throughout this section, that for all $d\geq 1$, $|2^{<d}|=2^d-1$. 

\begin{lemma}[Finding a dense subspace]\label{fact:densesubspace}
For each $d\geq 1$ and $\e\in (0,1)$, there are $N=N(\e,d)$, $m=m(d)\geq 1$ and $\mu=\mu(d,\e)>0$ with the following property.  Let $n\geq N$, let $A\subseteq G=\F_p^n$, and assume $|A|\geq \e |G|$.   Suppose that there are at most $\mu|G|^{2^d+2^d-1}$ encodings of $T(d)$ in $(G, A)$.  Then there is a subspace $H\leqslant G$ of codimension at most $\e^{-m}$   and $y\in G$ such that $|(A-y)\cap H|\geq (1-\e)|H|$.  
\end{lemma}
\begin{proof}
Without loss of generality, assume that $\e<1/2$.  Let $D_1=8(2^d+2^d-1)^2$.  Then set  $m=3D_1$, $\mu=p^{-\e^{-m}2^{d+1}}$, and $N=N(\e^{2D_1})$.  Suppose $n\geq N$, $A\subseteq G=\F_p^n$ satisfies $|A|\geq \e |G|$, and assume that there are at most $\mu |G|^{2^d+2^d-1}$ encodings of $T(d)$ in $(G,A)$.  By Fact \ref{lem:lemma7},  there are $H\leqslant G$ and $y\in G$ so that $H+y$ is $\e^{2D_1}$-uniform with respect to $A$ in $G$, so that $H$ has codimension at most $\lceil 2\e^{-2D_1}\rceil <\e^{-m}$, and so that $\rho:=|A\cap (H+y)|/|H+y|\geq \e |H|$.  

If $\rho\geq 1-\e$, then we are done, so we may assume that $\rho\in (\e,1-\e)$. Let $\calL$ be a linear factor in $H$ with $L(0)=H$ and let $a$ be such that $H+y=L(a)$.  Let $F=(U\cup V, E)$ be the bipartite graph where $U=2^{<d}$, $V=2^d$, and $E=\{\sigma \eta: \sigma \wedge 1\trianglelefteq \eta\}$.  We remind the reader that $|U|=2^d-1$ while $|V|=2^d$.  For each $\sigma\in 2^{<d}$, let $a_{\sigma}=a$ and for each $\eta\in 2^d$, let $b_{\eta}=0$.  We now have that for all $\sigma\in 2^{<d}$ and $\eta\in 2^d$, if $\alpha_{\sigma\eta}$ is the density of $A$ on $L(a_{\sigma}+b_{\eta})$, then $\e<\alpha_{\sigma\eta}<1-\e$, and moreover, by Observation \ref{ob:uniformu2},   $\|1_{A}-\alpha_{\sigma\eta}\|_{U^2(a_{\sigma}b_{\eta})}<\e^{2D_1}$.   Thus, by Proposition \ref{prop:countinggplinear}, the number of tuples $\xbar\ybar=(x_{\sigma})_{\sigma\in 2^{<d}}(y_{\eta})_{\eta\in 2^d}\in \prod_{\sigma\in 2^{<d}}L(a_{\sigma})\prod_{\eta\in 2^d}L(b_{\eta})$ satisfying $x_{\sigma}+y_{\eta}\in A$ if and only if $\sigma\eta\in E(F)$, is at least the following.
\begin{align*}
&\e^{2^d+2^d-1}\prod_{\sigma\in 2^{<d}}L(a_{\sigma})\prod_{\eta\in 2^d}L(b_{\eta}) - (2^d+2^d-1)^2\e^{2D_1}\prod_{\sigma\in 2^{<d}}L(a_{\sigma})\prod_{\eta\in 2^d}L(b_{\eta})\\
&=\e^{2^d+2^d-1}(1-(2^d+2^d-1)^2\e^{2D_1-2^d-2^d+1})|H|^{2^d+2^d-1}.
\end{align*}
By our choice of $D_1$ and because $\e<1/2$, this is at least  $\frac{1}{2}\e^{2^d+2^d-1}|H|^{2^d+2^d-1}$. But observe that
 $$
 \frac{1}{2}\e^{2^d+2^d-1}|H|^{2^d+2^d-1}\geq \frac{1}{2}\e^{2^d+2^d-1}p^{-\e^{-m}(2^d+2^d-1)}|G|^{2^d+2^d-1}>\mu |G|^{2^d+2^d-1},
 $$
 where the last inequality holds by definition of $\mu$ and since $\e\leq 1/2$.  Since all such tuples yield encodings of $T(d)$ in $(G, A)$, we have arrived at a contradiction. 
\end{proof}

We are now ready to prove Theorem \ref{thm:stableinh},  the base step for the proof of Theorem \ref{thm:almoststableinhstrong}.  Informally, Theorem \ref{thm:stableinh} tells us the following. Suppose $\e\in (0,1)$, $d\geq 1$, $A\subseteq G=\F_p^n$, and $H\leqslant \F_p^n$ satisfy the following two criteria: first, that there are few encodings of $T(d)$ in $(G,A)$ with leaves in $H$, and second, that for almost all cosets $H+g$,  there are few encodings of $T(d)$ in $(G,A)$ with leaves in $H$ and nodes in $H+g$.  Then there is a subspace $H'$ of $H$ which is almost $\e$-atomic with respect to $A$ and whose codimension in $H$ \emph{depends only on $\e$ and $d$}.  The form of this dependence is crucial for later applications.

\begin{theorem}\label{thm:stableinh}
For all $d\geq 1$, there is $\e^*>0$ such that for all $0<\e\leq \e^*$, there exist $\mu=\mu(d,\e)$, $L=L(d,\e)$, and $N=N(d,\e)$ such that the following holds.  Suppose that $A\subseteq G=\F_p^n$, $H\leqslant G$, and $\Omega\subseteq G/H$ satisfy
\begin{enumerate}[label=\normalfont(\roman*)]
\item $|H|\geq p^N$;
\item the number of encodings of $T(d)$ in $(G,A)$ with leaves in $H$ is at most $\mu  |G|^{2^{d}-1}|H|^{2^d}$;
\item $|\Omega|\leq \e |G|/2|H|$ and for all $K\in (G/H)\setminus \Omega$, the number of encodings of $T(d)$ in $(G,A)$ with leaves in $H$ and nodes in $K$ is at most $\mu |H|^{2^d+2^d-1}$.
\end{enumerate}
Then there exists $H'\leqslant H$, of codimension at most $L$ in $H$ such that $H'$ is almost $\e$-atomic with respect to $A$. 
\end{theorem}
\begin{proof}
We proceed by induction on $d$.  Suppose first that $d=1$.  Let $\e^*=\e^*(1)=1/2$.  Given  $0<\e\leq \e^*$, let $L=L(1,\e)=1$, $\mu=\mu(1,\e)=\e^3$, and  $N=N(1,\e)=\e^{-1}$.  Suppose that $A\subseteq  G=\F_p^n$, $H\leqslant \F_p^n$, and $\Omega\subseteq G/H$ satisfy the following.
\begin{enumerate}
\item  $|\Omega|\leq \e |G|/2|H|$ and $|H|\geq p^N$;
\item the number of encodings of $T(1)$ in $(G,A)$ with leaves in $H$ is at most $\mu|G||H|^2$;
\item for all $g+H\notin \Omega$, the number of encodings of $T(1)$ in $(G,A)$ with leaves in $H$ and nodes in $g+H$ is at most $\mu|H|^{3}$. 
\end{enumerate}
  Let $Y=\{K\in G/H: \min\{|A\cap K|, |(\neg A)\cap K|\}\geq \e |H|\}$.  We shall show that $Y\subseteq \Omega$.  Fix $K\in Y$.   We construct many encodings of $T(1)$ in $(G,A)$ with leaves in $H$ and nodes in $K$ as follows.
  \begin{itemize}
\item Choose $(g',g'')\in (A\cap K)\times ((\neg A)\cap K)$.  There are at least $\e^2 |H|^2$ choices since $K\in Y$.
  \item Choose $g\in K$ and set $(a_{<>},b_0,b_1):=(g,g'-g,g''-g)$.  Note that this triple is an encoding of $T(1)$ in $(G,A)$ with leaves in $H$ and nodes in $K$.  There are $|K|=|H|$ choices for this step. 
  \end{itemize}
 Thus  the number of encodings of $T(1)$ in $(G,A)$ with leaves in $H$ and nodes in $K$ is at least $\e^2|H|^3>\mu |H|^3$, so by (iii), we must have $K\in \Omega$.  We have shown that $Y\subseteq \Omega$, so $|Y|\leq |\Omega|\leq \e |G|/2|H|$.  Taking $H'=H$, this finishes case $d=1$.

Suppose now that $d>1$, and assume by induction that the conclusion of Theorem \ref{thm:stableinh} holds for $d-1$.  Define $\e'=\e^*(d-1)$ as in the inductive hypothesis, and set $\e^*=\e'/2^d$. Let $m=m(d)$ be as in Lemma \ref{fact:densesubspace}. Fix $0<\e\leq \e^*$. Let $\mu'=\mu'(d-1,\e)$ be from the inductive hypothesis, and set $\rho=\min\{\e,((\mu')^2-(\mu')^3)/4^d\}$.  Let $\mu_0=\mu_0(d,\rho)$ be as in Lemma \ref{fact:densesubspace} (recall from the proof of Lemma \ref{fact:densesubspace} that this means $\mu_0=p^{-\rho^{-m}2^{d+1}}$).  Then take $L'=L'(d-1,\rho)$ from the inductive hypothesis, and set  $L=2\rho^{-m}+L'$.  Let $\mu=\e(\mu')^{10}p^{-\rho^{-m}2^{d+2}}/4$.  Note that $\mu< \mu_0$.  Choose $N_1=N_1(d,\e)$ from Lemma \ref{fact:densesubspace}, and let $N\geq N_1$ be sufficiently large compared to $\mu^{-1}$.

Suppose now that $A\subseteq \F_p^n$, $H\leqslant \F_p^n$, and $\Omega\subseteq G/H$ satisfy 
\begin{enumerate}
\item $|\Omega|\leq \e |G|/2|H|$ and $|H|\geq p^N$;
\item the number of encodings of $T(d)$ in $(G,A)$ with leaves in $H$ is at most $\mu  |G|^{2^{d}-1}|H|^{2^d}$;
\item for all $g+H\notin \Omega$, the number of encodings of $T(d)$ in $(G,A)$ with leaves in $H$ and nodes in $g+H$ is at most $\mu |H|^{2^d+2^d-1}$.
\end{enumerate}
Define $Y=\{K\in G/H: \min\{|A\cap K|, |(\neg A)\cap K|\}\geq \e |H|\}$.  If $|Y|\leq \e |G|/|H|$, then taking $H'=H$ we are done.  So assume that $|Y|> \e |G|/|H|$.   Setting $Y'=Y\setminus \Omega$, we have that $|Y'|\geq \e |G|/2|H|$. 

Fix any $K\in Y'$ and $g\in K$.  Since $g+H=K\in Y$, both $(A-g)\cap H$ and $(\neg A-g)\cap H$ have size at least $\e|H|\geq \rho|H|$.  On the other hand, $K\notin \Omega$ implies that there are at most $\mu|H|^{2^d+2^d-1}<\mu_0|H|^{2^d+2^d-1}$ encodings of $T(d)$ in $(H,(A-g)\cap H)$.  Therefore, we may apply Lemma \ref{fact:densesubspace} in the group $H$ to both the sets $(A-g)\cap H$ and $(\neg A-g)\cap H=(\neg (A-g))\cap H$ with parameter $\rho>0$.  This yields subspaces $H_g(0),H_g(1)\leqslant H$, each of codimension at most $\rho^{-m}$ in $H$, and $y^g_0,y^g_1\in H$, such that $|(A-g-y^g_1)\cap H_1(g)|\geq (1-\rho)|H_1(g)|$ and  $|(\neg A-g-y^g_0)\cap H_0(g)|\geq (1-\rho)|H_0(g)|$.  Set $H_g=H_0(g)\cap H_1(g)$.  Note that $H_g$ has codimension at most $2\rho^{-m}$ in $H$.  By averaging, there are $z^g_1,z^g_0\in H$ such that for each $u\in \{0,1\}$, $|(A^u-g-y^g_u-z^g_u)\cap H_g|\geq (1-\rho)|H_g|$.  Setting $x_0^g=y_0^g+z_0^g$ and $x_1^g=y_1^g+z_1^g$, we have that
\begin{align}\label{align:stableinh1}
|(A-g-x^g_1)\cap (\neg A-g-x^g_0)\cap H_g|\geq (1-2\rho)|H_g|.
\end{align}

Define $\Omega_g$ to be the set of $K'\in G/H_g$ such that  there are at least $\mu'|H_g|^{2^{d-1}+2^{d-1}-1}$ encodings of $T(d-1)$ in $(G,A)$ with leaves in $H_g$ and nodes in $K'$. By the induction hypothesis and the definition of $\Omega_g$, one of the following holds.
\begin{enumerate}
\item[(a)] There is $H_g'\leqslant H_g$ of codimension at most $L'$ in $H_g$ which is almost $\e$-atomic with respect to $A$;
\item[(b)] $|\Omega_g|>\e|G|/2|H_g|$;
\item[(c)] The number of encodings of $T(d-1)$ in $(G,A)$ with leaves in $H_g$ is at least $\mu'|G|^{2^{d-1}-1}|H_g|^{2^{d-1}}$.
\end{enumerate}

We shall show that, in fact, (a) must hold.  First we eliminate (b) as a possibility.  In fact, we prove the following stronger statement.

\begin{claim}\label{cl:omega}
$\Omega_g=\emptyset$.
\end{claim}
\begin{proof}
Suppose towards a contradiction that there exists some $K'\in \Omega_g$.  Let $R$ be the set of encodings of $T(d-1)$ in $(G,A)$ with leaves in $H_g$ and nodes in $K'$.   Set 
\begin{align*}
R_1=&\{(a_{\eta})_{\eta\in 2^{d-1}}\in (H_g)^{2^{d-1}}:\\
&|\{ (b_{\sigma})_{ \sigma\in 2^{<d-1}}\in (K')^{2^{d-1}-1}: (a_{\eta})_{\eta\in 2^{d-1}}(b_{\sigma})_{ \sigma\in 2^{<d-1}}\in R\}|\geq (\mu')^2|H_g|^{2^{d-1}-1}\}.
\end{align*}
Note that 
$$
\mu'|H_g|^{2^{d-1}+2^{d-1}-1}\leq |R|\leq |R_1||H_g|^{2^{d-1}-1}+(\mu')^2|H_g|^{2^{d-1}-1}|(H_g)^{2^{d-1}}\setminus R_1|,
$$
where the first inequality holds by assumption, and the second holds by definition of $R_1$.  Rearranging, this yields 
$$
 \frac{\mu'-(\mu')^2}{1-(\mu')^2}|H_g|^{2^{d-1}}\leq |R_1|.
$$
Consequently, $|R_1|\geq (\mu')^2|H_g|^{2^{d-1}}$.  Now define
$$
W_g=\{\xbar\in (H_g)^{2^{d-1}}: \{x_1,\ldots, x_{2^{d-1}}\}\setminus ((\neg A-g-x_0^g)\cap (A-g-x_1^g))\neq \emptyset\}.
$$ 
By (\ref{align:stableinh1}), $|W_g|\leq 2^{d-1}2\rho|H_g|^{2^{d-1}}=\rho2^{d}|H_g|^{2^{d-1}}$.  Thus if $W_g':=(H_g)^{2^{d-1}}\setminus W_g$, then $|W_g'|\geq |H_g|^{2^{d-1}}(1-\rho2^d)$.   Consequently,
$$
|R_1\cap W_g'|\geq |H_g|^{2^{d-1}}((\mu')^2- \rho2^d)>(\mu')^3 |H_g|^{2^{d-1}},
$$
where the last inequality holds by definition of $\rho$.  We can now construct many encodings of $T(d)$ in $(G,A)$ with leaves in $H$ and nodes in $K$ as follows.
\begin{itemize}
\item Choose $(a^0_{\eta})_{\eta\in 2^{d-1}}$ and $(a^1_{\eta})_{\eta\in 2^{d-1}}$ in $R_1\cap W_g'$.  There are $|R_1\cap W_g'|^2\geq (\mu')^6 |H_g|^{2^{d}}$ choices.
\item Choose $(b^0_{\sigma})_{ \sigma\in 2^{<d-1}}\in (H_g)^{2^{d-1}}$ and $(b^1_{\sigma})_{ \sigma\in 2^{<d-1}}\in (H_g)^{2^{d-1}}$ with $(a^0_{\eta})_{\eta\in 2^{d-1}}(b^0_{\sigma})_{ \sigma\in 2^{<d-1}}\in R$ and $(a^1_{\eta})_{\eta\in 2^{d-1}}(b^1_{\sigma})_{ \sigma\in 2^{<d-1}}\in R$.  Since $(a^0_{\eta})_{\eta\in 2^{d-1}}$ and $(a^1_{\eta})_{\eta\in 2^{d-1}}$ are in $R_1$, there are at least $(\mu')^4|H_g|^{2^{d}-1}$ choices.
\item Set $c_{<>}=g$, and given $\eta\in 2^{d-1}$, $\sigma\in 2^{<d-1}$, and $u\in \{0,1\}$, set $e_{u\wedge \eta}=x^g_u+a_{\eta}^u$, and $c_{u\wedge \sigma}=b_{\sigma}^u-x_u^g$.  There is one way to do this.
\end{itemize}
We claim that every tuple $(c_{\alpha})_{\alpha\in 2^{<d}}(e_{\lambda})_{\lambda\in 2^d}$ constructed above is an encoding of $T(d)$ in $(G,A)$ with leaves in $H$ and nodes in $K$.  By construction, each $c_{\alpha}\in H$ and each $e_{\lambda}\in H+g=K$.

Now fix $\alpha\in 2^{<d}$ and $\lambda\in 2^d$ with $\mu\triangleleft \lambda$.    Let $u\in \{0,1\}$ and  $\eta\in 2^{d-1}$ be such that $\lambda=u\wedge\eta$.  Suppose first $\alpha=<>$, then
$$
c_{\alpha}+e_{\lambda}=g+e_{u\wedge\eta}=g+x_u^g+a_{\eta}^u.
$$
By construction, $a_{\eta}^u\in W_g'\subseteq A^u-g-x^g_u$, so $g+x_u^g+a_{\eta}^u\in A^u$.  Thus $c_{\alpha}+e_{\lambda}\in A$ if and only if $\alpha\wedge 1\trianglelefteq \lambda$.

Suppose now that $\alpha\neq <>$.  Then $\alpha\triangleleft \lambda$ implies there is a unique $\sigma\in 2^{<d-1}$ such that $\sigma\triangleleft \eta$ and $\alpha=u\wedge\sigma$.  Note that for each $v\in \{0,1\}$, $\sigma\wedge v\trianglelefteq \eta$ if and only if $\alpha\wedge v\trianglelefteq \lambda$.   In this case,
$$
c_{\alpha}+e_{\lambda}=c_{u\wedge \sigma}+e_{u\wedge \eta}=b_{\sigma}^u+a_{\eta}^u.
$$
By assumption, $b_{\sigma}^u+a_{\eta}^u\in A$ if and only if $\sigma\wedge 1\trianglelefteq \eta$.  By above, this holds if and only if $\alpha\wedge 1\trianglelefteq \lambda$.  Thus each $(c_{\alpha})_{\alpha\in 2^{<d}}(e_{\lambda})_{\lambda\in 2^d}$ constructed above is an encoding of $T(d)$ in $(G,A)$ with leaves in $H$ and nodes in $K$, so we have constructed at least 
\[(\mu')^{6} |H_g|^{2^{d}}(\mu')^{4}|H_g|^{2^{d}-1}=(\mu')^{10}|H_g|^{2^d+2^d-1}\geq (\mu')^{10}(p^{-2\rho^{-m}})^{2^d+2^d-1}|H|^{2^d+2^d-1}>\mu|H|^{2^d+2^d-1}\]
such encodings, all with $c_{<>}=g$ (here the first inequality holds because $H_g$ has codimension at most $2\rho^{-m}$ in $H$, and the last inequality holds by definition of $\mu$).  However, this contradicts the assumption that $K\in \Omega$.
\end{proof}

Thus, by Claim \ref{cl:omega}, (b) cannot hold.  We now show that (c) cannot hold. Suppose towards a contradiction that (c) holds.  We will show that there are too many encodings of $T(d)$ in $(G,A)$ with leaves in $H$ and nodes in $K$.  Let $S\subseteq G^{2^{d-1}-1}H^{2^{d-1}}$ be the set of encodings of $T(d-1)$ in $(G,A)$ with leaves in $H_g$.  Then define
\begin{align*}
S_1=&\{(a_{\eta})_{\eta\in 2^{d-1}}\in (H_g)^{2^{d-1}}:\\
&|\{ (b_{\sigma})_{ \sigma\in 2^{<d-1}}\in G^{2^{d-1}-1}: 
(a_{\eta})_{\eta\in 2^{d-1}}(b_{\sigma})_{ \sigma\in 2^{<d-1}}\in S\}|\geq (\mu')^2 |G|^{2^{d-1}-1}\}.
\end{align*}
Observe that 
$$
|S|\leq |S_1||G|^{2^{d-1}-1}+(\mu')^2 |G|^{2^{d-1}-1}|(H_g)^{2^{d-1}}\setminus S_1|.
$$
By (a), we have that $|S|\geq  \mu'|G|^{2^{d-1}-1}|H_g|^{2^{d-1}}$, and therefore the equation above implies that $|S_1|\geq (\mu')^2|H_g|^{2^{d-1}}$.  Recall that
$$
W'_g=\{\xbar\in (H_g)^{2^{d-1}}: \{x_1,\ldots, x_{2^{d-1}}\}\subseteq ((\neg A-g-y_0-z_0^g)\cap (A-g-y_1-z_1^g))\},
$$ 
and that we showed in the proof of Claim \ref{cl:omega} that $|W_g'|\geq |H_g|^{2^{d-1}}(1-\rho2^{d})$.  Combining our lower bounds on $S_1$ and $W_g'$, we have that
$$
|W_g'\cap S_1|\geq |H_g|^{2^{d-1}}((\mu')^2-\rho 2^d) \geq |H_g|^{2^{d-1}} (\mu')^3,
$$
where the second inequality holds by our choice of $\rho$.  We now give a procedure for generating many encodings of $T(d)$ in $(G,A)$ with leaves in $H$ and nodes in $K$, where $c_{<>}=g$.
\begin{itemize}
\item Fix any two elements, $\abar^1=(a^1_{\eta}: \eta\in 2^{d-1}), \abar^0=(a^0_{\eta}: \eta\in 2^{d-1})\in W_g'\cap S_1$.  There are at least $|W_g'\cap S_1|^2\geq  |H_g|^{2^{d}} (\mu')^6$ choices.
\item  Choose $\bbar^1=(b^1_{\sigma}:\sigma\in 2^{<d-1})$ and $\bbar^0=(b^0_{\sigma}:\sigma\in 2^{<d-1})$ such that $\bbar^1\abar^1$ and $\bbar^0\abar^0\in S$. Since $\abar^1,\abar^0\in S_1$, there are at least $(\mu' |G|^{2^{d-1}-1})^2\geq (\mu')^4 |G|^{2^{d}-2}$ choices.
\item  Set $c_{<>}=g$, and for each $\eta\in 2^{d-1}$, $\sigma\in 2^{<d-1}$ and $u\in \{0,1\}$, set $e_{u\wedge \eta}=x_u^g+a^u_{\eta}$ and $c_{u\wedge \sigma}=b_{\sigma}^u-x^g_u$.  
\end{itemize}
We now check that every $(c_{\mu})_{\mu\in 2^d}(e_{\lambda})_{\lambda\in 2^d}$ is an encoding of $T(d)$ with leaves in $H$ and nodes in $K$.  First, by construction, for each $\eta\in 2^d$, $e_{\eta}\in H$, and for each $\sigma\in 2^{<d}$, $c_{\sigma}\in K$.  Now fix $\eta\in 2^d$ and $\sigma\trianglelefteq \eta$.  There is some $u\in \{0,1\}$ and $\rho\in 2^{d-1}$ so that $\eta=u\wedge \rho$, and some $w\in \{0,1\}$ so that $\sigma \wedge w\trianglelefteq \eta$.  If $\sigma=<>$, then we must have $u=w$, and in this case,
$$
c_{\sigma}+e_{\eta}=c_{<>}+e_{u\wedge\rho}=g+x^g_u+a^u_{\rho}\in A^u.
$$
Assume now that $\sigma \neq <>$.  Since $\sigma\trianglelefteq \eta$ and $\eta=u\wedge \rho$, there is $\lambda\in 2^{<d-1}$ such that $\lambda\trianglelefteq \rho$, and $\sigma=u\wedge \lambda$.  Note we also have $\lambda\wedge w\trianglelefteq \rho$.  Then
$$
c_{\sigma}+e_{\eta}=c_{u\wedge\lambda}+e_{u\wedge\rho}= x^g_u+a^u_{\rho}+b^u_{\lambda}-x^g_u=a^u_{\rho}+b^u_{\lambda}\in A^w,
$$
where the inclusion in $A^w$ holds by construction.   So if $\sigma\wedge 1\trianglelefteq \eta$, then $c_{\sigma}+e_{\eta}\in A$ and if $\sigma\wedge 0\trianglelefteq \eta$, then $c_{\sigma}+e_{\eta}\notin A$.  This shows that $(c_{\sigma})_{\sigma\in 2^{<d}}(e_{\eta})_{\eta\in 2^d}$ is an encoding of $T(d)$ in $(G,A)$ with leaves in $H$ and nodes in $K$, and by construction $c_{<>}=g$.

Consequently, we can construct at least 
\[|H_g|^{2^{d}} (\mu')^{10} |G|^{2^{d}-2}\geq (\mu')^{10}(p^{-2\rho^{-m}})^{2^d}|H|^{2^{d}}|G|^{2^{d}-2}\]
encodings of $T(d)$ in $(G,A)$ with leaves in $H$ and nodes in $K$ with  $c_{<>}=g$. Clearly the encodings produced are distinct for difference choices of $K\in Y'$ and $g\in K$.  Consequently, the total number encodings of $T(d)$ in $(G,A)$ with leaves in $H$ is at least 
\begin{align*}
|H| (\e|G|/2|H|)(\mu')^{10}(p^{-2\rho^{-m}})^{2^d}|H|^{2^{d}}|G|^{2^{d}-2}&=(\e/2) (\mu')^{10}(p^{-2\rho^{-m}})^{2^d}|H|^{2^{d}}|G|^{2^{d}-1}\\
&> \mu|H|^{2^d}|G|^{2^d-1},
\end{align*}
where the final inequality follows from the definition of $\mu$ and because $H\geq N$. But this contradicts our original assumption about the number of such encodings. 

Thus (c) cannot hold.  We are left with the conclusion that (a) holds.  Consequently, there is some $H_g'\leqslant H$ of codimension at most $L'$ in $H_g$ which is almost $\e$-atomic with respect to $A$.  In this case we are done, since the codimension of $H_g'$ in $H$ is at most $2\rho^{-m}+L'\leq L$.  
\end{proof}

\subsection{Iterating to obtain a stronger version}\label{subsec:iterate}

In this section we prove Theorem \ref{thm:almoststableinhstrong}, which is a stronger version of Theorem \ref{thm:stableinh}, in which the error on the ``good'' atoms can be made arbitrarily small with respect to the codimension of $H'$ in $H$.  We will do this by iteratively applying Theorem \ref{thm:stableinh} to produce a ``factor chain'', as defined below.

\begin{definition}[Linear factor chain]\label{def:chainlin}
Let $\ell_0,T,D\geq 1$ be integers, and let $f,g:\mathbb{N}\rightarrow\mathbb{N}$ be functions.  Given $A\subseteq \mathbb{F}_p^n$, an \emph{$(T,\ell_0,D,f,g)$-linear factor chain for $A$ in $\mathbb{F}_p^n$} consists of 
\begin{itemize}
\item for each $0\leq i\leq T$, a linear factor $\calL_i$ of some complexity $\ell_i$,
\item for each $1\leq i\leq T$, a partition $\At(\calL_i)=\Gamma_i^0\cup \Gamma_i^1\cup \Gamma^{err}_i$,
\end{itemize}
satisfying
\begin{enumerate}[label=\normalfont(\arabic*)]
\item $\calL_T\prec_{syn}\ldots \prec_{syn}\calL_1\prec_{syn}\calL_0$,
\item for each $1\leq i\leq T$, $f(\ell_{i-1}-\ell_0) \leq \ell_i-\ell_{i-1} \leq D$,
\item  for each $1\leq i\leq T$, $u\in \{0,1\}$, and $L\in \Gamma_i^u$, $|A^u\cap L|/|L|\geq 1- p^{-g(\ell_{i-1}-\ell_0)}$.
\item for each $1\leq i\leq T$ and $L\in \At(\calL_{i-1})$, $|\{L'\in \Gamma^{err}_i: L'\subseteq L\}|\leq   p^{-g(\ell_{i-1}-\ell_0)}p^{\ell_i-\ell_{i-1}}$.
\end{enumerate}
\end{definition}

Note that in a linear factor chain, one has some control over the growth of the complexities from one stage to another, and the error parameter at step $i$ is allowed to depend on the growth in complexity from previous steps.   The goal of Lemma \ref{lem:chainlinearnew} below is to show that linear factor chains can be built to refine a given factor $\calL$ whenever hypotheses similar to the ones in Theorem \ref{thm:stableinh} are satisfied.

\begin{lemma}\label{lem:chainlinearnew}
For all $d\geq 1$, $T\geq 1$, and increasing functions $f,g:\mathbb{N}\rightarrow \mathbb{N}$, there are reals $\e'=\e'(d,T,f,g)$ and  $\mu=\mu(d, T,f,g)$ in $(0,1)$ and positive integers  $D=D(d, T,f,g)$ and $N=N(d, T,f,g)$ so that the following holds.

Let $\ell\geq 0$, let $n-\ell\geq N$,  and let $A\subseteq G=\F_p^n$. Suppose that  $\calL$ is a linear factor of complexity $\ell$ on $\F_p^n$, and that $H=L(0)$ and $\Omega\subseteq \At(\calL)$ satisfy the following.
\begin{enumerate}[label=\normalfont(\roman*)]
\item $|\Omega|\leq \e' p^{\ell}/2$;
\item  the number of encodings of $T(d)$ in $(G,A)$ with leaves in $H$ is at most $\mu |H|^{2^d}|G|^{2^d-1}$;
\item for all $K\in (G/H)\setminus \Omega$, the number of encodings of $T(d)$ in $(G,A)$ with leaves in $H$ and nodes in $K$ is at most $\mu |H|^{2^d+2^d-1}$.
\end{enumerate}
Then there exists a  $(T,\ell, D, f,g)$-factor chain for $A$ with $\calL_0=\calL$.
\end{lemma}
\begin{proof}
Fix $d,T\geq 1$ and increasing functions $f,g:\mathbb{N}\rightarrow \mathbb{N}$.  We inductively define parameters $\rho_i,\e_i, \mu_i,\mu_i', N_i,N_i', D_i,D_i'$ for each $0\leq i\leq T$ as follows.

\underline{Step $i=0$}. Let $\e_0=\e_0(d)$ be from Theorem \ref{thm:stableinh}.  Set $\rho_0=\mu_0=\mu_0'=1$,  $N_0=N'_0=1$, and $D_0=D_0'=0$. 

\underline{Step $i+1$}.  Suppose that $0\leq i<T$ and assume that we have defined, for each $0\leq s\leq i$, parameters $\rho_s, \e_s, \mu_s,\mu_s',N_s,N_s',D_s,D_s'$, so that $D_0\leq D_1\leq \ldots \leq D_i$, and so that for each $0\leq s\leq i$, the choice of parameters $\e_s$, $\mu_s$, $\mu_s'$, $\rho_s$, $N_s'$, $D_s$, $D_s'$ depends only on $d$, $f$, $g$, $d$, and the choices of $D_j,\mu_j$ for $0\leq j<s$.

Set $\e_{i+1}=\e_0 p^{-g(iD_i)}$.  Choose $N'_{i+1}=N(d,\e_{i+1})$, $D'_{i+1}=L(d,\e_{i+1})$, and $\mu'_{i+1}=\mu(d,\e_{i+1})$ as in Theorem \ref{thm:stableinh}.  Then set $\mu_{i+1}=\min(\{\mu_{i+1}'\}\cup \{\mu_j: 1\leq j\leq i\})$, and define 
$$
\rho_{i+1}=p^{-f(iD_i)-iD_i}\text{ and }D_{i+1}=f(iD_i)+D_{i+1}'.
$$
Observe that since $f$ is increasing, $D_{i+1}>D_i$.  Finally, choose some $N_{i+1}$ sufficiently large compared to $N_{i+1}'+\rho_{i+1}^{-1}$.  Note that the choices of $\rho_{i+1},\e_{i+1}$, $\mu_{i+1}$, $\mu_{i+1}'$, $N_{i+1}$, $N_{i+1}'$, $D_{i+1}'$, $D_{i+1}$, depend only on $D_i$, $\mu_0,\ldots,\mu_i$,  $d$, $f$, and $g$.   This finishes the construction of the parameter sequences. 

Set $N=N_T$ and $D=\sum_{j=0}^TD_j$, and define $\mu=\min\{\mu_{i+1}\rho_{i+1}^{2^{d+1}}: i\in [T]\}$ and $\e'=\e_T$.    Note that the choice of these parameters depends only of $\e,d,T,f,g$.

Now suppose $\ell\geq 0$, $n-\ell\geq N$, $A\subseteq \F_p^n$, $\calL$ is a linear factor in $\F_p^{\ell}$ of complexity $\ell$, and $\Omega\subseteq \At(\calL)$ are such that the following hold with $H=L(0)$. 
\begin{enumerate}[label=\normalfont(\roman*)]
\item $|\Omega|\leq \e' p^{\ell}/2$;
\item the number of encodings of $T(d)$ in $(G,A)$ with leaves in $H$ is at most $\mu |H|^{2^d}|G|^{2^d-1}$;
\item for every $K\in (G/H)\setminus \Omega$, the number of encodings of $T(d)$ in $(G,A)$ with leaves in $H$ and nodes in $K$ is at most $\mu |H|^{2^d+2^d-1}$.  
\end{enumerate}
We construct a  $( T, \ell, D, f, g)$-factor chain for $A$ by induction. 

\noindent \underline{Step 0:} Let $\calL_0=\calL$ and $\ell_0=\ell$.

\noindent\underline{Step $i+1$:} Suppose that $0\leq i<T$, and assume by induction that we have constructed an $(  i, \ell_{0},D_i,f,g)$-factor chain $\calL_0\supseteq \ldots \supseteq \calL_i$, where each $|\calL_j|=\ell_j$.   By definition, this implies that for each $1\leq j\leq i$, $f(\ell_{j-1}-\ell_0)\leq \ell_j-\ell_{j-1}\leq D_i$.  

Since $n\geq N$, there exists a set $X\subseteq \F_p^n$ such that $\calL'_i=X\cup \calL_i$ is linearly independent, and has complexity $\ell_i'=  \ell_i+ f(iD_i)$.  Set $\tau_i=|L_i'(0)|/|H|$, and observe that 
\begin{align*}
\tau_i=p^{-(\ell'_i-\ell_{i-1})-(\ell_{i-1}-\ell_{i-2})-\ldots-(\ell_1-\ell_0)}&=p^{-f(iD_i)-(\ell_i-\ell_{i-1})-(\ell_{i-1}-\ell_{i-2})-\ldots-(\ell_1-\ell_0)}\geq \rho_{i+1},
\end{align*}
where the last inequality is by definition of $\rho_{i+1}$ and, since $\ell_j-\ell_{j-1}\leq D_i$ for each $1\leq j\leq i$.  Let $M_i$ be the number of encodings of $T(d)$ in $(G,A)$ with leaves in $L_i(0)$.  Since $n\geq N$, we have $n-\ell_i'\geq N'_{i+1}$, and since any encoding of $T(d)$ in $(G,A)$ with leaves in $L'_i(0)$ is also an encoding of $T(d)$ in $(G,A)$ with leaves in $H$, we have
$$
M_i\leq \mu |H|^{2^d}|G|^{2^d-1}=\mu \tau_i^{-2^d}|L_i'(0)|^{2^d}|G|^{2^d-1}\leq \mu \rho_{i+1}^{-2^d}|L_i'(0)|^{2^d}|G|^{2^d-1}\leq \mu_{i+1} |L_i'(0)|^{2^d}|G|^{2^d-1},
$$
where the second inequality holds since $\tau_i\geq \rho_{i+1}$, and the third inequality holds by definition of $\mu$.  Now let $\Omega'_i$ be the set of cosets $K$ of $L'_i(0)$ such that the number of encodings of $T(d)$ in $(G,A)$ with leaves in $L'_i(0)$ and nodes in $K$ is at least $\mu_{i+1} |L_i'(0)|^{2^d+2^d-1}$.  We claim that for all $H+k\in (G/H)\setminus \Omega$ and $L'_i(0)+h+k\subseteq L(0)+k$, $L'_i(0)+h+k\notin \Omega'_i$.  Indeed, suppose towards a contradiction that there exists some $H+k\in (G/H)\setminus \Omega$ and  $L'_i(0)+h+k\subseteq H+k$ with $L'_i(0)+h+k\in \Omega'_i$.  Then the number of encodings of $T(d)$ in $(G,A)$ with leaves in $H$ and nodes in $L_i'(0)+h+k$ is at least
\begin{align*}
\mu_{i+1} |L_i'(0)|^{2^d+2^d-1}=\mu_{i+1} \tau_i^{2^d+2^d-1}|H|^{2^d+2^d-1}\geq \mu_{i+1}\rho_{i+1}^{2^d+2^d-1}|H|^{2^d+2^d-1}>\mu |H|^{2^d+2^d-1},
\end{align*}
where the first inequality is because $\tau_i\geq \rho_{i+1}$, and the second inequality holds by definition of $\mu$. Since $L_i'(0)+h+k\subseteq H+k$, this implies that there are more than $\mu|H|^{2^d+2^d-1}$ encodings of $T(d)$ in $(G,A)$ with leaves in $H$ and nodes in $H+k$, contradicting that $H+k\notin \Omega$.  Therefore, we have that for all $H+k\notin \Omega$ and all $L'_i(0)+h+k\subseteq L(0)+k$, $L'_i(0)+h+k\notin \Omega'_i$.  Consequently, using (i) and the definition of $\e'$, we have
$$
|\Omega_i'|\leq |\Omega||L(0)|/|L_i'(0)|\leq \e' p^{\ell} |L(0)|/2|L_i'(0)|=\e' p^{\ell_i'}/2\leq \e_{i+1}p^{\ell_i'}/2.
$$
Thus, the hypotheses of Theorem \ref{thm:stableinh} are satisfied by $A$ and $\calL'_i$ for $d$ and $\e_{i+1}$.

By Theorem \ref{thm:stableinh}, there is $\calL_{i+1}\supseteq \calL'_i$ of complexity $\ell_{i+1}\leq\ell'_i+D_{i+1}'$, such that $\calL_{i+1}$ is almost $\e_{i+1}$-atomic with respect to $A$. Note $\ell_{i+1}\leq \ell_i+f(iD_i)+D_{i+1}'\leq \ell_i+D_{i+1}$.  Clearly $\ell_{i+1}-\ell_i\geq f(iD_i)$ since $\calL_i'\subseteq \calL_{i+1}$, and thus, since $\ell_j-\ell_{j-1}\leq D_i$ holds for each $1\leq j\leq i$, we have that $\ell_{i+1}-\ell_i\geq f(\ell_i-\ell_0)$.  Combining these yields that $f(\ell_i-\ell_0)\leq \ell_{i+1}-\ell_i\leq D_{i+1}$.  Further, since $g$ is increasing, $\e_0\leq 1$, and again because $\ell_j-\ell_{j-1}\leq D_i$ holds for each $1\leq j\leq i$, we have that
$$
\e_{i+1}=\e_0 p^{-g(iD_i)}\leq p^{-g(\ell_i-\ell_{0})}.
$$
Thus  $\calL_0\supseteq \dots \supseteq \calL_{i+1}$ forms an $(i+1,\ell_0,D_{i+1},f,g)$-factor chain.

Upon completion of step $T$, $\calL_0\supseteq \ldots \supseteq \calL_T$ will be the desired $(T,\ell, D,f,g)$-factor chain for $A$.
\end{proof}

Our next goal is to use Lemma \ref{lem:chainlinearnew} to prove Theorem \ref{thm:almoststableinhstrong}. For this we require two facts.  The first is a well known relationship between the density of a set $A$ on a sum of cosets and the density of its sum-graph on said pair of cosets.  Its proof is a simple exercise.

\begin{fact}\label{fact:sumdensities}
Suppose $A\subseteq \F_p^n$ and $K, L$ are both cosets of subgroups of $\F_p^n$ (not necessarily of the same subgroup).  Then 
$$
\frac{|\{(x,y)\in K\times L: x+y\in A\}|}{|K||L|}=\frac{|A\cap (K+L)|}{|K+L|}.
$$
\end{fact}
We note that for  Fact \ref{fact:sumdensities} to hold,  it is crucial that the sets $K,L$ are both cosets of subspaces. The second fact we will use is the following   simple counting lemma for linear atoms on which a set has density near $0$ or $1$.  We note this is not a direct corollary of Proposition \ref{prop:countinggplinear} as there are several different linear factors involved.  

\begin{lemma}\label{lem:graphcounting1}
For all $t\in \mathbb{N}$ and $\e\in (0,1)$ there exists $\delta>0$ so that the following holds.  Let $U,W$ be sets satisfying $|U|+|W|=t$, and assume $E_0,E_1$ disjoint subsets of $U\times W$.  Let $A\subseteq \mathbb{F}_p^n$, and for each $x\in U\cup W$, suppose $C_x$ is a coset of some subspace of $\F_p^{n}$ (possibly the $C_x$ are cosets of different subspaces for different choices of $x$). Assume that for all $(u,w)\in E_1$, 
 $$
|A\cap (C_u+C_w)|\geq (1-\delta)|C_u+C_w|,
 $$
 and for $(u,w)\in E_0$, 
 $$
 |A\cap (C_u+C_w)|\leq \delta|C_u+C_w|.
 $$
 Then there are at least $(1-\e)\prod_{x\in U\cup W}|C_x|$ tuples $(v_u)_{u\in U}(v_w)_{w\in W}\in \prod_{x\in U\cup W}C_x$ such that $v_u+v_w\in A$ if $uw\in E_1$ and $v_u+v_w\notin A$ if $uw\in E_0$.   
\end{lemma}
\begin{proof}
Fix $t\geq 2$ and $\e\in (0,1)$.  Without loss of generality, assume that $0<\e<1/4$. Set $\delta=\e/t^2$.  Let $U,W$ be sets satisfying $|U|+|W|=t$, and let $E_0,E_1$ be disjoint subsets of $U\times W$. Suppose $A\subseteq \F_p^n$, and for each $x\in U\cup W$, $C_x$ is a subspace of $\F_p^n$.  Assume $(u,w)\in E_1$ implies $|A\cap (C_u+C_w)|\geq (1-\delta)|C_u+C_w|$ and $(u,w)\in E_0$ implies $|A\cap (C_u+C_w)|\leq \delta|C_u+C_w|$. 

 Let  $I\subseteq \prod_{u\in U}C_u\times \prod_{w\in W}C_w$ be the set of tuples $(v_{u})_{u\in U}(v_{w})_{w\in W}$, satisfying $v_{u}+v_{w}\in A$ for all $(u,w)\in E_1$ and  $v_{u}+v_{w}\notin A$ for all $(u,w)\in E_0$.  Using our assumptions and  Fact \ref{fact:sumdensities}, we can then compute 
\begin{align*}
|(\prod_{u\in U}C_u\times \prod_{w\in W}C_w)\setminus I|&\leq \sum_{uw\in E_1}|\{(v_u,v_w)\in C_u\times C_w: v_u+v_w\notin A\}|  \prod_{u'\in U\setminus \{u\}, w'\in W\setminus \{w\}}|C_{u'}||C_{w'}|\\
&+\sum_{uw\in E_0}|\{(v_u,v_w)\in C_u\times C_w: v_u+v_w\in A\}|  \prod_{u'\in U\setminus \{u\}, w'\in W\setminus \{w\}}|C_{u'}||C_{w'}|\\
&\leq 2 \delta|U||W| \prod_{u\in U}|C_u|\times \prod_{w\in W}|C_w|\\
&\leq \e \prod_{u\in U}|C_u|\times \prod_{w\in W}|C_w|,
\end{align*}
where the last inequality is by choice of $\delta$.
\end{proof}

We now prove Theorem \ref{thm:almoststableinhstrong}. 
\vspace{2mm}

\begin{proofof}{Theorem \ref{thm:almoststableinhstrong}} Fix an integer $d\geq 1$ and an increasing function $\psi:\mathbb{N}\rightarrow \mathbb{R}^+$. Choose $\delta=\delta(2^d+2^d-1,1/p)$  as in Lemma \ref{lem:graphcounting1}, and define  $f,g:\mathbb{N}\rightarrow \mathbb{N}$ by setting $g(x)=\lceil 9x\psi(dx)+\log_p(\delta^{-1})\rceil$ and $f(x)=10g(x)+2x+2$ for all $x\in \mathbb{N}$. Let $\e_1=\e_1(d,d,f,g)$, $D_1=D_1(d,d,f,g)$, $\mu_1=\mu_1(d, d,f,g)$, and $N_1=N_1(d, d,f,g)$ be as in Lemma \ref{lem:chainlinearnew}. Set $D=D_1d$, $C=(D_1(d+2)+2)2^d+2$, $\mu_0=\min\{\e_1/2^{2^d}, 1/p,\mu_1 \}$, and $N=(D_1d)^dN_0N_1$.

Now suppose $0<\mu\leq \mu_0$, $\ell\geq 0$, $n-\ell \geq N$, $A\subseteq G=\F_p^n$,  $\calL$ is a linear factor on $\F_p^n$ of complexity $\ell$, and $\Omega\subseteq G/L(0)$ satisfy the following with $H=L(0)$.
\begin{enumerate}
\item   $|\Omega|\leq \mu p^{\ell}$ and $|H|\geq p^N$;
\item the number of encodings of $T(d)$ in $(G,A)$ with leaves in $H$ is at most $\mu^C |H|^{2^d}|G|^{2^d-1}$;
\item for all $K\in (G/H)\setminus \Omega$, the number of encodings of $T(d)$ in $(G,A)$ with leaves in $H$ and nodes in $K$ is at most $\mu |H|^{2^d+2^d-1}$. 
\end{enumerate}

Since $\mu<\mu_0<\min\{\e_1,\mu_1\}$, Lemma \ref{lem:chainlinearnew} implies there exists a  $( d, \ell, D_1, f, g)$-linear factor chain for $A$ where $\calL_0=\calL$.  Say the chain consists of, for each $0\leq i\leq d$, a linear factor $\calL_i$ of complexity $\ell_i$, and a partition $\At(\calL_i)=\Gamma_i^1\cup \Gamma_i^0\cup \Gamma_i^{err}$.  Note that by the definition of a factor chain, $\ell_d-\ell_0\leq dD_1$.  We make note of a few inequalities here which will we will use later on in the proof.  By definition of $f$, we have that for all $x$,
\begin{align}\label{in1}
f(x)>\max\{2\psi(x),x\}.
\end{align}
Note that for all $1\leq i\leq d$, 
\begin{align}\label{in0}
\ell_i-\ell_{i-1}\geq f(\ell_{i-1}-\ell_{0})>1,
\end{align}
 where the last inequality is by definition of $f$. For all $0\leq i\leq d$, we have that 
\begin{align}\label{in2}
\ell_i-\ell_0=\sum_{j=1}^i\ell_j-\ell_{j-1}\geq \sum_{j=1}^{i-1} f(\ell_j-\ell_0)\geq \sum_{j=1}^{i-1} 2(\ell_j-\ell_0)+2>2(i-1)+2> i,
\end{align}
where the first inequality holds by definition of a factor chain, the second  by definition of $f$, and the third by (\ref{in0}).  Finally, by definition of $g$, we have that for all $0\leq i\leq d$,
\begin{align}\label{in3}
p^{-g(\ell_i-\ell_0)/2}<p^{-4\psi(\ell_i-\ell_0)}.
\end{align}

We will show that for some $0\leq i\leq d-1$, $\calL_i$ is $p^{-\psi(\ell_i-\ell_0)}$-atomic with respect to $A$.  For each $0\leq i\leq d-1$, define 
\begin{align*}
\Sigma_i^{err}=\{L_i(\abar)\in \At(\calL_i): & \text{ for each }u\in \{0,1\},\\
& |\{L\in \At(\calL_{i+1}): L\subseteq L_i(\abar)\}\cap \Gamma_{i+1}^u\}|
\geq  p^{\ell_{i+1}-\ell_i-2\psi(\ell_i-\ell_0)} \}.
\end{align*}
We claim that there is $0\leq i\leq d-1$ such that $|\Sigma^{err}_i|\leq \mu p^{\ell_i}$. Suppose towards a contradiction that this is false, i.e. that for every $0\leq i\leq d-1$, $|\Sigma_i^{err}|\geq \mu p^{\ell_i}$.  Let $X$ be the set of encodings of $T(d)$ in $(G,A)$ with leaves in $H$.  We construct many distinct elements of $X$ as follows.
\begin{enumerate}
\item[(a)] For each $0\leq i\leq d-1$, choose a set $I_i=\{L_i(\abar_{\sigma}): \sigma\in 2^i\}\in {\Sigma_i^{err}\choose 2^i}$.  By assumption there are at least $\prod_{i=0}^{d-1}{\mu p^{\ell_i}\choose 2^i}\geq \prod_{i=0}^{d-1}(\mu/2^i)^{2^i} p^{2^i \ell_i}$ ways to do this.
\item[(b)] For each $0\leq i\leq d-1$ and $\sigma\in 2^i$, choose any  $L_{i+1}(\abar_{\sigma}\cbar_{\sigma \wedge 1})\in \Gamma_{i+1}^1$, and  $L_{i+1}(\abar_{\sigma}\cbar_{\sigma\wedge 0})\in \Gamma_{i+1}^0$.   For each $0\leq i\leq d-1$ and $\sigma\in 2^i$, since $L_i(\abar_{\sigma})\in \Sigma_i^{err}$, there are at least $(p^{\ell_{i+1}-\ell_i-2\psi(\ell_i-\ell_0)})^2$ ways to do this.  Therefore the number of choices for this step is at least $\prod_{i=0}^{d-1}(p^{\ell_{i+1}-\ell_i-2\psi(\ell_i-\ell_0)})^{2^i}$.
\item[(c)] For each $0\leq i\leq d-1$ and $\sigma\in 2^i$, let $h^{\sigma}_i$ be such that $L_i(\abar_{\sigma})=h^{\sigma}_i+L_i(0)$, and for each $u\in \{0,1\}$, let $g_i^{\sigma \wedge u}\in L_i(0)$ be such that $L_{i+1}(\abar_{\sigma}\cbar_{\sigma\wedge u})=h^{\sigma}_i+g_i^{\sigma \wedge u}+L_{i+1}(0)$.  Now for each $0\leq i\leq d-1$ and $\sigma\in 2^i$, and each $\eta\in 2^d$, set 
\[V_{\sigma}=L_{i+1}(0)+h^{\sigma}_{i}-\sum_{1\leq j\leq i-1}g_j^{\sigma(1)\wedge \ldots \wedge \sigma(j+1)}\]
and
\[W_{\eta}=L_d(0)+\sum_{j=1}^{d-1} g_j^{\eta(1)\wedge \cdots \wedge \eta(j+1)}.\]
Let $J_i$ be the set of tuples $(v_{\sigma})_{\sigma\in 2^{<d}}(w_{\eta})_{\eta \in 2^d}\in \prod_{\sigma\in 2^{<d}}V_{\sigma}\times \prod_{\eta\in 2^d}W_{\eta}$ where $v_{\sigma}+w_{\eta}\notin A$ holds for all  $\sigma \wedge 0\trianglelefteq \eta$, and $v_{\sigma}+w_{\eta}\in A$ holds for all  $\sigma \wedge 1\trianglelefteq \eta$.   Note that any element of $J_i$  is an encoding of $T(d)$ in $(G,A)$  with the property that each $v_{\sigma}\in V_{\sigma}$ and each $w_{\eta}\in W_{\eta}$.  We give a lower bound for the size of $J_i$.   Note that each $W_{\eta}\subseteq H$, so each such encoding will have leaves in $H$.

Observe that for all $\sigma \in 2^i$ and $\eta\in 2^d$, if $\sigma \wedge 1\triangleleft \eta$, then  we have
$$
\frac{|A\cap (V_{\sigma}+W_{\eta})|}{| V_{\sigma}+W_{\eta}|}=\frac{|A\cap (L_{i+1}(0)+h^{\sigma}_{i}+g_i^{\sigma \wedge 1})|}{|(L_{i+1}(0)+h^{\sigma}_{i}+g_i^{\sigma \wedge 1})|}\geq 1-p^{-g(\ell_i-\ell_0)}
$$
and if $\sigma \wedge 0\triangleleft \eta$ then we have
$$
\frac{|A\cap (V_{\sigma}+W_{\eta})|}{| V_{\sigma}+W_{\eta}|}=\frac{|A\cap (L_{i+1}(0)+h^{\sigma}_{i}+g_i^{\sigma \wedge 0})|}{|(L_{i+1}(0)+h^{\sigma}_{i}+g_i^{\sigma\wedge 0})|}\leq p^{-g(\ell_i-\ell_0)}.
$$
Since $ p^{-g(\ell_i-\ell_0)}<\delta$ (by definition of $g$), we have by Lemma \ref{lem:graphcounting1} that
$$
|J_i|\geq (1-p^{-1})\prod_{\sigma\in 2^{<d}}|V_{\sigma}|\prod_{\eta\in 2^d}|W_{\eta}|=(1-p^{-1})|L_d(0)|^{2^d} \prod_{i=0}^{d-1}|L_{i+1}(0)|^{2^{i}}.
$$
\end{enumerate}

Combining these, (a)-(c), we have that 
\begin{align*}
|X|&\geq \Big(\prod_{i=1}^{d-1}(\mu/2^i)^{2^i} p^{2^i \ell_i}\Big)\Big(\prod_{i=0}^{d-1}( p^{\ell_{i+1}-\ell_i-2\psi(\ell_i-\ell_0)})^{2^i}\Big)\Big((1-p^{-1})|L_d(0)|^{2^d} \prod_{i=0}^{d-1}|L_{i+1}(0)|^{2^{i}}\Big)\\
&\geq (1-p^{-1})|L_d(0)|^{2^d}\prod_{i=0}^{d-1}(\mu 2^{-i}p^{\ell_i} p^{\ell_{i+1}-\ell_i-2\psi(\ell_i-\ell_0)}|L_{i+1}(0)|)^{2^{i}}\\
&\geq (1-p^{-1})|L_d(0)|^{2^d}\prod_{i=0}^{d-1}(\mu 2^{-i} p^{\ell_i}|L_{i+1}(0)|)^{2^{i}},
\end{align*}
where the last inequality holds because, for each $0\leq i\leq d-1$, $\ell_{i+1}-\ell_i\geq f(\ell_i-\ell_0)>2\psi(\ell_i-\ell_0)$ (here the first inequality is by (\ref{in0}) and the second by (\ref{in1})).  The last line above is equal to 
\begin{align*}
(1-p^{-1})|L_d(0)|^{2^d}\prod_{i=0}^{d-1}( \mu 2^{-i} p^{\ell_i}p^{-\ell_{i+1}}|G|)^{2^{i}}&\geq (1-p^{-1})|L_d(0)|^{2^d}\prod_{i=0}^{d-1}( \mu  p^{\ell_i-\ell_{i+1}-i}|G|)^{2^{i}}\\
&\geq (1-p^{-1})|L_d(0)|^{2^d}\prod_{i=0}^{d-1}( \mu  p^{2(\ell_i-\ell_{i+1})}|G|)^{2^{i}}\\
&\geq (1-p^{-1})|L_d(0)|^{2^d}\prod_{i=0}^{d-1}( \mu p^{-2D_1}|G|)^{2^{i}}\\
&= (1-p^{-1})|L_d(0)|^{2^d}(\mu p^{-2D_1})^{2^d-1}|G|^{2^d-1},
\end{align*}
where the second inequality holds because, by (\ref{in0}), (\ref{in1}) and $(\ref{in2})$, we have that for each $0\leq i\leq d-1$, $\ell_{i+1}-\ell_i \geq f(\ell_i-\ell_0)>\ell_i-\ell_0\geq i$, and the third inequality  is because, by the definition of a factor chain, $\ell_{i+1}-\ell_i\leq D_1$ holds for each $0\leq i\leq d-1$.  Combining the above with the fact that $|L_d(0)|=p^{\ell_0-\ell_d}|L_0(0)|=p^{\ell_0-\ell_d}|H|$ yields that $|X|$ is at least
\begin{align*}
(1-p^{-1})(p^{\ell_0-\ell_d})^{2^d}|H|^{2^d}(\mu  p^{-2D_1})^{2^d-1} |G|^{2^d-1}&\geq (1-p^{-1})(p^{-dD_1})^{2^d}(\mu   p^{-2D_1})^{2^d-1}|H|^{2^d} |G|^{2^d-1}\\
&\geq (1-p^{-1})(\mu  p^{-(d+2)D_1})^{2^d}|H|^{2^d}|G|^{2^d-1}\\
&\geq (\mu)(\mu \mu  \mu ^{(d+2)D_1})^{2^d}|H|^{2^d}|G|^{2^d-1}\\
&>\mu^C|H|^{2^d}|G|^{2^d-1},
\end{align*}
where the first inequality uses the fact that $\ell_d-\ell_0\leq dD_1$, the third inequality is because, by definition, $\mu<\min\{1-p^{-1}, p^{-1}\}$, and the last inequality holds by definition of $C$.  But this contradicts our assumption (ii), which says that $|X|\leq \mu^C|H|^{2^d}|G|^{2^d-1}$.

Thus there is some $0\leq i\leq d-1$ such that $|\Sigma^{err}_i|\leq \mu p^{\ell_i}$. Let 
$$
\Omega'=\Sigma_i^{err}\cup \{L\in \At(\calL_i): |\{L'\subseteq L: L'\in \Gamma^{err}_i\}|\geq   p^{\ell_{i+1}-\ell_i-g(\ell_i-\ell_0)/2}\}.
$$
Then, using the definition of a factor chain and the bound on the size of $\Sigma_i^{err}$, we have that 
$$
|\Omega'|\leq |\Sigma^{err}_i|+  p^{\ell_i-g(\ell_i-\ell_0)/2}\leq \mu p^{\ell_i}+  p^{\ell_i-g(\ell_i-\ell_0)/2}=p^{\ell_i}(\mu+p^{-g(\ell_i-\ell_0)/2})\leq p^{\ell_i}(\mu+p^{-\psi(\ell_i-\ell_0)}),
$$
where the last inequality holds by (\ref{in3}).   Then for all $L\in \At(\calL_i)\setminus \Omega'$, it is straightforward to see there is $u\in \{0,1\}$ such that 
    $$
    |A^u\cap L|/|L|\geq 1-  p^{-g(\ell_i-\ell_0)/2}-  p^{-2\psi(\ell_i-\ell_0)}\geq 1- 2  p^{-2\psi(\ell_i-\ell_0)}\geq 1- p^{-\psi(\ell_i-\ell_0)},
    $$
where the second inequality uses (\ref{in3}), and the final inequality holds because $\psi$ is a growth function, so $p^{-\psi(\ell_i-\ell_0)}\leq 1/2$.  This finishes the proof, taking $H'=L_i(0)$, $m=\ell_i-\ell_0$, and the $\Omega'$ above.
\end{proofof}

\appendix

\section{Hypergraph formulations of the main theorems}\label{ss:hypergraphs}  
In this section we demonstrate how Theorems \ref{thm:vc2finite} and \ref{thm:fop} translate into regular decompositions for ternary sum-graphs of subsets of $\F_p^n$, allowing us to relate these results to those of \cite{Terry.2021b}, which deal with hypergraph regularity.   The types of graphs and hypergraphs of interest in this section are as follows.

\begin{definition}
Suppose $A\subseteq \F_p^n$. The \emph{binary sum-graph associated to $A$} is the pair $(\F_p^n, E_A^{(2)})$, where 
$$
E_A^{(2)}=\{xy\in {\F_p^n\choose 2}: x+y\in A\}.
$$
 Similarly, the \emph{ternary sum-graph associated to $A$} is the pair $(\F_p^n, E_A^{(3)})$, where
 $$
 E_A^{(3)}=\{xyz\in {\F_p^n\choose 3}: x+y+z\in A\}.
 $$
\end{definition}

It was shown by Green that the binary sum-graph associated to a set $A\subseteq \F_p^n$ admits a special version of Szemer\'{e}di's regularity lemma, where all parts in the regular partition are cosets of a single subgroup (see \cite[Section 9]{Green.2005}). In other words, the regular partition is exactly $\At(\calL)$ for some linear factor $\calL$.  One goal of this section is to prove a hypergraph analogue of this. In particular, we will show that the ternary sum-graph associated to a subset $A\subseteq \F_p^n$ admits a special version of the hypergraph regularity lemma, where the components of the decomposition are structured according to a high rank quadratic factor (Theorem \ref{thm:hypergraphreg}).  As a corollary, we show that Theorems \ref{thm:vc2finite} and \ref{thm:fop} imply algebraically structured versions of the corresponding $\NIP_2$ and $\NFOP_2$ hypergraph regularity lemmas from  \cite{Terry.2021b}, in the case where the hypergraph arises as the ternary sum-graph of subset of $\F_p^n$ (see Corollaries \ref{cor:hypergraphvc2} and \ref{cor:hypergraphfop}).  

\subsection{Regularity for  hypergraphs}
We provide in this subsection an abridged introduction to hypergraph regularity.   To simplify the presentation, we will use a slightly different version of hypergraph regularity than that used in \cite{Terry.2021b}.  This is no loss of generality, as the two versions are known to be equivalent. We refer the reader to \cite{Nagle.2013} and \cite{Nagle.2024} for a thorough history of the subject, and the equivalence of the version used here with that from \cite{Terry.2021b}.

We begin by defining quasirandomenss for a  bipartite graph. This particular notation is taken from \cite{Nagle.2013}, but has older roots (see e.g. \cite{Chung.1990}). 

\begin{definition}[$\dev_2$]\label{def:dev2}
Let $\e\in (0,1)$ and let $G=(U\cup V,E)$ be a bipartite graph. Then \emph{$G$ has $\dev_2(\e)$} if the following holds for $d_2=|E|/|U||V|$.
$$
\sum_{u_0,u_1\in U}\sum_{v_0,v_1\in V} \prod_{i\in \{0,1\}}\prod_{j\in \{0,1\}}g(u_i,v_j)\leq \e |U|^2|V|^2,
$$
where $g(u,v)=1-d_2$ if $uv\in E$ and $g(u,v)=-d_2$ if $uv\notin E$.

We say $G$ has $\dev_2(\e; d)$ if it has $\dev_2(\e)$ and $|E|/|U||V|\in (d-\e, d+\e)$.
\end{definition}

We next work towards defining when a 3-uniform hypergraph (or 3-graph, for short), is quasirandom relative to an underlying graph.  We will use the following notation. Given sets $V_1,V_2\subseteq V$, set $K_2[V_1,V_2]=\{xy: x\in V_1, y\in V_2\}$.  If  $G=(V,E)$ is a graph, we define 
$$
E[V_1,V_2]=E\cap K_2[V_1,V_2]\text{ and }G[V_1,V_2]:=(V_1\cup V_2,E[V_1,V_2]).  
$$
We will also use the following notation for the set of triangles in a graph.

\begin{definition}
Given a graph $G=(V,E)$ be a graph, define
\[K_3^{(2)}(G)=\left\{xyz\in {V\choose 3}: xy, yz, xz \in E\right\}.\]
\end{definition}

We need the following terminology for certain graph/hypergraph pairs.

\begin{definition}
Let $G=(V,E)$ be a graph and let  $H=(V,E')$ be a $3$-graph on the same vertex set as $G$.  We say $G$ \emph{underlies} $H$ if $E'\subseteq K_3^{(2)}(G)$.
\end{definition}

In other words,  $G$ underlies $H$ if every edge in $H$ sits atop a triangle from $G$. We next define what it means for a hypergraph to be quasi-random, relative to an underlying graph.

\begin{definition}[$\dev_{2,3}(\e_1,\e_2)$]\label{def:regtrip}
Let $\e_1,\e_2>0$, and let $H=(V,R)$ be a $3$-graph. Suppose $G=(V,E)$ is a $3$-partite graph underlying $H$ with vertex partition $V=V_1\cup V_2\cup V_3$.  Let $d_3$ be such that $|R|=d_3|K_3^{(2)}(G)|$.  We say that \emph{$(H,G)$ has $\dev_{2,3}(\e_1,\e_2)$} if there is $d_2\in (0,1]$ such that for each $1\leq i<j\leq 3$, the graph $G[V_i,V_j]$ has $\dev_2(\e_2;d_2)$, and 
$$
\sum_{u_0,u_1\in V_1}\sum_{v_0,v_1\in V_2}\sum_{w_0,w_1\in V_3}\prod_{i,j,k\in \{0,1\}}h_{H,G}(u_i,v_j,w_k)\leq \e_1d_2^{12}|V_1|^2|V_2|^2|V_3|^2,
$$
where 
\[h_{H,G}(u,v,w)=\begin{cases}
1-d_3 &\text{ if }uvw\in R\\
-d_3 &\text{ if }uvw\in K_3^{(2)}(G)\setminus R\\
0&\text{ if }uvw\notin K_3^{(2)}(G).
\end{cases}\]
\end{definition}

Our next goals are to define regular decompositions and to state the regularity lemma for $\dev_{2,3}$-quasirandomness.  For this, we need the following definition of a decomposition.

\begin{definition}[$(t,s)$-decomposition]\label{def:partn}
Let $V$ be a set. A \emph{$(t,s)$-decomposition} $\calP$ of $V$ is a partition $V=V_1\cup \ldots \cup V_t$, along with, for each $ij\in {[t]\choose 2}$, a partition $K_2[V_i, V_j]=\bigcup_{\alpha\leq s}P_{ij}^{\alpha}$.

A \emph{triad of $\calP$} is a tripartite graph of the form $G_{ijk}^{\alpha\beta\gamma}:=(V_i\cup V_j\cup V_k, P_{ij}^{\alpha}\cup P_{jk}^{\beta}\cup P_{ik}^{\gamma})$ where $ijk\in {[t]\choose 3}$.
\end{definition}

Given a $(t,s)$-decomposition $\calP$ as in Definition \ref{def:partn}, we let $\calP_{vert}=\{V_i: i\in [t]\}$ and $\calP_{pairs}=\{P_{ij}^{\alpha}: ij\in {[t]\choose 2}, \alpha\in [s]\}$, and write $\calP=(\calP_{vert},\calP_{pairs})$.  Continuing with this notation, our next definition describes decompositions where most of the bipartite graphs arising from $\calP_{pairs}$ are quasirandom.

\begin{definition}[$(t,s,\e_1,\e_2)$-decomposition]\label{def:partn2}
 A \emph{$(t,s,\e_1,\e_2)$-decomposition} of a set $V$ is a $(t,s)$-decomposition $\calP$ of $V$ such that
$$
|\bigcup_{G_{ijk}^{\alpha\beta\gamma}\in \Sigma}K_3^{(2)}(G_{ijk}^{\alpha\beta\gamma})|\geq (1-\e_1){|V|\choose 3},
$$ 
where $\Sigma$ is the set of triads $G_{ijk}^{\alpha\beta\gamma}$ of $\calP$ with the property that each of $(V_i\cup V_j, P_{ij}^{\alpha})$, $(V_i\cup V_k; P_{ik}^{\beta})$, and $(V_j\cup V_k; P_{jk}^{\gamma})$ have $\dev_2(\e_2)$.
\end{definition}

We now define a regular triad.

\begin{definition}[Regular triad]\label{def:regtriad}
Let $H=(V,E)$ be a $3$-graph and let $\calP$ be a $(t,s,\e_1,\e_2)$-decomposition of $V$.  A triad of $\calP$ is called \emph{$\dev_{2,3}(\e_1,\e_2)$-regular} with respect to $H$ if $(H_{ijk}^{\alpha\beta\gamma}, G_{ijk}^{\alpha\beta\gamma})$ has $\dev_{2,3}(\e_1,\e_2)$, where $H_{ijk}^{\alpha\beta\gamma}=(V_i\cup V_j\cup V_k, E\cap K_3^{(2)}(G_{ijk}^{\alpha\beta\gamma}))$.
\end{definition}

We now define a regular decomposition of a $3$-graph.

\begin{definition}[Regular decomposition]\label{def:regptn}
Given a $3$-graph $H=(V,E)$, a  $(t,s,\e_1,\e_2)$-decomposition $\calP$ of $V$ is called  \emph{$\dev_{2,3}(\e_1,\e_2)$-regular with respect to $H$} if
$$
|\bigcup_{G_{ijk}^{\alpha,\beta,\gamma}\in \Sigma}K_3^{(2)}(G_{ijk}^{\alpha\beta\gamma})|\geq (1-\e_1){|V|\choose 3},
$$
where $\Sigma$ is the set of triads of $\calP$ which are $\dev_{2,3}(\e_1,\e_2 )$-regular with respect to $H$. 
\end{definition}

We now state a regularity lemma for $3$-graphs. This particular version is essentially due to Gowers \cite{Gowers.20063gk} (see also  \cite[Theorem 1.1]{Nagle.2013}).

\begin{theorem}\label{thm:hgreg}
For all $\e_1>0$, every function $\e_2:\mathbb{R}^+\rightarrow (0,1]$, and every $s_0, t_0\in \mathbb{N}$, there exist positive integers $T_0,S_0$ such that for any sufficiently large $3$-graph $H=(V,E)$, there exists a $(t,s, \e_1,\e_2(s))$-decomposition $\calP$ of $V$ which is $\dev_{2,3}(\e_1,\e_2(s))$-regular with respect to $H$, and where $t_0\leq t \leq T_0$ and $s_0\leq s\leq S_0$.
\end{theorem}

We will soon state the first goal of this section, namely a special version of Theorem \ref{thm:hgreg} in the case where $H$ is the ternary sum-graph of a subset of $\F_p^n$.   First, we recall the following notation from earlier in the paper.  
 
 \begin{definition}\label{def:betaref}
 Suppose $\calQ=\{M_1,\ldots, M_q\}$ is a purely quadratic factor on $\F_p^n$. Given $b=(b_1,\ldots, b_q)\in \mathbb{F}_p^q$, let
$$
\beta_\calQ(b)=\{(x,y)\in \mathbb{F}_p^{2n}: \text{for each }i\in [q], x^TMy=b_i\}.
$$
\end{definition}

We will throughout this section abuse notation slightly by also writing $\beta_{\calQ}(b)$ to denote the set of unordered pairs $xy\in {\F_p^n\choose 2}$ with the property that $(x,y)\in \beta_{\calQ}(b)$.  With this in mind, given any quadratic factor $\calB=(\calL,\calQ)$ on $\F_p^n$ of complexity $(\ell,q)$, we naturally obtain a decomposition of the set $\F_p^n$ in the sense of Definition \ref{def:partn} by taking $\calP_{vert}$ to be $\At(\calB)$ and letting $\calP_{pairs}$ be 
$$
\{ K_2[B_1,B_2]\cap \beta_{\calQ}(b):B_1,B_2\in \At(\calB), b\in \F_p^q\}.  
$$
The first main result of the section will say that any ternary sum-graph associated to a subset of $\F_p^n$ has a regular decomposition in the sense of Theorem \ref{thm:hgreg}, which moreover arises from a quadratic factor in the manner just described.

\begin{theorem}[Hypergraph regularity with quadratic decompositions]\label{thm:hypergraphreg}
For all reals $\e_1\in (0,1)$, all functions $\e_2:\mathbb{R}^+\rightarrow (0,1]$, all integers $t_0,s_0\geq 1$, and all growth functions $\tau$, there are $S_0$ and $T_0$ such that for all sufficiently large $n$, the following holds.  For all $A\subseteq \F_p^n$, there exist  $t_0\leq t\leq T_0$, $s_0\leq s\leq S_0$, and a $(t,s,\e_1,\e_2(s))$-decomposition $\calP$ of $\F_p^n$ which is $\dev_{2,3}(\e_1,\e_2(s))$-regular with respect to $(\F_p^n, E_A^{(3)})$.

Moreover, there is a quadratic factor $\calB=(\calL,\calQ)$ on $\F_p^n$ of complexity $(\ell,q)$ and rank at least $\tau(\ell+q)$, such that 
$$
\calP_{vert}=\At(\calB)\text{ and }\calP_{edge}= \{K_2[B_1,B_2]\cap \beta_{\calQ}(b): b\in \F_p^{q}, B_1,B_2\in \At(\calB)\}.
$$
 \end{theorem}
 
 The main difference between this theorem and what one would obtain by simply applying Theorem \ref{thm:hgreg} to $(\F_p^n, E_A^{(3)})$ is of course the ``moreover" statement, which says that the regular decomposition exhibits extremely strong algebraic structure. The next three subsections will focus on proving this result.

 \subsection{Graphs arising from quadratic factors}
In this subsection, we focus on proving that the bipartite graphs appearing in the statement of Theorem \ref{thm:hypergraphreg} are sufficiently quasirandom.  This is needed to show that the decompositions appearing there will satisfy Definition \ref{def:partn2}.  It will  be useful to have the following compact notation for bipartite graphs arising from restricting  $\beta_{\calQ}(b)$ to pairs of atoms from the factor $\calB$.
 
 \begin{definition}\label{def:triads}
Suppose $\calB=(\calL,\calQ)$ is a quadratic factor on $\mathbb{F}_p^n$ of complexity $(\ell,q)$.   Given $e=(a_1,a_2,b_{12})\in \mathbb{F}_p^{\ell+q}\times \F_p^{\ell+q}\times \F_p^{	q}$, define \emph{bipartite graph corresponding to $e$} to be
$$
\Delta_{\calB}^{(2)}(e):=(B(a_1)\cup B(a_2), \beta_{\calQ}(b_{12})\cap K_2[B(a_1), B(a_2)]).
$$
 Note that if $(x,y)\in \Delta_{\calB}^{(2)}(e)$, then $x+y\in B(\Sigma_{\calB}(e))$, for 
\begin{align*}
\Sigma_{\calB}(e):=a_1+a_2+2(0b_{12}),
\end{align*}
where $0b_{12}$ is obtained by adding a string of $\ell$ zeros to the left of $b_{12}$.   
\end{definition}

In practice, we will suppress the subscripts in Definition \ref{def:triads} when the factor is clear from context. Our next goal is to show that  if a quadratic factor has high rank, then the bipartite graphs appearing Definition \ref{def:triads} are quasirandom.  This will require the following size estimates, which follow from \cite{Terry.2021d}.\footnote{See also Lemma B.2 in \cite{Terry.2021av2}.}  

\begin{lemma}\label{lem:sizeofbetaref}
Let $\calB=(\calL,\calQ)$ a quadratic factor on $\F_p^n$ of complexity $(\ell,q)$ and rank $\tau$. Then for any $b\in \F_p^q$, the following hold.
\begin{enumerate}
\item $|\beta_{\calQ}(b)|=p^{2n-q}(1+O(p^{q-\tau}))$.
\item For any $B_1,B_2\in \At(\calB)$, 
\begin{align*}
|\beta_{\calQ}(b)\cap K_2[B_1,B_2]|&=p^{-q}|B_1||B_2|(1+O(p^{O(\ell+q)-\tau/2})).
\end{align*} 
\end{enumerate}
\end{lemma}
\begin{proof}
Part (i) is Corollary 2.24 of \cite{Terry.2021d}.  For (ii), apply Lemma A.2 in \cite{Terry.2021d} with $m=1$, $I=J=K=\{1\}$, with $a_1,a_2$ chosen so that $B_1=B(a_1)$ and $B_2=B_2(a_2)$, with any choice of $a_3$, with $b_{12}=b$,  with any choices for $b_{13}$ and $b_{23}$, and with all instances of $\mu_{\beta_{23}}$ and $\mu_{\beta_{13}}$ replaced by $1$.  This implies 
$$
|\beta_{\calQ}(b)\cap K_2[B(a_1),B(a_2)]|=p^{-2n}|\beta_{\calQ}(b)||B(a_1)||B(a_2)|(1+O(p^{O(\ell+q)-\tau/2})).
$$
Combining this with the estimate for $|\beta_{\calQ}(b)|$ from part (i) yields the desired equality.  
\end{proof}

We will also use the following lemma, which says the balanced indicator function of a quadratic atom is locally $U^2$-uniform. This is Lemma 3.6 in \cite{Terry.2021d}.

\begin{lemma}\label{lem:quadatuni}
Let $\rho>0$, let $(\calL,\calQ)$ be a quadratic factor on $\F_p^n$ of complexity $(\ell,q)$ and rank at least $\rho$, and let $A=B(c;d)$ be an atom of $(\calL,\calQ)$. Then, for any tuple $d=(a_1,a_2)\in \F_p^{\ell}\times \F_p^{\ell}$ with $a_1+a_2=c$,
\[\|1_A-\alpha\|_{U^2(d)}=O(p^{-(\rho-2\ell)/2}),\]
where $\alpha=|A\cap L(c)|/|L(c)|$.
\end{lemma}

We now deduce that bipartite graphs of the form $\Delta^{(2)}_{\calB}(e)$ are quasirandom, whenever the factor $\calB$ has sufficiently high rank.\footnote{See also Lemma 3.13 in \cite{Terry.2021av2}.} We will use the well known fact that for any sets $U,V$ and function $f:U\times V\rightarrow \mathbb{R}$, the quantity $\sum_{u_0,u_1\in U}\sum_{v_0,v_1\in V}\prod_{i\in \{0,1\}}\prod_{j\in \{0,1\}}f(u_i,v_i)$ is non-negative.  Indeed, this follows from its equality with $\sum_{u_0,u_1}(\sum_{v\in V}f(u_0,v)f(u_1,v))^2$.

\begin{lemma}\label{lem:formqr}
For any $\e_2:\mathbb{R}^+\rightarrow (0,1]$, there is a growth function $\tau$ such that the following holds.  Suppose $\calB=(\calL,\calQ)$ is a quadratic factor on $\F_p^n$ of complexity $(\ell,q)$ and rank at least $\tau(\ell+q)$.  Then for any  $e\in \F_p^{\ell+q}\times \F_p^{\ell+q}\times \F_p^q$, the bipartite graph $\Delta_{\calB}^{(2)}(e)$ has $\dev_2(\e_2(p^{\ell+q}); p^{-q})$.  

Moreover,  it suffices to choose $\tau(x)$ a sufficiently fast growing polynomial in $x$ and $\e_2(x)$.  
\end{lemma}
\begin{proof}
Let $\tau$ grow sufficiently fast compared to $\e_2$ (we will see from the proof below $\tau(x)$ can be taken to be polynomial in $x$ and $\e_2(x)$). Suppose $\calB=(\calL,\calQ)$ is a quadratic factor on $\F_p^n$ of complexity $(\ell,q)$ and rank $\rho\geq \tau(\ell+q)$.  Fix 
$$
e=(a_1,b_1,a_2,b_2,b_{12})\in \F_p^{\ell}\times \F_p^{q}\times  \F_p^{\ell}\times \F_p^{q}\times \F_p^q.
$$
  Let $B_1=B(a_1;b_1)$ and $B_2=B(a_2;b_2)$, and let $c,d$ be such that $B(\Sigma(e))=B(c;d)$.  Setting 
  $$
  \alpha=\frac{|\beta_{\calQ}(b_{12})\cap K_2[B_1,B_2]|}{|B_1||B_2|},
  $$
   our goal is to show that $|\alpha-p^{-q}|\leq \e_2(p^{\ell+q})$ and that (\ref{bgoal}) below holds, where $g$ is defined by setting $g(x,y)=1-\alpha$ if $(x,y)\in \beta_{\calQ}(b)$ and otherwise $g(x,y)=-\alpha$.  
 \begin{align}\label{bgoal}
 \sum_{u_0,u_1\in B_1}\sum_{v_0,v_1\in B_2}\prod_{i\in \{0,1\}}\prod_{j\in \{0,1\}}g(u_i,v_j) \leq \e_2(p^{\ell+q})|B_1|^2|B_2|^2.
 \end{align}
 We already know $|\alpha-p^{-q}|\leq O(p^{O(\ell+q)-\rho/2})\leq \e_2(p^{\ell+q})$ by Lemma \ref{lem:sizeofbetaref} and our choice of $\tau$, so we just need to show (\ref{bgoal}).  

To ease notation, let $L_1=L(a_1)$ and $L_2=L(a_2)$.  Define $A$ to be the atom $B(c;d)$, and let $\alpha'=|A\cap L(c)|/|L(c)|$. Note that by Lemma \ref{lem:sizeofatoms}, $|\alpha'-p^{-q}|\leq O(p^{O(\ell+q)-\rho/2})$, and thus, combining with the above estimate on $\alpha$ we have 
\begin{align}\label{alphaprime}
|\alpha'-\alpha|\leq |\alpha-p^{-q}|+|\alpha'-p^{-q}|=O(p^{O(\ell+q)-\rho/2}).
\end{align}
  By Lemma \ref{lem:quadatuni}, 
\begin{align}\label{u2d}
\|1_A-\alpha'\|_{U^2(c)}=O(p^{-(\rho-2\ell)/2}).
\end{align}
Since $\|\cdot \|_{U^2(d)}$ satisfies the triangle inequality (see \cite{Terry.2021d}), we can conclude that 
\begin{align}\label{u2dd}
\|1_A-\alpha\|_{U^2(c)}\leq \|1_A-\alpha'\|_{U^2(c)}+\|\alpha-\alpha'\|_{U^2(c)}= \|1_A-\alpha'\|_{U^2(c)}+|\alpha-\alpha'|\leq O(p^{O(\ell+q)-\rho/2}),
\end{align}
where the second inequality uses  the definition of $\| \cdot \|_{U^2(c)}$, and the last inequality uses (\ref{alphaprime}) and (\ref{u2d}). Using the definition $\|\cdot \|_{U^2(d)}$, we can rewrite (\ref{u2d}) as
\begin{align}\label{u2ddd}
\sum_{f_1,f_2\in \F_p^q}\Big(\sum_{u_0,u_1\in B(a_1;f_1)}\sum_{v_0,v_1\in B(a_2;f_2)}\prod_{i\in \{0,1\}}\prod_{j\in \{0,1\}}1_A(u_i+v_j)-\alpha\Big)\leq O(p^{O(\ell+q)-\rho(\ell+q)/2})|L_1|^2|L_2|^2.
\end{align}
By the comment before Lemma \ref{lem:formqr}, we know that for each $f_1,f_2\in \F_p^q$, the above sum over $u_0,u_1\in B(a_1;f_1)$ and $v_0,v_1\in B(a_2;f_2)$ is non-negative.  Consequently, we can conclude that all such sums are at most the right hand side of (\ref{u2ddd}), so in particular,
$$
\sum_{u_0,u_1\in B(a_1;b_1)}\sum_{v_0,v_1\in B(a_2;b_2)}\prod_{i\in \{0,1\}}\prod_{j\in \{0,1\}}1_A(u_i+v_j)-\alpha\leq  O(p^{O(\ell+q)-\rho/2})|L_1|^2|L_2|^2.
$$
By definition of $A$, we have that for all $u\in B(a_1;b_1)$ and $v\in B(a_2;b_2)$, $1_A(u+v)=1$ if and only if $(u,v)\in \beta_{\calQ}(b_{12})$.  Thus, we have in fact shown  
\begin{align*}
\sum_{u_0,u_1\in B(a_1;b_1)}\sum_{v_0,v_1\in B(a_2;b_2)}\prod_{i\in \{0,1\}}\prod_{j\in \{0,1\}}g(u_i,v_j)&\leq O(p^{O(\ell+q)-\rho/2})|L_1|^2|L_2|^2\\
&\leq O(p^{O(\ell+q)-\rho/2})|B_1|^2|B_2|^2\\
&\leq \e_2(p^{\ell+q})|B_1|^2|B_2|^2,
\end{align*}
where the last two inequalities use Lemma \ref{lem:sizeofatoms} and the fact $\tau$ grows sufficiently fast.
\end{proof}

\subsection{Triads arising from quadratic factors}

In this section we prove several important properties regarding the behavior of ternary sum-graphs on the types of triads appearing in Theorem \ref{thm:hypergraphreg}.  We begin by setting some compact notation for such triads.   

\begin{definition}\label{def:triads2}
Let $\calB=(\calL,\calQ)$ be a quadratic factor  on $\F_p^n$ of complexity $(\ell,q)$. Given $d=(a_1,a_2,a_3,b_{12},b_{13},b_{23})\in \mathbb{F}_p^{\ell+q}\times \mathbb{F}_p^{\ell+q}\times \mathbb{F}_p^{\ell+q}\times \mathbb{F}_p^{q}\times \mathbb{F}_p^{q}\times \mathbb{F}_p^{q}$, define the \emph{triad corresponding to $d$} to be the $3$-partite graph
$$
\Delta_{\calB}^{(3)}(d):=(B(a_1)\cup B(a_2)\cup B(a_3), \beta_{12}\cup \beta_{13}\cup \beta_{23}), 
$$
where for each $1\leq i,j\leq 3$, $\beta_{ij}=\beta_{\calQ}(b_{ij})\cap K_2[B(a_i)\cup B(a_j)]$.  Note that if $xyz\in \Delta_{\calB}^{(3)}(d)$, then $x+y+z$ is in $B(\Sigma_{\calB}(d))$, for
$$
\Sigma_{\calB}(d):=a_1+a_2+a_3+2(0b_{12})+2(0b_{13})+2(0b_{23}),
$$
where $0b_{ij}$ denotes the vector obtained by adding $\ell$ many $0$'s to the left of $b_{ij}$.
\end{definition}
We will suppress the subscripts in Definition \ref{def:triads2} when ambiguity does not arise. It will be useful to have the following estimate on the number of triangles in such a triad.

\begin{lemma}\label{lem:qrcount}
Let $\calB=(\calL,\calQ)$ be a quadratic factor  on $\F_p^n$ of complexity $(\ell,q)$ and rank $\rho$. For any $d\in \mathbb{F}_p^{\ell+q}\times \mathbb{F}_p^{\ell+q}\times \mathbb{F}_p^{\ell+q}\times \mathbb{F}_p^{q}\times \mathbb{F}_p^{q}\times \mathbb{F}_p^{q}$,
$$
|K_3^{(2)}(\Delta^{(3)}(d))|=(1+O(p^{O(\ell+q)-\rho/2}))p^{n-3\ell-3q}.
$$
\end{lemma}
\begin{proof}
This is immediate from combining Lemma A.2 in \cite{Terry.2021d} with Lemmas \ref{lem:sizeofatoms} and  \ref{lem:sizeofbetaref}.  
\end{proof}

We will need the following lemma which allows us to essentially interchange the density of $E_A^{(3)}$-edges on a triad from Definition \ref{def:triads} with the density of the set $A$ on certain atoms of the factor in question.  It is a corollary of results in \cite{Terry.2021d}.\footnote{See also Lemma 3.15 in \cite{Terry.2021av2}.}

\begin{lemma}\label{lem:formqr4}
Let $\calB=(\calL,\calQ)$ be a quadratic factor  on $\F_p^n$ of complexity $(\ell,q)$ and rank $\rho$.  Fix $d\in \mathbb{F}_p^{\ell+q}\times \mathbb{F}_p^{\ell+q}\times \mathbb{F}_p^{\ell+q}\times \mathbb{F}_p^{q}\times \mathbb{F}_p^{q}\times \mathbb{F}_p^{q}$ and $A\subseteq \F_p^n$, and let $\alpha$ and $\alpha_{B(\Sigma(d))}$ be such that 
$$
 |E_A^{(3)}\cap K_3^{(2)}(\Delta^{(3)}(d))|=\alpha |K_3^{(2)}(\Delta^{(3)}(d))|\text{ and } |A\cap B(\Sigma(d))|=\alpha_{B(\Sigma(d))}|B(\Sigma(d))|.
$$
Then 
\begin{enumerate}
\item  $|\alpha-\alpha_{B(\Sigma(d))}| =O(p^{O(\ell+q)-\rho/8})$;
\item $\|1_A-\alpha\|_{U^3(d)}=\|1_A-\alpha_{B(\Sigma(d))}\|_{U^3(d)}+O(p^{O(\ell+q)-\rho/8})$.
\end{enumerate}
\end{lemma}
\begin{proof}
Part (i) is immediate from \cite[Lemma A.2 and Proposition C.1]{Terry.2021d}. For part (ii), since $\|\cdot \|_{U^3(d)}$ satisfies the triangle inequality (see \cite{Terry.2021d}), we have 
\begin{align}\label{u3dtriangle}
\|1_A-\alpha\|_{U^3(d)}= \|1_A-\alpha_{\calB(\Sigma(d))}\|_{U^3(d)}+ \|\alpha-\alpha_{\calB(\Sigma(d))}\|_{U^3(d)}.
\end{align}
By definition of $\|\cdot \|_{U^3(d)}$, Lemma \ref{lem:sizeofatoms}, and Lemma \ref{lem:sizeofbetaref}, 
$$
\|\alpha-\alpha_{\calB(\Sigma(d))}\|_{U^3(d)}=|\alpha-\alpha_{\calB(\Sigma(d))}|(1+O(p^{O(\ell+q)-\rho/2})).
$$
Combining this with (i) and (\ref{u3dtriangle}) yields (ii).
\end{proof}

Using Lemma \ref{lem:formqr4}, we obtain the following translation between the local $U^3$-norm on an atom, and the $\dev_{2,3}$-regularity of the corresponding triad.

\begin{lemma}\label{lem:qrtranslate}
For all $\e_1>0$ and $\e_2:\mathbb{R}^+\rightarrow (0,1]$, there is a polynomial growth function $\tau$ with the following property. Suppose $A\subseteq \mathbb{F}_p^n$, $\calB=(\calL,\calQ)$ be a quadratic factor on $\F_p^n$ of complexity $(\ell,q)$ and rank at least $\tau(\ell+q)$, and  
$$
d=(a_1,a_2,a_3,b_{12},b_{23},b_{13})\in \mathbb{F}_p^{\ell+q}\times \F_p^{\ell+q}\times \F_p^{\ell+q}\times \F_p^q\times \F_p^q\times \F_p^q.
$$
Then the following implications hold for the triad  $G=\Delta^{(3)}(d)$  and the $3$-graph $H=(B(a_1)\cup B(a_2)\cup B(a_3), E_A^{(3)}\cap K_3^{(2)}(G))$.
\begin{enumerate}[label=\normalfont(\roman*)]
\item If $\|1_A-\alpha_{B(\Sigma(d))}\|_{U^3(d)}<\e_1 $, then $(H, G)$ has $\dev_{2,3}(2\e_1,\e_2(p^{-q}))$.
\item If $(H, G)$ has $\dev_{2,3}(\e_1,\e_2(p^{-q}))$, then $\|1_A-\alpha_{B(\Sigma(d))}\|_{U^3(d)}< 2 \e_1$.
\end{enumerate}
Moreover, $\tau(x)$ can be taken as a polynomial in $x$, $\e_2(x)$, and $\e_1$.
\end{lemma}
\begin{proof}
Fix $\e_1>0$ and  $\e_2:\mathbb{R}^+\rightarrow (0,1]$.  Let $\tau$ grow sufficiently fast compared to the function $\e_2$ and the constant $\e_1$.   Suppose $A\subseteq \mathbb{F}_p^n$, $\calB=(\calL,\calQ)$ is a quadratic factor on $\F_p^n$ of complexity $(\ell,q)$ and rank $\rho\geq \tau(\ell+q)$, and fix
$$
d=(a_1,a_2,a_3,b_{12},b_{23},b_{13})\in \mathbb{F}_p^{\ell+q}\times \F_p^{\ell+q}\times \F_p^{\ell+q}\times \F_p^q\times \F_p^q\times \F_p^q.
$$
To ease notation, let $B_i=B(a_i)$ for each $i\in [3]$, and let $\beta_{ij}=\beta_{\calQ}(b_{ij})$ for each $1\leq i<j\leq 3$. Let $\alpha$ and $\alpha_{\Sigma(d)}$ be such that 
$$
|E_A^{(3)}\cap K_3^{(2)}(G)|=\alpha|K_3^{(2)}(G)|\text{ and } |A\cap B(\Sigma(d))|=\alpha_{B(\Sigma(d))}|B(\Sigma(d))|.
$$ 
By Lemma \ref{lem:formqr4}, 
\begin{align}\label{u3d1}
\Big|\|1_A-\alpha\|_{U^3(d)}- \|1_A-\alpha_{\Sigma(d)}\|_{U^3(d)}\Big|\leq O(p^{O(\ell+q)-\rho/8}).
\end{align}
Let $h(x,y,z)=1-\alpha$ for $xyz\in K_3^{(2)}(G)$ and $h(x,y,z)=-\alpha$ for $xyz\notin K_3^{(2)}(G)$.  Note that for all $xyz\in K_3^{(2)}(G)$, $h(x,y,z)=1-\alpha$ if $x+y+z\in A$ and otherwise, $h(x,y,z)=-\alpha$.  Thus, combining with the definition of $G$, we can rewrite 
  \begin{align*} 
 &\sum_{u_0,u_1\in B_1}\sum_{v_0,v_1\in B_2}\sum_{w_0,w_1\in B_3}\prod_{i,j,k\in \{0,1\}}h(u_i,v_j,w_k)
 \end{align*}
 as 
 \begin{align*}
 \sum_{u_0,u_1\in B_1}\sum_{v_0,v_1\in B_2}\sum_{w_0,w_1\in B_3} \prod_{(i,j)\in \{0,1\}^2} 1_{\beta_{12}}(u_i,v_j)\prod_{(i,k)\in \{0,1\}^2}1_{\beta_{13}}(u_i,w_k) \prod_{(j,k)\in \{0,1\}^2}1_{\beta_{23}}(v_j,w_k)\cdot &\\
  \prod_{i,j,k\in \{0,1\}}(1_A(u_i+v_j+w_k)-\alpha)&.
\end{align*}  
Combining with the definition of $\|\cdot \|_{U^3(d)}$,  we can conclude that 
\begin{align*}
\sum_{u_0,u_1\in B_1}\sum_{v_0,v_1\in B_2}\sum_{w_0,w_1\in B_3}\prod_{i,j,k\in \{0,1\}}h(u_i,v_j,w_k)=\frac{|B_1|^2|B_2|^2|B_3|^2|\beta_{12}|^4|\beta_{13}|^4|\beta_{23}|^4}{p^{24 n}} \|1_A-\alpha\|_{U^3(d)}&\\
=\frac{|B_1|^2|B_2|^2|B_3|^2|\beta_{12}|^4|\beta_{13}|^4|\beta_{23}|^4}{p^{24 n}} \Big(\|1_A -\alpha_{\Sigma(d)}\|_{U^3(d)}+O(p^{O(\ell+q)-\rho/8})\Big)&,
\end{align*}
where the final inequality uses (\ref{u3d1}).  Combining with Lemmas \ref{lem:sizeofatoms} and \ref{lem:sizeofbetaref}, this implies 
\begin{align}\label{u3dd}
\sum_{u_0,u_1\in B_1}\sum_{v_0,v_1\in B_2}\sum_{w_0,w_1\in B_3}\prod_{i,j,k\in \{0,1\}}h(u_i,v_j,w_k)&=\frac{|B_1|^2|B_2|^2|B_3|^2}{p^{12q}}\Big(\|1_A -\alpha_{\Sigma(d)}\|_{U^3(d)}+O(p^{O(\ell+q)-\rho/8})\Big).
\end{align}
 
We now show (i).  Assume $\|1_A-\alpha_{B(\Sigma(d))}\|_{U^3(d)}<\e_1$.  By our choice of $\tau$ and Lemma \ref{lem:formqr}, we know that for each $ij\in {[3]\choose 2}$, $G^{ij}:=\Delta^{(2)}(a_ia_jb_{ij})$ has $\dev_2(\e_2(p^{-q}); p^{-q})$.  Therefore, we can conclude that $(H, G)$ has $\dev_{2,3}(2\e_1,\e_2(p^{-q}))$, since our choice of $\tau$ and (\ref{u3dd}) imply
\begin{align*}
\sum_{u_0,u_1\in B_1}\sum_{v_0,v_1\in B_2}\sum_{w_0,w_1\in B_3}\prod_{i,j,k\in \{0,1\}}h(u_i,v_j,w_k)&\leq  (\|1_A-\alpha_{B(\Sigma(d))}\|_{U^3(d)}+\e_1) p^{-12q}|B_1|^2|B_2|^2|B_3|^2\\
&\leq 2\e_1  p^{-12q}|B_1|^2|B_2|^2|B_3|^2,
\end{align*} 
where the final inequality uses the assumption  $\|1_A-\alpha_{B(\Sigma(d))}\|_{U^3(d)}<\e_1$.

We now prove (ii).  To this end, assume $(H,G)$ has $\dev_{2,3}(\e_1,\e_2(p^{-q}))$.  Then combining the definition of $\dev_{2,3}(\e_1,\e_2(p^{-q}))$ with (\ref{u3dd}), we have 
\begin{align*}
 \|1_A-\alpha_{\Sigma(d)}\|_{U^3(d)}\leq \e_1+O(p^{O(\ell+q)-\rho/8})\leq 2\e_1,
\end{align*} 
where the last inequality is by our choice of $\tau$.  The moreover statement regarding our choice of $\tau$ is straightforward to see from the above argument, along with the moreover statement from Lemma \ref{lem:formqr}.
\end{proof}

We end this subsection by showing that given a set $A$ which is locally uniform on an atom $B$ of a high rank factor, the density of $A$ on any high rank sub-atom $B'\subseteq B$ is approximately the same as the density of $A$ on $B$ (Proposition \ref{prop:preservedensity} below). This is not needed for the main results of this section, but will be used in analyzing a lower bound construction in a forthcoming paper of the authors \cite{Terry.2025}, and seems generally useful for such purposes.  It is a straightforward corollary of the results in this section along with the analogous result in the hypergraph setting, which shows that in a sufficiently $\dev_{2,3}$-regular triad, hyergraph edge densities are approximately preserved when passing to large sub-triads. This was first shown several years ago \cite{Nagle.2013}, with the following a quantitative improvement appearing more recently in \cite{Nagle.2024}.  

\begin{proposition}[Proposition 1.4 in \cite{Nagle.2024}]\label{prop:disclem}
There is a polynomial $\mu(x)$ so that for all $0<\e_1,d_1,d_2<1$ and all $0<\e_2<\mu(d_2)$, the following holds.  Let $H=(V_1\cup V_2\cup V_3, F)$ be $3$-graph underlied by a tripartite graph $G=(V_1\cup V_2\cup V_3,E)$ with  $|F|=d_1|K_3^{(2)}(G)|$. Assume $(H, G)$ satisfies $\dev_{2,3}(\e_1,\e_2)$, and for each $1\leq i<j\leq 3$,   $|E[V_i,V_j]|=d_{ij}|V_i||V_j|$ for some  $d_{ij}\geq d_2$. The for any subgraph $G'$ of $G$,   
$$
\Big||E(H)\cap K^{(2)}_3(G')|-d_1|K^{(2)}_3(G')|\Big|\leq (2\e_1^{1/8})d_{12}d_{13}d_{23}|V_1||V_2||V_3|.
$$
\end{proposition}

 We note the statement above is slightly more general than what appears in \cite{Nagle.2024}, but clearly follows from the proof there.  We now prove an analogue for sets which are locally uniform on high rank quadratic atoms.

\begin{proposition}\label{prop:preservedensity}
 For all $\e\in (0,1)$, there is a polynomial growth function $\rho$ such that the following holds.  Assume $A\subseteq \F_p^n$, $\calB=(\calL,\calQ)$ is a quadratic factor on $\F_p^n$ of complexity $(\ell, q)$ and rank at least $\rho(\ell+q)$, and 
 $$
 d=(a_1,a_2,a_3,b_{12},b_{13},b_{23})\in \F_p^{\ell+q}\times \F_p^{\ell+q}\times\F_p^{\ell+q}\times\F_p^{q}\times\F_p^{q}\times\F_p^{q}
 $$
  is such that $\|1_A-\alpha_{B(\Sigma(d))}\|_{U^3(d)}<\e$.  
 
 Suppose $\calB'=(\calL',\calQ')$ is a quadratic factor of complexity $(\ell',q')$ and rank at least $\rho(\ell+q)$  satisfying $\calL'\supseteq \calL$, $\calQ'\supseteq \calQ$, and $\ell'+q'\leq \ell+q+\frac{1}{3}\log_p(\e^{-1/32})$, and suppose $(u,v)\in \F_p^{\ell}\times \F_p^{q}$ is such that $\Sigma_{\calB}(d)=(u,v)$. Then for any $(u',v')\in \F_p^{\ell'-\ell}\times \F_p^{q'-q}$, 
 $$
 \Bigg|\frac{|A\cap B(u;v)|}{|B(u;v)|}-\frac{|A\cap  B'(uu';vv')|}{|B'(uu';vv')|}\Bigg|\leq 5\e^{1/32}.
 $$ 
\end{proposition}
\begin{proof}
Fix $\e\in (0,1)$. Let $\mu(x)$ be the polynomial is Proposition \ref{prop:disclem}, and let $\tau$ grow sufficiently fast compared to the function $\mu(x)$ and the constant $\e$ (we will see in the proof below that $\tau$ may be taken to be polynomial).  Fix $A\subseteq \F_p^n$ and $\calB=(\calL,\calQ)$ a quadratic factor  on $\F_p^n$ of complexity $(\ell, q)$ and rank $\rho\geq \tau(\ell+q)$, and suppose $d$ is such that $\|1_A-\alpha_{B(\Sigma(d))}\|_{U^3(d)}<\e$.  

Let $\calB'=(\calL',\calQ')$ be a quadratic factor of complexity $(\ell',q')$ and rank at least $\rho$, satisfying $\calL'\supseteq \calL$, $\calQ'\supseteq \calQ$, and $\ell'+q'\leq \ell+q+\frac{1}{3}\log_p(\e^{-1/32})$.  Now fix any $(u,v)\in \F_p^{\ell}\times \F_p^{q}$ such that $\Sigma_{\calB}(d)=(u,v)$, and any $(u',v')\in \F_p^{\ell'-\ell}\times \F_p^{q'-q}$. Our goal is to show that 
\[\frac{|A\cap B(u;v)|}{|B(u;v)|}=\frac{|A\cap  B'(uu';vv')|}{|B'(uu';vv')|}\pm 5\e^{1/32}.\]

Fix any $u_1,u_2,u_3\in \F_p^{\ell}$ such that $u_1+u_2+u_3=u$, any $v_1,v_2,v_3\in \F_p^q$ such that $v_1+v_2+v_3=v$, and any $w_1,w_2,w_3\in \F_p^{\ell'-\ell}$ with the property that $w_1+w_2+w_3=u'$. Define $d'$ to be the following element of $\F_p^{\ell'+q'}\times \F_p^{\ell'+q'}\times \F_p^{\ell'+q'}\times \F_p^{ q'}\times \F_p^{ q'}\times \F_p^{q'}$.
$$
d':=(u_1w_1v_10, u_2w_3v_20, u_3w_3 v_3v', b_{12}0, b_{13}0, b_{23}0).
$$
Note $\Sigma_{\calB'}(d')=(uu'; vv')$, and thus $B'(\Sigma_{\calB'}(d'))=B'(uu';vv')$.  To ease notation, let $B_i=B(u_i;v_i)$ for each $i\in [3]$.

Set $G=\Delta_{\calB}^{(3)}(d)$ and $H=(B_1\cup B_2\cup B_3, E_A^{(3)}\cap K_3^{(2)}(G))$.  By  Lemma \ref{lem:qrtranslate} and our choice of $\rho$, the assumption $\|1_A-\alpha_{B(\Sigma(d))}\|_{U^3(d)}<\e$ implies $(H, G)$ has $\dev_{2,3}(2\e, \mu(p^q))$.  Now set $G'=\Delta_{\calB'}^{(3)}(d')$, and observe that by construction, $G'$ is a sub-graph of $G$.  Let $\alpha$ and $\alpha'$ be such that 
$$
|E(H)\cap K_3(G)|=\alpha |K_3(G)|\text{ and }|E(H)\cap K_3(G')|=\alpha' |K_3(G')|.
$$
  By Lemma \ref{lem:formqr4} applied first to $\Delta^{(3)}_{\calB}(d)$ and then to $\Delta^{(3)}_{\calB'}(d')$, we have
  \begin{align*} 
   |\alpha-\alpha_{\Sigma_{\calB}(d)}|\leq O(p^{O(\ell+q)-\rho/8})\text{ and }|\alpha'-\alpha_{\Sigma_{\calB'}(d')}|\leq O(p^{O(\ell'+q')-\rho/8}).
   \end{align*}
Further, by Proposition \ref{prop:disclem} and Lemma \ref{lem:qrcount},
\begin{align*}
|\alpha'-\alpha |\leq \frac{(4\e^{1/8})p^{-3q}|B_1||B_2||B_3|}{|K_3(G')|}&\leq (4\e^{1/8})|B_1|^2|B_2|^2|B_3|^3\\
&=(4\e^{1/8})p^{3(\ell'-\ell+q'-q)}(1+O(p^{O(\ell'+q')-\rho/2})),
\end{align*}
where the last equality uses Lemma \ref{lem:sizeofatoms} and the fact $\ell'\geq \ell$ and $q'\geq q$.  Combining these bounds with the triangle inequality, we have  
 \begin{align*}
 |\alpha_{\Sigma_{\calB}(d)}-\alpha_{\Sigma_{\calB'}(d')}|&\leq  |\alpha-\alpha_{\Sigma_{\calB}(d)}|+|\alpha'-\alpha_{\Sigma_{\calB'}(d')}|+|\alpha'-\alpha|\\
 &\leq O(p^{O(\ell+q)-\rho/8})+O(p^{O(\ell'+q')-\rho/8}) +(4\e^{1/8})p^{3(\ell'-\ell +q'-q)}(1+O(p^{O(\ell'+q')-\rho/2}))\\
 &\leq O(p^{O(\ell'+q')-\rho/8}) +(4\e^{1/8})\e^{-3/32} (1+O(p^{O(\ell'+q')-\rho/2}))\\
 &\leq 5\e^{1/32},
 \end{align*}
where the second to last inequality uses our assumption that $\ell'-\ell +q'-q\leq \frac{1}{3}\log_p(\e^{-1/32})$, and the last inequality uses that $\rho$ is sufficiently large.
\end{proof}

\subsection{Special hypergraph regularity for ternary sum-graphs}

In this section we prove our special regularity lemma for ternary sum-graphs, Theorem \ref{thm:hypergraphreg}, after which we consider hypergraph analogues of Theorems \ref{thm:vc2finite} and \ref{thm:fop}.   The main ingredient, in addition to the lemmas of the preceding subsections, is  the following general quadratic arithmetic regularity lemma for subsets of $\F_p^n$, which first appeared in earlier versions of this paper, and has now been moved to \cite[Proposition 4.2]{Terry.2021d}.

 \begin{proposition}\label{prop:unifatom2}
 For all $\e\in (0,1)$, there is a growth function $\omega_1=\omega_1(\e)$ such that for all growth functions $\omega\geq\omega_1$ and all $\ell_0,q_0\geq 0$, there is $D=D(\e,\omega,\ell_0,q_0)$ such that the following holds.  Suppose $n$ is sufficiently large, $A\subseteq \F_p^n$, and $\calB_0=(\calL_0,\calQ_0)$ is a quadratic factor on $\F_p^n$ of complexity $(\ell_0,q_0)$. Then there exist $\ell,q\leq D$, a quadratic factor $\calB=(\calL,\calQ)\preceq \calB_0$ of complexity $(\ell, q)$, and a set $\Gamma\subseteq \mathbb{F}_p^{\ell+q}\times \mathbb{F}_p^{\ell+q}\times\mathbb{F}_p^{\ell+q}\times\mathbb{F}_p^{q}\times\mathbb{F}_p^{q}\times\mathbb{F}_p^{q}$ such that 
\begin{enumerate}[label=\normalfont(\roman*)]
\item  $\calB$ has rank at least $\omega(\ell+q)$;
\item $|\Gamma|\geq (1-\e^{2}/2^{29})p^{3\ell+6q}$;
\item for all $d\in \Gamma$, $\|1_A-\alpha_{B(\Sigma(d))}\|_{U^3(d)}<\e$.
\end{enumerate}
\end{proposition}

We now prove Theorem \ref{thm:hypergraphreg}.\footnote{See Corollary 3.27 in \cite{Terry.2021av2}.}  

\vspace{3mm}

\begin{proofof}{Theorem \ref{thm:hypergraphreg}} 
Fix $\e_1>0$,  $\e_2:\mathbb{R}^+\rightarrow (0,1]$, $s_0,t_0\geq 1$, and a growth function $\omega$.  To ease notation, let $\e=\e_1^2/4$.  Set $\ell_0=\max\{s_0, (4\e)^{-2}\}$, and choose $\rho$ to be a sufficiently fast growth function, depending on the constants $\e, \ell_0$ and the functions $\omega$ and $\e_2$.  Let $D=D(\e,\omega,\ell_0,q_0)$ be from Proposition \ref{prop:unifatom2}, and set $T_0=S_0=p^{2D}$.   

Assume $n$ is sufficiently large and $A\subseteq \F_p^n$, and let $H=(\F_p^n, E_A^{(3)})$ be the ternary sum-graph associated to $A$.  Choose any quadratic factor $\calB_0$ of complexity $(\ell_0,q_0)$.  By Proposition \ref{prop:unifatom2}, there are $\ell,q\leq D$, a quadratic factor $\calB=(\calL,\calQ)\preceq \calB_0$ of complexity $(\ell, q)$ and rank at least $\omega(\ell+q)$, and $\Gamma\subseteq \mathbb{F}_p^{3\ell+6q}$ with $|\Gamma|\geq (1-\e^2/2^{29})p^{3\ell+6q}$, such that for all $d\in \Gamma$,  
$$
\|1_A-\alpha_{B(\Sigma(d))}\|_{U^3(d)}<\e.
$$

Let $s=p^{q}$ and $t=p^{\ell+q}$.  Let $V_1,\ldots, V_t$ be any enumeration of $\At(\calB)$, and for each $ij\in {[t]\choose 2}$, let $\{P_{ij}^{\alpha}: \alpha\in [\ell]\}$ be any enumeration of $\{K_2[V_i,V_j]\cap \beta_{\calQ}(b): b\in \F_p^{q}\}$.  Note the triads of $\calP$ are exactly the triads of the form $\Delta^{(3)}(d)$ for $d\in \F_p^{3\ell+6q}$ from Definition \ref{def:triads2}.  By Lemma \ref{lem:formqr}, for every $ij$, $(V_i\cup V_j, P_{ij}^\alpha)$ has $\dev_2(\e_2(s);1/s)$, and by  Lemma \ref{lem:qrtranslate}, for all $d\in \Gamma$, $\Delta^{(3)}(d)$ has $\dev_{2,3}(\e_1,\e_2(s))$ with respect to $H$.  Thus, it suffices to show 
$$
\sum_{d\in \Gamma'}|K_3^{(2)}(d)|\geq (1-\e_1){p^n\choose 3},
$$
for some $\Gamma'\subseteq \Gamma$ satisfying the following:  for all $d=(a_1,a_2,a_3,b_{12},b_{23},b_{13})\in \Gamma'$ the  $a_1,a_2,a_3$ are pairwise distinct, and for all $d\neq d'$ in $\Gamma$,  $\Delta^{(3)}(d)$ and $\Delta^{(3)}(d')$ are distinct graphs.  

To construct such a set $\Gamma'$, we first observe that given $d\in \Gamma$, there are at most $3!$ many $d'\in \Gamma$ so that $\Delta^{(3)}(d)$ and $\Delta^{(3)}(d')$ are the same graphs.  Thus we can choose a subset $\Gamma_1\subseteq \Gamma$ with  $|\Gamma_1|\geq |\Gamma|/3!$ so that $\Gamma_1$ contains one representative from each such equivalence class.  We then obtain $\Gamma'$ from $\Gamma_1$ by deleting any tuples $(a_1,a_2,a_3,b_{12},b_{23},b_{13})$  where $a_1,a_2,a_3$ are not pairwise distinct. Observe 
\begin{align*}
|\Gamma'|\geq \frac{|\Gamma|}{3!}-|\{(a_1,a_2, a_3,b_{12},b_{23},b_{13})\in p^{2\ell+3q}: |\{a_1,a_2,a_3\}|\neq 3\}|&\geq \frac{|\Gamma|}{6}-3p^{2\ell+5q} \\
&\geq (1-\frac{\e^2}{2^{29}}- 18p^{-\ell-q})\frac{p^{3\ell+6q}}{6}\\
&\geq (1-\e^2)\frac{p^{3\ell+6q}}{6},
\end{align*}
where the last inequality uses our choice of $\ell_0$ and the fact that $\calB\preceq\calB_0$ implies $\ell+q\geq \ell_0$.  We can now compute, using Lemma \ref{lem:qrcount} and our choice of $\tau$, that
\begin{align*}
\sum_{d\in \Gamma'}|K_3^{(2)}(\Delta^{(3)}(d))|\geq (1-\e^2)|\Gamma'|p^{3n-3\ell-3q}\geq (1-\e^2)^2\frac{p^{3n}}{6}\geq (1-\e_1){p^n\choose 3}.
\end{align*}
This finishes the proof.
\end{proofof}

\vspace{2mm}

We next  obtain a stronger version of hypergraph regularity for tenary sum-graphs associated to $\NIP_2$ subsets of $\F_p^n$.  Given a $(t,s,\e_1,\e_2(s))$-decomposition $\calP$ of a set $V$, and a $3$-graph $H=(V,E)$, we say $\calP$ is \emph{$\e_1$-homegeneous with respect to $H$} if 
$$
|\bigcup_{G\in \Sigma}K_3^{(2)}(G)|\geq (1-\e_1){|V|\choose 3},
$$
where $\Sigma$ is the set of triads of $\calP$ where the density of $H$ on $K_3^{(2)}(G)$ is in  $[0,\e_1)\cup (1-\e_1,1]$.  In \cite[Theorem 2.23]{Terry.2021b} we showed that $3$-graphs of bounded $\VC_2$-dimension admit $\dev_{2,3}$-regular, homogeneous decompositions. We easily obtain the same result here for ternary sum-graphs of sets with bounded $\VC_2$-dimension, with additional algebraic structure built into the decompositions.

\begin{corollary}[$\NIP_2$ hypergraph regularity with quadratic partitions]\label{cor:hypergraphvc2}
For all $k\geq 1$, $\e_1>0$, $\e_2:\mathbb{R}^+\rightarrow [0,1]$, $t_0,s_0\geq 1$ and growth functions $\tau$, there are $S_0$ and $T_0$ such that the following holds. Let $n$ be sufficiently large, and suppose that $A\subseteq \F_p^n$ has no $k$-$\IP_2$.   Then there exists $t_0\leq t\leq T_0$, $s_0\leq s\leq S_0$, and a $(t,s,\e_1,\e_2(s))$-decomposition $\calP$ of $\F_p^n$ which is $\e_1$-homogeneous and $(\e_1,\e_2(s))$-regular with respect to $(\F_p^n, E_A^{(3)})$.  

Moreover, there is a quadratic factor $\calB=(\calL,\calQ)$ on $\F_p^n$ of complexity $(\ell,q)$ and rank at least $\tau(\ell+q)$, such that $\calP_{vert}=\At(\calB)$ and $\calP_{edge}= \{K_2[V_i,V_j]\cap \beta_{\calQ}(b): b\in \F_p^{q}\}$. 
\end{corollary}
\begin{proof}
Mimic the proof of Theorem \ref{thm:hypergraphreg} above, using Theorem \ref{thm:vc2finite} in place of Proposition \ref{prop:unifatom2}.  Note  Lemma \ref{lem:formqr4} will tell us the decomposition is $\e_1$-homogeneous, and  Proposition \ref{prop:locsparseuni} (combined with the same arguments from the proof of Theorem \ref{thm:hypergraphreg}) will tell us it is regular.
\end{proof}

In \cite{Terry.2021b} we also consider decompositions whose irregular triads can be controlled in a special way.  Specifically, given a $(t,s,\e_1,\e_2(s))$-decomposition $\calP$ of a set $V$ and a $3$-uniform hypergraph $H=(V,E)$, we say that $\calP$ \emph{has linear $\dev_{2,3}$-error} if there is a set $\Omega\subseteq {\calP_{vert}\choose 3}$ such that $|\Omega|\leq \e_1 t^3$, and every irregular triad of $\calP$ is contained in a triple $V_iV_jV_k$ from $\Omega$.   This nomenclature is no coincidence. If $V=\F_p^n$ and $\calP$ comes from a high rank factor as in Theorem \ref{thm:hypergraphreg}, then such sets $\Omega$ arise naturally from small sets of linear atoms.

\begin{fact}
For all $\e\in (0,1)$ there is a growth function $\omega$ so that if $\calB=(\calL,\calQ)$ is a quadratic factor on $\F_p^n$ of complexity $(\ell,q)$ and rank at least $\omega(\ell+q)$, then the following holds.   If $\Sigma \subseteq \At(\calL)$ satisfies   $|\Sigma|\leq \e p^{\ell}$ and 
$$
\Omega=\{B(a_1;a_2)B(b_1;b_2)B(c_1;c_2)\in {\At(\calB)\choose 3}: L(a_1+b_1+c_1)\in \Sigma\},
$$
then $|\Omega|\leq \e p^{3(\ell+q)}$.
\end{fact}
\begin{proof}
This is immediate from Lemma \ref{lem:sizeofatoms}.
\end{proof}

It is now straightforward to see that this paper's main result  about $\NFOP_2$ sets gives rise to regular decompositions with linear error for the corresponding ternary sum-graphs.

\begin{corollary}[$\NFOP_2$ hypergraph regularity with quadratic partition]\label{cor:hypergraphfop}
For all $k\geq 1$, $\e_1>0$, $\e_2:\mathbb{R}^+\rightarrow (0,1]$, $t_0,s_0\geq 1$, and growth functions $\omega$, there are $S_0$ and $T_0$ such that the following holds. Let $n$ be sufficiently large, and suppose that $A\subseteq \F_p^n$ has no $k$-$\FOP_2$.   Then there exists $t_0\leq t\leq T_0$, $s_0\leq s\leq S_0$, and a $(t,s,\e_1,\e_2(s))$-decomposition $\calP$ of $\F_p^n$ which is $\e_1$-homogeneous and $(\e_1,\e_2(s))$-regular with respect to $H=(\F_p^n, E_A^{(3)})$ and which has linear error.    

Moreover, there is a quadratic factor $\calB=(\calL,\calQ)$ on $\F_p^n$ of complexity $(\ell,q)$ and rank at least $\omega(\ell+q)$, such that $\calP_{vert}= \At(\calB)$ and $\calP_{edge}= \{K_2[V_i,V_j]\cap \beta_{\calQ}(b): b\in \F_p^{q}\}$. 

Further,  there is a set $\Sigma\subseteq \At(\calL)$ such that $|\Sigma|\leq \e p^{\ell}$ and such that any triad $G$ of $\calP$ which is irregular with respect to $H$ has vertex set of the form 
$$
B(a_1;a_2)\cup B(b_1;b_2)\cup B(c_1;c_2),
$$
 for some labels $a_1a_2,b_1b_2,c_1c_2\in \F_p^{\ell}\times \F_p^q$ satisfying $L(a_1+b_1+c_1)\in \Sigma$.
\end{corollary}
\begin{proof}
Like the proof of Corollary \ref{cor:hypergraphvc2} with  Theorem \ref{thm:fop} in place of Theorem \ref{thm:vc2finite}.
\end{proof}

In \cite[Theorem 2.50]{Terry.2021b}, we show that $\NFOP_2$ in fact \emph{characterizes} the hereditary properties of $3$-graphs which admit $\dev_{2,3}$-regular decompositions with linear error.  Thus Corollary \ref{cor:hypergraphfop} lines up closely with what occurs in the hypergraph setting.

\section{Lemma from Section \ref{sec:FOP}}\label{app:easy}

This section contains the proofs of some elementary lemmas used in Section \ref{sec:FOP}.  We begin with Fact \ref{fact:aveofavecor}, which shows that the property of being almost $\e$-atomic with respect to a set $A$ is approximately inherited when passing to refinements.  We include a restatement of this for the convenience of the reader. 

\begin{fact}\label{fact:aveofavecorrestate}
Suppose $A\subseteq \F_p^n$ and $\calY,\calY'$ are partitions of $\F_p^n$ with $\calY'$ refining $\calY$.  If $\calY$ is almost $\e$-atomic with respect to $A$, then $\calY'$ is almost $\sqrt{2\e}$-atomic with respect to $A$. 
\end{fact}
\begin{proof}
Set $\calY_1=\{Y\in \calY: |A\cap Y|\geq (1-\e)|Y|\}$,  $\calY_0=\{Y\in \calY: |A\cap Y|\leq \e |Y|\}$, and $\Sigma=\calY\setminus (\calY_0\cup \calY_1)$.  Now let 
$$
A'=\Big(\bigcup_{Y\in \calY_1}Y\setminus A\Big)\cup \Big(\bigcup_{Y\in \calY_0}Y\cap A\Big)\cup \Big(\bigcup_{Y\in\Sigma}  Y \Big).
$$
By assumption, 
$$
|A'|\leq \sum_{Y\in \calY\setminus \Sigma}\e|Y|+\sum_{Y\in \Sigma}|Y|\leq \e|\F_p^n|+\e|\F_p^n|=2\e|\F_p^n|.
$$
Now let $\Sigma'=\{Y'\in \calY': |Y'\cap A'|\geq \sqrt{2\e}|Y'|\}$.  Then   $2\e|\F_p^n|\geq |A'|\geq \sqrt{2\e}|\bigcup_{Y'\in \Sigma'}Y'|$, so  $|\bigcup_{Y'\in \Sigma'}Y'|\leq \sqrt{2\e}|\F_p^n|$.  On the other hand, given $Y'\in \calY'\setminus \Sigma'$, we have by definition that either $|A\cap Y'|\geq (1-\sqrt{\e})|Y'|$ or $|A\cap Y'|\leq \sqrt{\e}|Y'|$. This finishes the proof.  
\end{proof}

We now prove Lemma \ref{lem:padding}, restated as Lemma \ref{lem:paddingrestate} below, which allows us to add matrices to a high rank factor while maintaining the rank assumption.  In fact we will prove the following more detailed statement, which immediately implies Lemma \ref{lem:padding}.

\begin{lemma}\label{lem:paddingrestate}
For all integers $q,q'\geq 1$ and reals $r>0$, there is $N=N(q,q')$ such that the following holds.  Let $n\geq N+1$ be sufficiently large compared to $q,q'$ and $r$, let $\calS=\{M_1,\ldots, M_{N}\}$ be a purely quadratic factor on $\F_p^n$ of complexity $N$ and rank at least $n-1$, and let $\calQ$ be a purely quadratic factor on $\F_p^n$ of complexity $q$ and rank at least $r$.  Then there is a set $\calS'\subseteq \calS$ so that $\calQ\cup \calS'$ has complexity $q+q'$ and rank at least $r$.
\end{lemma}
\begin{proof}
We induct on $q\geq 1$.

\underline{Case $q=1$:}  Fix an integers $q'\geq 1$ and a real $r>0$, and set $N=N(1,q')=2q'$.  Suppose $n\geq \max\{2q'+1,r+1\}$,  $\calS=\{M_1,\ldots, M_{N}\}$ is a purely quadratic factor on $\F_p^n$ of rank at least $n-1$, and $\calQ=\{R\}$ is a purely quadratic factor on $\F_p^n$ of rank $\rho\geq r$. Since $N=2q'$, we can define $\calN_1=\{M_1,\ldots, M_{q'}\}\cup \calQ$ and $\calN_2=\{M_{q'+1},\ldots, M_{2q'}\}\cup \calQ$. We claim one of $\calN_1$ or $\calN_2$ has rank at least $r$.    Suppose towards a contradiction that both $\calN_1$ and $\calN_2$ have rank less than $r$.  Then there are $\lambda_1,\ldots, \lambda_{q'},\lambda\in \F_p$  not all zero, and $\lambda_{q'+1},\ldots, \lambda_{2q'},\lambda'\in \F_p$  not all zero, such that 
\begin{align*}
\rk(\lambda_1 M_1+\ldots +\lambda_{q'}M_{q'}+\lambda R)<r\text{ and }\rk(\lambda'_1 M_{q'+1}+\ldots +\lambda'_{2q'}M_{q'}+\lambda' R)<r.
\end{align*}
Since $\calQ$ has rank $\rho\geq r$, at least one of the $\lambda_i$ is nonzero,  and at least one of the $\lambda'_j$ is nonzero.  Since $\{M_1,\ldots, M_{N}\}$ has rank at least $n-1\geq r$, we must have that both $\lambda$ and $\lambda'$ are non-zero.  But now, by subadditivity, we have that
\begin{align*}
\lambda_1 M_1+\ldots +\lambda_{q'}M_{q'}+\lambda R-\frac{\lambda}{\lambda'}(\lambda'_1 M_{q'+1}+\ldots +\lambda_{2q'}M_{q'}+\lambda' R)&\\
=\lambda_1 M_1+\ldots +\lambda_{q'}M_{q'}+\frac{\lambda}{\lambda'}(\lambda'_1 M_{q'+1}+\ldots +\lambda_{2q'}M_{q'})&
\end{align*}
has rank less than $r$. Since at least one of the $\lambda_i$ is nonzero  and one of the $\lambda'_j$ is nonzero,  this contradicts that $\calQ$ has rank at least $r$.

\underline{Case $q+1$:}  Suppose $q\geq 1$ and assume by induction we have defined $N(q,q')$ for all $q'\geq 1$ so that the claim holds, and moreover, so that $N(q,q')\geq 2q'$.   Fix $q'\geq 1$.  Define  $h:\mathbb{Z}^{\geq 1}\rightarrow \mathbb{Z}^{\geq 1}$ by setting $h(x)=N(q,x)$, and given $t\geq 0$, let $h^{(t)}$ denotes the $t$-fold iterate of $h$.  We then set $N(q+1,q')=h^{(q+1)}(q')$.  

Fix $r>0$, and assume $n\geq \max\{N(q+1,q')+1, r+1\}$. Let $\calS=\{M_1,\ldots, M_{N}\}$ be a purely quadratic factor on $\F_p^n$ of complexity $N$ and rank at least $n-1$, and let $\calQ=\{R_1,\ldots, R_{q+1}\}$ be a purely quadratic factor on $\F_p^n$ of complexity $q+1$ and rank $\rho \geq r$.  For each $1\leq i\leq q+1$, set $\calQ_i=\calQ\setminus \{R_i\}$, and note each $\calQ_i$ has size $q$.

By our definition, $N=h(h^{(q)}(q'))=N(q,h^{(q)}(q'))$, so the induction hypothesis applied to $\calQ_1$ and $\calS$ implies there exists $\calS_1\subseteq \calS$ of size $h^{(q)}(q')$ so that $\calQ_1\cup \calS_1$ has rank at least $r$.  

Suppose $1\leq i<q$ and assume we have constructed $\calS_i\subseteq \calS$ so that $\calS_i\cup \calQ_i$ has rank   $r$ and $|\calS_i|=h^{(q-i+1)}(q')$.  Since $h^{(q-i+1)}(q')=N(q,h^{(q-i)}(q'))$, our induction hypothesis implies there is $\calS_{i+1}\subseteq \calS_i$ of size $h^{(q-i)}(q')$ so that $\calS_{i+1}\cup \calQ_{i+1}$ has rank at least $r$.

After $q$ steps, we will have constructed  $\calS_q \subseteq  \calS$  so that for each $i\in [q]$, $\calS_q\cup \calQ_i$ has rank at least $r$ and $|\calS_q|=h(q')=N(q,q')\geq 2q'$.  Choose two disjoint subsets $\calS',  \calS''$ of $\calS_q$, each of size $q'$.  Enumerate them as $\calS'=\{M_1',\ldots, M_{q'}'\}$ and $\calS''=\{M_1'',\ldots, M_{q'}''\}$.  We claim that one of $\calS'\cup \calQ$ or $\calS''\cup \calQ$  has rank at least $r$.    Suppose towards a contradiction this is not the case.  Then there are $\lambda'_1,\ldots, \lambda'_{q'},\rho'_1,\ldots, \rho'_q$ not all zero and $\lambda''_1,\ldots, \lambda''_{q'},\rho''_1,\ldots, \rho''_q$ not all zero such that
\begin{align*}
&\rk(\lambda'_1 M'_1+\ldots +\lambda'_{q'}M'_{q'}+\rho_1'R_1+\ldots+\rho_q'R_q)<r\text{ and }\\
&\rk(\lambda''_1 M''_1+\ldots +\lambda''_{q'}M''_{q'}+\rho_1''R_1+\ldots+\rho_q''R_q)<r.
\end{align*}
Since $\calS_q\cup \calQ_i$ has rank $r$ for each $i\in [q]$,  we must have that all of $\rho_1',\ldots, \rho_q'$ are non-zero and all of $\rho_1'',\ldots, \rho_q''$ are non-zero.  Further, at least one $\lambda'_i$ and at least one $\lambda_j''$ are non-zero, since $\calQ$ has rank at least $r$.  Therefore, by subadditivity,
\begin{align*}
&\lambda'_1 M'_1+\ldots +\lambda'_{q'}M'_{q'}+\rho_1'R_1+\ldots+\rho_q'R_q-\frac{\rho_1'}{\rho_1''}(\lambda''_1 M''_1+\ldots +\lambda''_{q'}M''_{q'}+\rho_1''R_1+\ldots+\rho_q''R_q)\\
&=\lambda'_1 M'_1+\ldots +\lambda'_{q'}M'_{q'}+\rho_2'R_2+\ldots+\rho_q'R_q-\frac{\rho_1'}{\rho_1''}(\lambda''_1 M''_1+\ldots +\lambda''_{q'}M''_{q'}+\rho_2''R_2+\ldots+\rho_q''R_q)
\end{align*}
has rank less than $r$.  This is a contradiction, since by assumption, $\calS_q\cup \calQ_1$ has rank $r$.  This finishes the proof.
\end{proof}

Finally, we prove Fact \ref{fact:basis}, restated below, which shows that we can always find enough matrices whose linear combinations have high rank. 

\begin{fact}\label{fact:basisrestate}
For all $n\geq 1$, there exists a purely quadratic factor $\{M_1,\ldots, M_{n-1}\}$ on $\F_p^n$ of complexity $n-1$ and rank at least $n-1$.
\end{fact}
\begin{proof}
If $n\geq 1$ is odd, then by  \cite[Lemma 3]{Dumas.2011}, there exists a purely quadratic factor $\{M_1,\ldots, M_n\}$ on $\F_p^n$ of rank $n$.  Taking $\{M_1,\ldots, M_{n-1}\}$ then fulfills the desired conclusion.

Suppose now $n\geq 2$ is even. By above, there exists a purely quadratic factor $\{N_1,\ldots, N_{n-1}\}$  in $\F_p^{n-1}$ of rank $n-1$.  Let $\psi_1,\ldots,\psi_{n-1}$ the linear transformations from $\F_p^{n-1}$ to $\F_p^{n-1}$ corresponding to $N_1,\ldots, N_{n-1}$.  Identify $\F_p^{n-1}$ with the subspace $H=\langle \ebar_n\rangle^{\perp}$ of $\F_p^n$ in the obvious way, and let $g_1,\ldots, g_{n-1}$ be the linear transformations from $H$ into $H$ corresponding to $\psi_1,\ldots, \psi_{n-1}$ under this identification.  For each $1\leq i\leq n-1$, let $f_i:\F_p^n\rightarrow \F_p^n$ be defined by setting
$$
f_i(a_1,\ldots,a_{n-1},a_n)=g_i(a_1,\ldots,a_{n-1},0).
$$
Each $f_i$ is clearly a linear transformation, so we can let $M_i$ be the $n\times n$ matrix satisfying $M_i\xbar=f_i(\xbar)$. Given $\lambda_1,\ldots, \lambda_{n-1}$ not all $0$, it is not difficult to see that the dimension of $\ker(\lambda_1M_1+\ldots+\lambda_{n-1}M_{n-1})$ in $\F_p^n$ is one plus the dimension of $\ker(\lambda_1N_1+\ldots+\lambda_{n-1}N_{n-1})$ in $\F_p^{n-1}$.  Since the dimension of  $\ker(\lambda_1N_1+\ldots+\lambda_{n-1}N_{n-1})$ in $\F_p^{n-1}$ is $0$, this implies the dimension of $\ker(\lambda_1M_1+\ldots+\lambda_{n-1}M_{n-1})$ is $1$.  Thus $\{M_1,\ldots, M_{n-1}\}$ has rank $n-1$.
\end{proof}

\section{Companion results to Section \ref{sec:context}}\label{app:generalprops}

This appendix contains proofs of two results complementary to the statements appearing in Section \ref{sec:context}.  Specifically we prove in this section that unions of quadratic atoms have no $\HOP_2$, and that the linear Green-Sanders example has no $4$-$\HOP_2$.

\subsection{Unions of quadratic atoms have no $k$-$\HOP_2$} \label{app:boolean}

In this section we prove Proposition \ref{prop:twotameunion}, which shows that a union of sets with no $2$-$\HOP_2$ has no $k$-$\HOP_2$ for some $k$ depending only on the number of sets in the union.  We will in fact prove this for general ternary relation symbols, not only those arising as ternary sum-graphs.

\begin{definition}[$\HOP_2$ for ternary relation symbols]\label{def:hoprel}
Suppose $R$ is a ternary relation symbol and $\calM$ is an $\{R\}$-structure.  We say $R$ \emph{has $\ell$-$\HOP_2$ in $\calM$} if there are $\{a_i,b_i,c_i: i\in [\ell]\}\subseteq M$ such that $\calM\models R(a_i,b_j,c_k)$ if and only if $i<j+k$.
\end{definition}

We will use throughout the following re-indexing of $\ell$-$\HOP_2$, which makes it clear that the role of all three variables is symmetric.

\begin{fact}\label{fact:reindex}
Suppose $R$ is a ternary relation and $\calM$ is an $\{R\}$-structure.  Then the following are equivalent.
\begin{enumerate}[label=\normalfont(\arabic*)]
\item There is a set $\{a_u,b_u,c_u:u\in [k]\}\subseteq M$ such that $\calM\models R(a_u, b_v,c_w)$ holds if and only if $u< v+w$.  
\item There is a set $\{a_u,b_u,c_u:u\in [k]\}\subseteq M$ such that $\calM\models R(a_u, b_v,c_w)$ holds if and only if $k+2\leq u+ v+w$.  
\end{enumerate}
\end{fact}
\begin{proof}
Suppose first there are $a_u,b_u,c_u\in M$ for each $u\in [k]$ such that $\calM\models R(a_u, b_v,c_w)$ holds if and only if $u\leq v+w$.  For each $u\in [k]$, let $a_u'=a_{k-u+1}$, $b_u'=b_{u}$, and $c'_u=c_{u}$.  Then $R(a_u',b_v',c_w')$ holds if and only if 
$$
k-u+1< v+w \Leftrightarrow k+2\leq u+v+w.
$$
Suppose on the other hand that there are $a_u,b_u,c_u\in M$ for each $u\in [k]$ such that $\calM\models R(a_u, b_v,c_w)$ holds if and only if $u\leq v+w$. For each $u\in [k]$, let $a_u'=a_{k-u+1}$, $b_u'=b_{u}$, and $c'_u=c_{u}$.  Then $R(a_u',b_v',c_w')$ holds if and only if 
$$
k-u+1+v+w \geq k+2 \Leftrightarrow v+w\geq u+1\Leftrightarrow v+w>u,
$$
which concludes the proof.
\end{proof}

To ease notation, we will use the indexing from Fact \ref{fact:reindex} (2) above for the remainder of this section.  We will first show that if $R_1\vee R_2$ has a large $n$-$\HOP_2$, and further $R_1$ is very dense, then $R_2$ has a large instance of $\HOP_2$.  For this, will use the following theorem of Gowers \cite[Theorem 10.3]{Gowers.2007}.

\begin{theorem}\label{thm:G} For every  $\delta> 0$, every positive integer $r$ and every finite subset $X \subseteq \mathbb{Z}^r$, there is a positive integer $N=N(\delta,r,X)$ such that every subset $A$ of
the grid $[N]^r$ of size at least $\delta N^r$ has a subset of the form $a + dX$ for some positive integer $d$.
\end{theorem}

\begin{lemma}\label{lem:almostallramsey}
For all $t\geq 1$ there is $\delta>0$ and $N\geq 1$ such that the following holds.  Suppose $n\geq N$, $R_1\cup R_2=\{(i,j,k)\in [n]^3: i+j+k\geq n+2\}$, and $|R_1|<\delta n^3$.  Then there are subsequences $1\leq i_1<\ldots<i_t\leq n$, $1\leq j_1<\ldots<j_t\leq n$, $1\leq k_1<\ldots<k_t\leq n$ such that for each $1\leq u,v,w\leq t$, $(i_u,j_v,k_w)\in R_2$ if and only if $u+v+w\geq t+2$.
\end{lemma}
\begin{proof}
If $t=1$, this is trivial. 

Assume now $t\geq 2$.  If $t=2$, let $\alpha=1/8$, and if $t>2$, let $\alpha:=(2t^3-6t+9)/(8t^3(t-1)(t-2))$.  Note that for all $t\geq 2$, $\alpha>0$. Choose any $\delta=\delta(t)<\alpha/(130 t^{100})$ and $N$ sufficiently large so that $N>12\alpha^{-1}$ and $(N/4t^3-2)^3>N^3/65 t^9$.

Suppose $n\geq N$, $R_1\cup R_2=\{(i,j,k)\in [n]^3: i+j+k\geq n+2\}$, and $|R_1|<\delta n^3$.  We say that $(\ibar,\jbar,\kbar)\in ([n]^t)^3$ is an \emph{good tuple} if the following hold.
\begin{enumerate}
\item There is some $d\geq 1$ such that for each $1\leq u\leq t$, $i_u=i_1+(u-1)d$, $j_u=j_1+(u-1)d$, and $k_u=k_1+(u-1)d$, and 
\item For each $1\leq u,v,w\leq t$, $i_u+j_v+k_w\geq n+2$ if and only if $u+v+w\geq t+2$.  
\end{enumerate}
We first give a lower bound on the number of good tuples.  Below is a procedure for generating good tuples.
\begin{itemize}
\item Fix $(i_1,j_1,k_1)\in [n]^3$ such that $\lceil n/4-n/8t^2\rceil \leq i_1,j_1,k_1\leq \lfloor n/4+n/8t^2\rfloor $.   The number of choices is at least 
$$
(\lfloor n/4+n/8t^3\rfloor -\lceil n/4-n/8t^3\rceil)^3 \geq (\lceil n/4t^3-2\rceil )^3\geq \lceil n^3/65t^9\rceil, 
$$
where the last inequality is by assumption on $n\geq N$.
\item Choose an integer $d$ such that 
$$
(t-2)d+i_1+j_1+k_1 \leq n+1 \text{ and } n+2\leq (t-1)d+i_1+j_1+k_1.
$$
We show there are at least $\lfloor \alpha n/2\rfloor $ many choices for this.  If $t=2$, we want to choose $d$ such that $\lceil \frac{n+2-i_1-j_1-k_1}{t-1}\rceil\leq d\leq n$, and if $t>2$, we want to choose $d$ such that $\lceil \frac{n+2-i_1-j_1-k_1}{t-1}\rceil\leq d\leq \lfloor \frac{n+1-i_1-j_1-k_1}{t-2}\rfloor$.  Observe that 
\begin{align*}
n+2-i_1-j_1-k_1&\geq n+2-3\lfloor n/4+n/8t^3\rfloor\geq n\Big(\frac{1}{4} -\frac{3}{8t^3}\Big)+2\text{ and }\\
n+1-i_1-j_1-k_1&\leq n+1-3\lceil n/4-n/8t^3\rceil\leq n\Big(\frac{1}{4} +\frac{3}{8t^3}\Big)+1.
\end{align*}
Therefore, if $t=2$, the number of choices is at least
\begin{align*}
n- \Big\lceil \frac{n+2-i_1-j_1-k_1}{t-1}\Big\rceil&\geq n- \Big(\frac{1}{t-1}\Big(n\Big(\frac{1}{4} +\frac{3}{8t^3}\Big)+2\Big)+1\Big)\\
&=n\Big(1-\frac{1/4+3/8t^3}{t-1}\Big)-\frac{2}{t-1}-1\\
&\geq n(1-1/4-3/8)-2-1\\
&\geq \lfloor \alpha n/2\rfloor ,
\end{align*}
where the last inequality is because $n$ is large, and $\alpha=1/8$ when $t=2$.  On the other hand, if $t>2$, the number of choices is at least
\begin{align*}
&\Big\lfloor \frac{n+1-i_1-j_1-k_1}{t-2}\Big\rfloor-\Big\lceil \frac{n+2-i_1-j_1-k_1}{t-1}\Big\rceil\\
&\geq \Big(\frac{n(\frac{1}{4} -\frac{3}{8t^3})+2}{t-2}-1\Big) - \Big(\frac{n(\frac{1}{4} +\frac{3}{8t^3})+1}{t-1}+1\Big)\\
&=n\Big(\frac{1}{4(t-2)}-\frac{3}{8t^3(t-2)}-\frac{1}{4(t-1)}-\frac{3}{8t^3(t-1)}\Big)+\frac{2}{t-2}-\frac{1}{t-1}-2\\
&=\alpha n +\frac{2}{t-2}-\frac{1}{t-1}-2\\
&\geq \lfloor \alpha n/2\rfloor,
\end{align*}
where the last inequality is because $n$ is large, and definition of $\alpha$ when $t>2$.
\item For each $1\leq u\leq t$, let $i_u=i_1+(u-1)d$, $j_u=j_1+(u-1)d$, and $k_u=k_1+(u-1)d$.  We claim the resulting tuple $(\ibar,\jbar,\kbar)$ is good.  Suppose first that $u+v+w<t+2$.  We show $i_u+j_v+k_w<n+2$.  If $t=2$, then $u+v+w<4$ implies $u=v=w=1$, so $i_u+j_v+k_w=i_1+j_1+k_1\leq 3\lfloor n/4+n/8t^2\rfloor<n+2$.    If $t>2$, then
\begin{align*}
i_u+j_v+k_w&=(i_1+j_1+k_1)+(u+v+w-3)d\\
&\leq i_1+j_1+k_1+(t-2)d\\
&\leq i_1+j_1+k_1+(t-2)\Big\lfloor \frac{n+1-i_1-j_1-k_1}{t-2}\Big\rfloor\\
&\leq i_1+j_1+k_1+(t-2)\Big(\frac{n+1-i_1-j_1-k_1}{t-2}\Big)\\
&<n+2.
\end{align*}  
On the other hand, assume $u+v+2\geq t+2$.  We show $i_u+j_v+k_w\geq n+2$.   Observe
\begin{align*}
i_u+j_v+k_w&=(i_1+j_1+k_1)+(u+v+w-3)d\\
&\geq i_1+j_1+k_1+(t-1)d\\
&\geq i_1+j_1+k_1+(t-1)\Big\lceil \frac{n+2-i_1-j_1-k_1}{t-1}\Big\rceil\\
&\geq i_1+j_1+k_1+(t-1)\Big( \frac{n+2-i_1-j_1-k_1}{t-1}\Big)\\
&=n+2.
\end{align*}
This completes the verification that $(\ibar,\jbar,\kbar)$ is a good tuple.
\end{itemize}
This shows the number of good triples is at least $\lfloor \alpha n/2\rfloor \lfloor n^3/65t^9\rfloor \geq \alpha n^4 /131 t^9$.    On the other hand, we claim that for any point $(i,j,k)\in [n]^3$, there are at most $nt^3$ many good triples containing $(i,j,k)$.  Clearly any good triple $(\ibar,\jbar,\kbar)$ containing $(i,j,k)$ is determined by its common distance $d$ (for which there are at most $n$ choices), along with the position $(u,v,w)\in [t]^3$ such that $(i,j,k)=(i_u,j_v,k_w)$ (for which there are at most $t^3$ choices).  Thus the number good triples containing $(i,j,k)$ is at most $t^3 n$.

Consequently, the number of good triples $(\ibar,\jbar,\kbar)$ for which there is $(u,v,w)\in [t]^3$ with $(i_u,j_v,k_w)\in R_1$ is at most $|R_1|nt^3<\delta t^3n^4$.  Thus if $\delta<\alpha/(130t^{10})$, there is a good triple $(\ibar,\jbar,\kbar)$ such that for all $(u,v,w)\in [t]^3$, $(i_u,j_v,k_w)\notin R_1$.  By assumption, and because it is a good tuple, this implies that for all $1\leq u,v,w\leq t$, $(i_u,j_v,k_w)\in R_2$ if and only if $u+v+w\geq t+2$. This finishes the proof.
\end{proof}

We now prove our desired result for a disjunction of one relation with no $2$-$\HOP_2$ and another with no $t$-$\HOP_2$.  This will provide the inductive step in the proof of our main result.

\begin{lemma}\label{lem:3dim2col}
For all $t\geq 1$ there is an $k=k(t)$ such that the following holds. Suppose $\calL=\{R_1,R_2\}$ where each $R_i$ is a ternary relation symbol.  If $\calM$ is a sufficiently large $\calL$-structure such that $R_1(x,y,z)$ does not have $2$-$\HOP_2$  in $\calM$, and $R_2$ does not have $t$-$\HOP_2$ in $\calM$.  Then $R_1\vee R_2$ does not have $k$-$\HOP_2$ in $\calM$.
\end{lemma}
\begin{proof}

Choose first $\delta=\delta(t)$ and $N_1=N_1(t)$ as in Lemma \ref{lem:almostallramsey}.  Given $m\in \mathbb{N}$, let $X_m=[m]^2\subseteq \mathbb{Z}^2$.  Now let $N_2$ be the smallest even integer which is greater than $\max\{2N(\delta,2,X_{N_1}),2N(\delta,2,X_{4t})\}$, where $N$ is as in Theorem \ref{thm:G}. Finally, set $k=2N_2$.  

Suppose towards a contradiction $\calM$ is an $\calL$-structure where $R_1$ has no $2$-$\HOP_2$, where $R_2$ has no $t$-$\HOP_2$, but where $R_1\vee R_2$ has $k$-$\HOP_2$. Let $\{a_i,b_i,c_i:i\in [k]\}\subseteq M$ be such that $\calM\models (R_1\vee R_2)(a_u,b_v,c_w)$ if and only if $u+v+w\geq k+2$.  Recall that given an $\calL$-formula $\varphi(x,y,z)$, $\varphi^\calM=\{(a,b,c)\in M^3: \calM\models \varphi(a,b,c)\}$.

\underline{Case 1}: Suppose for every $ k/2 +1\leq  i\leq k$, $|R^\calM_1\cap \{(a_i,b_v,c_w): (v,w)\in [ k/2]^2\}|<\delta (k/2)^2$.  In this case, for each $1\leq i\leq k/2$, let $a'_i=a_{k/2 +i}$, $b'_i=b_{i}$, and $c_i'=c_i$.  Then 
$$
|R^\calM_1\cap \{(a_i',b_j',c_t'): (i,j,t)\in [ k/2]^3\}|\leq \delta (k/2)^3.
$$
On the other hand, we claim that $(R_1\vee R_2)(a_i',b_j',c_t')$ holds if and only if $i+j+t\geq k/2+2$.  Indeed, given $i,j,t\in [\lfloor k/2\rfloor]$,
$$
a_i'+b_j'+c_t'=a_{k/2 +i}+b_{j}+c_{j}.
$$
Thus by assumption $(R_1\vee R_2)(a_i',b_j',c_t')$ holds if and only if $k/2+i+j+t\geq k+2$, i.e. if and only if $i+j+t\geq k/2+2$.  Since $k/2\geq N_2>N_1$, Lemma \ref{lem:almostallramsey} implies that $R_2$ has $t$-$\HOP_2$, a contradiction.

\underline{Case 2}: Suppose we are not in Case 1.  Then there is $k/2\leq i_*\leq k$ such that 
$$
|R^\calM_1\cap \{(a_{i_*},b_v,c_w): (v,w)\in [k/2]^2\}|\geq \delta (k/2)^2.
$$
Since $k/2\geq N_2$, Theorem \ref{thm:G} implies there is $(\beta,\gamma)\in [k/2]^2$ and $d\geq 1$ such that for each $1\leq s,s'\leq N_1$, $(a_{i_*},b_{\beta+sd}, c_{\gamma+s'd})\in R^\calM_1\cap \{(a_{i_*},b_v,c_w): (v,w)\in [k/2]^2\}$.  

Now for each $1\leq s\leq N_1$, define $b'_s=b_{\beta+sd}$ and $c'_s=c_{\gamma+sd}$.  Set $a'_{N_1}=a_{i_*}$ and for each $1\leq s\leq N_1-1$, let $a'_s= a_{k+2-\beta-\gamma -(N_1+2)d+sd}$.  Note this last definition makes sense because $\beta+N_1d\leq k/2$ and $\gamma+N_1d\leq k/2$, so for each $1\leq s\leq N_1-1$, 
$$
k\geq k+2-\beta-\gamma-(N_1+2)d+sd\geq 1.
$$
  We claim that given $u,v,w\in [N_1]$, $(R_1\vee R_2)(a_u',b_v',c_w')$ holds if and only if $u+v+w\geq N_1+2$.

Suppose $1\leq u,v,w\leq N_1$.  If $u=N_1$ then we already know that $R_1(a_u',b_v',c_w')$.  So assume $u\neq N_1$.  Then 
\begin{align*}
(a'_u,b'_v,c'_w)&=(a_{k+1-\beta-\gamma-(N_1+2)d+ud},b_{\beta+vd},c_{\gamma+wd}). 
\end{align*}
By assumption, $(R_1\vee R_2)(a_{k-\beta-\gamma-(N_1+2)d+ud},b_{\beta+vd},c_{\gamma+wd})$ holds if and only if 
\begin{align*}
&k+2-\beta-\gamma-(N_1+2)d+ud+\beta+vd+\gamma+wd=k+2+d(u+v+w-N_1-2)\geq k+2.
\end{align*}
Note this holds if and only if $u+v+w-N_1-2\geq 0$, which holds if and only if $u+v+w\geq N_1+2$, as desired.  Thus we have   $\{(a'_u,b'_u,c'_u): u\in [N_1]\}$ such that that $(R_1\vee R_2)(a'_u,b'_v,c'_w)$ holds if and only if $u+v+w\geq N_1+2$, and further, for all $v,w\in [N_1]$, $R_1(a'_{N_1},b'_v,c'_w)$. 

\underline{Subcase 2.1}: Suppose that for all $N_1/2+1\leq j\leq N_1$, 
$$
|R^\calM_1\cap \{(a_u',b_j',c_w'): (u,w)\in [N_1/2]^2\}|<\delta (N_1/2)^2.
$$
In this case, for each $1\leq j\leq N_1/2$, let $a_j''=a_j'$, $c_j''=c_j'$ and $b_j''=b'_{N_1/2+j}$. Then 
$$
|R^\calM_1\cap \{(a_i'',b_j'',c_t''): (i,j,t)\in [N_1/2]^3\}|<\delta (N_1/2)^3.
$$
Further, given $u,v,w\in [N_1/2]$, $R_1\vee R_2(a_u'',b_v'',c_w'')$ holds if and only if $u+w+N_1/2+v\geq N_1+2$, which holds if and only if $u+v+w\geq N_1/2+2$.  Consequently, we may apply Lemma \ref{lem:almostallramsey} to $\{(a_u'',b_v'',c_w''): u,v,w\in [N_1/2]\}$, to obtain that $R_2$ must have $t$-$\HOP_2$, a contradiction.

\underline{Subcase 2.2}: There exists $N_1/2+1\leq j_*\leq N_2$, such that 
$$
|R^\calM_1\cap \{(a_u',b_{j_*}',c_w'): (u,w)\in [N_1/2]^2\}|\geq \delta (N_1/2)^2.
$$
By Theorem \ref{thm:G}, there are $(\sigma,\tau)\in [N_1/2]^2$ and $d'\geq 1$ such that for each $1\leq s,s'\leq 4t$, $(a'_{\sigma+sd'},b_{j_*}', c_{\tau+s'd'}')\in R^\calM_1\cap \{(a_v',b_{j_*}',c_w'): (v,w)\in [4t]^2\}$.  

Given $1\leq v\leq t$, let $x_v=a'_{\sigma+(v+t)d'}$, $z_v=c'_{\tau+(v+t)d'}$, and $y_v=b'_{N_1+1-\sigma-\tau+d'(v-3t-1)}$.  Note the last definition makes sense because $\sigma+4td'\leq N_1/2$ and $\tau+4td'\leq N_1/2$  so for each $1\leq v\leq t$,
$$
N_1\geq N_1+1-\sigma-\tau+d'(v-3t-1)\geq 1.
$$

Now for $u,v,w\in [t]$, $R_1\vee R_2(x_u,y_v,z_w)$ holds if and only if
\begin{align*}
N_1+1-\sigma-\tau+d'(u-3t-1)+\sigma+(v+t)d'+\tau+(w+t)d'&=N_1+1+d'(-1-t+u+v+w)\\
&\geq N_1+2.
\end{align*}
This holds if and only if $-t-1+u+v+w\geq 1$, which holds if and only if $u+v+w\geq t+2$.  Now I claim that in fact, in each case where $u+v+w\geq t+2$, we have that $R_2(x_u,y_v,z_w)$ holds.  Indeed, fix $1\leq u,v,w\leq t$ and assume $u+v+w\geq t+2$.  By definition, $x_u=a'_{\sigma+(u+t)d'}$, $y_v=b'_{N_1+1-\sigma-\tau+d'(v-3t-1)}$, and $z_w=c'_{\tau+(w+t)d'}$.

Let $x_1'=x_u$, $y_1'=y_v$ and $z_2'=z_w$.  Now let $x_2'=a'_{N_2}=a_{i_*}$, $y_2'=b_{j_*}'$, and $z_1'=c_1'$.  Clearly $R_1(x_i',y_2',z_j')\wedge R_1(x_2',y_i',z_j')$ holds for each $i,j\in \{1,2\}$.  Consider now $(x_1', y_1', z_1')=(a'_{\sigma+(u+t)d'},b'_{N_1+1-\sigma-\tau+d'(v-3t-1)}, c_1')$.  Note
\begin{align*}
(\sigma+(u+t)d')+(N_1+1-\sigma-\tau+d'(v-3t-1))+1&=N_1+2+d'(u+v-2t-1)\\
&\leq N_1+2-d',
\end{align*}
where the last inequality is because $1\leq u,v\leq t$.  Thus since the above sum is less than $N_1+2$, we have $\neg R_1(x_1',y_1',z_1')$.  Consider now $(x_1',y_1',z_2')=(a'_{\sigma+(u+t)d'},b'_{N_1+1-\sigma-\tau+d'(v-3t-1)}, c_{\tau+d'(w+t)})$.  Observe that
\begin{align*}
&(\sigma+(u+t)d')+(N_1+1-\sigma-\tau+d'(v-3t-1))+(\tau+d'(w+t))\\
&=N_1+1+d'(-t+u+v+w-1)\\
&\geq N_1+1+d'\\
&\geq N_1+2,
\end{align*}
where the second to last inequality is because $u+v+w\geq t+2$ by assumption.  Thus $(R_1\vee R_2)(x_1',y_1',z_2')$.  Observe that if $R_1(x_1',y'_1,z_2')$ held, then $R_1$ would have $2$-$\HOP_2$  (as witnessed by $\{x_1',x_2',y_1',y_2',z_1',z_2'\}$), a contradiction.  Thus we must have that $R_2(x_1',y_1',z_2')$, i.e. $R_2(x_u,y_v,z_w)$ holds.

Thus we have shown that for all $1\leq u,v,w\leq t$, $R_2(x_u,y_v,z_w)$ holds if and only if $u+v+w\geq t+2$, and consequently, $R_2$ has $t$-$\HOP_2$, a contradiction.
\end{proof}

It is now easy to deduce that an arbitrary disjunction of ternary relations with no $2$-$\HOP_2$ has no $k$-$\HOP_2$ for some $k$.

\begin{corollary}
For all $t\geq 1$, there exists $s=s(t)$ such that the following holds. Let $\calL=\{R_1,\ldots, R_t\}$, where each $R_i$ is a ternary relation symbol, and let $\calM$ be an $\calL$-structure such that each $R_i$ does not have $2$-$\HOP_2$ in $\calM$. Then $R_1\vee \ldots \vee R_t$ does not have $s$-$\HOP_2$ in $\calM$.
\end{corollary}
\begin{proof}
Let $s(1)=2$.  Now suppose by induction $t\geq 1$ and that we have defined $s(t)$. Let $s(t+1)=k(s(t))$, where $k(x)$ is as in Lemma \ref{lem:3dim2col} above.  Suppose $R_1,\ldots, R_{t+1}$ are each ternary relations which each have no $2$-$\HOP_2$ in $\calM$.  By induction $R_1\vee \ldots \vee R_t$ has no $s(t)$-$\HOP_2$, and by assumption $R_{t+1}$ has no $2$-$\HOP_2$.  Therefore, by Lemma \ref{lem:3dim2col}, $R_1\vee \ldots \vee R_t \vee R_{t+1}$ has no $k(s(t))=s(t+1)$-$\HOP_2$, as desired.
\end{proof}

\vspace{2mm}

\subsection{$A_{\GS}(3,n)$ has no $4$-$\HOP_2$}\label{app:no4hop} In this section, we show that $A_{\GS}(3,n)$ has no $4$-$\HOP_2$. In fact, the proof will show that $A_{\GS}(3,n)$ does have $3$-$\HOP_2$. Taken together, these facts will prove Proposition \ref{prop:GSnoHOP}.

Firstly, observe that the elements of $\F_3^4$ given by
\[\begin{array}{ccc}
x_1=2220 & y_1=2220 &z_1=2221\\
x_2=2210 & y_2=2200 & z_2=2011\\
x_3=2120 & y_3=0220 & z_3=2021
\end{array}\]
determine a $3$-$\HOP_2$ in $A_{\GS}(3,4)$. Indeed, it is easy to check by hand that $x_i+y_j+z_k\in A_{\GS}(3,n)$ if and only if $i<j+k$. This immediately gives rise to a $3$-$\HOP_2$ in $A_{\GS(3,n)}$ for any $n\geq 4$.

We now prove that $A_{\GS}(3,n)$ has no $4$-$\HOP_2$.  We begin with some lemmas about how certain small configurations can arise in $A_{\GS}(3,n)$. Throughout the remainder of this section, $A=A_{\GS}(3,n)\subseteq \F_3^n$.  Recall the following definition.

\begin{definition}
Given $x\in \F_3^n$, define
$$
f(x):=\max\Big( \{1\}\cup \{2\leq j\leq n: \text{ for all }1\leq i< j, x_i=0\}\Big).
$$
\end{definition}

Our proofs will rely heavily on proving restrictions on the behavior of $f$ on sums arising from various combinatorial configurations.  We begin with such a lemma regarding additive quadruples. 

\begin{lemma}[Necessary conditions for an additive quadruple]\label{lem:c4}
Let $a^1,a^2,b^1,b^2\in \mathbb{F}_3^n$ and suppose that for each $(i,j)\in \{0,1\}^2$, $a^i+b^j\in A$.  Then one of the following holds, where $u_{ij}=f(a^i+b^j)$.
\begin{enumerate}[label=\normalfont(\roman*)]
\item $u_{11}=u_{12}=u_{21}=u_{12}$,
\item $u_{11}=u_{12}<u_{21},u_{22}$,
\item $u_{11}=u_{21}<u_{12},u_{22}$,
\item $u_{22}=u_{12}<u_{21},u_{11}$,
\item $u_{22}=u_{21}<u_{21},u_{11}$.
\end{enumerate}
\end{lemma}
\begin{proof}
Let $y=(a^1+b^2)+(a^2+b^1)=(a^1+b^1)+(a^2+b^2)$, and each $a^i+b^j\in A$, we have the following.
\begin{align*}
&e^{u_{12}}+ e^{21}+\sum_{s=u_{12}+1}^n (a^1_s+b^2_s) e^s +\sum_{s=u_{21}+1}^n(a^2_s+b^1_s)e^s\\
&=e^{u_{11}}+ e^{u_{22}}+\sum_{s=u_{11}+1}^n (a^1_s+b^1_s) e^s + \sum_{s=u_{22}+1}^n(a^2_s+b^2_s)e^s.
\end{align*}
This shows that $f(y)=\min\{u_{12},u_{21}\}=\min\{u_{11},u_{22}\}$.  Suppose first $u_{12}=u_{21}$.  Then 
$$
y=2e^{u_{12}}+\sum_{s=u_{21}+1}^n(a^2_s+b^1_s+a^1_s+b^2_s)e^s=e^{u_{11}}+ e^{u_{22}}+\sum_{s=u_{11}+1}^n (a^1_s+b^1_s) e^s + \sum_{s=u_{22}+1}^n(a^2_s+b^2_s)e^s.
$$
This shows $u_{11}=u_{22}=u_{12}=u_{21}$.  Similarly if $u_{11}=u_{22}$, then $u_{11}=u_{22}=u_{12}=u_{21}$.   Suppose now $u_{12}<u_{21}$.  By what we just showed, we must have $u_{11}\neq u_{22}$ (since $u_{11}=u_{22}$ would imply $u_{12}=u_{21}$).  Since $\min\{u_{12},u_{21}\}=\min\{u_{11},u_{22}\}$, we have that either $u_{12}=u_{11}<u_{21},u_{22}$ or $u_{12}=u_{22}<u_{21},u_{11}$.  Similarly, if $u_{21}<u_{12}$, then either $u_{21}=u_{11}<u_{12},u_{22}$ or $u_{21}=u_{22}<u_{12},u_{11}$.  
\end{proof}

\begin{lemma}[Necessary conditions for a $2$-$\OP$]
Let $a^1,a^2,b^1,b^2\in \mathbb{F}_3^n$ and suppose that for each $(i,j)\in \{1,2\}^2\setminus \{(2,2)\}$, $a^i+b^j\in A$, but $a^2+b^2\notin A$.  Then one of the following holds, where for each $i,j\in \{1,2\}$, $u_{ij}=f(a^i+b^j)$.
\begin{enumerate}[label=\normalfont(\roman*)]
\item $u_{12}=u_{21}=u_{22}<u_{11}$,
\item $u_{11}=u_{22}<u_{12},u_{21}$,
\item $u_{11}=u_{12}<u_{22},u_{21}$,
\item $u_{11}=u_{21}<u_{22},u_{12}$.
\end{enumerate}
\end{lemma}
\begin{proof}
Let $y=a^1+b^1+a^2+b^2$.  For ease of notation, let $a=u_{11}$, $b=u_{12}$, $c=u_{21}$, and $d=u_{22}$.  Then we can write $a^1+b^1=\sum_{i=1}^nr_ie_i$, $a^1+b^2=\sum_{i=1}^ns_ie_i$, $a^2+b^1=\sum_{i=1}^nt_ie_i$, and $a^2+b^2=\sum_{i=1}^nv_ie_i$, where $r_i=0$ for all $i<a$, $s_i=0$ for all $i<b$, $t_i=0$ for all $i<c$, and $v_i=0$ for all $i<d$.  By our assumption, that $(i,j)\in \{1,2\}^2\setminus \{(2,2)\}$, $a^i+b^j\in A$, and $a^2+b^2\notin A$, we have $r_a=s_b=t_c=1$ and $v_d=2$.   Therefore,
$$
y=2e^d+e^a+\sum_{i=a+1}^nr_ie_i+\sum_{j=d+1}^nv_je^j=e^b+e^c+\sum_{i=b+1}^ns_ie_i+\sum_{j=c+1}^nt_je^j.
$$
Suppose first that $b=c$, this implies that either $b=c=d<a$ (so (i)) holds, or $a=d<b=c$ (so (ii) holds).

Suppose now that $b<c$.  Then the above equation for $y$ implies that either $a=b<c,d$ (so (iii) holds), or $a=d<b<c$ (so (ii) holds).  The case $c<b$ is proceeds in an analogous way to show that either (iv) or (ii) holds, so we are done.
\end{proof}

Our next step is a crucial lemma about any instance of $2$-$\HOP_2$.

\begin{lemma}[Necessary conditions for a $2$-$\HOP_2$]\label{lem:2hop}
Let $a^1,a^2,b^1,b^2,c^1,c^2$ be such that $a^i+b^j+c^k\in A$ if and only if $i+j+k\geq 4$.  Then one of the following holds, where $u_{ijk}=f(a^i+b^j+c^k)$. 
\begin{enumerate}[label=\normalfont(\arabic*)]
\item $u_{112}=u_{122}=u_{222}=u_{212}<u_{111}=u_{121}=u_{211}<u_{221}$,
\item $u_{112}=u_{122}=u_{222}=u_{212}<u_{111}=u_{221}<u_{121}, u_{211}$,
\item $u_{112}=u_{122}=u_{222}=u_{212}<u_{121}=u_{221}<u_{111}, u_{211}$,
\item $u_{112}=u_{122}=u_{222}=u_{212}<u_{211}=u_{221}<u_{111}, u_{121}$,

\item $u_{212}=u_{211}=u_{221}=u_{222}<u_{111}=u_{121}=u_{112}<u_{122}$,
\item $u_{212}=u_{211}=u_{221}=u_{222}<u_{111}=u_{122}<u_{121}, u_{112}$,
\item $u_{212}=u_{211}=u_{221}=u_{222}<u_{122}=u_{121}<u_{111}, u_{112}$,
\item $u_{212}=u_{211}=u_{221}=u_{222}<u_{122}=u_{112}<u_{111}, u_{121}$,

\item $u_{122}=u_{222}=u_{221}=u_{121}<u_{111}=u_{121}=u_{211}<u_{212}$,
\item $u_{122}=u_{222}=u_{221}=u_{121}<u_{111}=u_{212}<u_{112}, u_{211}$,
\item $u_{122}=u_{222}=u_{221}=u_{121}<u_{212}=u_{112}<u_{111}, u_{211}$,
\item $u_{122}=u_{222}=u_{221}=u_{121}<u_{212}=u_{211}<u_{111}, u_{112}$.
\end{enumerate}
\end{lemma}

In order to interpret the conclusions of Lemma \ref{lem:2hop}, the reader may benefit from visualizing a $2$-$\HOP_2$ as a cube whose corner in the $(a^i, b^j,c^k)$-coordinate is labelled by the value $u_{ijk}=f(a^i+b^j+c^k)$ as in the diagram below. A filled-in circle represents a sum in $A$ and a void circle represents a sum not in $A$.
 
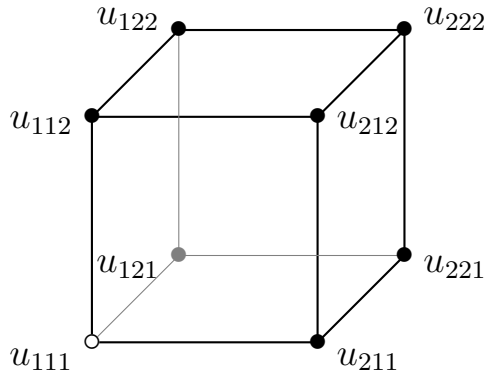
 \begin{figure}[H]
\begin{tikzpicture}[scale=1.5, every node/.style={scale=1.2}]
  \draw[thick](2,2,0)--(0,2,0)--(0,2,2)--(2,2,2)--(2,2,0)--(2,0,0)--(2,0,2)--(0,0,2)--(0,2,2);
  \draw[thick](2,2,2)--(2,0,2);
  \draw[gray](2,0,0)--(0,0,0)--(0,2,0);
  \draw[gray](0,0,0)--(0,0,2);
  \draw[gray](0,0,0) node{$\bullet$};
  \draw(-0.45,-0.1,0) node{$u_{121}$};
  \draw(0,0,2) node{\textcolor{white}{$\bullet$}};
  \draw(0,0,2) node{$\circ$};
  \draw(-0.45,-0.15,2) node{$u_{111}$};
  \draw(0,2,0) node{$\bullet$};
  \draw(-0.45,2.1,0) node{$u_{122}$};
  \draw(2,0,0) node{$\bullet$};
  \draw(2.45,-0.1,0) node{$u_{221}$};
  \draw(2,2,0) node{$\bullet$};
  \draw(2.45,2.1,0) node{$u_{222}$};
  \draw(2,0,2) node{$\bullet$};
  \draw(2.45,-0.15,2) node{$u_{211}$};
  \draw(0,2,2) node{$\bullet$};
  \draw(-0.45,1.95,2) node{$u_{112}$};
  \draw(2,2,2) node{$\bullet$};
  \draw(2.45,1.95,2) node{$u_{212}$};
\end{tikzpicture}
\end{figure}

Each of (1)-(4), (5)-(8), and (9)-(12) in Lemma \ref{lem:2hop} are equivalent up to relabeling, and correspond to a cube in which one face has a constant $f$-value, and this face does not touch the corner whose sum is not in $A$.

\vspace{2mm}

\begin{proofof}{Lemma \ref{lem:2hop}}
For ease of notation, let $u_{111}=a$, $u_{112}=b$, $u_{122}=c$, $u_{121}=d$, $u_{221}=e$, $u_{211}=f$, $u_{212}=g$, $u_{222}=h$.  The diagram below is intended to be helpful to the reader as we proceed.

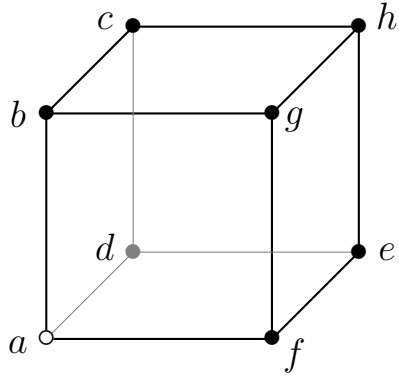
\begin{figure}[H]
\begin{tikzpicture}[scale=1.5, every node/.style={scale=1.2}]
  \draw[thick](2,2,0)--(0,2,0)--(0,2,2)--(2,2,2)--(2,2,0)--(2,0,0)--(2,0,2)--(0,0,2)--(0,2,2);
  \draw[thick](2,2,2)--(2,0,2);
  \draw[gray](2,0,0)--(0,0,0)--(0,2,0);
  \draw[gray](0,0,0)--(0,0,2);
  \draw[gray](0,0,0) node{$\bullet$};
  \draw(-0.25,0.05,0) node{$d$};
  \draw(0,0,2) node{\textcolor{white}{$\bullet$}};
  \draw(0,0,2) node{$\circ$};
  \draw(-0.25,-0.05,2) node{$a$};
  \draw(0,2,0) node{$\bullet$};
  \draw(-0.25,2.05,0) node{$c$};
  \draw(2,0,0) node{$\bullet$};
  \draw(2.25,0,0) node{$e$};
  \draw(2,2,0) node{$\bullet$};
  \draw(2.25,2.1,0) node{$h$};
  \draw(2,0,2) node{$\bullet$};
  \draw(2.2,-0.15,2) node{$f$};
  \draw(0,2,2) node{$\bullet$};
  \draw(-0.25,2,2) node{$b$};
  \draw(2,2,2) node{$\bullet$};
  \draw(2.2,1.95,2) node{$g$};
\end{tikzpicture}
\end{figure}

Note that $a^1+c^1, a^2+c^1,b^1,b^2$ form a $2$-$\OP$, and $f((a^1+c^1)+b^1)=a$, $f((a^1+c^1)+b^2)=d$, $f((a^2+c^1)+b^1)=f$, and $f((a^2+c^1)+b^2)=e$.   Applying Lemma \ref{lem:2hop}, we have that one of the following holds.
\begin{enumerate}
\item[(a)] $d=f=a<e$
\item[(b)] $e=d<f,a$
\item[(c)] $e=a<d,f$
\item[(d)] $e=f<a,d$
\end{enumerate}
Similarly, $a^1+b^1,a^2+b^1, c^1,c^2$ form a $2$-$\OP$, so one of the following holds.
\begin{enumerate}
\item[(a')] $b=f=a<g$
\item[(b')] $a=g<b,f$
\item[(c')] $b=g<a,f$
\item[(d')] $f=g<b,a$
\end{enumerate}
Similarly, $c^1+a^1,c^2+a^1,b^1,b^2$ form a $2$-$\OP$,  which implies one of the following holds.
\begin{enumerate}
\item[(a'')] $b=d=a<c$
\item[(b'')] $c=a<b,d$
\item[(c'')] $c=b<d,a$
\item[(d'')] $c=d<a,b$
\end{enumerate}
Consider now $a^1,a^2,b^1+c^2,b^2+c^2$.  These form an additive quadruple, so by Lemma \ref{lem:c4}, one of the following holds.

\begin{enumerate}
\item $b=c=g=h$
\item $b=g<c,h$
\item $b=c<g,h$
\item $c=h<g,b$
\item $h=g<b,c$
 \end{enumerate}
Observe that $a^1,a^2,c^1+b^2,c^2+b^2$ form an additive quadruple, so one of the following holds.
\begin{enumerate}
\item[(i')] $c=e=d=h$
\item[(ii')] $c=h<d,e$
\item[(iii')] $c=d<h,e$
\item[(iv')] $d=e<h,c$
\item[(v')] $e=h<c,d$
\end{enumerate}
Consider $b^1,b^2,c^1+a^2,c^2+a^2$.  These form an additive quadruple, so one of the following holds.
\begin{enumerate}
\item[(i'')] $g=h=e=f$
\item[(ii'')] $g=h<e,f$
\item[(iii'')] $g=f<h,e$
\item[(iv'')] $e=h<g,f$
\item[(v'')] $e=f<g,h$
\end{enumerate}

Suppose first that (a) holds, i.e. $d=f=a<e$. Then we have to have that either (a') or (c') hold, and either (a'') hold or (c'') hold.  Suppose (a') holds.  So we have that $d=f=a=b<g,c$.  But none of (i)-(v) are consistent with this, so (a') cannot hold.  A symmetric argument shows (a'') also cannot hold.  So we must have (c') and (c''), i.e. $b=g<a,f$ and $c=b<a,d$.  Combining these we have $e>d=f=a>b=g=c$.  This is (1).  By symmetry, it is easy to see that if (a') holds, then (9) holds, and if (a'') holds, then (5) holds.   This takes care of the cases where one of (a), (a'), or (a'') hold.

Assume now none of (a), (a'), or (a'') hold. 

Suppose (b) holds, so $e=d<f,a$.  Note we cannot have (b''), which leaves us that (c'') or (d'') hold, i.e. either $c=b<d,a$ or $c=d<a,b$.  

Suppose first (d'') holds, i.e. $c=d<a,b$.  So we have $e=c=d<a,b,f$.  The only possibility from (i')-(v') consistent with $e=c=d$ is (i'), so we have $e=c=d=h$.  Combining with (b'),(c'), or (d'), we have the following possibilities.
\begin{itemize}
\item $c=d=e=h<a=g<b,f$
\item $c=d=e=h<b=g<a,f$
\item $c=d=e=h<f=g<a,b$.
\end{itemize}
These correspond to (10), (11), and (12) respectively.

Suppose now (c'') holds, i.e. $c=b<d,a$.  Combining (c'') and (b), we have $c=b<d=e<f,a$.  Since one of (i')-(iv') holds, and $c<d$, the only option is $c=h<d,e$.  Overall we now have $c=b=h<d=e<f,a$.   The only option among (i)-(iv) which is consistent with  $c=b=h<d=e<f,a$ is (i), so we have that $c=b=h=g<d=e<f,a$. This is (10).

This deals with all the cases where (a), (a'), (a'') do not hold and (b) holds.  By symmetry, it also takes care of the cases where   (a), (a'), (a'') do not hold and one of (b') or (b'') hold.  So assume now that none of (a), (a'), (a'') or (b), (b') or (b'') hold.

Suppose (c) holds, i.e. $e=a<d,f$.  Then neither (d') nor (d'') can hold, so we are left with (c') and (c'') holding.  Combining these we have that $c=b=g<e=a<f,d$.  The only of of (i)-(iv) hold consistent with   $c=b=g$ is (i), so we have $c=b=g=h<e=a<f,d$.  This is (8).

This deals with the case where none of (a), (a'), (a'') or (b), (b') or (b'') hold, and where (c) holds.  By symmetry, this finishes all cases where none of (a), (a'), (a'') or (b), (b') or (b''), or (c), (c'), (c'') hold.  We are left with the case where (d), (d'), and (d'') hold.  Combining these together we have that $e=f=g<c=d<b,a$.  The only one of (i'')-(v'') consistent with $e=f=g$ is (i''), so we have that $e=f=g=h<c=d<b,a$.   This is (6).			
\end{proofof}

\vspace{2mm}

\begin{lemma}[Necessary conditions for a $3$-$\HOP_2$]\label{lem:3hop}
Suppose $\{a^i,b^i,c^i: i\in [3]\}$ are such that $a^i+b^j+c^k\in A$ if and only if $i+j+k\geq 5$. Then one of the following holds, where $u_{ijk}=f(a^i+b^j+c^k)$ and $\gamma:=\min\{u_{ijk}: i,j,k\in [3]\}$.
\begin{enumerate}[label=\normalfont(\roman*)]
\item For all $i,j\in [3]$, $u_{ij3}=\gamma$ and $u_{ij1},u_{ij2}>\gamma$.
\item For all $i,k\in [3]$, $u_{i3k}=\gamma$ and $u_{i1k},u_{i2k}>\gamma$.
\item For all $j,k\in [3]$, $u_{3jk}=\gamma$ and $u_{1jk},u_{2jk}>\gamma$.
\end{enumerate}
\end{lemma}
\begin{proof}
Consider $a^1, b^1, c^1, a^3,b^3,c^3$.  This makes a $2$-$\HOP_2$, so by Lemma \ref{lem:2hop}, after relabeling, we may assume that $u_{131}=u_{331}=u_{133}=u_{333}<u_{111},u_{311},u_{113},u_{313}$.

Observe that $a^2,b^1,c^1, a^3,b^3,c^3$ also forms a $2$-$\HOP_2$.  We already know that $u_{333}=u_{331}<u_{313}, u_{311}$, so we cannot have $u_{333}=u_{313}=u_{311}=u_{331}$ or $u_{333}=u_{233}=u_{213}=u_{313}$. The only possibility for the four equal terms in Lemma \ref{lem:2hop} is then $u_{333}=u_{331}=u_{231}=u_{233}$, so we have $u_{333}=u_{331}=u_{231}=u_{233}<u_{211},u_{311},u_{213},u_{313}$.

Now consider $a^1,a^2,b^1,b^3,c^2,c^3$.  This also forms a $2$-$\HOP_2$.  We already know that $u_{233}=u_{133}<u_{213}, u_{113}$ so it is not possible to have either $u_{233}=u_{213}=u_{212}=u_{232}$ or $u_{233}=u_{213}=u_{113}=u_{133}$.  So this leaves only one possibility for the face of four equal terms from Lemma \ref{lem:2hop}, namely $u_{233}=u_{133}=u_{132}=u_{232}$.  Thus we have $u_{233}=u_{133}=u_{132}=u_{232}<u_{112},u_{213},u_{212},u_{113}$.

Now consider $a^1,a^3,b^1,b^3,c^2,c^3$.  Again this makes a $2$-$\HOP_2$.  Note we already know that $u_{133}=u_{333}=u_{132}<u_{113},u_{313}, u_{312}$.  This leaves only one possibility for face of four equal terms from Lemma \ref{lem:2hop}, namely $u_{333}=u_{133}=u_{332}=u_{132}$ so $u_{333}=u_{133}=u_{332}=u_{132}<u_{112},u_{313},u_{312},u_{113}$.

Observe we now know that there is some $a$ such that for all $i,k\in [3]$, $u_{i3k}=\gamma$, and for all $i,k\in [3]$, $\gamma<u_{i1k}$.  The following sets also form $2$-$\HOP_2$s.
\begin{itemize}
\item $\{a^2,a^3,b^2,b^3,c^1,c^2\}$.
\item   $\{a^2,a^3,b^2,b^3,c^1,c^3\}$. 
\item $\{a^1,a^2,b^2,b^3,c^1,c^3\}$.
\item $\{a^1,a^2,b^2,b^3,c^1,c^2\}$.
\end{itemize}
Based on what we already know, and Lemma \ref{lem:2hop}, there is only one compatible choice for how these all behave, so we can conclude that for all $i,k\in [3]$, we also have $\gamma<u_{i2k}$.  In conclusion, we have that $u_{i3k}=\gamma$ for all $i,k\in [3]$ and $u_{ijk}>\gamma$ for all $i,k\in [3]$ and $j\in [2]$.
\end{proof}

The following lemma provides some restrictions on how things can align in a $4$-$\HOP_2$.

\begin{lemma}\label{lem:notallowed}
Suppose $\{x_i,y_i:i\in [3]\}\subseteq \F_p^n$ and $z,z',w\in \F_p^n$ are such that
\begin{itemize}
\item $w+x^i+y^j\in A$ if and only if $i+j\geq 3$,  
\item $z+x^i+y^j\in A$ if and only if $i+j\geq 5$, and 
\item $z'+x^i+y^j\in A$ if and only if $i+j\geq 4$.  
\end{itemize}
Then none of the following can hold.  
\begin{enumerate}[label=\normalfont(\arabic*)]
\item $f(w+y^1+x^3)=f(w+y^2+x^3)<f( w+y^i+x^j)$ for all $i\in [3]$ and $j\in [2]$,
\item $f(w+x^3+y^2)=f(w+x^2+y^2)<f(w+x^2+y^3)=f(w+x^1+y^3)$,
\item $f(w+x^3+y^3)=f(w+x^3+y^2)=f(w+x^2+y^3)=f(w+x^2+y^2)>f(w+x^1+y^3)$,
\item $f(w+x^3+y^3)=f(w+x^3+y^2)=f(w+x^2+y^3)=f(w+x^2+y^2)=f(w+x^1+y^1)$,
\item $f( w+x^1+y^3)=f(w+x^2+y^3)=f(w+x^3+y^3)>f(w+x^3+y^2)=f(w+x^2+y^2)=f(w+x^1+y^2)$,
\item $f( w+x^1+y^3)=f(w+x^2+y^3)=f(w+x^3+\ybar_3)=f(w+x^3+y^2)=f(w+x^2+y^2)=f(w+x^1+y^2)<f(w+x^3+y^1)$.
\end{enumerate}

\end{lemma}
\begin{proof}
Suppose first that (1) holds.  Note
$$
(z-w)+(w+y^1+x^3)\notin A\text{ and }(z-w)+(w+y^2+x^3)\in A.
$$
Since $f(w+y^1+x^3)=f(w+y^2+x^3):=\gamma$, this implies that 
$$
f(z-w)=f(w+y^1+x^3)=f(w+y^2+x^3)=\gamma.
$$ 
Further since $w+y^1+x^3,w+y^2+x^3\in A$ and have $f(w+y^1+x^3)=f(w+y^2+x^3)=\gamma$, they are both in $H_\gamma+e^\gamma$.  So for the above to happen, we have to have that $z-w\in H_\gamma+2e^\gamma$.    But now since $f(w+y^3+x^2)>\gamma$, $(z-w)+w+y^3+x^2\in H_\gamma+2e^\gamma$ so $(z-w)+w+y^3+x^2=z+y^3+x^2\notin A$, a contradiction.

Suppose now that (2) holds.  Since $(z-w)+w+x^3+y^2\in A$ while $(z-w)+w+ x^2+y^2\notin A$, and $f(w+x^3+y^2)=f(w+x^2+y^2):=\gamma$, we must have that $f(z-w)=\gamma$ as well.  And moreover, we have to have $(z-w)_\gamma=2$.  But then since $\gamma<f(w+x^2+y^3)=f(w+x^1+y^3)$, we have to have $(z-w)+w+x^2+y^3$, and $(z-w)+w+x^1+y^3\in H_\gamma+2e^\gamma$ so both $z+x^2+y^3, z+x^1+y^3\notin A$, a contradiction.

Suppose now that (3) holds.  Since $a:=f(w+x^2+y^2)=f(w+x^2+y^3)$ and $(z-w)+w+x^2+y^2\notin A$ while $(z-w)+w+x^2+y^3 \in A$, we must have that $f(z-w)=\gamma$ and further, $(z-w)_\gamma=2$ (to cancel the 1).  But then $f(w+x^1+y^3)<\gamma$ implies $(z-w)+w+x^1+y^3=z+x^1+y^3 \in A$, a contradiction.

Suppose now that (4) holds.  As above, since $\gamma:=f(w+x^2+y^2)=f(w+x^2+y^3)$ and $(z-w)+w+x^2+y^2\notin A$ while $(z-w)+w+x^2+y^3 \in A$, so we must have that $f(z-w)=\gamma$ and further, $(z-w)_\gamma=2$ (to cancel the 1).  But since $f(w+x^1+y^1)=\gamma$ and $w+x^1+y^1\notin A$ we know that $f(w+x^1+y^1)=\gamma$ and $(w+x^1+y^1)_\gamma=2$.  But then $(z-w)+(w+x^1+y^1)=z+x^1+y^1$ is in $A$ (since the first nonzero coordinate will be $2+2=1$).

Suppose that (5) holds.  Since $\gamma:=f(w+x^2+y^3)=f(w+x^3+y^3)$ and $w+x^3+y^3\in A$ and $w+x^2+y^3\in A$, but $(z-w)+w+x^2+y^3\in A$ while $(z-w)+w+x^2+y^3\notin A$, we must have that $f(z-w)=\gamma$ and $(z-w)_\gamma=2$.  Similarly, since $\gamma':=f(w+x^3+y^2)=f(w+x^2+y^2)$ and $w+x^2+y^2\in A$ and $w+x^3+y^2\in A$ and since $(z-w)+w+x^3+y^2\in A$ while $(z-w)+w+x^2+y^2\notin A$, we must have $\gamma'=(z-w)$.  Thus $\gamma=\gamma'$, a contradiction to (5).

Suppose that (6) holds.  Since $\gamma:=f(w+x^2+y^3)=f(w+x^3+y^3)$ and both $w+x^2+y^3\in A$, $w+x^3+y^3\in A$, and $(z-w)+w+x^2+y^3\in A$ while $(z-w)+w+x^2+y^3\notin A$, we must have that $f(z-w)=\gamma$ and $(z-w)_\gamma=2$.  By (6), $\gamma=f(w+x^2+y^2)=f(w+x^3+y^2)$, both $w+x^2+y^2\in A$, $w+x^3+y^2\in A$, but $(z'-w)+w+x^2+y^2\notin A$ while $(z'-w)+w+x^3+y^2\in A$, we must have that $\gamma=f(z'-w)$ and $(z'-w)_\gamma=(z-w)_\gamma=2$ (to cancel the $1$).  But (6) implies $\gamma<f(w+x^3+y^1)$, which implies that $(z'-w)+w+x^3+y^1=z'+x^3+y^1\notin A$, a contradiction.
\end{proof}

\begin{lemma}[Necessary conditions on $3$-$\HOP_2$ extensions]\label{lem:3hopfact}
Suppose we have $a^i, b^i, c^i$ for $i\in [3]$ forming an instance of a $3$-$\HOP_2$.  Then none of the following hold.
\begin{enumerate}[label=\normalfont(\arabic*)]
\item There exist $\mathbf{c}^1,\mathbf{a}^4,\mathbf{b}^4\in \F_3^n$, such that if we set $\mathbf{c}^2=c^1,\mathbf{c}^3=c^2, \mathbf{c}^4=c^3$ and for each $i\in [3]$, we set $\mathbf{a}^i=a^i$, and $\mathbf{b}_i=b^i$, then $\{\mathbf{a}^i,\mathbf{b}^i,\mathbf{c}^i: i\in [4]\}$ is a $4$-$\HOP_2$.
\item There exist $\mathbf{a}^1,\mathbf{c}^4,\mathbf{b}^4\in \F_3^n$, such that if we set $\mathbf{a}^2=a^1,\mathbf{a}^3=a^2, \mathbf{a}^4=a^3$ and for each $i\in [3]$, we set $\mathbf{b}^i=b^i$, and $\mathbf{c}^i=c^i$, then $\{\mathbf{a}^i,\mathbf{b}^i,\mathbf{c}^i: i\in [4]\}$ is a $4$-$\HOP_2$.
\item There exist $\mathbf{b}^1,\mathbf{c}^4,\mathbf{a}^4\in \F_3^n$, such that if we set $\mathbf{b}^2=b^1,\mathbf{b}^3=b^2, \mathbf{b}^4=b^3$ and for each $i\in [3]$, we set $\mathbf{a}^i=a^i$, and $\mathbf{c}^i=c^i$, then $\{\mathbf{a}^i,\mathbf{b}^i,\mathbf{c}^i: i\in [4]\}$ is a $4$-$\HOP_2$.
\end{enumerate} 
\end{lemma}
\begin{proof}
By Lemma \ref{lem:3hop}, one face of the $3$-$\HOP_2$ has the same $f$-value.  Without loss of generality, we may assume that $a=u_{i3k}$ each $i,k\in [3]$ and $a<u_{ijk}$ for all $i,k\in [3]$ and $j\in \{1,2\}$.  Let $d=u_{113}$, $s=u_{123}$, $p=u_{122}$, $t=u_{223}$, $h=u_{222}$, $f=u_{221}$, $v=u_{323}$, $u=u_{322}$, $r=u_{321}$, $m=u_{213}$, $c=u_{313}$, $e=u_{312}$, $b=u_{311}$, $g=u_{212}$, $\alpha=u_{211}$, $\beta=u_{112}$, $\delta=u_{121}$, and $\gamma=u_{111}$.  See the diagram below for a geometric representation of the cube with these labels.

\begin{figure}[H]
\begin{tikzpicture}[scale=1, every node/.style={scale=1}]
  \draw[thick](4,4,0)--(0,4,0)--(0,4,4)--(4,4,4)--(4,4,0)--(4,0,0)--(4,0,4)--(0,0,4)--(0,4,4);
  \draw[thick](4,4,4)--(4,0,4);
  \draw[thick](0,4,2)--(4,4,2);
  \draw[thick](2,4,0)--(2,4,4);
  \draw[gray](4,0,0)--(0,0,0)--(0,4,0);
  \draw[gray](0,0,0)--(0,0,4);
  \draw[gray](0,4,2)--(0,0,2);
  \draw[gray](0,2,4)--(0,2,0);
  \draw[gray](2,4,0)--(2,0,0);
  \draw[gray](0,2,0)--(4,2,0);
  \draw[thick](4,4,2)--(4,0,2);
  \draw[thick](4,2,0)--(4,2,4);
  \draw[gray](2,0,0)--(2,0,4);
  \draw[gray](0,0,2)--(4,0,2);
  \draw[thick](0,2,4)--(4,2,4);
  \draw[thick](2,4,4)--(2,0,4);
  \draw[gray,dashed](4,2,2)--(0,2,2);
  \draw[gray,dashed](2,0,2)--(2,4,2);
  \draw[gray,dashed](2,2,4)--(2,2,0);
  
  \draw(0,4,0) node{$\bullet$};
  \draw(0,4.25,0) node{$a$};
  
  \draw(2,4,0) node{$\bullet$};
  \draw(2,4.25,0) node{$a$};
  
  \draw(4,4,0) node{$\bullet$};
  \draw(4,4.25,0) node{$a$};
  
  \draw[gray](0,0,0) node{$\bullet$};
  \draw(-0.25,0.1,0) node{$a$};
  
  \draw(4,0,0) node{$\bullet$};
  \draw(4,0,-0.5) node{$a$};
  
  \draw(0,0,4) node{\textcolor{white}{$\bullet$}};
  \draw(0,0,4) node{$\circ$};
  \draw(0,0,4.5) node{$\gamma$};
  
  \draw(4,0,4) node{$\bullet$};
  \draw(4.35,0,4.4) node{$b$};
  
  \draw(0,4,4) node{$\bullet$};
  \draw(-0.25,4,4) node{$d$};
  
  \draw(4,4,4) node{$\bullet$};
  \draw(4.2,3.95,4) node{$c$};
  
  \draw(0,4,2) node{$\bullet$};
  \draw(-0.2,4.1,2) node{$s$};
  
  \draw(4,4,2) node{$\bullet$};
  \draw(4.2,4,2) node{$v$};
  
  \draw(2,4,2) node{$\bullet$};
  \draw(1.9,4.25,2) node{$t$};
  
  \draw(2,4,4) node{$\bullet$};
  \draw(1.85,4.25,4) node{$m$};
  
  \draw(0,2,4) node{\textcolor{white}{$\bullet$}};
  \draw(0,2,4) node{$\circ$};
  \draw(-0.25,2,4) node{$\beta$};
  
  \draw[gray](0,2,0) node{$\bullet$};
  \draw(-0.2,2.2,0) node{$a$};
  
  \draw[gray](0,2,2) node{$\bullet$};
  \draw(-0.2,2.05,2) node{$p$};
  
  \draw(0,0,2) node{\textcolor{white}{$\bullet$}};
  \draw[gray](0,0,2) node{$\circ$};
  \draw(-0.2,0.1,2) node{$\delta$};
  
  \draw[gray](2,2,0) node{$\bullet$};
  \draw(1.8,2.2,0) node{$a$};
  
  \draw(4,2,0) node{$\bullet$};
  \draw(4.23,2,0) node{$a$};
  
  \draw[gray](2,0,0) node{$\bullet$};
  \draw(1.85,0.2,0) node{$a$};
  
  \draw(4,2,4) node{$\bullet$};
  \draw(4.2,1.8,4) node{$e$};
  
  \draw(4,2,2) node{$\bullet$};
  \draw(4.25,1.9,2) node{$u$};
  
  \draw(4,0,2) node{$\bullet$};
  \draw(4.2,-0.15,2) node{$r$};
  
  \draw(2,2,4) node{$\bullet$};
  \draw(1.8,2.2,4) node{$g$};
  
  \draw[gray](2,0,2) node{$\bullet$};
  \draw(1.82,0.25,2) node{$f$};
  
  \draw(2,0,4) node{\textcolor{white}{$\bullet$}};
  \draw(2,0,4) node{$\circ$};
  \draw(2,-0.3,4) node{$\alpha$};
  
  \draw[gray](2,2,2) node{$\bullet$};
  \draw(1.82,2.25,2) node{$h$};
\end{tikzpicture}

\end{figure}

Suppose towards a contradiction that (1) holds.  So we have new elements $\mathbf{c}^1,\mathbf{a}^4,\mathbf{b}^4\in \F_p^n$ and setting $\mathbf{c}^2=c^1,\mathbf{c}_3=c^2, \mathbf{c}^4=c^3$ and $\mathbf{a}^i=a^i$, $\mathbf{b}^i=b^i$ for $i\in [3]$, we have $\mathbf{a}^i+\mathbf{b}^j+\mathbf{c}^k\in A$ if and only if $i+j+k\geq 6$.

Consider now $a^1,a^2,a^3,\mathbf{a}^4, b^1,b^2,b^3,\mathbf{b}^4$.  Let $y^i=a^i$ and $x^i=b^i$, and let $w=c^2$ and $z=\mathbf{c}_1$.  It is straightforward to check that $w+x^i+y^j\in A$ if and only if $i+j\geq 3$ and $z+x^i+y^j\in A$ if and only if $i+j\geq 5$.  However, we also have that $f(w+y^1+x^3)=a=f(w+y^2+x^3)< f(w+y^i+x^j)$ for all $i\in [3]$ and $j\in [2]$.  Thus condition (1) of Lemma \ref{lem:notallowed} holds, a contradiction.  Thus (1) cannot hold.

Suppose towards a contradiction (2) holds.  So we have new elements $\mathbf{a}^1,\mathbf{b}^4,\mathbf{c}^4\in \F_p^n$ such that, setting, $\mathbf{a}^2=a^1,\mathbf{a}^3=a^2, \mathbf{a}^4=a^3$ and for each $i\in [3]$, $\mathbf{b}^i=b^i$, and $\mathbf{c}^i=c^i$, we have $\mathbf{a}^i+\mathbf{b}^j+\mathbf{c}^k\in A$ if and only if $i+j+k\geq 6$.  Consider $b^1,b^2,b^3,\mathbf{b}^4, c^1,c^2,c^3, \mathbf{c}^4$.  Let $x^i=b^i$, $y^i=c^i$ for $i\in [3]$, $w=a^2$ and let $z=\mathbf{a}^1$.  Then $w+x^i+y^j\in A$ if and only if $i+j\geq 3$ while $z+x^i+y^j\in A$ if and only if $i+j\geq 5$.  However, we also have that $f(w+y^1+x^3)=a=f(w+y^2+x^3)< f(w+y^i+x^j)$ for all $i\in [3]$ and $j\in [2]$.  Thus condition (1) of Lemma \ref{lem:notallowed} holds, a contradiction.  Thus (2) cannot hold.

Suppose towards a contradiction (3) holds.  So we have new elements  $\mathbf{b}^1,\mathbf{a}^4,\mathbf{c}^4\in \F_p^n$, and setting $\mathbf{b}_2=b^1,\mathbf{b}^3=b^2, \mathbf{b}^4=b^3$, and for each $i\in [3]$, setting  $\mathbf{a}^i=a^i$, $c^i=\mathbf{c}^i$, we have $\mathbf{a}^i+\mathbf{b}^j+\mathbf{c}^k\in A$ if and only if $i+j+k\geq 6$.   Setting $x^i=a^i$ and $y^i=c^i$ for $i\in [3]$, $w=b^2$ and $z=\mathbf{b}^4$, we have that $w+x^i+y^j\in A$ if and only if $i+j\geq 3$ while $z+x^i+y^j\in A$ if and only if $i+j\geq 5$.  Consequently, Lemma \ref{lem:notallowed} implies that \emph{none} of the following hold. 

\begin{enumerate}
\item $r=u<t$
\item $t=s<u$.
\item $u=h<t=s$
\item $t=s<u=h$
\item $v=t=u=h>s$
\item $v=t=u=h=\delta$.
\item $s=t=v>u=h=p$
\item $v=t=s=h=p<r$
\end{enumerate}

Suppose first that $t\neq s$ and $r\neq u$.  Note that $\{a^2,a^1,b^1,b^2,c^2,c^3\}$ form a $2$-$\HOP_2$. Consequently, since $t\neq s$, Lemma \ref{lem:2hop} implies we must have $t=m=g=h<s,d,p,\beta$.  Note the following is also a $2$-$\HOP_2$: $\{a^2,a^3,b^1,b^2,c^1,c^2\}$. Since $r\neq u$, Lemma \ref{lem:2hop} implies we must have $e=u=g=h<b,r,f,\alpha$.  So in this case $t=m=g=h=u=e$. Note the following are additive quadruples: $\{a^2+b^1,a^3+b^1, c^2,c^3\}$ and $\{ b^1+a^3,b^2+a^3, c^2,c^3\}$, which correspond to $e,c,g,m$ and $e,c,v,u$.  Therefore, Lemma \ref{lem:c4} along with the equalities we have already established implies that $e=u=g=h=e=c=m=t=v$.  Since (v) and (vi) cannot hold, and $t\neq s$, we have that $v=t=u=h<s$ and $v=t=u=h\neq \delta$.  Note that $\{a^1,a^3,b^2,b^3,c^1,c^3\}$ also form a $2$-$\HOP_2$.  Consequently, Lemma \ref{lem:2hop} along with the fact we know $a<\delta,s,v,r$, and $t=v=u\neq s,r$, we either have $\delta=r=s<v$ or $\delta=v<s,r$.  The first is incompatible with  $v=t=u=h<s$, so we must have $\delta=v$.  But we now have $\delta=t=v=u=h$, contradicting the fact that (vi) cannot hold.

Suppose now $t=s$ and $u\neq r$.  Since (ii) cannot hold, we have $t=s\geq u\neq r$. Note $\{a^2+b^2,a^3+b^2,c^2,c^3\}$ forms an additive quadruple, and  $t\geq u=h$, Lemma \ref{lem:c4} implies that either $t,v>u=h$ or $t=v=u=h$.  In the latter case we would have $t=s>u=h$, contradicting that (iii) cannot hold.  Thus $t=v=u=h$. Similarly since $\{a^1+b^2,a^2+b^2,c^2,c^3\}$ is an additive quadruple and $t=s\geq h$, we must have that either $t=s=h=p$ or $t=s>h=p$. In the latter case we would have $t=s>h=u$, again contradicting that (iii) cannot hold.  Thus $t=s=h=p$.  Overall we have that $t=s=h=p=u=v$.  Since (vi) cannot hold, $t=s=h=p=u=v\neq \delta$.

Note $\{a^1,a^2,b^2,b^3,c^1,c^2 \}$ forms a $2$-$\HOP_2$.  Since we already know $a<h,p,f,\delta$, and $p,h,\neq \delta$, Lemma \ref{lem:2hop} implies that $h=p<\delta,f$.  So we have $h=p=u=h=t=s=v<\delta,f$.

We also have a $2$-$\HOP_2$ consisting of $\{a^1,a^2,b^2,b^3,c^1,c^3\}$. So by Lemma \ref{lem:2hop}, along with the fact that $a<v,s,r,\delta$ and $s=v<\delta$, we have $s=v<\delta,r$.  Thus we have $h=p=u=h=t=s=v<\delta,f,r$.  This contradicts the fact that (viii) cannot hold.

A symmetric argument deals with the case where $t\neq s$ and $u=r$.

So we are left with the case where $t=s$ and $u=r$.  Since (i) and (ii) cannot hold, we must have that $r=u\geq t=s\geq u=r$, so $u=r=s=t$.  Since $\{a^2+b^2,a^3+b^2,c^2,c^3\}$ forms an additive quadruple, and $t=u$, Lemma \ref{lem:c4} implies $v=t=u=h$, so $u=r=s=t=u=h=v$.  Since $\{a^2+b^2,a^3+b^2,c^1,c^2\}$ forms an additive quadruple and $r=h$, Lemma \ref{lem:c4} implies $u=h=r=f$, so we have that $u=r=s=t=h=v=f$.  Similarly, since $\{a^1+b^2,a^2+b^2,c^2,c^3\}$ is an additive quadruple and $s=h$, Lemma \ref{lem:c4} implies $t=s=p=h$.  Consequently, we have that $u=r=s=t=h=v=f=p$.  But now the $2$-$\HOP_2$ formed by $\{a^1,a^2,b^2,b^3,c^1,c^2\}$ fits none of the possibilities in Lemma \ref{lem:2hop}.  In particular, it has $a<h=p=f$ while Lemma \ref{lem:2hop} implies that if $a<h=p$, then $h=p<f,\delta$.  This is a contradiction. 
\end{proof}

We are now able to prove Proposition \ref{prop:GSnoHOP}.

\vspace{2mm}

\begin{proofof}{Proposition \ref{prop:GSnoHOP}}
Suppose towards a contradiction that there are $\{\mathbf{a}^i, \mathbf{b}^i, \mathbf{c}^i: i\in [4]\}$ so that $\mathbf{a}^i+\mathbf{b}^j+\mathbf{c}^k\in A=A_{\GS}(p,n)$ if and only $i+j+k\geq 6$.  Setting $c^1=\mathbf{c}^2$, $c^2=\mathbf{c}^3$, $c^3=\mathbf{c}^4$ and for each $i\in [3]$, $a^i=\mathbf{a}^i$ and $b^i=\mathbf{b}^i$, we have a contradiction to Lemma \ref{lem:3hopfact}.  Thus no such $\{\mathbf{a}^i, \mathbf{b}^i, \mathbf{c}^i: i\in [4]\}$ exists. This shows that $A$ has no $4$-$\HOP_2$.
\end{proofof}




\providecommand{\bysame}{\leavevmode\hbox to3em{\hrulefill}\thinspace}
\providecommand{\MR}{\relax\ifhmode\unskip\space\fi MR }
\providecommand{\MRhref}[2]{%
  \href{http://www.ams.org/mathscinet-getitem?mr=#1}{#2}
}
\providecommand{\href}[2]{#2}

\end{document}